\pdfoutput=1
 \documentclass[10pt]{article}
\usepackage[utf8]{inputenc} 
\usepackage[T1]{fontenc}    
\usepackage{comment}      
\usepackage{url}            
\usepackage{nicefrac}       
\usepackage{microtype}     
\usepackage{fullpage}

\usepackage{pifont}
\usepackage{xcolor}
\usepackage{framed}
\colorlet{shadecolor}{pink} 
\usepackage{authblk}
\usepackage{adjustbox}
\usepackage{comment}
\usepackage{multirow}
\usepackage{amsfonts} 	
\usepackage{color, colortbl}
\usepackage{graphicx}
\usepackage{soul}
\usepackage{subcaption}
\usepackage{booktabs}
\usepackage{tablefootnote}

\usepackage[colorlinks,linkcolor=magenta,citecolor=blue, pagebackref=true,backref=true]{hyperref}

\usepackage{amsmath,amsthm,amssymb,amsfonts}
\usepackage{algorithm}
\usepackage{algorithmic}
\usepackage{enumerate}
\usepackage{cleveref}
\usepackage{comment}
\usepackage{bm}
\usepackage{pifont}
\theoremstyle{plain} 
\newcommand{\vertiii}[1]{{\left\vert\kern-0.25ex\left\vert\kern-0.25ex\left\vert #1 
\right\vert\kern-0.25ex\right\vert\kern-0.25ex\right\vert}}

\usepackage{algorithm}
\usepackage[ruled,algo2e]{algorithm2e}

\definecolor{antiquewhite}{rgb}{0.98, 0.92, 0.84} 
\definecolor{blizzardblue}{rgb}{0.67, 0.9, 0.93}

\usepackage{epsf}
\usepackage{graphics}
\usepackage{wrapfig}
\usepackage{psfrag}

\usepackage{color}

\usepackage{mathtools}
\usepackage{amsfonts}
\usepackage{amsthm}
\usepackage{amsmath}
\usepackage{amssymb}

\newcommand{\pr}[1]{\left( #1 \right)} 
\newcommand{\cbr}[1]{\left\{ #1 \right\}}

\newcommand{\inprod}[2]{\ensuremath{\left\langle #1 , \, #2 \right\rangle}}
\newcommand{\pare}[1]{\left( #1 \right)}

\newtheorem{theorem}{Theorem}
\newtheorem{proposition}{Proposition}
\newtheorem{lemma}{Lemma}

\newtheorem{definition}{Definition}

\newtheorem{assumption}{Assumption}

\long\def\comment#1{}

\newcommand{\norm}[1]{\left\| #1 \right\|}

\newcommand{\E}{\ensuremath{{\mathbb{E}}}}

\DeclareMathOperator{\diag}{diag}

\newcommand{\R}{\mathbb{R}}

\newcommand{\cB}{\mathcal{B}}  
\newcommand{\cD}{\mathcal{D}}  
\newcommand{\cF}{\mathcal{F}}  
\newcommand{\cH}{\mathcal{H}}  
\newcommand{\cI}{\mathcal{I}}

\newcommand{\cM}{\mathcal{M}}  
\newcommand{\cS}{\mathcal{S}}
\newcommand{\cT}{\mathcal{T}}   
\newcommand{\cW}{\mathcal{W}}   
  
\newcommand{\cP}{\mathcal{P}}  
\newcommand{\cL}{\mathcal{L}} 
   
\newcommand{\cX}{\mathcal{X}}  
\newcommand{\cY}{\mathcal{Y}}  

\newcommand{\balpha}{\boldsymbol{\alpha}}

\newcommand{\ba}{\mathbf{a}}
\newcommand{\bb}{\mathbf{b}}

\newcommand{\bx}{\mathbf{x}}
\newcommand{\bu}{\mathbf{u}}

\newcommand{\bv}{\mathbf{v}}
\newcommand{\w}{\mathbf{w}}
\newcommand{\bw}{\mathbf{w}}
\newcommand{\bq}{\mathbf{q}}

\newcommand{\bz}{\mathbf{z}}

\newcommand{\by}{\mathbf{y}}
\newcommand{\bg}{\mathbf{g}}

\newcommand{\mI}{{\bf I}}

\newcommand{\mM}{{\bf M}}

\usepackage[colorinlistoftodos,bordercolor=orange,backgroundcolor=orange!20,linecolor=orange,textsize=scriptsize]{todonotes}

\title{\bf Stochastic Compositional Minimax Optimization with Provable  Convergence Guarantees}

\author{Yuyang Deng \qquad Fuli Qiao \qquad Mehrdad Mahdavi \vspace*{.2em} \\ 
The Pennsylvania State University \vspace*{.2em} \\ 
\texttt{ \{yzd82,fvq5015,mzm616\}@psu.edu}
}
\date{}

\begin{document}

\maketitle

\begin{abstract}
Stochastic compositional minimax problems are prevalent in machine learning, yet there are only limited established on the convergence of  this class of problems. In this paper, we propose a formal definition of the stochastic compositional minimax problem, which involves optimizing a minimax loss with a compositional structure either in primal , dual, or both primal and dual variables. We introduce a simple yet effective algorithm, stochastically Corrected stOchastic gradient Descent Ascent (CODA), which is a descent ascent type algorithm with compositional correction steps, and establish its convergence rate in aforementioned three settings. In the presence of the compositional structure in primal, the objective function typically becomes nonconvex in primal due to function composition. Thus, we consider the nonconvex-strongly-concave and nonconvex-concave settings and show that CODA can efficiently converge to a stationary point. In the case of composition on the dual, the objective function becomes nonconcave in the dual variable, and we demonstrate convergence in the strongly-convex-nonconcave and convex-nonconcave setting. In the case of composition on both variables, the primal and dual variables may lose convexity and concavity, respectively. Therefore, we anaylze the convergence in weakly-convex-weakly-concave setting. We also give a variance reduction version algorithm, CODA+, which achieves the best known rate on nonconvex-strongly-concave and nonconvex-concave compositional minimax problem.  This work initiates the theoretical study of the stochastic compositional minimax problem on various settings and may inform modern machine learning scenarios such as domain adaptation or robust model-agnostic meta-learning.    
\end{abstract}
  \section{Introduction}

We consider the following minimax optimization problem:
\begin{align}
    \min_{\bx \in \cX} \max_{\by \in \cY} F(\bx,\by) := h(\bx) +  f(\bx,\by)  - r(\by) \label{eq:obj},
\end{align}
where $f(\bx,\by)$ is a stochastic compositional function. The compositional part can be defined over $\bx$, $\by$ or both $\bx$ and $\by$, which correspond to the following three scenarios: 
\begin{itemize}
    \item Composition on $\bx$:   $ f(\bx,\by):=  \E_{\zeta \sim \cS_f}[f( \E_{\xi\sim \cS_g}[g( \bx;\xi)],\by;\zeta) ]$
    \item Composition on $\by$: 
    $f(\bx,\by) :=  \E_{\zeta\sim \cS_f}[ f(\bx,\E_{\xi\sim \cS_g}[g( \by;\xi)];\zeta) ] $ 
    \item Composition on $\bx$ and $\by$:  
    $ f(\bx,\by) := \E_{\zeta\sim \cS_f}[ f(\E_{\xi\sim S_g}[g(\bx, \by;\xi)];\zeta) ]$,
\end{itemize}
where $\cS_f$ and $\cS_g$ capture inner and outer stochasticity in composition, respectively. 
  
Many practical machine learning problems can be cast as an instance of the above optimization problem, such as robust meta learning~\cite{collins2020task}, AUC maximization~\cite{zhang2023federated}, statistical discrepancy minimization~\cite{mansour2009domain}, multiple source domain adaptation~\cite{konstantinov2019robust,mansour2021theory,deng2023mixture}, and stochastic compositional constrained programming~\cite{yang2022stochastic}. Here we briefly describe several important applications in machine learning that motivate the proposed compositional minimax problem (\ref{eq:obj}).

\paragraph{AUC Optimization~\cite{yuan2021compositional}}~AUC optimization has been a powerful training paradigm that can mitigate the training data imbalance issue in classification problems.
Consider a binary classification problem with label $y \in \{-1,+1\}$, and model parameter $\bw \in \R^d$.  Our goal is to solve the following problem:
\begin{align*}
   \min_{\bw, a, b} \max_{\theta \in \Omega}  &\Phi(\bw - \alpha \nabla L(\bw), a,b,\theta)  = \frac{1}{n} \sum_{i=1}^n \phi(\bw - \alpha \nabla L(\bw), a,b,\theta;\bx_i,y_i),
\end{align*}
where $L(\bw) = \frac{1}{n}\sum_{i=1}^n \ell(\bw;\bx_i,y_i)$. Usually, the objective has a negative quadratic dependency on dual variable $\theta$, which renders this objective a nonconvex-strongly-concave minimax problem.\\
\paragraph{Robust MAML~\cite{collins2020task}}~In task robust MAML, our goal is to learn a model that can perform well over all distributions of the observed tasks, hence making it robust to task shifts in the testing time. The objective can be formulated as the following minimax problem:
\begin{align*}
    \min_{\bw \in \R^d}\max_{\balpha \in \Delta^N}  \sum_{i=1}^N \alpha(i) L_i(\bw - \eta \nabla L_i(\bw))
\end{align*}
where $\Delta^N$ denotes simplex, $L_i(\bw) := \frac{1}{n_i} \sum_{j=1}^{n_i}\ell(\bw;\xi_j)$, $N$ is the number of tasks with task $i$ having $n_i$ samples. It is a nonconvex-concave compositional minimax problem, with composition on the primal variable.

\paragraph{Stochastic Minimization with Compositional Constraint~\cite{yang2022stochastic}} 
In many risk management problems, such as conditional value-at-risk (CVaR)~\cite{rockafellar2000optimization}, risk-averse mean-deviation constraint~\cite{10.1093/rfs/hhw080}, or portfolio optimization with high-moment constraint~\cite{harvey2010portfolio}, we are tasked with solving the following programming with compositional constraint: 
\begin{align*}
    \min_{\bx \in \R^d} f(\bx)  \ \quad    \text{s.t.} \ \quad & g_1(\E[g_2(\bx;\xi)])\leq 0.
\end{align*}
To solve the above constraint problem, one may  consider solving the unconstrained variant:
\begin{align*}
    \min_{\bx} \max_{\lambda} f(\bx) +  \lambda g_1(\E[g_2(\bx;\xi)]) 
\end{align*}
which turns the original objective into a nonconvex-concave compositional minimax problem, with composition on the primal variable.\\
\paragraph{Mixture Weights Estimation~\cite{konstantinov2019robust}}
In multiple source domain adaptation, the goal is learn a predictive model for a target domain by  finding a mixture weights $\balpha$ to optimally aggregate empirical risks over $N$ data sources, which can be achieved by solving the following minimax problem:
\begin{align*}
    \min_{\balpha \in \Delta^N} \max_{\bw \in \cW }\sum_{i=1}^N \alpha(i) f(L_T(\bw) - L_i(\bw)) + \lambda \norm{\balpha}_{\mM} 
\end{align*} 
where $\balpha$ is the mixture weight, $\bx$ is the model parameter, $L_T(\bw)$ is the empirical risk realized by target domain data, and $L_i(\bw)$ is the empirical risk realized by $i$th source domain data. $f(z) $ is the smooth approximation of the absolute value function. $\norm{\cdot}^2_\mM = \balpha^\top \mM \balpha$ where $\mM = \diag\pr{m_1,...,m_N}$. It is a convex-nonconcave minimax problem with composition on the dual variable, and if we relax the norm of $\balpha$ to the squared norm, it becomes a strongly-convex-nonconcave minimax problem.\\
\paragraph{Discrepancy Minimization~\cite{mansour2009domain,sriperumbudur2009integral}}~In the domain adaptation and  empirical estimation of integral probability metrics (IPMs), the goal is  to estimate the statistical discrepancy among two distributions or to manipulate a source distribution to align it well with a target distribution, leading to the following  problem:
\begin{align*}
    \min_{  \cS}  \sup_{\w,\w'\in\cW} |L_{ \cS}(\w,\w') - L_{ \cT}(\w,\w')|  
\end{align*}
where $L_\cS$ and $L_{\cT}$ denote empirical risks over source and target domains, respectively. Here the goal is to find a (empirical) distribution $\cS$, such that its discrepancy with target distribution $\cT$ is minimized, leading to a nonconvex-nonconcave compositional minimax problem.\\

Despite the well-established convergence theory on non-compositional minimax,  there is few literature focusing on its compositional counterpart. The prior works on compositional minimax~\cite{gao2021convergence,liu2023breaking} only consider the case of function composition on the primal variable, and hence they only establish the convergence rate in the nonconvex-strongly-concave setting. However,  a significant portion of ML applications including aforementioned applications do not satisfy this assumption. The function composition can also be on the dual variable, or on both primal and dual variables, which makes the problem convex-nonconcave or even nonconvex-nonconcave. To fill this gap, we propose a simple yet effective algorithm dubbed as stochastically Corrected stOchastic gradient Descent Ascent (CODA). The proposed algorithm is a variant of the celebrated Stochastic Gradient Descent-Ascent (SGDA) for minimax optimization, with a compositional correction step  inspired from~\cite{chen2021solving}. Depending on the objective structure, the correction step is either applied on the primal update, dual update, or both updates, which yields three algorithm instances: CODA-Primal, CODA-Dual, and CODA-PD (Primal-Dual), respectively. We establish their convergence in nonconvex-(strongly)-concave, (strongly)-convex-nonconcave and weakly-convex-weakly-concave settings correspondingly. Our results demonstrate that  intricate compositional structure inherent to these minimax problems does not impose an  asymptotic increase to their computational complexity compared to non-compositional variants. This coincides with the findings in compositional minimization setting~\cite{chen2021solving}.

\begin{table*}[t!]
\centering
\footnotesize
\resizebox{ \columnwidth}{!}{%
\begin{tabular}{|c|c|c|c|} 
\hline
Objective type & Assumption &     {Algorithm}&  {Complexity}  \\ 
\hline \hline  \\
 [-1.1em] 
\multirow{2}{*}{Composition on $\bx$} &  \multirow{2}{*}{Nonconvex-strongly-concave}
  &  SCGDA~\cite{gao2021convergence} &$O(\frac{\kappa^4}{\epsilon^4})$ \\
  \cline{3-4} 
  &  &   \cellcolor{antiquewhite} CODA-Primal (Theorem~\ref{thm:NCSC})  & \cellcolor{antiquewhite} $O(\frac{\kappa^4}{\epsilon^4})$  \\
 \cline{3-4}
   &  & $^\dagger$NSTORM~\cite{liu2024faster} &$O(\frac{\kappa^3}{\epsilon^3})$ \\
  \cline{3-4}
   &  &  \cellcolor{antiquewhite} $^\dagger$CODA-Primal+ (Theorem~\ref{thm:NCSC+}) & 
 \cellcolor{antiquewhite} $O(\frac{\kappa^2}{\epsilon^3})$ \\
 \cline{2-4}
  &  \multirow{2}{*}{Nonconvex-concave}&  \cellcolor{antiquewhite} CODA-Primal (Theorem~\ref{thm:ncc}) & \cellcolor{antiquewhite} $O(\frac{1}{\epsilon^8})$\\
   \cline{3-4}
  &  &   \cellcolor{antiquewhite} CODA-Primal+ (Theorem~\ref{thm:ncc+}) & \cellcolor{antiquewhite} $O(\frac{1}{\epsilon^7})$\\
 \hline
 \multirow{2}{*}{Composition on $\by$}  & {Strongly-convex-nonconcave}&  \cellcolor{antiquewhite} CODA-Dual (Theorem~\ref{thm:SCNC}) &\cellcolor{antiquewhite} $O(\frac{\kappa^5}{\epsilon^4})$ \\
  \cline{2-4}
    & {Convex-nonconcave}&  \cellcolor{antiquewhite} CODA-Dual (Theorem~\ref{thm: CNC}) & \cellcolor{antiquewhite}$ O(\frac{1}{\epsilon^8})$  \\
 \hline
 {Composition on $\bx$ and $\by$}   & {Weakly-convex-weakly-concave}&  \cellcolor{antiquewhite} CODA-PD (Theorem~\ref{thm:wcwc}) &\cellcolor{antiquewhite} $ O(\frac{1}{\epsilon^4})$ \\
\hline 
\end{tabular}}\label{tab:results} 
\caption{A summary of prior works and our results on compositional minimax optimization, in different objective settings.  $^\dagger$: These algorithms use variance reduction techniques.  }\vspace{-0.6cm}
\end{table*}

\paragraph{Contributions} The main contributions of the present work are summarized as follows:
\begin{itemize} 
    \item We introduce a general algorithmic framework, dubbed CODA,  for stochastic compositional minimax optimization that encompasses three fundamental scenarios of objective composition: composition on the primal variable (Section~\ref{sec:com:primal}), composition on the dual variable (Section~\ref{sec:com:dual}), and composition on both variables (Section~\ref{sec:com:both}). We provide a rigorous convergence analysis for various CODA variants tailored to these settings. \vspace{1mm}
    \item In composition on the primal variable setting, we examine both nonconvex-strongly-concave and nonconvex-concave cases. Remarkably, our analysis demonstrates that the CODA-Primal algorithm exhibits the same convergence rate as the widely acclaimed SGDA algorithm in their respective non-compositional counterparts. \vspace{1mm}
    \item  In composition on the dual variable setting, we study strongly-convex-nonconcave and convex-nonconcave cases, and show that our CODA-Dual algorithm can provably converge. This is also the first work that establishes the convergence guarantee on stochastic (compositional) minimax optimization on the general convex-nonconcave objective. \vspace{1mm}
    \item In composition on the both variable setting,  we study weakly-convex-weakly-concave cases and also provide the convergence guarantee of our CODA-PD algorithm. This work represents the first instance of a compositional minimax algorithm that guarantees convergence in scenarios featuring weakly-convex-weakly-concave objectives. \vspace{1mm}
    \item We further propose a variance reduced variant, CODA+, which achieves the best-known convergence rates on nonconvex-strongly-concave and nonconvex-concave compositional minimax problems (Section~\ref{sec:VR}). \vspace{1mm}
    \item  Lastly, we conduct  experiments across various applications that corroborate the theoretical analysis and demonstrate the effectiveness of the proposed algorithms (Section~\ref{sec:exp}).
\end{itemize}

\section{Related Works}\label{app:related}
\paragraph{Compositional Optimization}
Stochastic compositional optimization has wide applications in risk management and machine learning. Wang et al~\cite{wang2017stochastic} is the pioneering work in this field where they proposed the first compositional optimization algorithm, SCGD, and achieved $O(\epsilon^{-4})$ rate. An accelerated version of SCGD is then proposed by~\cite{wang2017accelerating}, and achieves faster rate $O(\epsilon^{-2.25})$. Tutunov et al~\cite{tutunov2020compositional} also proposed an accelerated compositional optimization algorithm and achieved the same rate as~\cite{wang2017accelerating}. The SOTA rate is achieved by~\cite{chen2021solving} where their rate matches with the non-compositional setting, demonstrating that solving the compositional minimization problem can be as easy as the non-compositional counterpart.

A line of works also study the compositional minimax problem, particularly in the AUC maximization area. Yuan et al~\cite{yuan2021compositional} proposed to optimize a compositional minimax objective, which improves the AUC training results. In their setting, the objective is nonconvex-strongly-concave, and they achieve $O( {\epsilon^{-4}})$ gradient complexity. Deng et al~\cite{deng2023mixture} studied the mixture weight estimation problem and cast it as a strongly-convex-nonconcave compositional minimax problem. They also proposed an SGDA-based compositional algorithm and achieved $O( \epsilon^{-4})$ gradient complexity to find the stationary point. Gao et al~\cite{gao2021convergence} for the first time study the general compositional minimax problem. They consider nonconvex-strongly-concave setting, and achieve $O( {\kappa^4}{\epsilon^{-4}})$ rate. Liu et al~\cite{liu2023breaking} proposed the first variance reduction algorithm for stochastic compositional minimax problem and achieved $O( {\kappa^3}{\epsilon^{-3}})$ rate.

\paragraph{Nonconvex/Nonconcave Minimax optimization}
Nonconvex minimax optimization has received increasing attention over the past decade. For nonconvex-strongly-concave minimax optimization,
 Lin et al~\cite{lin2019gradient} proved the first non-asymptotic convergence of GDA to $\epsilon$-stationary point of primal function, with the gradient complexity of $O({\kappa^2}{\epsilon^{-2}})$. Mahdavinia et al~\cite{mahdavinia2022tight} proved that OGDA and EG method enjoy the same rate in solving nonconvex-strongly-concave problem.  Yang et al~\cite{yang2021faster} developed a single loop smoothed AGDA, which achieved gradient complexity of $O( \kappa \epsilon^{-2})$. Lin et al~\cite{lin2020near} proposed an accelerated gradient method based algorithm, and achieved the best known rate $O( \sqrt{\kappa}\epsilon^{-2})$. Zhang et al~\cite{zhang2021complexity} proposed triple loop algorithms which also achieved gradient complexity of $O( \sqrt{\kappa}\epsilon^{-2})$ using catalyst idea, and inexact proximal point method. These two results match the existing lower bound for deterministic setting~\cite{Li2021ComplexityLB,zhang2021complexity,Han2021LowerCB}.  For stochastic setting, Lin et al~\cite{lin2019gradient} proved that  SGDA finds $\epsilon$-stationary point of primal function, with the gradient complexity of  $O( \kappa^3 \epsilon^{-4})$ respectively. Chen et al~\cite{chen2021tighter} proposed a double loop algorithm and using this idea they achieved gradient complexity of $O( \kappa^3 \epsilon^{-4})$ with fixed batch size.  Yang et al~\cite{yang2021faster} introduced stochastic smoothed AGDA with gradient complexity of $O( \kappa^2 \epsilon^{-4})$ to find stationary point. The significant improvement over SGDA is that stochastic smoothed AGDA uses fixed batch size. Zhang et al~\cite{zhang2022sapd+} proposed SAPD+, which achieve the SOTA rates of $O( \kappa \epsilon^{-4} )$ and $O(\kappa^2\epsilon^{-3})$ under average smoothness and point-wise smoothness assumptions respectively.

For nonconvex-concave optimization, Rafique et al~\cite{rafique2018non} proposed the pioneering work Proximally Guided Stochastic Mirror Descent Method, which achieves $O({\epsilon^{-6}})$ gradient complexity to find the stationary point. Nouiehed et al~\cite{nouiehed2019solving} presented a double-loop algorithm with $O(\epsilon^{-7})$ rate. Lin et al~\cite{lin2020gradient} gave the first proof of GDA and SGDA on nonconvex-concave functions, with $O({\epsilon^{-6}})$ rate and $O({\epsilon^{-8}})$ rate respectively. Zhang et al~\cite{zhang2020single} proposed smoothed-GDA and also achieve $O(\epsilon^{-4})$ rate. Thekumparampil et al~\cite{thekumparampil2019efficient} proposed Proximal Dual Implicit Accelerated Gradient method and achieved the best known rate $O({\epsilon^{-3}})$ for nonconvex-concave problem. Kong et al~\cite{kong2019accelerated} proposed an accelerated inexact proximal point method and also achieve $O({\epsilon^{-3}})$ rate. Lin et al~\cite{lin2020near} designed 
an algorithm which uses accelerated gradient method as a sub-problem solver, with $O(\epsilon^{-3})$ rate. Zhang et al~\cite{zhang2022sapd+} proposed a variance reduction algorithm SAPD+, which achieve the $O({\epsilon^{-6}})$ convergence rate in stochastic setting. 

For convex-nonconcave optimization, Xu et al~\cite{xu2023unified} gave the first convergence proof in this setting, where they consider alternating gradient descent ascent. They achieved $O( \epsilon^{-2})$ and $O( \epsilon^{-4})$ rate for strongly-convex-nonconcave and convex-nonconcave setting respectively. Later on Deng et al~\cite{deng2023mixture} gave a stochastic algorithm for a special case of compositional strongly-convex-nonconcave minimax problem and achieved $O( \epsilon^{-4} )$ rate.

For nonconvex-nonconcave optimization, first-order stationary points may not exist if we do not post any further assumptions. Hence, a line of work~\cite{liu2021first,diakonikolas2021efficient} made extra assumptions and developed convergent algorithms. Liu et al~\cite{liu2021first} assumed there is a solution of the MVI induced by objective, and gave the first provably convergent algorithm framework for nonconvex-nonconcave minimax problem. Diakonikolas et al~\cite{diakonikolas2021efficient} considered the relaxed assumption that only weak MVI solution exists and proposed an EG-based algorithm with convergence guarantee. Zheng et al~\cite{zheng2023universal} consider nonconvex-nonconcave minimax problem with one-side KL condition, and proposed Doubly Smoothed GDA which enjoys $O( \epsilon^{-4})$ convergence rate.
Some work also explores to design new convergence measures in general nonconvex-nonconcave minimax problems. Mangoubi et al~\cite{mangoubi2021greedy} proposed a greedy adversarial equilibrium and developed an algorithm with polynomial gradient and hessian complexity to find it.

\section{Composition on Primal Variable}\label{sec:com:primal}
We begin by considering the  composition in the primal variable:
\begin{align*} 
    F(\bx,\by) := h(\bx) + \E_{\zeta\sim \cS_f}[ f( \E_{\xi\sim \cS_g}[g( \bx);\xi],\by;\zeta) ] - r(\by) . 
\end{align*}
Due to the  compositional structure, the objective will become nonconvex in $\bx$, and hence we consider two settings: nonconvex-strongly-concave and nonconvex-concave. Our algorithm builds on top of classic SGDA algorithm for minimax optimization. However, due to the compositional structure, we cannot have an unbiased estimation of the gradient with respect to $\bx$. To see this, notice the fact: $ 
   \E_{\xi,\zeta}[ \nabla_1 f( g(\bx;\xi),\by ;\zeta ) \nabla g(\bx;\xi)] \neq \nabla_1 f( g(\bx ),\by ) \nabla g(\bx ). $ 
Hence, we borrow the technique from compositional minimization problem~\cite{chen2021solving} and maintain an auxiliary variable $\bz$ to approximate inner function $g(\bx)$. Given a mini-batch  $\cM^t$ with size $M$, we update $\bz$ as:
\begin{align*}
  \bz^{t+1} =  (1-\beta) \pare{\bz^t + g(\bx^t; \cM^t)  -  g(\bx^{t-1}; \cM^t)  }  + \beta  g(\bx^t; \cM^t) , 
\end{align*}
where $g(\bx;\cM^t):= \frac{1}{M} \sum_{\xi \in \cM^t} g(\bx;\xi)$.
Then, the proposed CODA-Primal algorithm computes the following primal-dual gradients using a mini-batch $\cB^t$ with size $B$:
\begin{align*} 
          & \bg_{\bx}^t = \frac{1}{B}\sum_{(\zeta,\xi) \in \cB^t}\nabla_1 f( \bz^{t+1},\by^t;\zeta )  \nabla g(\bx^t ; \xi) + \nabla h(\bx^t)  , \bg_{\by}^t =  \frac{1}{B}\sum_{(\zeta,\xi) \in \cB^t} \nabla_2 f( \bz^{t+1},\by^t;\zeta )  - \nabla r(\by^t),
\end{align*}
and then performs primal-dual updates on $\bx$ and $\by$: ${\bx}^{t+1} =  \bx^{t} - \eta_{\bx} \bg_{\bx}^t  ,  \by^{t+1} =\cP_{\cY}\pare{ \by^{t} + \eta_y \bg_{\by}^t}$.

\noindent\textbf{\ding{110}~Nonconvex-Strongly-Concave Setting.} We start by providing   the convergence rate of CODA-Primal (Algorithm~\ref{algorithm: CODA-Primal}), in the nonconvex-strongly-concave setting. From an algorithmic standpoint, the proposed method is reminiscent of existing methods, but the convergence analysis is more involved due to the minimax, compositional structure, and non-convexity of the objective. We make the following assumptions on the objective:
\begin{assumption}\label{ass:ncsc}
   $F(\bx,\by)$ is $\mu$-strongly-concave in $\by$,  $\forall \bx \in \R^d$.
\end{assumption}
\begin{assumption}\label{ass:bounded Y}
    Domain $\cY$ has bounded diameter, i.e., $\norm{\by-\by'} \leq D_{\cY}$ $\forall  \by,\by'\in\cY$. 
\end{assumption}

\begin{assumption}\label{ass:primal standard}
    We make the following assumptions:
    \begin{enumerate}  
        \item  $f(\bz,\by )$ is $L_f$ smooth,  $g(\bx)$ is $L_g$ smooth , and  $F(\bx,\by)$ is $L $ smooth.\vspace{-1mm}
        \item $\E[\norm{\nabla g(\bx;\xi)}^2] \leq G_g^2$ and $\E[\norm{\nabla f(\bz,\by;\zeta)}^2] \leq G_f^2$.\vspace{-1mm}
        \item $ \max\cbr{\E\norm{  g(\bx;\xi) -    g(\bx ) }^2, \E\norm{  \nabla g(\bx;\xi) -   \nabla g(\bx ) }^2,\E\norm{ \nabla f(\bz,\by;\zeta)   - \nabla f(\bz,\by ) }^2  }\leq \sigma^2$.\vspace{-1mm}
    \end{enumerate}
\end{assumption}

\begin{algorithm2e}[t] 
	\DontPrintSemicolon
    \caption{\sffamily{CODA-Primal}}
    	\label{algorithm: CODA-Primal}
    	\textbf{Input:} $({\bx}^{0},{\by}^{0})$,  stepsizes $\eta_{\bx},\eta_{\by}$,   $\bz^0$ such that $\E\norm{\bz^0 - g(\bx^0)}^2 \leq \delta$.
	\\

   {
    {

	\For{$t = 0, \ldots,T$}
	 {   Sample a mini-batch $\cM^t$ of $M$ data from $\cS_g$.\\ 
            $\bz^{t+1} = (1-\beta)$
            $ \pare{\bz^t + \frac{1}{M}\sum_{\xi^t \in \cM^t} \pr{ g(\bx^t; \xi^t)  -  g(\bx^{t-1}; \xi^{t })}  }  + \beta \frac{1}{M}\sum_{\xi^t \in \cM^t} g(\bx^t; \xi^t)  $. \\
         Sample a mini-batch $\cB^t $ of $B$ data pair $(\zeta^t,\xi^t)$ from $\cS_g \times \cS_f$.\\
            Compute $\bg_{\bx}^t = \frac{1}{B}\sum_{(\zeta,\xi) \in \cB^t}\nabla_1 f( \bz^{t+1},\by^t;\zeta )  \nabla  g(\bx^t ; \xi) + \nabla h(\bx^t) $, \ $\bg_{\by}^t =  \frac{1}{B}\sum_{(\zeta,\xi) \in \cB^t} \nabla_2 f( \bz^{t+1},\by^t;\zeta )  - \nabla r(\by^t) .  $\\
 	    ${\bx}^{t+1} =  \bx^{t} - \eta_{\bx} \bg_{\bx}^t $, $\by^{t+1} =\cP_{\cY}\pare{ \by^{t} + \eta_\by \bg_{\by}^t} $.\\ 
 	 }
      
 	      }
 	     }  
       \textbf{Output:}  $(\hat{\bx}, \hat{\by})$ uniformly sampled from $\cbr{ ( {\bx}^t,  {\by}^t)   }_{t=1}^{T}$.
\end{algorithm2e}
Since the objective is nonconvex, we are unable to show the convergence to the global saddle point. Thus, following the standard machinery in nonconvex-concave analysis~\cite{lin2019gradient,thekumparampil2019efficient,rafique2018non}, we introduce the following primal function to measure the convergence.
\begin{definition} We define  primal function $  \Phi(\bx) = \max_{\by\in\cY} F(\bx,\by )$ to facilitate our analysis.
\end{definition}
\noindent We consider the convergence rate to the first order stationary point of $ \Phi(\bm{x})$, as advocated in seminal nonconvex-concave minimax literature~\cite{lin2019gradient,rafique2018non,thekumparampil2019efficient}.  
\begin{theorem}\label{thm:NCSC}
Under Assumptions~\ref{ass:ncsc},~\ref{ass:bounded Y},~\ref{ass:primal standard}, defining $\kappa:= L/\mu$, for  Algorithm~\ref{algorithm: CODA-Primal}, if we choose $\delta = \frac{1}{\kappa}$, $M=B = \Theta\pr{\max\cbr{\frac{\kappa^2 L   \sigma^2}{\epsilon^2},1 }}$,  $\beta =  \frac{1}{2} $, $\eta_{\bx} = \Theta\pare{\frac{1}{\kappa^2 L}}$,  $\eta_{\by} = \Theta\pare{\frac{1}{L}}$, then Algorithm~\ref{algorithm: CODA-Primal} guarantees that $   \frac{1}{T}\sum_{t=1}^T \E\norm{\nabla \Phi(\bx^t )}^2 \leq \epsilon^2 $ 
with the gradient complexity bounded by:
 \begin{align*}
  O\pr{\frac{\kappa^2( \Delta_{\Phi} +  L^2 D_{\cY} )}{\epsilon^2}\max\cbr{\frac{\kappa^2 L   \sigma^2}{\epsilon^2},1 } },
 \end{align*} \vspace{-2mm}
where $\Delta_{\Phi}:=\Phi(\bx^0) - \min_{\bx\in\cX}\Phi(\bx)$.
\end{theorem} 
\noindent The proof of Theorem~\ref{thm:NCSC} is deferred to Appendix~\ref{app: proof NCSC}.   Here we achieve an $O(\frac{\kappa^4}{\epsilon^4})$ rate, which is the same as the rate of SGDA on non-compositional objective~\cite{lin2019gradient}. This suggests that we can solve nonconvex-strongly-concave compositional minimax problems as readily as non-compositional cases. The most relevant work is~\cite{gao2021convergence}, where they also achieve an $O(\frac{\kappa^4}{\epsilon^4})$ rate.  An even faster rate is attainable through the use of variance reduction techniques, as will be shown and discussed in  Section~\ref{sec:VR}. As a final remark, we need to choose a good initialization of $\bz^0$ such that it is close to $g(\bx^0)$, i.e., $\E\norm{ \bz^0 - g(\bx^0) }^2 \leq \frac{1}{\kappa}$. This will take $\kappa\sigma^2$ samples drawn from $\cS_g$.

\noindent\textbf{\ding{110}~Nonconvex-Concave Setting.} We now turn to the analysis in  nonconvex-concave setting. 
\begin{assumption}\label{ass:ncc}
  $F(\bx,\by)$ is concave in $\by$,   $\forall \bx \in \R^d$. $h(\bx )$ is $G_h$ Lipschitz.
\end{assumption} 
Since $F$ is merely concave in $\by$, primal function $\Phi(\bx)$ is not smooth anymore, and its gradient is not a viable quantity to measure the convergence. Following the convention~\cite{lin2019gradient}, we consider the gradient of Moreau Envelope of primal function as the convergence measure.

\begin{definition}[Moreau Envelope] A function $\Phi_{\rho} (\bw)$ is the $\rho$-Moreau envelope of a function $\Phi$ if $    \Phi_{\rho} (\bw) := \min_{\bw'\in \cW} \{ \Phi  (\bw') + \frac{1}{2\rho}\|\bw'-\bw\|^2\}$. 
 \end{definition}
\begin{theorem}\label{thm:ncc}
Under Assumptions~\ref{ass:bounded Y},~\ref{ass:primal standard}, and~\ref{ass:ncc}, for  Algorithm~\ref{algorithm: CODA-Primal}, if we choose $\delta = O(1)$, $B = \Theta\pr{1}$,  $M = \Theta\pr{\frac{\sigma^2 L^4 D^2_{\cY}}{\epsilon^4}}$, $\beta = \frac{1}{2}$, $\eta_{\by} = \Theta\pr{  \frac{\epsilon^2}{L\sigma^2}  }$, $\eta_{\bx} = \Theta \pr{ \frac{\epsilon^6}{ LD_{\cY} \sigma^2 (G_h + G_f G_g )}}$,  then Algorithm~\ref{algorithm: CODA-Primal} guarantees that $ 
    \frac{1}{T}\sum_{t=1}^T \E\norm{\nabla \Phi_{1/2L}(\bx^t )}^2 \leq \epsilon^2 $ 
with gradient complexity bounded by 
\begin{align*}
    O\pare{  \max\left\{ \frac{ L^3 \sigma^2D^2_{\cY}(G_h + G_fG_g) \Delta_{\hat{\Phi}}}{\epsilon^8},  \frac{L\Delta_{\Phi}}{\epsilon^2} \right\} },
\end{align*}
where $\Delta_{\Phi}:=\Phi(\bx^0) - \min_{\bx\in\cX}\Phi(\bx)$,   $\Delta_{\hat\Phi}:=\Phi_{1/2L}(\bx^0) - \min_{\bx\in\cX}\Phi_{1/2L}(\bx)$.
\end{theorem}
\noindent The proof of Theorem~\ref{thm:ncc} is deferred to Appendix~\ref{app: proof NCC}. \textcolor{black}{Unlike the proof technique used in compositional minimization or minimax optimization~\cite{chen2021solving,gao2021convergence} that only considers potential function with second moment of inner function estimation error $\E\norm{\bz^t - g(\bx^{t-1})}^2$, in nonconvex-concave minimax setting we rely on the potential function with both first and second moment of : $ \Psi^t = \Phi_{1/2L} (\bx^t) + O(\frac{\eta_{\bx} }{ (2\beta-\beta^2)}G_g^2 L_f^2) \E\norm{\bz^t - g(\bx^{t-1})}^2 + O(\frac{ \eta_\bx L^2 D_{\by}}{\beta}) \E\norm{\bz^t - g(\bx^{t-1})}$.}
We recover $O(\epsilon^{-8})$ gradient complexity, the same as SGDA on the non-compositional nonconvex-concave objective. Notice that we also need $T\cdot M = O(\epsilon^{-12})$ zeroth order evaluations for maintaining $\bz$ variable, but zeroth order evaluation (i.e., inference) is usually considered much cheaper than gradient evaluation, and can be implemented efficiently in practice. This is also the first convergence result for nonconvex-concave compositional minimax optimization.

\section{Composition on Dual Variable }\label{sec:com:dual}
In this section, we consider the setting where the composition happens on the dual variable:  
\begin{align*}
    F(\bx,\by) := h(\bx) +  \E_{\zeta\sim\cS_f}[f(\bx, \E_{\xi\sim \cS_g}[g( \by);\xi] ;\zeta) ]  - r(\by) .
\end{align*}
In this setting  the objective is possibly nonconcave in $\by$. We hence consider strongly-convex-nonconcave and convex-nonconcave settings to establish the convergence rate. The algorithm employed is CODA-Dual, which is akin to CODA-Primal, but instead utilizes a compositional correction step on the dual variable as shown in Algorithm~\ref{algorithm: CODA-Dual}. The little difference is that, we have a $\alpha_t \bx^t$ term which serves as a regularizer. When the objective is strongly convex in $\bx$, we will set $\alpha_t$ to be $0$. When it is merely convex in $\bx$, we will set $\alpha_t$ to be a small value to make the objective strongly convex. This technique was used in vanilla convex-nonconcave minimax problem~\cite{xu2023unified}. 

\begin{algorithm2e}[t] 
	\DontPrintSemicolon
    \caption{\sffamily{CODA-Dual}}
    	\label{algorithm: CODA-Dual}
  	\textbf{Input:} $({\bx}^{0},{\by}^{0})$,   $\cbr{\alpha_t}$,  stepsizes $\eta_{\bx},\eta_{\by}$, $\bz^0$ such that $\E\norm{\bz^0 - g(\bx^0)}^2 \leq \delta$.

 { 
   
	\For{$t = 0,...,T-1$}
	 {
             Sample a mini-batch $\cM^t$ of $M$ data from $\cS_g$.\\  
            $\bz^{t+1} = (1-\beta) \pare{\bz^t + \frac{1}{M}\sum_{\xi \in \cM^t} \pr{ g( \by^t; \xi )  -  g(  \by^{t-1}; \xi )  }}  + \beta  \frac{1}{M}\sum_{\xi  \in \cM^t} g( \by^t; \xi )  $. \\
               Sample a mini-batch $\cB^t $ of $B$ data pair $(\zeta ,\xi )$ from $\cS_f \times \cS_g$,\\
            Compute $  \bg_{\bx}^t = \frac{1}{B}\sum_{(\zeta ,\xi ) \in \cB^t} \nabla_1 f( \bx^t, \bz^{t+1};\zeta  )   + \nabla h(\bx^t) + \alpha_t \bx^t$, $ {\bg}_{\by}^t =\frac{1}{B}\sum_{(\zeta ,\xi ) \in \cB^t} \nabla_2 f( \bx^t,  \bz^{t+1}; \zeta  )  \nabla g( \by^t; \xi ) - \nabla r(\by^t)$. \\
 	    ${\bx}^{t+1} = \cP_{\mathcal{X}}\pare{\bx^{t} - \eta_{\bx} \bg_{\bx}^t }$     , $\by^{t+1} =\cP_{\cY}\pare{ \by^{t} + \eta_{\by} \bg_{\by}^t} $. 
 	 } 
\textbf{Output:}  $(\hat{\bx}, \hat{\by})$ uniformly sampled from $\cbr{ ( {\bx}^t,  {\by}^t)   }_{t=1}^{T}$.  } 
\end{algorithm2e}

\noindent\textbf{\ding{110}~Strongly-convex-Nonconcave Setting.}
We first present the convergence result of CODA-Dual in the strongly-convex-nonconcave setting by making the following standard assumptions.
\begin{assumption}\label{ass:sc x}
    $F(\bx,\by)$ is $\mu$-strongly-convex in $\bx$,   $\forall \by \in \cY$ and  $L $ smooth.
\end{assumption}
\begin{assumption}\label{ass:dual standard}
    We make the following assumptions:
    \begin{enumerate}  
       \item $f(\bx,\bz)$ is $L_f$ smooth,    $g(\by)$ is $L_g$ smooth, and  $F(\bx,\by)$ is $L $ smooth.
        \item $\E[\norm{\nabla g(\by;\xi)}^2] \leq G^2_g$ and $\E[\norm{\nabla f(\bx,\bz;\zeta)}^2] \leq G^2_f$.
        \item $\max\cbr{\E\norm{  g(\by;\xi) -  g(\by ) }^2  ,  \E\norm{  \nabla g(\by;\xi) -   \nabla g(\by ) }^2 ,  \E\norm{ \nabla f(\bx,\bz;\zeta)   - \nabla f(\bx,\bz ) }^2} \leq \sigma^2$.
    \end{enumerate}
\end{assumption}
\begin{assumption}\label{ass:bounded domain}
    Domain $\cX$ and $\cY$ have bounded diameter, i.e., $\norm{\bx-\bx'} \leq D_{\cX}$ and $\norm{\by-\by'} \leq D_{\cY}$, $\forall \bx,\bx' \in \cX, \by,\by'\in\cY$. We use $D$ to denote $D = \max\cbr{D_{\cX}, D_{\cY}}$.
\end{assumption}
\begin{assumption} \label{ass:boundedness F}
    There exists a constant $F_{\max}$ such that $\max_{\bx \in \cX, \by \in \cY} F(\bx,\by) \leq F_{\max}$.
\end{assumption}
Next, we consider the following convergence measure:
 
\begin{definition}[Convergence Measure~\cite{xu2023unified}]\label{def: convergence measure} Given two parameters, $\bx$ and $\by$, we define the following quantity as a stationary gap  $ 
    \nabla G(\bx,\by)   = \begin{pmatrix} \frac{1}{\eta_{\bx}}  \pare{ \bx - \cP_{\cX}\pare{\bx-\eta_{\bx}\nabla_{\bx} F(\bx,\by)}  } \\
    \frac{1}{\eta_{\by}}  \pare{ \by - \cP_{\cY}\pare{\by+\eta_{\by} \nabla_{\by} F(\bx,\by)}  }\end{pmatrix}~. $ 
\end{definition}
\noindent Given the nonconcave nature of $F(\bx,\by)$, we are only able to show the convergence to a stationary point. Definition~\ref{def: convergence measure} measures the stationarity given parameter pair $(\bx,\by)$ by examining how much the parameter will change if we run one step projected gradient descent-ascent on them.
The widely employed convergence measure \emph{primal function}  $\norm{\nabla \Phi(\bx)}$ does not apply here due to non-concavity.
\begin{theorem}\label{thm:SCNC}
Under Assumptions~\ref{ass:sc x}-\ref{ass:boundedness F}, if we choose  $  M =\Theta \pare{ \max\cbr{ \frac{\kappa^3 L \sigma^2}{\epsilon^2}   1} }, B= \Theta \pare{ \max\cbr{   \frac{\kappa^2 L_f^2\sigma^2}{\epsilon^2}, 1} }, \beta = 0.1,  
     \eta_\bx  = \Theta\pare{\min \cbr{ \frac{1}{ L^2}, \frac{\mu}{ L^2_f}}},  \eta_{\by} = \Theta\pare{ \frac{
     \eta_\bx
     }{\kappa^2}}$, $\delta = O(1)$, and $\alpha_t = 0$, then Algorithm~\ref{algorithm: CODA-Dual} guarantees that $ 
    \frac{1}{T}\sum_{t=1}^T \E\norm{\nabla G(\bx^t, \by^t)}^2 \leq \epsilon^2 $ 
with the gradient complexity bounded by:
 \begin{align*}
 O\pr{ \max\cbr{   \frac{\kappa^2 L^2_f \sigma^2}{\epsilon^2}, 1} \cdot \frac{\kappa^3 F_{\max}}{\epsilon^2} }. 
 \end{align*} 
\end{theorem} 
\noindent The proof of Theorem~\ref{thm:SCNC} is deferred to Appendix~\ref{app: proof SCC}. The heart of our proof is constructing a two-level potential function, together with controlling the iterate difference of auxiliary variable $\norm{\bz^{t+1} - \bz^t}^2$, so that we can derive the proper descent inequality. Here we achieve $O(\epsilon^{-4})$ rate, worse than {\em deterministic} non-compositional setting~\cite{xu2023unified} by  $O(\epsilon^{-2})$. The most relevant work is~\cite{deng2023mixture}, where they also solve a strongly-convex-nonconcave compositional minimax problem. The main difference is that their objective has a linear coupling of primal and dual variables, while we consider a  more general objective. To guarantee the convergence, the ratio of primal stepsize to dual stepsize needs to be $\kappa^2$, due to the asymmetry nature of our objective.

\noindent\textbf{Reducing SC-NC to NC-SC Problem.} One may argue that finding stationary point of an strongly-convex-nonconcave (SC-NC) minimax problem can be solved by an algorithm that can find that of an nonconvex-strongly-concave (NC-SC) problem. However, in the SC-NC problem, we must consider the gradient norm of $F(\bx,\by)$ as the convergence measure, while all existing NC-SC compositional minimax algorithm only has convergence guarantee in terms of gradient norm of primal function $\Phi(\bx)$. It is well-known that (\cite{lin2020gradient}, Proposition 4.11), translating an $\epsilon$ stationary point in term of $\Phi(\bx)$ to $\epsilon$ stationary point in term of $F$ will need $O(\epsilon^{-2})$ more stochastic gradient evaluations, hence if we simply use the result from Theorem 1, we will need another sub-rountine to translate the stationary point, which is less preferred than our single loop algorithm CODA-Dual.

\noindent\textbf{\ding{110}~Convex-Nonconcave Setting}
In this section, we present the convergence result of CODA-Dual when the objective is merely convex in $\bx$.
\begin{assumption}\label{ass:convex x}
    $F(\bx,\by)$ is convex in $\bx$,  $\forall \by \in \cY$. 
\end{assumption}


\begin{theorem}\label{thm: CNC}
Under Assumptions~\ref{ass:dual standard}-\ref{ass:convex x}, if we choose $M = \Theta \pare{ \max\cbr{   \frac{ L^6 L_f^2 D_{\cX}^4 \sigma^2}{\epsilon^4}, 1} }  $, $ B = \Theta \pare{ \max\cbr{   \frac{ L^3 L_f^2 D_{\cX}^2 \sigma^2}{\epsilon^4}, 1} }$,  $\delta = O(1)$,  $\beta = 0.1$,  
    $ \eta_{\bx} = \Theta\pare{\frac{1}{L^2}}$,  $\eta_{\by} = \Theta\pare{\frac{\epsilon^2}{L^4 D_{\cX}}}$, and $\alpha_t = \Theta\pare{\frac{\epsilon}{L D_{\cX}}}$, then for Algorithm~\ref{algorithm: CODA-Dual} it is guaranteed that $ 
    \frac{1}{T}\sum_{t=1}^T \E\norm{\nabla G(\bx^t, \by^t)}^2 \leq \epsilon^2 $ holds 
with the gradient complexity of 
 \begin{align*}
     O\pr{  \max\cbr{   \frac{ L^3 L_f^2 D_{\cX}^2 \sigma^2}{\epsilon^4}, 1} \cdot  \frac{F_{\max} D_{\cX}^2 L^4}{\epsilon^4} }.
 \end{align*}

\end{theorem} 
\noindent The proof of Theorem~\ref{thm: CNC} is deferred to Appendix~\ref{app: proof CNC}. 
 Here we achieve an $O(\epsilon^{-8})$ convergence rate to the stationary point of $F$ in the convex-nonconcave setting. This is also the first convergence rate for stochastic convex-nonconcave minimax optimization.
\section{Composition on Primal and Dual  Variables}\label{sec:com:both}
In this section, we consider composition on both variables
\begin{align*} 
    F(\bx,\by) := h(\bx) +  \E_{\zeta\sim \cS_f}[f( \E_{\xi\sim \cS_g}[g( \bx,\by;\xi)];\zeta )] - r(\by).  
\end{align*}
In this setting, we may lose convexity and concavity on primal and dual variables, so the problem becomes nonconvex-nonconcave minimax optimization. In general, if we make no assumptions (except smoothness), the first-order stationarity may not be a viable quantity to control convergence, and a more sophisticated optimality measure  is necessary, e.g., greedy adversarial equilibrium~\cite{mangoubi2021greedy}. In this paper, following a branch of nonconvex-nonconcave optimization literature~\cite{liu2021first,diakonikolas2021efficient}, we make a mild assumption on our nonconvex-nonconcave objective $F(\bx,\by)$, that there exists a point that is the solution of Minty Variational Inequality induced by $F(\bx,\by)$. 

\begin{assumption}[Existence of Solution for Minty Variational Inequality]\label{ass:MVI}
    There exists a $(\bx^*,\by^*) \in \cX \times \cY $, such that the following Minty variational inequality induced by $\nabla F$ can hold:
    \begin{align*}
        \inprod{\nabla F(\bx,\by)}{ (\bx,\by) - (\bx^*,\by^*) } \geq 0, \forall (\bx,\by) \in \cX \times \cY.
    \end{align*}
\end{assumption}

\noindent\textbf{\ding{110}~Weakly-convex-weakly-concave Setting.}
Relying on reduction introduced in~\cite{liu2021first}, we can break solving a weakly-convex-weakly-concave minimax problem into solving a bunch of strongly monotone variational inequality problems. Hence we first derive the following convergence of CODA-SCSC on solving Strongly Monotone VI.
\begin{lemma}\label{lem:conv PD SMV}  Assume that $F_{k}(\bx,\by)$ is $\mu$ strongly convex in $\bx$ and $\mu$ strongly concave in $\by$, and $ L$ smooth, and $\max\cbr{D_\cX,D_{\cY}} \leq D$, then if we run Algorithm~\ref{algorithm: CODA-SCSC} on $F_{k}(\bx,\by)$ with $\eta_{\bx} = \eta_{\by} = \eta$, the solution $\hat{\bw} = (\hat{\bx},\hat{\by})$ returned by algorithm guarantees that $
      \max_{\bw \in \cW}   \E [ \nabla F_{k}(\hat\bw)^\top (\hat\bw - \bw) ] \leq \epsilon$
    with gradient complexity $ O\pare{ \pare{D+\frac{1}{\eta}}\frac{L^2\sigma^2}{\mu^2 \eta\epsilon} \ln \pare{\frac{D+1/\eta}{\epsilon}}  } .$
 
\end{lemma}
\begin{algorithm2e}[t]
    \caption{\sffamily{CODA-SCSC}($F,\bx^0,\by^0, T$)}
    	\label{algorithm: CODA-SCSC}
    	 
 \textbf{Input:} $({\bx}^{0},{\by}^{0})$,   $\bz^0$ such that $\E\norm{\bz^0 - g(\bx^0)}^2 \leq \delta$.

   {
    {

	\For{$t = 0,...,T-1$}
	 {
            Sample a mini-batch $\cM^t$ of $M$ data from $\cS_g$.\\ 
            $\bz^{t+1} = (1-\beta ) \pare{\bz^t +\frac{1}{M}\sum_{\xi \in \cM^t} g(\bx^t, \by^t; \xi)  -  g( \bx^{t-1},  \by^{t-1}; \xi)  }  + \beta  \frac{1}{M}\sum_{\xi \in \cM^t} g(\bx^t, \by^t; \xi)  $.\\
            Sample a mini-batch $\cB^t $ of $B$ data pair $(\zeta^t,\xi^t)$ from $\cS_f \times \cS_g$.\\
            Compute $\bg_{\bx}^t = \frac{1}{B}\sum_{(\zeta,\xi) \in \cB^t}\nabla f(\bz^{t+1}; \zeta)\nabla_{\bx} g( \bx^t, \by^t; \xi)   + \nabla h(\bx^t) + \frac{2}{\gamma} \pare{\bx^t-\bx^0} $,  $\bg_{\by}^t = \frac{1}{B}\sum_{(\zeta,\xi) \in \cB^t} \nabla f(  \bz^{t+1};\zeta )  \nabla_{\by} g( \bx^t, \by^t; \xi) - \nabla r(\by^t) - \frac{2}{\gamma} \pare{\by^t-\by^0}   $.\\
 	    ${\bx}^{t+1} = \cP_{\mathcal{X}}\pare{\bx^{t} - \eta_{\bx} \bg_{\bx}^t }$     , $\by^{t+1} =\cP_{\cY}\pare{ \by^{t} + \eta_{\by} \bg_{\by}^t} $.\\ 
 	 }
      
 	      }
 	     } 
       \textbf{Output:} $(\bx^T,\by^T)$.
\end{algorithm2e}
\begin{algorithm2e}[t] 
	\DontPrintSemicolon
    \caption{\sffamily{CODA-PD}}
    	\label{algorithm: CODA-PD}
    	\textbf{Input:}   $({\bx}_{0},{\by}_{0})$,  $\gamma$, stepsizes $\eta_{\bx},\eta_{\by}$, weights $\cbr{\theta^k}_{k=0}^{K-1}$. 
	\\

   {
    {
   
	\For{$k = 0,\ldots,K-1$}
	 {
           Construct $F_k = F(\bx,\by) + \frac{1}{2\gamma} \pare{\norm{\bx - \bx_k}^2 - \norm{\by - \by_k}^2  }$\\
           $(\bx_{k+1},\by_{k+1}) = \texttt{CODA-SCSC}(F_k, \bx_k, \by_k,  T_k )$
 	 }
      
 	      }
 	\textbf{Output:}  $(\hat{\bx}, \hat{\by})$ sampled from $\cbr{ ( {\bx}_k,  {\by}_k)   }_{k=0}^{K-1}$ with $ \Pr\pare{ (\hat{\bx}, \hat{\by}) = ( {\bx}_k,  {\by}_k)  }  = \frac{\theta^k}{\sum_{k=0}^{K-1} \theta^k}$.}  
\end{algorithm2e}
\begin{assumption}\label{ass:wcwc}
    We make the following assumptions: \vspace{-1mm}
    \begin{enumerate} 
        \item $F(\bx,\by)$ is is $L$ smooth, $\rho$-weakly-convex in $\bx$, and $\rho$-weakly-concave in $\by$. 
        \item $f(\bz)$ is $L_f$ smooth and $g(\bx,\by)$ is $L_g$ smooth. \vspace{-1mm}
        \item $\E[\norm{\nabla g(\bx,\by;\xi)}^2] \leq G^2_g$ and $\E[\norm{\nabla f(\bz ;\zeta)}^2] \leq G^2_f$.\vspace{-1mm}
        \item $\max\cbr{\E\norm{  g(\bx,\by;\xi) -  g(\bx,\by ) }^2 ,   \E\norm{  \nabla g(\bx,\by;\xi) -  \nabla g(\bx,\by ) }^2 ,  \E\norm{ \nabla f(\bz ;\zeta)   - \nabla f(\bz ) }^2} \leq \sigma^2 $.
    \end{enumerate}
\end{assumption}

\begin{definition}[Convergence Measure~\cite{liu2021first}]\label{def: convergence measure WCWC} A solution $\bw = (\bx,\by)$ is $\epsilon$ stationary point of problem~(\ref{eq:obj}) if there exists $\bar\bw = (\bar{\bx},\bar{\by})$ such that $
 \norm{\bw - \bar\bw} \leq c \epsilon,  dist\pare{0, \partial (F(\bar\bx,\bar\by) + \mathbf{1}_{\cX\times \cY}(\bar\bx,\bar\by) )} \leq \epsilon$ for some $c >0$.
\end{definition}

\begin{theorem} \label{thm:wcwc}
Under Assumptions~\ref{ass:bounded domain}, ~\ref{ass:MVI} and~\ref{ass:wcwc},  if we choose $B = \Theta \pare{ \max\cbr{   \frac{ \sigma^2}{\epsilon^2}, 1} },   \beta = 1-   \frac{\mu}{ 16 G^2_g  L_f} ,  
     \eta_\bx  =  \eta_{\by} = \Theta\pare{ \frac{
    1
     }{L^2}}$, $\delta = O(1)$, $\theta^k = (k+1)^\alpha$ for some $\alpha \in (0,1)$, $\gamma = \frac{1}{\rho}$, then   Algorithm~\ref{algorithm: CODA-Dual} returns $\epsilon$ stationary solution $(\hat\bx,\hat\by)$  with the gradient complexity bounded by: 
 \begin{align*}
      O\pare{ \pare{D+L^2}\frac{  \rho^2 D  L^4\sigma^2}{  \epsilon^4} \log \pare{\frac{1}{\epsilon}}   }.
 \end{align*} 
\end{theorem}
\noindent The proof is deferred to Appendix~\ref{app: proof PD}. Here we achieved an $\tilde O(\epsilon^{-4})$ rate. The most relevant work is~\cite{liu2021first}, demonstrating an  $\tilde O(\epsilon^{-2})$ rate using (deterministic) gradient descent ascent as inner problem solver.

\section{Faster Rates with Variance Reduction}\label{sec:VR}

In this section, our goal is to explore the possibility of achieving faster rates. If we assume individual or point-wise smoothness of $f$ and $g$, i.e., smoothness of $f(g(\bx),\by;\zeta)$ and $g(\bx;\xi)$, then the natural idea would be employing variance reduction minimax algorithms, such as VR-SAPD~\cite{zhang2022sapd+}, combined with our compositional correction steps. Unfortunately, this integration solely does not lead to  an improved rate since the main bottleneck comes from the estimation error of inner function $g(\bx)$, i.e., $\E\norm{\bz^{t+1} - g(\bx^t)}^2$. This observation raises the following question: can we accelerate the convergence of estimation error of inner function to entail a faster convergence rate overall?  Inspired by Spider~\cite{fang2018spider}, we propose to use the following  \textbf{variance reduced inner function estimation} when maintaining $\bz$ variable:
\begin{align*}
    \bz^{t+1} = (1-\beta)(\bz^t + g(\bx^t; \cI^t_{\bz}) - g(\bx^{t-1}; \cI^t_{\bz})) + \beta g^{t-1},
\end{align*}
where $g^{t} = g^{t-1} + g(\bx^t;\cI^t_{\bz}) - g(\bx^{t-1};\cI^t_{\bz})$ and every $\tau$ iterations, we use a large mini-batch to update $g^{t}$: $g^{t} =  g(\bx^t;\cB^t_{\bz})$. 

It can be shown, the gap $\E\norm{\bz^{t+1} - g(\bx^t)}^2$ converges faster than standard compositional optimization as in~\cite{chen2021solving}. Equipped with this technique, we propose our CODA+ algorithm  depicted in Algorithm~\ref{algorithm: CODA-Primal+}, Appendix~\ref{app:VR algorithm}.

\noindent\textbf{Nonconvex-strongly-concave Setting.} The following theorem establishes the convergence:

\begin{assumption}\label{ass:primal vr}
    We make the following assumptions:
    \begin{enumerate}  
        \item  $f(\bz,\by;\zeta )$ is $L_f$ smooth,  $g(\bx;\xi)$ is $L_g$ smooth for any $\zeta \in \cS_f, \xi \in \cS_g$ , and  $F(\bx,\by)$ is $L $ smooth.\vspace{-1mm}
        \item $\E[\norm{\nabla g(\bx;\xi)}^2] \leq G_g^2$ and $\E[\norm{\nabla f(\bz,\by;\zeta)}^2] \leq G_f^2$.\vspace{-1mm}
        \item $ \max\cbr{\E\norm{  g(\bx;\xi) -    g(\bx ) }^2, \E\norm{  \nabla g(\bx;\xi) -   \nabla g(\bx ) }^2,\E\norm{ \nabla f(\bz,\by;\zeta)   - \nabla f(\bz,\by ) }^2  }\leq \sigma^2$.\vspace{-1mm}
    \end{enumerate}
\end{assumption}

\begin{theorem}\label{thm:NCSC+}
Under Assumptions~\ref{ass:ncsc},~\ref{ass:bounded Y},~\ref{ass:primal vr}, for  Algorithm~\ref{algorithm: CODA-Primal+}, if we choose  $\delta = O(\frac{\epsilon^2}{L})$,   $B_{\tau} = \Theta\pare{\frac{L^2 \sigma^2}{\mu  \epsilon^2}}$, $\tau = \sqrt{\frac{B_{\tau}}{\kappa}}$, $B = \sqrt{ \kappa B_\tau}$,  and $\eta_{\bx} = \Theta\pare{\frac{1}{L}}$, $\eta_{\by} = \Theta\pare{\frac{1}{L}}$, $\mu_{\bx} = 2L$, then it guarantees that $   \frac{1}{T}\sum_{t=1}^T \E\norm{\nabla \Phi_{1/2L}(\bx^t )}^2 \leq \epsilon^2 $ 
with the gradient complexity bounded by: 
\begin{align*}
    O\pare{ \frac{\kappa^2  {L}^{1.5} \Delta_\Phi \sigma}{\epsilon^3}       }. 
\end{align*}

\end{theorem} 
\noindent The proof of Theorem~\ref{thm:NCSC+} is deferred to Appendix~\ref{app: proof ncsc+}. Here we achieve $O\pare{ \frac{\kappa^2  {L}^{1.5}  }{\epsilon^3}  }$ rate for finding the stationary point of Moreau envelope of primal function. Notice that once we find $\epsilon$ stationary point of Moreau envelope, there will be an efficient algorithm that takes only $\tilde O(\frac{1}{\epsilon})$ stochastic gradient complexity to find $\epsilon$ stationary point of primal function, as indicated in~\cite{zhang2022sapd+}. 
This rate is faster than previous SOTA~\cite{liu2024faster}, by $O(\kappa)$ factor. Compared to non-compositional minimax optimization, the best-known rate is $O\pare{ \frac{\kappa^2  {L}   }{\epsilon^3}  }$, so our result almost recovers the current non-compositional minimax SOTA rate. As a last final, we need to choose a good initialization of $\bz^0$ such that   $\E\norm{ \bz^0 - g(\bx^0) }^2 \leq O({\epsilon^2}/{L})$. This will take $O({L\sigma^2}/{\epsilon^2})$ samples drawn from $\cS_g$.

\noindent\textbf{\ding{110}~Nonconvex-concave Setting.}~In the nonconvex-concave setting, we consider CODA-Primal+  on the following augmented strongly concave objective with strongly concavity parameter $\mu = \frac{\epsilon^2}{LD^2_{\cY}}$:\vspace{-2mm}
\begin{align}
    \tilde F(\bx,\by) := F(\bx,\by) - \frac{\epsilon^2}{LD^2_{\cY}} \norm{\by - \by_0}^2 \label{eq:augmented ncc}.
\end{align}\vspace{-2mm}
\begin{theorem}\label{thm:ncc+}
Under Assumptions~\ref{ass:bounded Y},~\ref{ass:ncc} and~\ref{ass:primal vr}, if we run  Algorithm~\ref{algorithm: CODA-Primal+} on $ \tilde F(\bx,\by)$ defined in (\ref{eq:augmented ncc}), and use the same parameter choice as in Theorem~\ref{thm:NCSC+} with $\mu = \frac{\epsilon^2}{LD^2_{\cY}}$, then it guarantees that $ \frac{1}{T}\sum_{t=1}^T \E\norm{\nabla \Phi_{1/2L}(\bx^t )}^2 \leq \epsilon^2 $ 
with the gradient complexity bounded by  $
    O\pare{ \frac{ {L}^{5.5} \Delta_\Phi \sigma}{\epsilon^7} }.$
\end{theorem} 
\noindent The proof of Theorem~\ref{thm:ncc+} is deferred to Appendix~\ref{app: proof ncc+}. 
This result improves the original $O(\frac{1}{\epsilon^{8}})$ rate in Theorem~\ref{thm:ncc}. Faster rate is possible, if we can develop an algorithm which converges at the rate of $O(\frac{\kappa}{\epsilon^4})$ in the strongly concave setting, and choose $\mu = \frac{\epsilon^2}{LD^2_{\cY}}$ which will yield $O(\frac{1}{\epsilon^6})$ rate in the merely concave setting. Such an algorithm can be possibly achieved by modifying SAPD~\cite{zhang2022sapd+}, combined with our variance-reduced inner function estimation, which we leave as future work.

\begin{figure*}[t]
    \centering
    \includegraphics[width=0.8\textwidth]{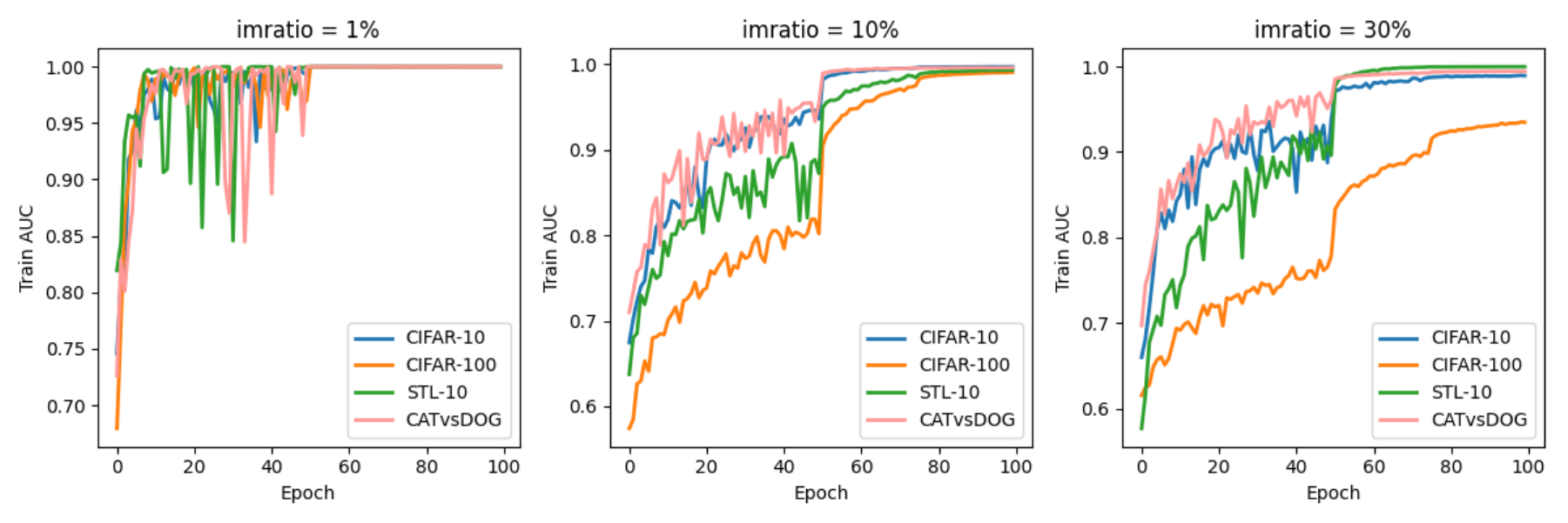}\vspace{-0.4cm}
    \caption{Convergence curves of CODA-Primal on AUC maximization task on four benchmark datasets with different imbalance ratios}
    \label{fig:Performance_imbalanced_ratios_CODA}
\end{figure*}

\begin{table}[t!]
    \centering
    \begin{tabular}{l c c c c}
    \hline
        \multirow{2}{*}{} & \multicolumn{2}{c}{\textbf{Meta-training Alphabets}} & \multicolumn{2}{c}{\textbf{Meta-testing Alphabets}}\\
        \cline{2-3}\cline{4-5}
        $(N,K)$ & Mean(Std. Dev.)& Worst(Std. Dev.) & Mean(Std. Dev.) & Worst(Std. Dev.)\\
        \hline
        (5,1) & 98.56(0.001) & 98.41(0.012) & 92.57(0.006) & 92.08(0.038)\\
        \hline
        (20,1) & 95.66(0.007) &  94.64(0.011) & 81.49(0.072)&75.28(0.086) \\
        \hline
    \end{tabular}
    \caption{Omniglot $N$-way, $K$-shot classification accuracies over three random runs}
    \label{tab:TR_MAML_Train_valid_test_5W1S_20W1S} 
\end{table}


\section{Experiment}\label{sec:exp}
In this section, we conduct experiemnts to  verify the effectiveness of our proposed algorithms. In primal composition setting, we run CODA-Primal on AUC optimization and robust MAML tasks, where we are tasked with solving nonconvex-strongly-concave and nonconvex-concave compositional minimax problems, respectively. In dual composition setting, we run CODA-Dual on mixture weights estimation problem, where we need to solve strongly-convex-nonconcave compositional minimax problem. The  additional results are provided in~\ref{app: exp}.

\begin{figure*}[t]
    \centering\vspace{-0.5cm}
    \includegraphics[width=0.9\textwidth]{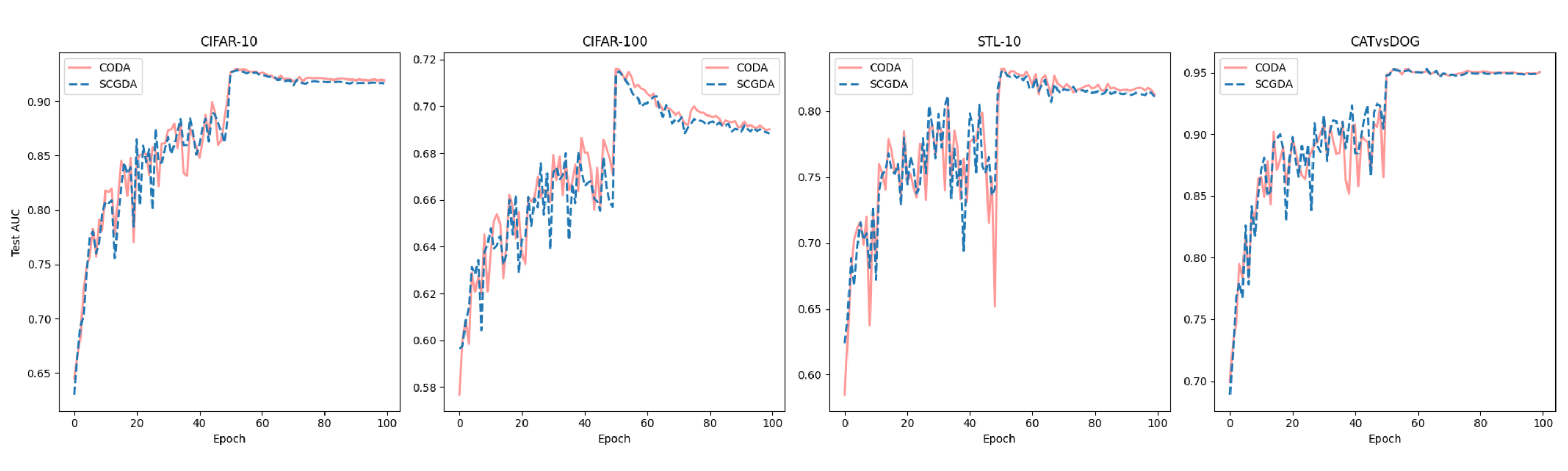}\vspace{-0.4cm}
    \caption{Testing AUC comparison of CODA and SCGDA on four benchmarks under imratio$=10\%$}
    \label{fig:coda_scgda_test_auc_0_1}
\end{figure*}

\paragraph{Deep AUC Maximization.} To evaluate our CODA-Primal algorithm on nonconvex-strongly-concave objective, we consider deep AUC maximization task~\cite{yuan2021compositional}. We choose four benchmark image classification datasets: CIFAR-10, CIFAR-100, CATvsDog, and STL-10. We generate imbalanced binary versions of these datasets with imbalanced ratios of $1\%, 10\%$, and $30\%$, following the positive examples' proportion to total training samples, following the work~\cite{yuan2021compositional}. Employing ResNet20 as our prediction network and setting the weight decay to 1e-4, we optimize for AUC over 100 epochs with a batch size of 128. For all benchmark datasets, we conduct experiments under three different random seeds, calculating mean values and standard deviations. Figure~\ref{fig:Performance_imbalanced_ratios_CODA} illustrates the convergence and AUC scores of CODA-Primal across these four datasets at different imbalance ratios. We also compare CODA with SCGDA~\cite{gao2021convergence} on testing AUC performance, and the results are shown in Figure \ref{fig:coda_scgda_test_auc_0_1}.

\paragraph{Task-agnostic Robust MAML.} In this experiment we apply  CODA-Primal algorithm to a nonconvex-concave objective. by considering task-agnostic robust MAML~\cite{collins2020task}. We use the Omniglot dataset~\cite{lake2015human}, 1623 handwritten characters from 50 alphabets, to learn a task-robust meta-model to enhance performance and adaptability across diverse task distributions. We define tasks as $N$-way, $K$-shot classification problems, each involving characters from a single alphabet, and use a 4-layer CNN~\cite{finn2017model}. Our approach involves meta-training on 25 alphabets and meta-testing on 20 different alphabets, following the same data splits in meta-learning evaluations~\cite{triantafillou2019meta}. After 60,000 meta-training iterations, we evaluate our model on 5,000 tasks from meta-test alphabets and 5,000 tasks from meta-training alphabets. The detailed accuracies, particularly when $N=5$ and $N=20$, are shown in Table~\ref{tab:TR_MAML_Train_valid_test_5W1S_20W1S},  demonstrating the effectiveness of  CODA  in nonconvex-concave compositional minimax problem.
\begin{table}[t]
    \centering
    \begin{tabular}{l c c c c}
    \hline
    \textbf{Algorithm} & \textbf{Group 1} & \textbf{Group 2} & \textbf{Group 3} & \textbf{Group 4} \\
    \hline
        Target-only &  86.97(0.01) & 69.43(0.02) & 90.20(0.01) & 65.30(0.02)\\
        Average(equal)-weight & 33.63(0.04) & 23.53(0.03) & 27.40(0.05) & 29.63(0.02)\\
        CODA-Dual & \textbf{92.90}(0.02) & \textbf{78.03}(0.02) & \textbf{96.17}(0.005) & \textbf{68.50}(0.02)\\
    \hline
    \end{tabular}
    \caption{Test accuracy of algorithms across various target domains, highlighting items where the algorithm outperforms others. Standard deviations are derived from ten random runs.}
    \label{tab:Test_Acc_for_diff_target_domains}
\end{table}

\paragraph{Mixture Weights Estimation in Multi-source Domain Adaptation.} In this experiment we test our CODA-Dual algorithm on strongly-convex-nonconcave objective. We consider mixture weight estimation task~\cite{konstantinov2019robust,deng2023mixture}. Given multiple source domains and a target domains, our goal is to find a good mixture of sources such that learning on it can yield a good target model. We use the CIFAR-10 dataset, divided into 4 different groups. For Groups 1, 2, 3, each group contains 5 domains drawing data from distinct classes, e.g., domains in Group 1 only draw data from first 3 classes of CIFAR10, and Group 4 contains 5 domains drawing data from classes in Groups 1 and 2 equally. In each test, we treat domains in one group as target and rest as sources. Hence the training involves 15 source domains (5 in each of three groups), with 1500 samples per domain. Following~\cite{mansour2021theory}, our evaluation includes two baselines: (i) target-only: training only on limited target domain data, (ii) average(equal)-weight: equal weight training across all source domains. Table ~\ref{tab:Test_Acc_for_diff_target_domains} shows CODA-Dual outperforms in test accuracy, effectively identifying source domains with mixture weight $\balpha$. Results are averaged over ten random runs.
\section{Conclusion}
In this paper, we study the compositional minimax problem, with three different settings: composition on primal variable, on dual variable, and on both variables. We give simple yet effective provable algorithms for three settings. Additionally,  we  give a variance reduction algorithm which achieves SOTA convergence rates. Extensive experimental results across different applications corroborate our theoretical findings and effectiveness of the proposed methods.

 \section*{Acknowledgement}

 This work was partially supported by NSF CAREER Award \#2239374 and NSF CNS   Award \#1956276. 

\clearpage
\appendix
 \clearpage
 \appendix
 \onecolumn
\paragraph{Organization} The appendix is organized as follows.  In Appendix~\ref{app: exp} we will provide more experimental results and setup details.  In Appendix~\ref{app:VR algorithm} we provide the detailed  algorithmic flow of the proposed variance reduction algorithm. In Appendix~\ref{app: proof primal} we provide proofs for CODA-Primal. In Appendix~\ref{app: proof dual} we provide proofs for CODA-Dual. In Appendix~\ref{app: proof PD} we provide proofs for CODA-PD. In Appendix~\ref{app: proof VR} we provide proofs for CODA-Primal+. 
\section{Experiments Details}\label{app: exp}

\begin{figure*}[b]
    \centering
    \includegraphics[width=0.8\textwidth]{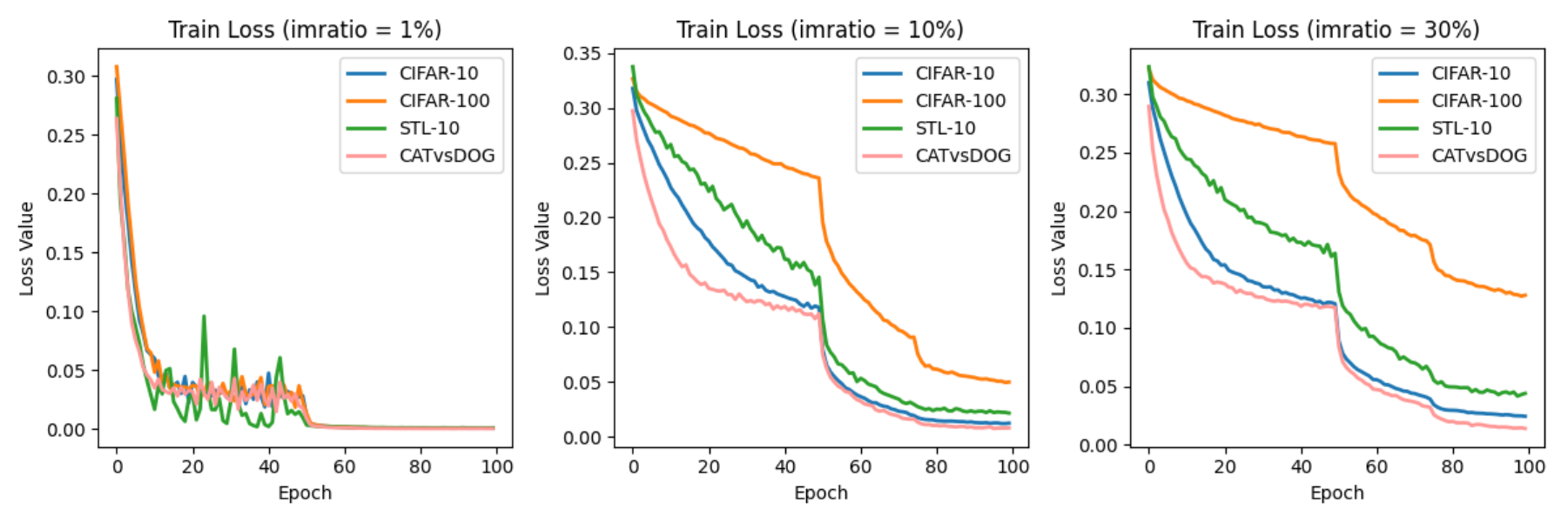}
    \caption{Training loss with different imbalance ratios on four benchmark datasets}
    \label{fig:train_loss_imbalanced_ratios_CODA}
\end{figure*}
\begin{figure} 
    \centering
    \begin{subfigure}{ \textwidth}
        \centering
        \includegraphics[width= 0.9\textwidth]{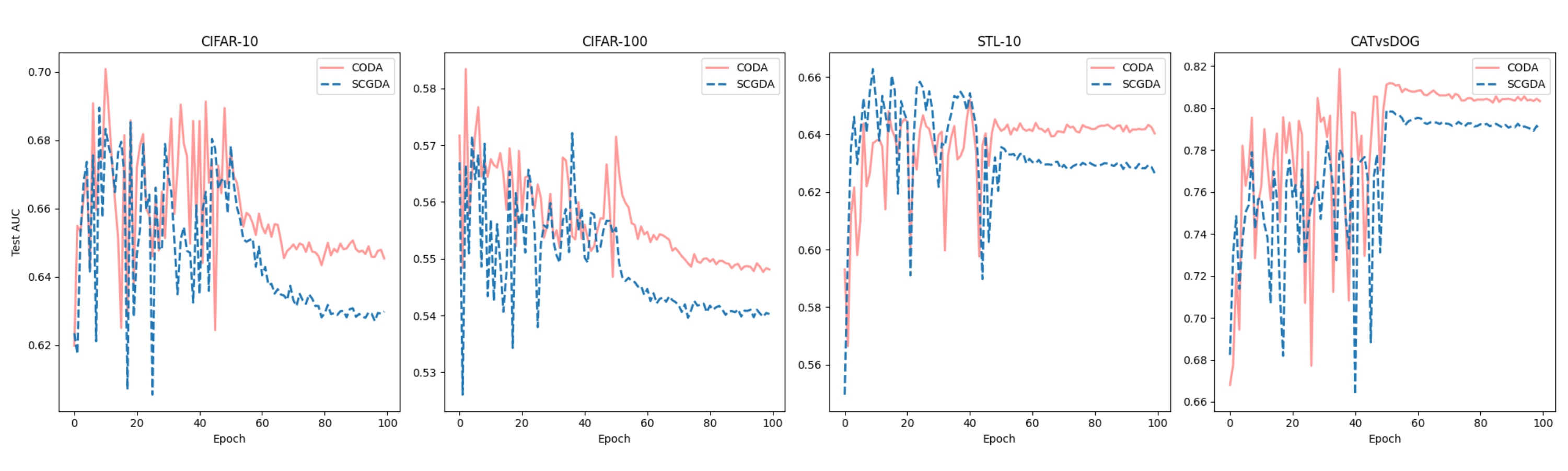}
        \label{fig:coda_scgda_test_auc_0_01}
    \end{subfigure} 
    \begin{subfigure}{ \textwidth}
        \centering
        \includegraphics[width= 0.9\textwidth]{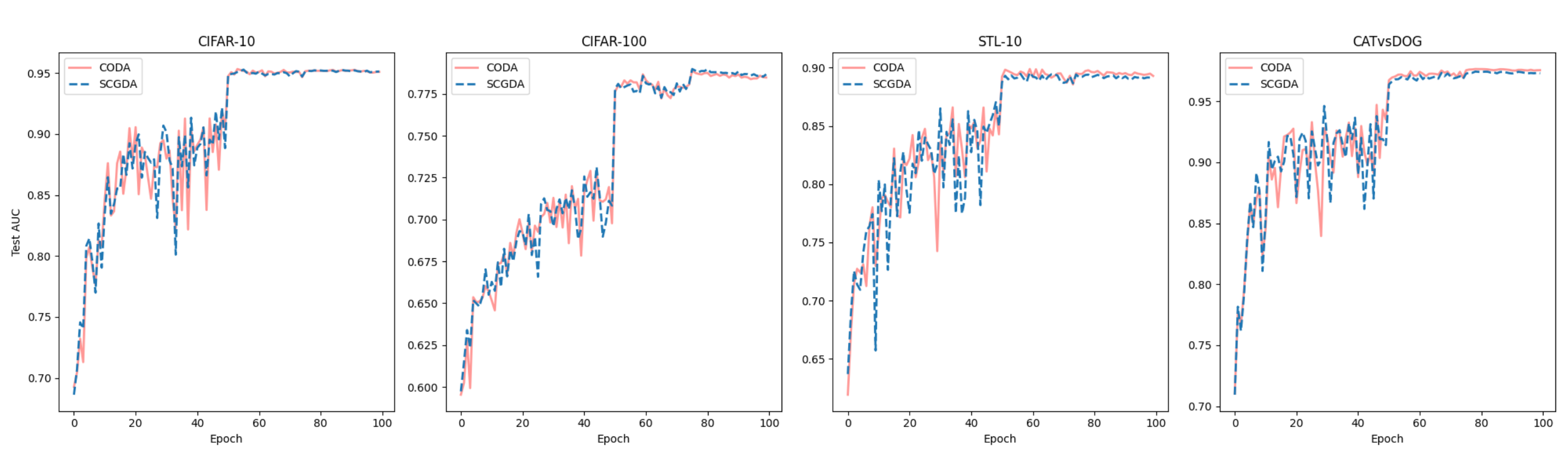}
        \label{fig:coda_scgda_test_auc_0_3}
    \end{subfigure}
    \caption{Testing AUC performance comparison of CODA and SCGDA on four benchmarks. The first row is under imratio$=1\%$ and the third row is under imratio$=30\%$}
    \label{fig:coda_scgda_compare}
\end{figure}
\paragraph{Deep AUC Maximization: }
In a binary classification task where the label $y \in \{-1,+1\}$ and the model parameter $\bw \in \R^d$, the model prediction on sample $\bx$ is $f(\bw; \bx_i)$. The objective is formulated as a minimax problem:
\begin{align*}
    \min_{\bw, a, b} \max_{\theta \in \Omega} \Phi(\bw - \alpha \nabla L(\bw), a,b,\theta) = \frac{1}{n} \sum_{i=1}^n \phi(\bw - \alpha \nabla L(\bw), a,b,\theta;\bx_i,y_i)
\end{align*}
where $L(\bw) = \frac{1}{n}\sum_{i=1}^n \ell(\bw;\bx_i,y_i)$ and
\begin{align*}
    \phi(\bw, a,b,\theta;\bx_i,y_i) = &(1-p) (f(\bw; \bx_i) -a)^2 \mathbb{I}[y_i=1] + p(f(\bw; \bx_i) - b)^2\mathbb{I}[y_i=-1]\\
    &+ 2\theta \left(p f(\bw;\bx_i)\mathbb{I}[y_i=-1] - (1-p)f(\bw; \bx_i)\mathbb{I}[y_i=1]\right)\\
    &+ 2 pf(\bw; \bx_i)\mathbb{I}[y_i=1] - 2(1-p)f(\bw; \bx_i)\mathbb{I}[y_i=-1] - p(1-p)\theta^2
\end{align*}
At time t, sample a minibatch $(\bx_1,y_1),...,(\bx_B,y_B)$ from the training set.\\
First, update the auxiliary variable:
\begin{align*}
    \bz^{t+1} = &(1-\beta^t) \pare{\bz^t + \pare{\bw^t- \alpha \frac{1}{B}\sum_{i=1}^B\nabla \ell(\w^t;\bx_i,y_i)} - \pare{\bw^{t-1}- \alpha  \frac{1}{B}\sum_{i=1}^B\nabla \ell(\w^{t-1};\bx_i,y_i) } } \\
    + &\beta^t\pare{\bw^t-  \alpha \frac{1}{B}\sum_{i=1}^B\nabla \ell(\w^t;\bx_i,y_i)}
\end{align*}
Then, sample another minibatch $(\tilde\bx_1,\tilde y_1),...,(\tilde\bx_B,\tilde y_B)$ and compute: 
   \begin{align*}
   \bg_{\bw}^t= \frac{1}{B}\sum_{i=1}^B \nabla \phi(\bz^{t+1},a^t,b^t,\theta^t;\tilde\bx_i, \tilde y_i),  & \
      g_a^t  = \frac{1}{B}\sum_{i=1}^B \nabla_a \phi(\bz^{t+1},a^t,b^t,\theta^t;\tilde\bx_i, \tilde y_i),\\
       g_b^t  = \frac{1}{B}\sum_{i=1}^B \nabla_b \phi(\bz^{t+1},a^t,b^t,\theta^t;\tilde\bx_i, \tilde y_i),  & \
        g_\theta^t  = \frac{1}{B}\sum_{i=1}^B \nabla_\theta \phi(\bz^{t+1},a^t,b^t,\theta^t;\tilde\bx_i, \tilde y_i)
\end{align*}
 Finally, update the model and $\balpha$:
 \begin{align*}
 \bw^{t+1} =  \bw^{t} - \eta_w \bg_{w}^t, \ a^{t+1}  = a^t - \eta_a g_a^t, \ b^{t+1}= b^t - \eta_b g_b^t, \  {\theta}^{t+1}  = \cP_{ {\Omega}}\pare{\theta^{t} + \eta_\theta \bg_{\theta}^t }   
 \end{align*}
 
 \noindent Table~\ref{tab:classification_task_data_description} is the dataset description in our experimental setting. 
For the benchmark datasets, the training and validation splits are set as follows: 19k/1k for CatvsDog, 45k/5k for CIFAR10, 45k/5k for CIFAR100, and 4k/1k for STL10. The learning rate is initially 0.1, with a reduction by a factor of 10 at $50\%$ and $75\%$ of total training time. Figure~\ref{fig:train_loss_imbalanced_ratios_CODA} shows the training loss under different imbalanced ratios, which illustrates the convergence of our CODA-Primal algorithm across the four datasets. The testing AUC scores and accuracy results are detailed in Table~\ref{tab:Performance_imbalanced_ratio_CODA}, where we compare our CODA-Primal algorithm with the PDSCA algorithm in the work~\cite{yuan2021compositional}. We also compare CODA with SCGDA ~\cite{gao2021convergence} on testing AUC performance, and the results are shown in Figure \ref{fig:coda_scgda_compare}.
 \begin{table*}[t]
\centering
\begin{tabular}{c c c c}
\hline
\textbf{Dataset} & \textbf{Number of samples} & \textbf{Number of classes} & \textbf{Imbalanced Ratio} \\
\hline
CIFAR-10 & 50,000 & 2 (binary) & 1\%,10\%,30\%\\
\hline
CIFAR-100 & 50,000 & 2 (binary) & 1\%,10\%,30\%\\
\hline
CATvsDOG & 20,000 & 2 (binary) & 1\%,10\%,30\%\\
\hline
STL10 & 5,000 & 2 (binary) &  1\%,10\%,30\%\\
\hline
\end{tabular}
\caption{Dataset description for classification tasks}
\label{tab:classification_task_data_description}
\end{table*}
\begin{table*}[t]
    \centering
    \footnotesize
    \setlength{\tabcolsep}{4pt}
    \renewcommand{\arraystretch}{1.2}
    \begin{tabular}{c | l c c c c c c}
    \hline
        \multirow{2}{*}{Datasets} & \multicolumn{1}{c}{} & \multicolumn{3}{c}{AUC score} & \multicolumn{3}{c}{Accuracy} \\
        \cline{3-8}
        & \textbf{imratio} & \textbf{1\%} & \textbf{10\%} & \textbf{30\%} & \textbf{1\%} & \textbf{10\%} & \textbf{30\%} \\
        \hline
       \multirow{2}{*}{CIFAR-10} & CODA-Primal & \textbf{63.6(0.003)} & 91.5(0.005) & \textbf{95.3(0.001)} & \textbf{50.7(0.002)} & 79.7(0.006) & \textbf{87.7(0.003)} \\
       \cline{2-8}
        & PDSCA & 62.7(0.009) & \textbf{92.3(0.016)} & 95.2(0.0003) & 50.6(0.0006) & \textbf{81.2(0.033)} & 86.9(0.0003) \\
        \hline
        \multirow{2}{*}{CIFAR-100} & CODA-Primal & \textbf{55.1(0.001)} & \textbf{68.1(0.001)} & \textbf{78.6(0.003)} & \textbf{50.1(0.0002)} & \textbf{57.6(0.003)} & \textbf{66.0(0.006)} \\
        \cline{2-8}
       & PDSCA & 54.3(0.005) & 68.5(0.001) & 78.6(0.002) & 50.1(0.0001) & 57.7(0.003) & 65.6(0.002) \\
       \hline
    \end{tabular}
    \caption{Testing Performance on datasets over three random runs}
    \label{tab:Performance_imbalanced_ratio_CODA}
\end{table*}
\begin{figure*}[b]
    \centering
    \includegraphics[width=0.8\textwidth]{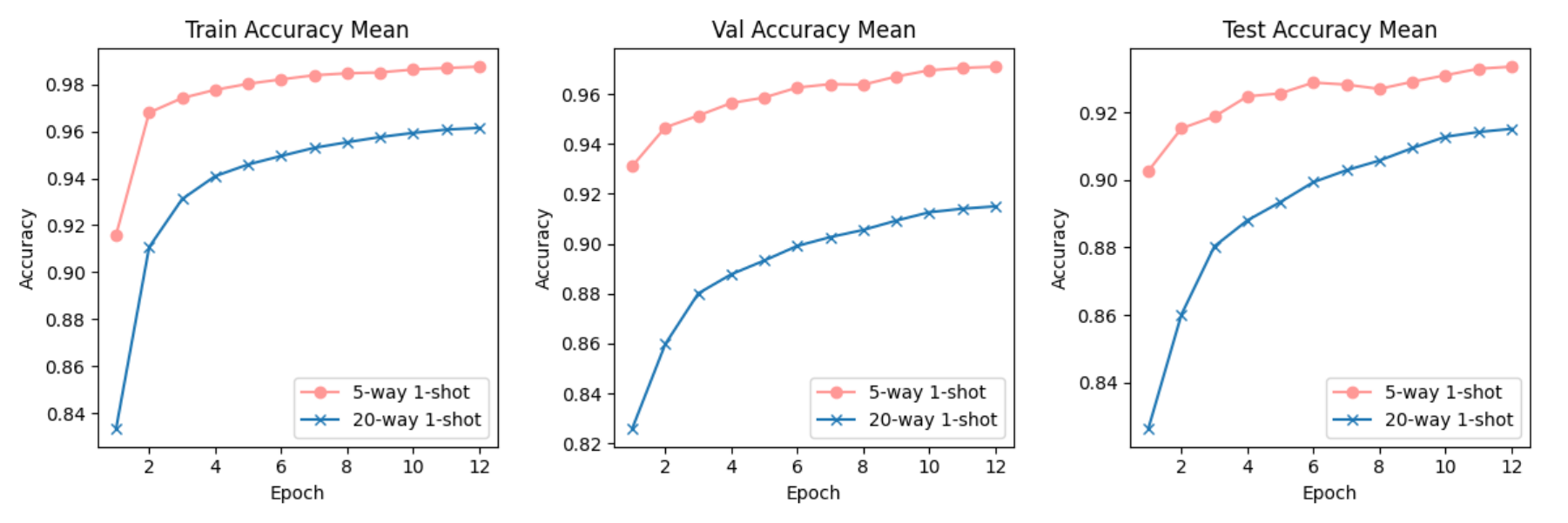}
    \caption{Meta-Training, Meta-Validation, and Meta-Testing accuracy over epochs for different task sizes ($K=5$ and $K=20$)}
    \label{fig:TRMAML_Train_Valid_Test_5W20W_1S_acc}
\end{figure*}
\begin{figure*}[t]
    \centering
    \includegraphics[width=0.8\textwidth]{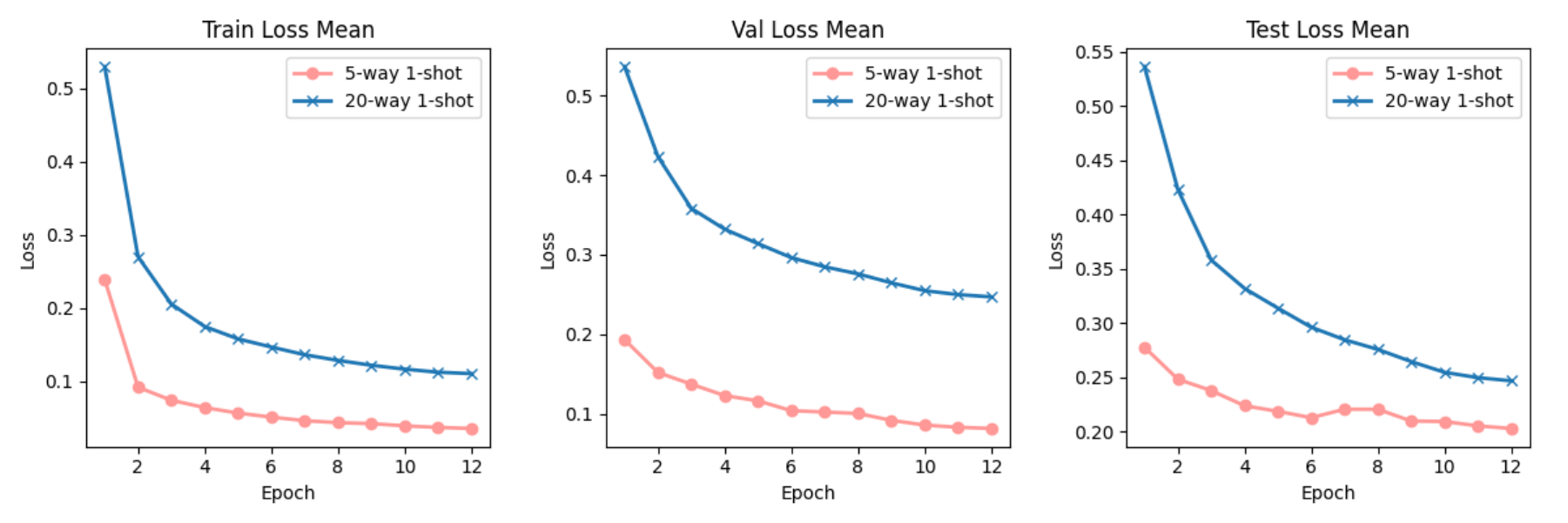}
    \caption{Meta-Training, Meta-Validation, and Meta-Testing losses over epochs for different task sizes ($K=5$ and $K=20$)}
    \label{fig:TRMAML_Train_Valid_Test_5W20W_1S_loss}
\end{figure*}
\paragraph{Task-agnostic Robust MAML:}
In task robust MAML, we aim to train a model that exhibits consistent performance across all observed task distributions, thereby ensuring robustness against task shifts at test time. The objective can be formulated as the following minimax problem:
\begin{align*}
    \min_{\bw \in \R^d}\max_{\balpha \in \Delta^N}  \sum_{i=1}^N \alpha(i) \ell_i^{\mathrm{test}}(\bw - \eta \nabla \ell_i^{\mathrm{ train}}(\bw))
\end{align*}
where $\ell_i^{\mathrm{train}}(\bw) := \frac{1}{|S_i^{\mathrm{train}}|} \sum_{\xi \in S_i^{\mathrm{train}}} \ell(\bw;\xi)$, $\ell_i^{\mathrm{test}}(\bw) := \frac{1}{|S_i^{\mathrm{test}}|} \sum_{\xi \in S_i^{\mathrm{test}}} \ell(\bw;\xi)$. This constitutes a nonconvex-concave compositional minimax problem involving the primal variable. At time t, sample minibatch $\xi^{\mathrm{train}}_1,...,\xi^{\mathrm{train}}_N$ from the training set of tasks $1,...,N$, and similarly for $\xi^{\mathrm{test}}_1,...,\xi^{\mathrm{test}}_N$.\\
First, update the auxiliary variable:
\begin{align*}
    \bz^{t+1}_i = (1-\beta^t) \pare{\bz_i^t + (\bw^t- \eta \nabla \ell_i(\w^t;\xi^{\mathrm{train}}_i)) - (\bw^{t-1}- \eta \nabla \ell_i(\w^{t-1};\xi^{\mathrm{train}}_i) ) }  + \beta^t(\bw^t-  \eta \nabla \ell_i(\w^t;\xi^{\mathrm{train}}_i)) \forall i \in [N]
\end{align*}
Next, compute the gradients:
\begin{align*}
    \bg_{\alpha}^t &= [\ell^{\mathrm{test}}_1(\bz_1^{t+1};\xi^{\mathrm{test}}_1), ..., \ell^{\mathrm{test}}_N(\bz_N^{t+1};\xi^{\mathrm{test}}_N)],\\  \bg_{\bw}^t &= \sum_{i=1}^N \alpha^t(i)\nabla \ell_i^{\mathrm{test}}(\bz_i^{t+1};\xi^{\mathrm{test}}_i) (\mI - \eta \nabla^2 \ell_i^{\mathrm{train}}(\bw^t;\xi_i^{\mathrm{train}})) 
\end{align*}
Finally, update the model and the mixture weights $\balpha$:
 \begin{align*}
    {\balpha}^{t+1} &= \cP_{ {\Delta^N}}\pare{\balpha^{t} + \eta_\alpha \bg_{\balpha}^t }   , \bw^{t+1}  =  \bw^{t} - \eta_w \bg_{w}^t
 \end{align*}
It's notable that $(\mI - \eta \nabla^2 \ell_i^{\mathrm{train}}(\bw^t;\xi_i^{\mathrm{train}}))$ involves the Hessian matrix $\nabla^2 \ell_i^{\mathrm{train}}(\bw^t;\xi_i^{\mathrm{train}})$. Inspired by practical implementations of MAML~\cite{finn2017model} that simply ignore the second-order term, we use the same trick in our experiments. Our experimental framework involves a range of tasks, each comprising several instances. Contrasting with the standard few-shot image classification including characters from different alphabets, we consider a 5-way 1-shot classification problem involving characters from 5 distinct classes (characters) which are from a specific alphabet. Note that during meta-training, we first sample an alphabet uniformly, then sample an $N$-way, $K$-shot problem uniformly from that alphabet. The convergence of our CODA-Primal algorithm for cases with $K = 5$ and $K = 20$ is illustrated in Figure~\ref{fig:TRMAML_Train_Valid_Test_5W20W_1S_loss}, which shows the meta-training, meta-validation, and meta-testing loss. Figure~\ref{fig:TRMAML_Train_Valid_Test_5W20W_1S_acc} compares the accuracies in different settings, which shows that the accuracies in 5-way 1-shot learning are higher than those in 20-way 1-shot learning.
\paragraph{Mixture Weights Estimation in Multi-source Domain Adaptation:}
Consider $N$ source domains $\cD_1,...,\cD_N$ and a single target domain $\cT$. The objective is to train a model on a mix of multiple source domains to achieve high accuracy on the target domain. The proposed algorithm consists of two stages.\\
\textbf{Phase I: mixing weight estimation}
We are aimed at solving:
\begin{align*}
    \min_{\alpha \in \Delta_N} \max_{\bw \in \mathcal{W}} F(\boldsymbol{\alpha},\mathbf{w}) &=\sum_{j=1}^N \alpha(j) g(f_\cT(\bw) - f_j(\bw)) + C {\sum\nolimits_{j=1}^N\frac{\alpha^2(j)}{m_j}}, 
    \label{eq:F}
\end{align*} 
where $g(x) = \sqrt{x^2+c}$ is the approximator of the absolute value and $f_\cT, f_1,...,f_N$ represent loss functions realized by target domain and source domains (training) data. Here, $m_1,...,m_N$ denote the number of training data from source $1,..., N$.\\
\begin{table*}[b]
\centering
\begin{tabular}{c|c|c|c}
\hline
\textbf{Group} & \textbf{Classes} & \textbf{Domains per Group} & \textbf{Samples per Domain} \\
\hline
1 & Airplanes, Cars, Birds & 5 & 1500 \\
\hline
2 & Cats, Deer, Dogs & 5 & 1500 \\
\hline
3 & Frogs, Horses, Ships, Trucks & 5 & 1500 \\
\hline
4 & Groups 1 \& 2 Combined & 5 & 1500 \\
\hline
\end{tabular}
\caption{Classes and Samples per Domain for Each CIFAR-10 Group}
\label{tab:CIFAR10_non_iid_groups}
\end{table*}
\textbf{Phase II: Weighted empirical risk minimization}
After we obtain the mixture weight $\balpha$, we perform empirical risk minimization on the mixed domain:
\begin{align*}
    \hat h = \arg\min_{h\in\cH} \sum_{i=1}^N \alpha(i) \cL_{\hat\cD_i}(h).
\end{align*}
Subsequently, we evaluate $\hat h$ on the target domain's test data and present the accuracy findings. The CIFAR-10 non-IID data generation has been shown in Table~\ref{tab:CIFAR10_non_iid_groups}.

\newpage

\section{Details of Variance Reduction Algorithm}\label{app:VR algorithm}
In this section we present whole algorithm pseduo-code of CODA-Primal+, as in Algorithm~\ref{algorithm: CODA-Primal+}. We define the following notations for the convenience of presenting:
\begin{align*}
    &G_k^{\bx}(\bx,\bz, \by; \cI) = \frac{1}{|\cI|} \sum_{(\zeta,\xi) \in \cI } 
 \nabla_1 f(\bz,\by;\zeta)  \nabla g(\bx;\xi)   + \mu_x  (\bx - \bx^0), G_k^{\by}(\bx,\bz, \by; \cI) = \frac{1}{|\cI|} \sum_{(\zeta,\xi) \in \cI }  \nabla_2 f(\bz,\by;\zeta)  .   
\end{align*}
Roughly speaking, we use large mini-batch every $\tau$ iterations to update gradients and inner function estimation variable, and use small mini-batch for the rest of iterations.

  \begin{algorithm2e}[H]  
    \caption{\sffamily{CODA-Primal+}}
    	\label{algorithm: CODA-Primal+}
    	\textbf{Input:}  ${\bx}^{0},{\by}^{0}$, $\mu_\bx$, stepsizes $\eta_{\bx},\eta_{\by}$. 
	\\

    {
   
	\For{$k = 0,...,K-1$}
	 {
           Construct $F_k(\bx,\by) = F( \bx ,\by)  + \frac{\mu_\bx}{2}  {\norm{\bx - \bx_k}^2   }$\\
           $(\bx_{k+1},\by_{k+1}) = \texttt{CODA-SCSC+}(F_k, \bx_k, \by_k, T_k )$
 	 }
      
 	      }
 
\end{algorithm2e}
 \vspace{\floatsep}
 \begin{algorithm2e}[H] 
	\DontPrintSemicolon
    \caption{\sffamily{CODA-SCSC+}$(F,\bx^0,\by^0,T)$}
    	\label{algorithm: CODA-SCSC+}
     	\textbf{Input:}  ${\bx}^{0},{\by}^{0}$, $\bz^0$ such that $\E \norm{\bz^0 - g(\bx^0)}^2 \leq \delta$.

   {
    {

	\For{$t = 0,...,T-1$}
	 {     
  \If{$t \mod \tau = 0$ }
            {
            Sample mini-batch $\cB^t_{\bz} = \cbr{\xi_1,...,\xi_{B_{\tau}}}$, and $\cB^t_{\by} = \cbr{\xi_1,...,\xi_{B_{\tau}}}$

            $g^{t} = g(\bx^t; \cB^t_{\bz})$, $\bz^{t+1} = (1-\beta)(\bz^t + g(\bx^t; \cB^t_{\bz}) - g(\bx^{t-1}; \cB^t_{\bz})) + \beta g^{t} $\\
                $\bq^{t} = G_k^{\by}(\bx^{t},\bz^{t+1}, \by^{t};\cB^t_{\by}) $

            }
            \Else
            {
             Sample mini-batch $\cI^t_{\bz} = \cbr{\xi_1,...,\xi_{B}}$, and $\cI^t_{\by} = \cbr{\xi_1,...,\xi_{B}}$

            $g^{t} = g^{t-1} + g(\bx^t;\cI^t_{\bz}) - g(\bx^{t-1};\cI^t_{\bz})$,  $\bz^{t+1} = (1-\beta)(\bz^t + g(\bx^t; \cI^t_{\bz}) - g(\bx^{t-1}; \cI^t_{\bz})) + \beta g^{t-1} $\\
                $\bq^{t} = \bq^{t-1} + G_k^{\by}(\bx^{t},\bz^{t+1}, \by^{t};\cI^t_{\by}) - G_k^{\by}(\bx^{t-1},\bz^{t}, \by^{t-1};\cI^t_{\by})$

            }
            $\bg_{\by}^{t} = 2 \bq^{t} -   \bq^{t-1} $\\
           $\by^{t+1}  =   \by^t  + \eta_{\by} ( \bg_{\by}^t - \nabla r(\by^{t+1}) ) $\\
           
           \If{$t \mod \tau = 0$ }
           {
            Sample mini-batch $\cB^t_{\bx} = (\cB^t_{\bx,f},\cB^t_{\bx,g}) = (\cbr{\zeta_1,...,\zeta_{B_{\tau}}}, \cbr{\xi_1,...,\xi_{B_{\tau}}}) $\\

            $\bg_{\bx}^{t} =  G_k^{\bx}(\bx^t, \bz^{t+1}, \by^{t+1};\cB^t_{\bx})$\\

           }
           
           \Else
           {
             Sample mini-batch $\cI^t_{\bx} = (\cI^t_{\bx,f},\cI^t_{\bx,g}) = (\cbr{\zeta_1,...,\zeta_{B}}, \cbr{\xi_1,...,\xi_{B}}) $\\

$\bg_{\bx}^{t} = \bg_{\bx}^{t-1} + G_k^{\bx}(\bx^t,\bz^{t+1}, \by^{t+1}; \cI^t_{\bx})- G_k^{\bx}(\bx^{t-1},\bz^{t}, \by^{t}; \cI^t_{\bx}) $
           
           }
        $\bx^{t+1} = \bx^t - \eta_{\bx} (\bg_{\bx}^{t} + \nabla h(\bx^t)) $
 
 	      }
 	 }

 	\textbf{Output:}  $(\hat{\bx}, \hat{\by}) = \frac{1}{T}\sum_{t=1}^T ( {\bx}^t, {\by}^t)$  }   
\end{algorithm2e}


\newpage

\section{Proof of Primal Composition Setting} \label{app: proof primal}
 
In this section we provide the proof of results in primal compostion setting (Theorem~\ref{thm:NCSC} and Theorem~\ref{thm:ncc}).

The following Lemma is from~\cite{chen2021solving}, which bound the tracking error of the stochastic correction algorithm.
\begin{lemma}[Tracking Error~\cite{chen2021solving}]\label{lem:tracking error}
For Algorithm~\ref{algorithm: CODA-Primal}, under the assumptions of Theorem~\ref{thm:NCSC}, the following statement holds true:
\begin{align*}
    \E \norm{\bz^{t+1}  - g(\bx^{t}) }^2 &\leq (1-\beta)^2 \E \norm{\bz^{t}  - g(\bx^{t-1}) }^2  + 4(1-\beta)^2 G_g^2 \norm{\bx^t - \bx^{t-1}}^2 + 2\beta^2 \frac{\sigma^2}{B}.
\end{align*} 
\end{lemma}
\subsection{Proof of Nonconvex-strongly-concave Setting} \label{app: proof NCSC}

\begin{lemma} \label{lem:NCSC one iteration}
For Algorithm~\ref{algorithm: CODA-Primal}, under the assumptions of Theorem~\ref{thm:NCSC}, the following statement holds true:
    \begin{align*}
        \E[ \Phi(\bx^{t+1})]  
       &\leq  \E[ \Phi(\bx^{t})] - \frac{\eta_{\bx}}{2} \E \norm{\nabla \Phi(\bx^t)}^2 - \pare{\frac{\eta_{\bx}}{2} -  {\eta_{\bx}^2 \kappa L}   }\E\norm{\bar{\bg}_{\bx}^t}^2  \\
       & \quad+ {\eta_{\bx}}G_g^2L^2_f\E\norm{ g(\bx^{t})  -  \bz^{t+1}  }^2   + {\eta_{\bx}} G_g^2L^2_f \E\norm{  \by^*(\bx^t) -  \by^t   }^2  + \frac{2\eta_{\bx}^2\kappa L( G_f^2\sigma^2 + G_g^2\sigma^2)}{B},
    \end{align*}
    where $\bar{\bg}_{\bx}^t =  \nabla_1 f(\bz^{t+1}, \by^t ) \nabla g(\bx^t )$.
\begin{proof}
According $\kappa L + L$ smoothness [Lemma 4.3~\cite{lin2019gradient}] of $\Phi$ and updating rule for $\bx$, we have:
\begin{align*}
    \Phi(\bx^{t+1}) &\leq  \Phi(\bx^{t}) + \inprod{\nabla \Phi(\bx^t) }{ \bx^{t+1} - \bx^t} +  {\kappa L}  \norm{ \bx^{t+1} - \bx^t }^2 \\
    &= \Phi(\bx^{t}) - \inprod{\nabla \Phi(\bx^t) }{ \eta_{\bx} \bg_{\bx}^t } +  {\kappa L} \eta_{\bx} \norm{ \bg^t_{\bx}}^2,
\end{align*}
where $\bg_{\bx}^t = \frac{1}{B} \sum_{(\zeta,\xi) \in \cB^t} \nabla_1 f(\bz^{t+1}, \by^t; \zeta) \nabla g(\bx^t; \xi)$. Notice the fact that
\begin{align*}
    \E\norm{  \nabla_1 f(\bz^{t+1}, \by^t; \zeta) \nabla g(\bx^t; \xi) - \nabla_1 f(\bz^{t+1}, \by^t ) \nabla g(\bx^t )}^2 \leq 2G_f^2\sigma^2 + 2G_g^2\sigma^2.
\end{align*}

Taking expectation over the randomness of $\cB^t$ yields:
 
\begin{align*}
   \E [\Phi(\bx^{t+1})] &\leq  \E[\Phi(\bx^{t})]  - \E \inprod{\nabla \Phi(\bx^t) }{ \eta_{\bx}  \nabla_1 f(\bz^{t+1}, \by^t ) \nabla g(\bx^t )  } \\
   &+ {\eta_{\bx}^2\kappa L}   \E\norm{\nabla_1 f(\bz^{t+1}, \by^t ) \nabla g(\bx^t ) + \nabla h(\bx^t)}^2  + \frac{2\eta_{\bx}^2\kappa L( G_f^2\sigma^2 + G_g^2\sigma^2)}{ B} .
\end{align*}

Due to the identity $\inprod{\ba}{\bb} = \frac{1}{2}\norm{\ba}^2 + \frac{1}{2}\norm{\bb}^2 - \frac{1}{2}\norm{\ba-\bb}^2 $,  we have:
    \begin{align*}
        \E[\nabla \Phi(\bx^{t+1})] &\leq \E[\nabla \Phi(\bx^{t})] - \frac{\eta_{\bx}}{2} \E \norm{\nabla \Phi(\bx^t)}^2 - \frac{\eta_{\bx}}{2} \norm{\bar{\bg}_{\bx}^t}^2 \\
        &\quad + \frac{\eta_{\bx}}{2}\norm{\nabla \Phi(\bx^{t}) - \nabla_1 f(\bz^{t+1},\by^t) \nabla g(\bx^t) - \nabla h(\bx^t) }^2 + {\eta_{\bx}^2 \kappa L}   \E\norm{ \bar{\bg}_{\bx}^t }^2 + \frac{2\eta_{\bx}^2\kappa L( G_f^2\sigma^2 + G_g^2\sigma^2)}{ B} \\
\end{align*}
Due to [Lemma 4.3~\cite{lin2019gradient}], we know $\nabla \Phi(\bx) = \nabla F(\bx,\by^*(\bx))$, we have:
  \begin{align*}
      \E[\nabla \Phi(\bx^{t+1})]    &\leq  \E[\nabla \Phi(\bx^{t})] - \frac{\eta_{\bx}}{2} \E \norm{\nabla \Phi(\bx^t)}^2 - \frac{\eta_{\bx}}{2} \norm{\bar{\bg}_{\bx}^t}^2\\
        &\quad + \frac{\eta_{\bx}}{2}\norm{\nabla_\bx f(g(\bx^{t}),\by^*(\bx^t))   -  \nabla_1 f(\bz^{t+1},\by^t) \nabla g(\bx^t)  }^2 + {\eta_{\bx}^2 \kappa L}   \E\norm{ \bar{\bg}_{\bx}^t }^2 + \frac{2\eta_{\bx}^2\kappa L( G_f^2\sigma^2 + G_g^2\sigma^2)}{ B} \\
        &\leq  \E[\nabla \Phi(\bx^{t})] - \frac{\eta_{\bx}}{2} \E \norm{\nabla \Phi(\bx^t)}^2 - \frac{\eta_{\bx}}{2} \norm{\bar{\bg}_{\bx}^t}^2  + {\eta_{\bx}^2 \kappa L}   \E\norm{ \bar{\bg}_{\bx}^t }^2 + \frac{2\eta_{\bx}^2\kappa L( G_f^2\sigma^2 + G_g^2\sigma^2)}{ B} \\
        &\quad + \frac{\eta_{\bx}}{2}\norm{\nabla_\bx f(g(\bx^{t}),\by^*(\bx^t)) -\nabla_\bx f(g(\bx^{t}),\by^t)+\nabla_\bx f(g(\bx^{t}),\by^t)  -  \nabla_1 f(\bz^{t+1},\by^t) \nabla g(\bx^t)  }^2 \\ 
        &\leq  \E[\nabla \Phi(\bx^{t})] - \frac{\eta_{\bx}}{2} \E \norm{\nabla \Phi(\bx^t)}^2 - \frac{\eta_{\bx}}{2} \norm{\bar{\bg}_{\bx}^t}^2+ {\eta_{\bx}^2 \kappa L}   \E\norm{ \bar{\bg}_{\bx}^t }^2 + \frac{2\eta_{\bx}^2\kappa L( G_f^2\sigma^2 + G_g^2\sigma^2)}{ B}\\
        &\quad + {\eta_{\bx}} \norm{\nabla_1 f(g(\bx^{t}),\by^*(\bx^t)) \nabla g(\bx^t) -  \nabla_1 f(g(\bx^{t}),\by^t) \nabla g(\bx^t)  }^2  \\
        &\quad + {\eta_{\bx}} \norm{ \nabla_1 f(g(\bx^{t}),\by^t) \nabla g(\bx^t)  -  \nabla_1 f(\bz^{t+1},\by^t) \nabla g(\bx^t)  }^2 \\
       &\leq  \E[\nabla \Phi(\bx^{t})] - \frac{\eta_{\bx}}{2} \E \norm{\nabla \Phi(\bx^t)}^2 - \pare{\frac{\eta_{\bx}}{2} - {\eta_{\bx}^2 \kappa L}   }\E\norm{\bar{\bg}_{\bx}^t}^2  \\
       & \quad+ {\eta_{\bx}}G_g^2L^2_f\E\norm{ g(\bx^{t})  -  \bz^{t+1}  }^2   + {\eta_{\bx}} G_g^2 L^2_f \E\norm{  \by^*(\bx^t) -  \by^t   }^2  + \frac{2\eta_{\bx}^2\kappa L( G_f^2\sigma^2 + G_g^2\sigma^2)}{B}.
    \end{align*}
        
\end{proof}
\end{lemma}

\begin{lemma} \label{lem:NCSC dual}
For Algorithm~\ref{algorithm: CODA-Primal}, under the assumptions of Theorem~\ref{thm:NCSC}, the following statement holds true:
   \begin{align*}
        \norm{\by^{t+1} - \by^*(\bx^{t+1})}^2  
         \leq  &\pare{1-\frac{1}{16\kappa} }\norm{\by^{t}- \by^*(\bx^t)}^2  +   2 \pare{\frac{4\eta_{\by}}{ \mu}+2\eta_{\by} ^2  } L^2_f\norm{ \bz^{t+1} -g(\bx^{t}) }^2 \\
         &  \quad + 16  \kappa^3 \norm{\bx^{t+1} - \bx^t}^2 + \frac{2\eta_{\by}^2\sigma^2}{B}.
\end{align*}
\begin{proof}
We notice the following canonical decomposition due to Young's inequality:
\begin{align*}
        \norm{\by^{t+1} - \by^*(\bx^{t+1})}^2 \leq \pare{1+\frac{1}{2(2\kappa-1)}} \underbrace{\norm{\by^{t+1} - \by^*(\bx^t)}^2}_{T_1} + \pare{1+2(2\kappa-1)} \underbrace{\norm{\by^*(\bx^t) - \by^*(\bx^{t+1})}^2}_{T_2}.
\end{align*}
We bound $T_1$ first. Define $\bar\bg_{\by}^t = \nabla_{\by} f(\bz^{t+1},\by^t) - \nabla r(\by^t)$. By the updating rule we have:
\begin{align*}
    \E\norm{\by^{t+1} - \by^*(\bx^t)}^2 &= \E\norm{\cP_{\cY}(\by^{t} + \eta_{\by} \bg_{\by}^t) - \cP_{\cY} \left(\by^*(\bx^t) +  \eta_{\by} \nabla_{\by} F(\bx^t, \by^*(\bx^t)) \right) }^2 \\
    &\leq \E\norm{\by^{t} + \eta_{\by} \bar{\bg}_{\by}^t - \by^*(\bx^t)  -   \eta_{\by} \nabla_{\by} F(\bx^t, \by^*(\bx^t))}^2 + \frac{\eta_{\by}^2\sigma^2}{B} \\
    & = \E\norm{\by^{t}- \by^*(\bx^t)}^2 + 2\eta_{\by} \inprod{\bar\bg_{\by}^t-   \nabla_{\by} F(\bx^t, \by^*(\bx^t))}{\by^{t}- \by^*(\bx^t)} + \frac{\eta_{\by}^2\sigma^2}{B} \\
    &\quad + \eta_{\by} ^2 \norm{\bar\bg_{\by}^t -    \nabla_{\by} F(\bx^t, \by^*(\bx^t))}^2\\
    & = \norm{\by^{t}- \by^*(\bx^t)}^2 + 2\eta_{\by} \inprod{ \nabla_{\by} F(\bx^t, \by^t)-   \nabla_{\by} F(\bx^t, \by^*(\bx^t))}{\by^{t}- \by^*(\bx^t)}\\
    &\quad +  2\eta_{\by} \inprod{\nabla_{\by} f( \bz^{t+1} , \by^t) - \nabla_{\by} f(g(\bx^{t}), \by^t)  }{\by^{t}- \by^*(\bx^t)} \\
    &\quad + 2\eta_{\by} ^2 \norm{ \nabla_{\by} F(\bx^{t}, \by^t) -\nabla_{\by} F(\bx^{t}, \by^*(\bx^t))}^2  +  2\eta_{\by} ^2 \norm{ \bar\bg_{\by}^t -\nabla_{\by} F(\bx^{t}, \by^t)}^2+ \frac{\eta_{\by}^2\sigma^2}{B}\\
   & \leq \norm{\by^{t}- \by^*(\bx^t)}^2 + (2\eta_{\by}-2\eta^2_{\by} L )\inprod{ \nabla_{\by} F(\bx^t, \by^t)-   \nabla_{\by} F(\bx^t, \by^*(\bx^t))}{\by^{t}- \by^*(\bx^t)}\\
    &\quad +  2\eta_{\by} \inprod{\nabla_{\by} f( \bz^{t+1} , \by^t) - \nabla_{\by} f(g(\bx^{t}), \by^t)  }{\by^{t}- \by^*(\bx^t)}   +  2\eta_{\by} ^2 \norm{ \bar\bg_{\by}^t -\nabla_{\by} F(\bx^{t}, \by^t)}^2+ \frac{\eta_{\by}^2\sigma^2}{B}\\
     & \leq \left(1- \frac{1}{2\kappa} \right)\norm{\by^{t}- \by^*(\bx^t)}^2  +  2\eta_{\by} \inprod{\nabla_{\by} f( \bz^{t+1} , \by^t) - \nabla_{\by} f(g(\bx^{t}), \by^t)  }{\by^{t}- \by^*(\bx^t)}  \\
     &\quad +  2\eta_{\by} ^2 \norm{\nabla_{\by} f( \bz^{t+1} , \by^t) - \nabla_{\by} f(g(\bx^{t}), \by^t)}^2+ \frac{\eta_{\by}^2\sigma^2}{B},
\end{align*}
where at first step, we use the fact that $\by^*(\bx^t) = \cP_{\cY}(\by^*(\bx^t) +  \eta_{\by} \nabla_{\by} F(\bx^t, \by^*(\bx^t)))$ for any $\eta_{\by} \geq 0$~\cite[Proposition 7.4]{recht-wright}; at fifth step, we use the classic result for smooth concave function: $\frac{1}{L}\norm{ \nabla_{\by} F(\bx^{t}, \by^t) -\nabla_{\by} F(\bx^{t}, \by^*(\bx^t))}^2 \leq \inprod{   \nabla_{\by} F(\bx^t, \by^*(\bx^t)) - \nabla_{\by} F(\bx^t, \by^t)}{\by^{t}- \by^*(\bx^t)}$, and at last step we use the $\mu$ strongly concavity of $F(\bx,\cdot)$, and $\eta_{\by} \leq \frac{1}{2L}$.

For the remaining inner product, we  apply Cauchy-schwartz inequality:
\begin{align*}
     \inprod{\nabla_{\by} f( \bz^{t+1} , \by^t) - \nabla_{\by} f(g(\bx^{t}), \by^t)  }{\by^{t}- \by^*(\bx^t)} &\leq \norm{ \nabla_{\by} f( \bz^{t+1} , \by^t) - \nabla_{\by} f(g(\bx^{t}), \by^t) }\norm{\by^{t}- \by^*(\bx^t)}\\
     &\leq \frac{2}{ \mu}   \norm{\nabla_{\by} f( \bz^{t+1} , \by^t) - \nabla_{\by} f(g(\bx^{t}), \by^t) }^2 + \frac{\mu}{8}\norm{\by^{t}- \by^*(\bx^t)}^2\\ 
     &\leq  \frac{ 2L^2_f}{ \mu} \norm{ \bz^{t+1} -g(\bx^{t}) }^2+ \frac{\mu}{8}\norm{\by^{t}- \by^*(\bx^t)}^2.
\end{align*}

Similarly for $\norm{\nabla_{\by} f( \bz^{t+1} , \by^t) - \nabla_{\by} f(g(\bx^{t}), \by^t)}^2$ we have:
\begin{align*}
    \norm{\nabla_{\by} f( \bz^{t+1} , \by^t) - \nabla_{\by} f(g(\bx^{t}), \by^t)}^2 \leq L_f^2 \norm{ \bz^{t+1} -g(\bx^{t})}^2
\end{align*}

Putting pieces together yields:
\begin{align*}
    \norm{\by^{t+1} - \by^*(\bx^t)}^2 & \leq \pare{1-\frac{1}{8\kappa} }\norm{\by^{t}- \by^*(\bx^t)}^2  +   \pare{\frac{4\eta_{\by}}{ \mu}+2\eta_{\by} ^2  } L^2_f\norm{ \bz^{t+1} -g(\bx^{t}) }^2 + \frac{\eta_{\by}^2\sigma^2}{B} .
\end{align*}
We then turn to bounding $T_2$.
\begin{align*}
    T_2 = \norm{\by^*(\bx^t) - \by^*(\bx^{t+1})}^2 \leq \kappa^2 \norm{\bx^{t+1} - \bx^t}^2
\end{align*}
Putting pieces together will conclude the proof:
\begin{align*}
        \norm{\by^{t+1} - \by^*(\bx^{t+1})}^2 & \leq \pare{1+\frac{1}{2(8\kappa-1)}}  \pare{\pare{1-\frac{1}{8\kappa} }\norm{\by^{t}- \by^*(\bx^t)}^2  +   \pare{\frac{4\eta_{\by}}{ \mu}+2\eta_{\by} ^2  } L^2_f\norm{ \bz^{t+1} -g(\bx^{t}) }^2+ \frac{\eta_{\by}^2\sigma^2}{B}     } \\
        & \quad + \pare{1+2(8\kappa-1)} \norm{\by^*(\bx^t) - \by^*(\bx^{t+1})}^2\\
        & \leq  \pare{1-\frac{1}{16\kappa} }\norm{\by^{t}- \by^*(\bx^t)}^2  +   \pare{1+\frac{1}{2(8\kappa-1)}} \pare{\frac{4\eta_{\by}}{ \mu}+2\eta_{\by} ^2  } L^2_f\norm{ \bz^{t+1} -g(\bx^{t}) }^2 \\
        & \quad + 16  \kappa^3  \norm{\bx^{t+1} - \bx^t}^2+ \frac{2\eta_{\by}^2\sigma^2}{B} \\ 
         & \leq  \pare{1-\frac{1}{16\kappa} }\norm{\by^{t}- \by^*(\bx^t)}^2  +   2 \pare{\frac{4\eta_{\by}}{ \mu}+2\eta_{\by} ^2  } L^2_f\norm{ \bz^{t+1} -g(\bx^{t}) }^2 \\
         &  \quad + 16  \kappa^3 \norm{\bx^{t+1} - \bx^t}^2 + \frac{2\eta_{\by}^2\sigma^2}{B}.
\end{align*}
\end{proof}

\end{lemma} 
\subsubsection{Proof of Theorem~\ref{thm:NCSC}}

\paragraph{Proof Sketch.} Before delving into the detailed proof, we first provide a sketch.  Unlike the classic nonconvex-strongly-concave analysis, we consider the following potential function: 
\begin{align*}
   \Psi(\bx^{t+1}) =&\E[  \Phi(\bx^{t+1})] +  O\pare{\kappa  \eta_{\bx}G_g^2 L_f^4   \pare{\frac{\eta_{\by}}{ \mu}+ \eta_{\by} ^2  }}  \E\norm{ g(\bx^{t})  -  \bz^{t+1}  }^2 \\
   &+  O(\kappa  \eta_{\bx}G_g^2 L_f^2) \E\norm{  \by^*(\bx^t) -  \by^t   }^2 +  O(\kappa^4  \eta_{\bx}G_g^2 L_f^2) \E\norm{\bx^{t+1} - \bx^{t}}^2.
\end{align*}
Then, we derive the convergence of dual gap $\norm{\by^{t+1} - \by^*(\bx^{t+1})}^2$:
   \begin{align*}
         \norm{\by^{t+1} - \by^*(\bx^{t+1})}^2  
         \leq  &\pare{1-\frac{1}{16\kappa} }\norm{\by^{t}- \by^*(\bx^t)}^2  +   2 \pare{\frac{4\eta_{\by}}{ \mu}+2\eta_{\by} ^2  } L^2_f\norm{ \bz^{t+1} -g(\bx^{t}) }^2 \\
         &  \quad + 16  \kappa^3 \norm{\bx^{t+1} - \bx^t}^2 + \frac{2\eta_{\by}^2\sigma^2}{B},
\end{align*}
where $\norm{ \bz^{t+1} -g(\bx^{t}) }^2$ is the error of inner function estimation, and can be bounded using standard result from~\cite{chen2021solving}. Combining above ingredients will yield
 \begin{align*}
       \Psi(\bx^{t+1})   &\leq\Psi(\bx^{t}) - \frac{\eta_{\bx}}{2} \E \norm{\nabla \Phi(\bx^t)}^2  \\
       & \quad +O\pare{ \frac{( \kappa^4  \eta_{\bx}^3G_g^2 L_f^2 + \eta_\bx^2 \kappa L)( G_f^2  + G_g^2)\sigma^2}{B} }+ O\pare{\kappa^2 L \eta_{\bx}  \frac{\sigma^2}{B}}+  O\pare{\kappa  \eta_{\bx}G_g^2 L_f^2\frac{2\eta_{\by}^2 \sigma^2}{ B}  }.
    \end{align*}
Unrolling above recursion will yield desired result.

\begin{proof}
First due to Lemma~\ref{lem:NCSC one iteration} we have:
      \begin{align*}
        \E[ \Phi(\bx^{t+1})]  
       &\leq  \E[ \Phi(\bx^{t})] - \frac{\eta_{\bx}}{2} \E \norm{\nabla \Phi(\bx^t)}^2 - \pare{\frac{\eta_{\bx}}{2} -  {\eta_{\bx}^2 \kappa L}   }\E\norm{\bar{\bg}_{\bx}^t}^2  \\
       & \quad+ {\eta_{\bx}}G_g^2L^2_f\E\norm{ g(\bx^{t})  -  \bz^{t+1}  }^2   + {\eta_{\bx}} G_g^2L^2_f \E\norm{  \by^*(\bx^t) -  \by^t   }^2  + \frac{2\eta_{\bx}^2\kappa L( G_f^2\sigma^2 + G_g^2\sigma^2)}{B}.
    \end{align*}
    Define potential function $\Psi(\bx^{t+1}) =\E[  \Phi(\bx^{t+1})] + C_1 \E\norm{ g(\bx^{t})  -  \bz^{t+1}  }^2 + C_2 \E\norm{  \by^*(\bx^t) -  \by^t   }^2 + C_3 \E\norm{\bx^{t+1} - \bx^{t}}^2 $, where $C_1,C_2$ and $C_3$ are determined later. Hence we have :
    \begin{align*}
        \Psi(\bx^{t+1})  &\leq  \Psi(\bx^{t+1}) - \frac{\eta_{\bx}}{2} \E \norm{\nabla \Phi(\bx^t)}^2 - \pare{\frac{\eta_{\bx}}{2} - {\eta_{\bx}^2 \kappa L}   }\E\norm{\bar{\bg}_{\bx}^t}^2  \\
       & \quad+ {\eta_{\bx}}G_g^2L^2_f\E\norm{ g(\bx^{t})  -  \bz^{t+1}  }^2   + {\eta_{\bx}} G_g^2L^2_f \E\norm{  \by^*(\bx^t) -  \by^t   }^2  +\frac{2\eta_{\bx}^2\kappa L( G_f^2\sigma^2 + G_g^2\sigma^2)}{B}\\
       & \quad +  C_1 \E\norm{ g(\bx^{t})  -  \bz^{t+1}  }^2 + C_2 \E\norm{  \by^*(\bx^t) -  \by^t   }^2 + C_3 \E\norm{\bx^{t+1} - \bx^{t}}^2  \\
       & \quad - C_1 \E\norm{ g(\bx^{t-1})  -  \bz^{t}  }^2 - C_2 \E\norm{  \by^*(\bx^{t-1}) -  \by^{t-1}   }^2 - C_3 \E\norm{\bx^{t} - \bx^{t-1}}^2.
    \end{align*}

Since $\E\norm{\bx^{t+1} - \bx^{t}}^2  \leq \eta_{\bx}^2 \E \norm{\bar\bg_{\bx}^t}^2 + \eta_{\bx}^2 \frac{2(G_f^2 + G_g^2)\sigma^2}{B} $, we have
    \begin{align*} 
     \Psi(\bx^{t+1})     &\leq  \Psi(\bx^{t+1}) - \frac{\eta_{\bx}}{2} \E \norm{\nabla \Phi(\bx^t)}^2 - \pare{\frac{\eta_{\bx}}{2} - {\eta_{\bx}^2 \kappa L}   -C_3 \eta_{\bx}^2}\E\norm{\bar{\bg}_{\bx}^t}^2+ \frac{2C_3\eta_{\bx}^2( G_f^2\sigma^2 + G_g^2\sigma^2)}{B}  \\
       & \quad+ \pare{{\eta_{\bx}}G_g^2L^2_f+C_1}\E\norm{ g(\bx^{t})  -  \bz^{t+1}  }^2   + \pare{{\eta_{\bx}} G_g^2L^2_f +C_2}\E\norm{  \by^*(\bx^t) -  \by^t   }^2  + \frac{2\eta_{\bx}^2\kappa L( G_f^2\sigma^2 + G_g^2\sigma^2)}{B} \\ 
       & \quad - C_1 \E\norm{ g(\bx^{t-1})  -  \bz^{t}  }^2 - C_2 \E\norm{  \by^*(\bx^{t-1}) -  \by^{t-1}   }^2 - C_3 \E\norm{\bx^{t} - \bx^{t-1}}^2.
    \end{align*}
    We plug in the bound for $\E\norm{ g(\bx^{t})  -  \bz^{t+1}  }^2 $ and $\E\norm{  \by^*(\bx^t) -  \by^t   }^2$:
   \begin{align*}
        \Psi(\bx^{t+1}) 
       &\leq  \Psi(\bx^{t}) - \frac{\eta_{\bx}}{2} \E \norm{\nabla \Phi(\bx^t)}^2 - \pare{\frac{\eta_{\bx}}{2} - {\eta_{\bx}^2 \kappa L}   -C_3 \eta_{\bx}^2}\E\norm{\bar{\bg}_{\bx}^t}^2   + \frac{(2C_3\eta_{\bx}^2 + 2\eta_\bx^2 \kappa L)( G_f^2\sigma^2 + G_g^2\sigma^2)}{B}\\
       & + \pare{{\eta_{\bx}}G_g^2L^2_f+C_1}\pare{ (1-\beta)^2 \E \norm{\bz^{t}  - g(\bx^{t-1}) }^2  + 4(1-\beta)^2 G_g^2 \norm{\bx^t - \bx^{t-1}}^2 + 2\beta^2 \frac{\sigma^2}{B}} \\
       &  + \pare{{\eta_{\bx}} G_g^2L^2_f +C_2}\\
       &  \times \pare{    \pare{1-\frac{1}{16\kappa} }\norm{\by^{t-1}- \by^*(\bx^{t-1})}^2  +   2 \pare{\frac{4\eta_{\by}}{ \mu}+2\eta_{\by} ^2  } L^2_f\norm{ \bz^{t} -g(\bx^{t-1}) }^2  + 16  \kappa^3 \norm{\bx^{t-1} - \bx^t}^2 + \frac{2\eta_{\by}^2\sigma^2}{B} } \\ 
       & - C_1 \E\norm{ g(\bx^{t-1})  -  \bz^{t}  }^2 - C_2 \E\norm{  \by^*(\bx^{t-1}) -  \by^{t-1}   }^2 - C_3 \E\norm{\bx^{t} - \bx^{t-1}}^2\\
       \leq & \Psi(\bx^{t+1}) - \frac{\eta_{\bx}}{2} \E \norm{\nabla \Phi(\bx^t)}^2 - \pare{\frac{\eta_{\bx}}{4}    -C_3 \eta_{\bx}^2}\E\norm{\bar{\bg}_{\bx}^t}^2\\
       & + \frac{(2C_3\eta_{\bx}^2 + \eta_\bx^2 \kappa L)( G_f^2  + G_g^2)\sigma^2}{B} +\pare{{\eta_{\bx}} G_g^2L^2 +C_1}2\beta^2 \frac{\sigma^2}{B}+ \pare{{\eta_{\bx}} G_g^2L^2_f +C_2}\frac{2\eta_{\by}^2 \sigma^2}{ B}\\
       &   + \left[\pare{{\eta_{\bx}}G_g^2L^2_f+C_1} (1-\beta)^2 +\pare{{\eta_{\bx}} G_g^2L^2 +C_2} 2 \pare{\frac{4\eta_{\by}}{ \mu}+ 2\eta_{\by} ^2  } L^2_f - C_1 \right]\E \norm{\bz^{t}  - g(\bx^{t-1}) }^2 \\
       & + \left[ \pare{{\eta_{\bx}} G_g^2L^2_f +C_2}  \pare{1-\frac{1}{16\kappa} } - C_2\right]\E\norm{\by^{t-1}- \by^*(\bx^{t-1})}^2 +\left[  16  \kappa^3 \pare{{\eta_{\bx}} G_g^2L^2_f +C_2}-C_3\right] \E\norm{\bx^{t} - \bx^{t-1}}^2, 
    \end{align*}
    where at second inequality we use the fact $\eta_{\bx} \leq \frac{1}{4\kappa L}$, so $\frac{\eta_{\bx}}{2}  - \eta_{\bx}\kappa L \geq \frac{\eta_{\bx}}{4}$.
    We choose $C_1,C_2,C_3$ such that the following conditions are satisfied:
    \begin{align*}
    \frac{\eta_{\bx}}{4}    -C_3 \eta_{\bx}^2 \geq 0,\\
        \pare{{\eta_{\bx}}G_g^2L^2_f+C_1} (1-\beta)  +\pare{{\eta_{\bx}} G_g^2L^2_f +C_2} 2 \pare{\frac{4\eta_{\by}}{ \mu}+ 2\eta_{\by} ^2  } L^2_f - C_1  \leq 0,\\
         \pare{{\eta_{\bx}} G_g^2L^2_f +C_2}  \pare{1-\frac{1}{16\kappa} } - C_2 \leq 0,\\
         16  \kappa^3 \pare{{\eta_{\bx}} G_g^2L^2_f +C_2}-C_3 \leq 0.
    \end{align*}
    It suffices to guarantee the following inequality holding:
    \begin{align*}
     & 16\kappa   \pare{1-\frac{1}{16\kappa}} \eta_{\bx} G_g^2 L_f^2 \leq C_2 \leq \frac{1}{16\kappa^3} C_3 - \eta_{\bx} G_g^2 L_f^2,\\
     & C_1 \geq \frac{1}{\beta}\pare{{\eta_{\bx}} G_g^2L^2_f +C_2} 2 \pare{\frac{4\eta_{\by}}{ \mu}+ 2\eta_{\by} ^2  } L^2_f  + \frac{1-\beta}{\beta}{\eta_{\bx}} G_g^2L^2_f. 
    \end{align*}
    A set of choices that can satisfy above conditions are the following:
\begin{align*}
  &\beta = \frac{1}{2},  C_1 = 80  \kappa  \eta_{\bx}G_g^2 L_f^4   \pare{\frac{4\eta_{\by}}{ \mu}+ 2\eta_{\by} ^2  }  + {\eta_{\bx}} G_g^2L^2_f,\\
  &C_2 = 16\kappa  \eta_{\bx}G_g^2 L_f^2, \ C_3 = 320\kappa^4  \eta_{\bx}G_g^2 L_f^2.
\end{align*}
We also need $C_3 = 320\kappa^4  \eta_{\bx}G_g^2 L_f^2 \geq \frac{1}{4\eta_{\bx}}$, which can be satisfied by choosing $\eta_{\bx} = \frac{1}{36 \kappa^2 L  }$.
    
    Hence we have:
       \begin{align*}
       \Psi(\bx^{t+1})   &\leq\Psi(\bx^{t}) - \frac{\eta_{\bx}}{2} \E \norm{\nabla \Phi(\bx^t)}^2  \\
       & \quad + \frac{(640\kappa^4  \eta_{\bx}^3G_g^2 L_f^2 + \eta_\bx^2 \kappa L)( G_f^2  + G_g^2)\sigma^2}{B} +\pare{{\eta_{\bx}} G_g^2L^2_f +C_1} \frac{\sigma^2}{B}+ 17\kappa  \eta_{\bx}G_g^2 L_f^2\frac{2\eta_{\by}^2 \sigma^2}{ B}  .
    \end{align*}

    Notice that $L \geq \max\cbr{ G_g L_f, L_f}$, and we choose $\eta_{\bx} \leq \frac{1}{L}$, we know $C_1 \leq 320 \kappa \eta_{\bx} \kappa L^2 + \eta_{\bx} L^2 \leq 321\kappa^2 L^2 \eta_{\bx} $. 
    Summing over $t=0$ to $T$ and re-arranging terms yield:
    \begin{align*}
      \frac{1}{T}\sum_{t=0}^{T-1} &\E \norm{\nabla \Phi(\bx^t)}^2 \leq \frac{2(\Psi(\bx^{0}) -  \Psi(\bx^{T}))}{\eta_{\bx} T}   + \frac{2(640\kappa^4  \eta_{\bx}^3G_g^2 L_f^2 + \eta_\bx^2 \kappa L)( G_f^2  + G_g^2)\sigma^2}{\eta_{\bx}B} + 642\kappa^2 L^2  \frac{\sigma^2}{B}+ 17\kappa  \frac{   \sigma^2}{  B}\\
      & \leq 2\frac{E[  \Phi(\bx^{0})] + 321\kappa^2 L^2 \eta_{\bx} \E\norm{ g(\bx^{-1})  -  \bz^{0}  }^2 + 16\kappa  \eta_{\bx} L^2 \E\norm{  \by^*(\bx^{-1}) -  \by^{-1}   }^2 + 320\kappa^4  \eta_{\bx}L^2 \E\norm{\bx^{0} - \bx^{-1}}^2 }{\eta_{\bx} T}  \\
      &\quad   + \frac{2(640\kappa^4  \eta_{\bx}^3G_g^2 L_f^2 + \eta_\bx^2 \kappa L)( G_f^2  + G_g^2)\sigma^2}{\eta_{\bx}B} + 642\kappa^2 L^2  \frac{\sigma^2}{B}+ 17\kappa  \frac{   \sigma^2}{  B}.
    \end{align*}
    By convention, we set $\bx^{-1} = \bx^0$, and with our choice $\eta_{\bx} = \frac{1}{36 \kappa^2 L  }$  we have
    \begin{align*}
      \frac{1}{T}\sum_{t=0}^{T-1} \E \norm{\nabla \Phi(\bx^t)}^2  
      & \leq O\pare{\frac{E[  \Phi(\bx^{0}) - \min_{\bx}\Phi(\bx)] +   \kappa    \E\norm{ g(\bx^{0})  -  \bz^{0}  }^2 +      L  D_{\cY}^2   }{\eta_{\bx} T} } \\
      &\quad + O\pare{\frac{ ( \kappa^4  \eta_{\bx}^3 L^2 + \eta_\bx^2 \kappa L)( G_f^2  + G_g^2)\sigma^2}{\eta_{\bx}B} +  \kappa^2 L  \frac{\sigma^2}{B}+  \kappa  \frac{   \sigma^2}{  B}}.
    \end{align*}
    Due to our initialization, $\E\norm{ g(\bx^{0})  -  \bz^{0}  }^2 \leq  \frac{1}{\kappa}$, and this can be done by $O(\kappa^2)$ zero order sampling.
    We choose  $B = \Theta\pr{\max\cbr{\frac{\kappa^2 L   \sigma^2}{\epsilon^2},1 }}$,   $T = O\pr{\max\cbr{\frac{\kappa^2 L\Delta_{\Phi}}{\epsilon^2},\frac{\kappa^2 L^2 D_{\cY}}{\epsilon^2}}  }$, which yields final gradient complexity of 
    \begin{align*}
        O\pr{\frac{\kappa^2( \Delta_{\Phi} +  L^2 D_{\cY} )}{\epsilon^2}\max\cbr{\frac{\kappa^2 L   \sigma^2}{\epsilon^2},1 } },
    \end{align*}
as desired.
\end{proof}
 
\subsection{Proof of Nonconvex-concave setting}\label{app: proof NCC}

\paragraph{Proof Sketch:} Before delving into the detailed proof, we first provide a sketch.
In nonconvex-concave setting, we consider the following potential function:
\begin{align}
    \Psi^t = \Phi_{1/2L} (\bx^t) + C_1 \E\norm{\bz^t - g(\bx^{t-1})}^2 + C_2 \E\norm{\bz^t - g(\bx^{t-1})}
\end{align}
where $C_1 =\frac{\eta_{\bx} }{2(2\beta-\beta^2)}G_g^2 L_f^2 $ and $C_2 = \frac{2\eta_\bx L^2 D_{\by}}{\beta}$. Another key step in nonconvex-concave analysis is to bound primal function gap. During the dynamic of our algorithm, the following statement holds:
\begin{align*}
    &\mathbb{E}\left(\Phi(\bx^{t-1}) - F(\bx^{t-1},\by^{t-1})\right)  \\
    \leq &(2t-1-2s) \eta_{\bx}\sqrt{2G_h^2+ 2G_f^2G_g^2}+ \frac{1}{2\eta_{\by}}( \E \|\by^*(\bx^s) - \by^{t-1}  \|^2 -  \E\|\by^*(\bx^s) - \by^t\|^2 ) +  \eta_{\by} \frac{\sigma^2}{M}  \\
    & +   L D_{\cY} \E\| \bz^{t}  -g(\bx^{t-1 })\|+\E\left[F(\bx^{t},\by^{t}) - F(\bx^{t-1},\by^{t-1})\right].
\end{align*}

Then we can derive that our algorithm satisfies the following descent inequality:
 \begin{align*}
    \E[\Psi^t - \Psi^{t-1}]  &\leq   2\eta_{\bx} L \left( \frac{ \mathbb{E}\|\by^*(\bx^s) - \by^{t-1}  \|^2 -   \mathbb{E}\|\by^*(\bx^s) - \by^t\|^2 }{2\eta_{\by}}  +  \eta_{\by} \frac{\sigma^2}{M}  + \mathbb{E}[F(\bx^{t},\by^{t}) - F(\bx^{t-1},\by^{t-1}) ]\right)  \\
    &\quad- \frac{\eta_{\bx} }{8} \mathbb{E} \|\nabla \Phi_{1/2L}(\bx^{t-1})\|^2+  2L\eta_{\bx}^2(G_h^2+ G_f^2G_g^2) + 4\eta_{\bx} L(2t-1-2s) \eta_{\bx} (G_h+ G_f G_g)\\
    &\quad + C_1\pare{ 8(1-\beta)^2 G_g^2  \eta_{\bx}^2(G_h^2 + G_f^2G_g^2) + 2\beta^2 \frac{\sigma^2}{B}}   + C_2 \sqrt{   8(1-\beta)^2 G_g^2  \eta_{\bx}^2(G_h^2 + G_f^2G_g^2) + 2\beta^2 \frac{\sigma^2}{B} }.
\end{align*}

\begin{lemma}
\label{lemma: coda-primal descent}
For CODA-Primal (Algorithm~\ref{algorithm: CODA-Primal}), under same assumptions as in Theorem~\ref{thm:ncc}, the following statement holds :
\begin{align*} 
    \mathbb{E}[\Phi_{1/2L}(\bx^t)]  
    & \leq \mathbb{E}[\Phi_{1/2L}({\bx}^{t-1})] + 2\eta_{\bx} L \mathbb{E}\left(\Phi(\bx^{t-1}) - F(\bx^{t-1},\by^{t-1})\right) - \frac{\eta_{\bx} }{8} \mathbb{E} \|\nabla \Phi_{1/2L}(\bx^{t-1})\|^2 \\
    &\quad+  L  \eta_{\bx}^2 G_f^2G_g^2 + \frac{\eta_{\bx} }{2}G_g^2 L_f^2\mathbb{E} \| \bz^t - g(\bx^{t-1})\|^2. 
\end{align*}
 \begin{proof}
 Let $\hat{\bx}^{t-1} = \arg\min_{\bx\in \mathbb{R}^d} \Phi(\bx)+L \|\bx-\bx^{t-1}\|^2$. Notice that:
\begin{align*}
   \mathbb{E}[\Phi_{1/2L}(\bx^t)] &\le \mathbb{E}[\Phi_{1/2L}(\hat{\bx}^{t-1})]+L\mathbb{E}\|\hat{\bx}^{t-1}-\bx^t\|^2\\
    &\leq \mathbb{E}[\Phi_{1/2L}(\hat{\bx}^{t-1})] +L ( \mathbb{E}\|\bx^{t-1}-\hat{\bx}^{t-1} \|^2 + 2\eta_{\bx}\langle \bg_{\bx}^t, \hat{\bx}^{t-1}-\bx^{t-1}  \rangle   +  \eta_{\bx}^2 \E\norm{\bg_{\bx}^t}^2 ).
\end{align*}
Due to the boundedness of gradients, we have
\begin{align*}
    \mathbb{E}[\Phi_{1/2L}(\bx^t)]  =  & \mathbb{E}[\Phi_{1/2L}(\hat{\bx}^{t-1})]\\
   & +L ( \mathbb{E}\|\bx^{t-1}-\hat{\bx}^{t-1} \|^2 + 2\eta_{\bx}\langle \nabla_1 f(\bz^t,\by^{t-1}) \nabla g(\bx^{t-1}) + \nabla h(\bx^{t-1}), \hat{\bx}^{t-1}-\bx^{t-1}\rangle   + 2\eta_{\bx}^2(G_h^2+ G_f^2G_g^2) )\\
     = &  \mathbb{E}[\Phi_{1/2L}(\hat{\bx}^{t-1})] +L ( \mathbb{E}\|\bx^{t-1}-\hat{\bx}^{t-1} \|^2 + 2\eta_{\bx}\langle   \nabla_x F(\bx^{t-1},\by^{t-1}),\hat{\bx}^{t-1}-\bx^{t-1} \rangle   + 2\eta_{\bx}^2(G_h^2+ G_f^2G_g^2) )\\
     &\quad + 2\eta_{\bx}L \langle \nabla_1 f(\bz^t,\by^{t-1}) \nabla g(\bx^{t-1}) - \nabla_1 f(\bx^{t-1},\by^{t-1}) \nabla g(\bx^{t-1}), \hat{\bx}^{t-1}-\bx^{t-1} \rangle .
\end{align*}
According to smoothness of $F(\cdot,\by)$, we have:
\begin{align*}
    \langle \hat{\bx}^{t-1} - {\bx}_{t-1}, \nabla_x F(\bx^{t-1},\by^{t-1})  \rangle &\leq F(\hat{\bx}^{t-1},\by^{t-1}) - F(\bx^{t-1},\by^{t-1}) + \frac{L}{2}\|\hat{\bx}^{t-1}-\bx^{t-1}\|^2 \\
    & = F(\hat{\bx}^{t-1},\by^{t-1})+ L\|\hat{\bx}^{t-1}-\bx^{t-1}\|^2 - F(\bx^{t-1},\by^{t-1}) - \frac{L}{2}\|\hat{\bx}^{t-1}-\bx^{t-1}\|^2 \\
    & \leq \Phi(\hat{\bx}^{t-1})+ L\|\hat{\bx}^{t-1}-\bx^{t-1}\|^2 - F(\bx^{t-1},\by^{t-1}) - \frac{L}{2}\|\hat{\bx}^{t-1}-\bx^{t-1}\|^2\\
    & \leq \Phi(\bx^{t-1}) - F(\bx^{t-1},\by^{t-1}) - \frac{L}{2}\|\hat{\bx}^{t-1}-\bx^{t-1}\|^2,
\end{align*}
where last step is due to that $\hat{\bx}^{t-1}$ is the minimizer of $\Phi(\cdot) + L\norm{\cdot - \bx^{t-1}}^2$.
So we have
\begin{align*}
     \mathbb{E}[\Phi_{1/2L}(\bx^t)] &\le  \mathbb{E}[\Phi_{1/2L}(\hat{\bx}^{t-1})]+L \mathbb{E}\|\hat{\bx}^{t-1}-\bx^t\|^2\\
    &\leq  \mathbb{E}[\Phi_{1/2L}(\hat{\bx}^{t-1})]+L  \mathbb{E}\|\bx^{t-1}-\hat{\bx}^{t-1} \|^2 \\
    &\quad+ 2\eta_{\bx} L  \mathbb{E}\left(\Phi(\bx^{t-1}) - f(\bx^{t-1},\by^{t-1}) - \frac{L}{2} \mathbb{E}\|\hat{\bx}^{t-1}-\bx^{t-1}\|^2\right)+2L\eta_{\bx}^2(G_h^2+ G_f^2G_g^2) \\
    &\quad+ \eta_{\bx} L \left(\frac{1}{2L}  \mathbb{E}\|\nabla_1 f(\bz^t,\by^{t-1}) \nabla g(\bx^{t-1}) - \nabla_1 f(\bx^{t-1},\by^{t-1}) \nabla g(\bx^{t-1})\|^2 +\frac{L}{2} \mathbb{E} \| \bx^{t-1}-\hat{\bx}^{t-1} \|^2   \right )\\
    & \leq  \mathbb{E}[\Phi_{1/2L}({\bx}^{t-1})] + 2\eta_{\bx} L \mathbb{E}\left(\Phi(\bx^{t-1}) - f(\bx^{t-1},\by^{t-1})\right) - \frac{\eta_{\bx}L^2}{2} \mathbb{E} \|\hat{\bx}^{t-1}-\bx^{t-1}\|^2 \\
    &\quad+  2L\eta_{\bx}^2(G_h^2+ G_f^2G_g^2) + \frac{\eta_{\bx} }{2}G_g^2 L_f^2\mathbb{E} \| \bz^t - g(\bx^{t-1})\|^2.  \\
\end{align*}
Finally using the fact that $\norm{\hat{\bx}^{t-1} - \bx^{t-1}} = \norm{\nabla \Phi_{1/2L}(\bx^{t-1})}/2L$ concludes the proof.
\end{proof}
\end{lemma}

\begin{lemma}\label{lemma: coda-primal dual descent}
For CODA-Primal (Algorithm~\ref{algorithm: CODA-Primal}), under the same assumptions made in Theorem~\ref{thm:ncc}, if we choose $\eta \leq 1/4L$ the following statement holds  for any $\by \in \mathcal{Y}$:
 \begin{align*}
 \E[F(\bx^{t-1},\by)- F(\bx^{t-1},\by^{t}) ] & \leq \frac{1}{2\eta_{\by}}( \E \|\by - \by^{t-1}  \|^2 -  \E\|\by - \by^t\|^2 )- \pare{\frac{1}{4\eta_{\by}}- \frac{L}{2}  }  \E\|\by^t - \by^{t-1}\|^2+  \eta_{\by}  {\sigma^2} \\
  &\quad +   L_f D_{\cY} \E\| \bz^{t}  -g(\bx^{t-1 })\|. 
\end{align*}

\begin{proof} 
According to updating rule of $\by$:
\begin{align*}
    \by^t &= \mathcal{P}_{\mathcal{Y}}\left(\by^{t-1}+ \eta_{\by} \bg_y^{t-1}\right) = \mathcal{P}_{\mathcal{Y}}\left(\by^{t-1}+ \eta_{\by} ( \nabla_2 f(\bz^{t },\by^{t-1};\zeta^{t-1}) - \nabla r(\by^{t-1})  )\right).
\end{align*}
We define
\begin{align*} 
    {\varepsilon}_{t-1} &=   ( \nabla_2 f(\bz^{t },\by^{t-1};\zeta^{t-1}) - \nabla r(\by^{t-1})  ) -   \nabla_{\by} F(\bx^{t-1},\by^{t-1})\\
  & =  \nabla_2 f(\bz^{t },\by^{t-1};\zeta^{t-1}) - \nabla_2 f(\bx^{t-1 },\by^{t-1}),
\end{align*}
 and re-write the updating rule as:
\begin{align*}
    \by^t = \mathcal{P}_{\mathcal{Y}}\left(\by^{t-1}+ \eta_{\by} \nabla_{\by} F(\bx^{t-1},\by^{t-1})+ \eta{\varepsilon}_{t-1}\right)  .
\end{align*}

Due to the property of projection we have:
\begin{align*}
    (\by-\by^t)^\top ( \by^t - \by^{t-1} - \eta_{\by} \nabla_{\by} F(\bx^{t-1},\by^{t-1})-\eta_{\by}  {\varepsilon}_{t-1} ) \geq 0.
\end{align*}
 Using the identity that $\langle \bm{a}, \bm{b} \rangle = \frac{1}{2}(\|\bm{a}+\bm{b}\|^2 - \|\bm{a}\|^2 - \|\bm{b}\|^2)$ we have:
\begin{align*}
   0 &\leq  \|\by - \by^{t-1} - \eta_{\by}\nabla_{\by} F(\bx^{t-1},\by^{t-1})-\eta_{\by} {\varepsilon}_{t-1} \|^2 - \|\by - \by^t\|^2 -  \|\by^t - \by^{t-1}- \eta_{\by} \nabla_{\by} F(\bx^{t-1},\by^{t-1})- \eta_{\by}{\varepsilon}_{t-1}\|^2   \\
   & =   \|\by - \by^{t-1}  \|^2 -  \|\by - \by^t\|^2 -  \|\by^t - \by^{t-1}\|^2 + 2\langle \by^{t-1} - \by,\eta_{\by} \nabla_{\by} F(\bx^{t-1}, \by^{t-1}) \rangle \\
  &\quad  + 2\langle \by^{t} - \by^{t-1},\eta_{\by} \nabla_{\by} F(\bx^{t-1}, \by^{t-1}) \rangle
   -2 \eta_{\by}\langle \by - \by^{t-1}, {\varepsilon}_{t-1} \rangle + 2\eta_{\by} \langle \by^t - \by^{t-1}, {\varepsilon}_{t-1} \rangle . 
\end{align*}
Due to $L$ smoothness and concavity, we have
\begin{align*}
    &2\langle \by^{t-1} - \by,\eta_{\by} \nabla_{\by} F(\bx^{t-1}, \by^{t-1}) \rangle+ 2\langle \by^{t} - \by^{t-1},\eta_{\by} \nabla_{\by} F(\bx^{t-1}, \by^{t-1}) \rangle \\
    \leq & 2\eta_{\by}(  f(\bx^{t-1},\by^{t}) -f(\bx^{t-1},\by)) +  {\eta_{\by}L}  \norm{\by^{t} - \by^{t-1}}^2.
\end{align*}
So we have:
\begin{align*}
 f(\bx^{t-1},\by)- f(\bx^{t-1},\by^{t})  & \leq \frac{1}{2\eta_{\by}}(  \|\by - \by^{t-1}  \|^2 -  \|\by - \by^t\|^2 )- \pare{\frac{1}{2\eta_{\by}}- \frac{L}{2} } \|\by^t - \by^{t-1}\|^2\\
  &\quad - \langle \by - \by^{t-1}, {\varepsilon}_{t-1} \rangle + \langle \by^t - \by^{t-1}, {\varepsilon}_{t-1} \rangle \\
  & \leq \frac{1}{2\eta_{\by}}(  \|\by - \by^{t-1}  \|^2 -  \|\by - \by^t\|^2 )- \pare{\frac{1}{2\eta_{\by}}- \frac{L}{2} } \|\by^t - \by^{t-1}\|^2\\
  &\quad -  \langle \by - \by^{t-1}, \nabla_2 f(\bz^{t },\by^{t-1};\zeta^{t-1}) - \nabla_2 f(\bx^{t-1 },\by^{t-1}) \rangle \\
  &\quad +  \langle \by^t - \by^{t-1}, \nabla_2 f(\bz^{t },\by^{t-1};\zeta^{t-1}) - \nabla_2 f(\bx^{t-1 },\by^{t-1}) \rangle\\
  & \leq \frac{1}{2\eta_{\by}}(  \|\by - \by^{t-1}  \|^2 -  \|\by - \by^t\|^2 )- \pare{\frac{1}{2\eta_{\by}}- \frac{L}{2} } \|\by^t - \by^{t-1}\|^2\\
  &\quad -  \langle \by - \by^{t-1}, \nabla_2 f(\bz^{t },\by^{t-1};\zeta^{t-1}) - \nabla_2 f(\bx^{t-1 },\by^{t-1}) \rangle \\
  &\quad +   \langle \by^t - \by^{t-1}, \nabla_2 f(\bz^{t },\by^{t-1};\zeta^{t-1}) - \nabla_2 f(\bz^{t },\by^{t-1} ) \rangle\\
   &\quad +  \langle \by^t - \by^{t-1},  \nabla_2 f(\bz^{t },\by^{t-1} )  - \nabla_2 f(\bx^{t-1 },\by^{t-1}) \rangle.
\end{align*}
Now we take expectation over the randomness of $\zeta^{t-1}$
\begin{align*}
\E[f(\bx^{t-1},\by)- f(\bx^{t-1},\by^{t})]  
  & \leq \frac{1}{2\eta_{\by}}( \E \|\by - \by^{t-1}  \|^2 -  \E\|\by - \by^t\|^2 )- \pare{\frac{1}{2\eta_{\by}}- \frac{L}{2} }  \E\|\by^t - \by^{t-1}\|^2\\
  &\quad -  \E\langle \by - \by^{t-1}, \nabla_2 f(\bz^{t },\by^{t-1} ) - \nabla_2 f(\bx^{t-1 },\by^{t-1}) \rangle \\
  &\quad +   \pare{ \frac{1}{4\eta_{\by}}\E \norm{\by^t - \by^{t-1}}^2 + \eta_{\by} \E\norm{\nabla_2 f(\bz^{t },\by^{t-1};\zeta^{t-1}) - \nabla_2 f(\bz^{t },\by^{t-1} )}^2}\\
   &\quad +  \E \langle \by^t - \by^{t-1},  \nabla_2 f(\bz^{t },\by^{t-1} )  - \nabla_2 f(\bx^{t-1 },\by^{t-1}) \rangle\\
    & \leq \frac{1}{2\eta_{\by}}( \E \|\by - \by^{t-1}  \|^2 -  \E\|\by - \by^t\|^2 )- \pare{\frac{1}{4\eta_{\by}}- \frac{L}{2}  }  \E\|\by^t - \by^{t-1}\|^2+  \eta_{\by}  {\sigma^2} \\
  &\quad +   L_f D_{\cY} \E\| \bz^{t}  -g(\bx^{t-1 })\|.  
\end{align*}

\end{proof}

\end{lemma}

\begin{lemma}\label{lemma: potential NC-C}
Define potential function:
\begin{align}
    \Psi^t = \Phi_{1/2L} (\bx^t) + C_1 \E\norm{\bz^t - g(\bx^{t-1})}^2 + C_2 \E\norm{\bz^t - g(\bx^{t-1})}
\end{align}
where $C_1 =\frac{\eta_{\bx} }{2(2\beta-\beta^2)}G_g^2 L_f^2 $ and $C_2 = \frac{2\eta_\bx L^2 D_{\by}}{\beta}$. Then the following statement holds:
 \begin{align*}
    \E[\Psi^t - \Psi^{t-1}]  &\leq   2\eta_{\bx} L \left( \frac{ \mathbb{E}\|\by^*(\bx^s) - \by^{t-1}  \|^2 -   \mathbb{E}\|\by^*(\bx^s) - \by^t\|^2 }{2\eta_{\by}}  +  \eta_{\by} \frac{\sigma^2}{M}  + \mathbb{E}[F(\bx^{t},\by^{t}) - F(\bx^{t-1},\by^{t-1}) ]\right)  \\
    &\quad- \frac{\eta_{\bx} }{8} \mathbb{E} \|\nabla \Phi_{1/2L}(\bx^{t-1})\|^2+  2L\eta_{\bx}^2(G_h^2+ G_f^2G_g^2) + 4\eta_{\bx} L(2t-1-2s) \eta_{\bx} (G_h+ G_f G_g)\\
    &\quad + C_1\pare{ 8(1-\beta)^2 G_g^2  \eta_{\bx}^2(G_h^2 + G_f^2G_g^2) + 2\beta^2 \frac{\sigma^2}{B}}   + C_2 \sqrt{   8(1-\beta)^2 G_g^2  \eta_{\bx}^2(G_h^2 + G_f^2G_g^2) + 2\beta^2 \frac{\sigma^2}{B} }.
\end{align*}
    \begin{proof}
    Recalling Lemma~\ref{lemma: coda-primal descent}, we have
        \begin{align*} 
    \mathbb{E}[\Phi_{1/2L}(\bx^t)]  
    & \leq \mathbb{E}[\Phi_{1/2L}({\bx}^{t-1})] + 2\eta_{\bx} L \mathbb{E}\left(\Phi(\bx^{t-1}) - F(\bx^{t-1},\by^{t-1})\right) - \frac{\eta_{\bx} }{8} \mathbb{E} \|\nabla \Phi_{1/2L}(\bx^{t-1})\|^2\\
    &\quad+  2L\eta_{\bx}^2(G_h^2+ G_f^2G_g^2) + \frac{\eta_{\bx} }{2}G_g^2 L_f^2\mathbb{E} \| \bz^t - g(\bx^{t-1})\|^2. 
\end{align*}
Following~\cite{lin2019gradient}, we split $\mathbb{E}\left(\Phi(\bx^{t-1}) - F(\bx^{t-1},\by^{t-1})\right) $ as follows:
\begin{align*}
    &\mathbb{E}\left(\Phi(\bx^{t-1}) - F(\bx^{t-1},\by^{t-1})\right)  \\
    = &\mathbb{E}\left[\underbrace{F(\bx^{t-1},\by^*(\bx^{t-1})) -F(\bx^{t-1},\by^*(\bx^{s}))}_{A} +\underbrace{F(\bx^{t-1},\by^*(\bx^{s}))-F(\bx^{t-1},\by^t)}_{B} +\underbrace{F(\bx^{t-1},\by^{t})-F(\bx^{t},\by^{t})}_{C}\right]\\
    &+\E\left[F(\bx^{t},\by^{t}) - F(\bx^{t-1},\by^{t-1})\right]
\end{align*}
for some iteration $s$.

For A:
\begin{align*}
    F(\bx^{t-1},\by^*(\bx^{t-1})) -F(\bx^{t-1},\by^*(\bx^{s})) &\leq F(\bx^{t-1},\by^*(\bx^{t-1})) - F(\bx^{s},\by^*(\bx^{t-1}))+ F(\bx^{s},\by^*(\bx^{s}))-F(\bx^{t-1},\by^*(\bx^{s}))\\
    &\leq 2(t-1-s) \eta_{\bx}\sqrt{2G_h^2+ 2G_f^2G_g^2}.
\end{align*}
For B, we evoke Lemma~\ref{lemma: coda-primal dual descent}:
\begin{align*}
    \E[F(\bx^{t-1},\by^*(\bx^{s}))-F(\bx^{t-1},\by^t)] \leq   \frac{1}{2\eta_{\by}}( \E \|\by - \by^{t-1}  \|^2 -  \E\|\by - \by^t\|^2 ) +  \eta_{\by} \frac{\sigma^2}{M}  +   L D_{\cY} \E\| \bz^{t}  -g(\bx^{t-1 })\|.
\end{align*}
For C:
\begin{align*}
    F(\bx^{t-1},\by^t) -F(\bx^{t},\by^t) \leq   \eta_{\bx}\sqrt{2G_h^2+ 2G_f^2G_g^2}.
\end{align*}
Putting pieces together yields:
\begin{align*}
    &\mathbb{E}\left(\Phi(\bx^{t-1}) - F(\bx^{t-1},\by^{t-1})\right)  \\
    \leq &(2t-1-2s) \eta_{\bx}\sqrt{2G_h^2+ 2G_f^2G_g^2}+ \frac{1}{2\eta_{\by}}( \E \|\by^*(\bx^s) - \by^{t-1}  \|^2 -  \E\|\by^*(\bx^s) - \by^t\|^2 ) +  \eta_{\by} \frac{\sigma^2}{M}  +   L D_{\cY} \E\| \bz^{t}  -g(\bx^{t-1 })\| \\
    &+\E\left[F(\bx^{t},\by^{t}) - F(\bx^{t-1},\by^{t-1})\right].
\end{align*}
Notice that $\sqrt{2G_h^2+ 2G_f^2G_g^2} \leq 2(G_h+ G_f G_g)$
Plugging above bound back to Lemma~\ref{lemma: coda-primal descent} yields:

{\small \begin{align*} 
   &  \mathbb{E}[\Phi_{1/2L}(\bx^t)]  
    \leq \mathbb{E}[\Phi_{1/2L}({\bx}^{t-1})] \\
    &+ 2\eta_{\bx} L \mathbb{E}\left( (2t-1-2s) \eta_{\bx}2(G_h+ G_f G_g)+ \frac{ \E \|\by^*(\bx^s) - \by^{t-1}  \|^2 -  \E\|\by^*(\bx^s) - \by^t\|^2 }{2\eta_{\by}}  +  \eta_{\by} \frac{\sigma^2}{B}  +   L D_{\cY} \E\| \bz^{t}  -g(\bx^{t-1 })\| \right. \\
    &\qquad \qquad \left.+\E\left[F(\bx^{t},\by^{t}) - F(\bx^{t-1},\by^{t-1})\right]\right)  \\
    &\quad- \frac{\eta_{\bx} }{8} \mathbb{E} \|\nabla \Phi_{1/2L}(\bx^{t-1})\|^2+  2L\eta_{\bx}^2(G_h^2+ G_f^2G_g^2) + \frac{\eta_{\bx} }{2}G_g^2 L_f^2\mathbb{E} \| \bz^t - g(\bx^{t-1})\|^2\\
    &\leq \mathbb{E}[\Phi_{1/2L}({\bx}^{t-1})]+ \frac{\eta_{\bx} }{2}G_g^2 L_f^2\mathbb{E} \| \bz^t - g(\bx^{t-1})\|^2   +  2\eta_{\bx} L^2 D_{\cY} \sqrt{\E\| \bz^{t}  -g(\bx^{t-1 })\|^2 } \\
    &+ 2\eta_{\bx} L \mathbb{E}\left( (2t-1-2s) \eta_{\bx}2(G_h+ G_f G_g)+ \frac{   \|\by^*(\bx^s) - \by^{t-1}  \|^2 -   \|\by^*(\bx^s) - \by^t\|^2 }{2\eta_{\by}}  +  \eta_{\by} \frac{\sigma^2}{M}  + F(\bx^{t},\by^{t}) - F(\bx^{t-1},\by^{t-1}) \right)  \\
    &\quad- \frac{\eta_{\bx} }{8} \mathbb{E} \|\nabla \Phi_{1/2L}(\bx^{t-1})\|^2+  2L\eta_{\bx}^2(G_h^2+ G_f^2G_g^2) .
\end{align*}}
Recalling the definition of $\Psi^t$
{\small\begin{align*}
    &\E[\Psi^t - \Psi^{t-1}] = \mathbb{E}[\Phi_{1/2L}(\bx^t) - \Phi_{1/2L}(\bx^{t-1}) ]\\
    &\quad + C_1 \mathbb{E} \| \bz^{t+1} - g(\bx^{t})\|^2 + C_2 \sqrt{\mathbb{E} \| \bz^{t+1} - g(\bx^{t})\|^2 } - \pare{C_1 \mathbb{E} \| \bz^{t} - g(\bx^{t-1})\|^2 + C_2 \sqrt{\mathbb{E} \| \bz^{t} - g(\bx^{t-1})\|^2 }}\\
    &\leq  \frac{\eta_{\bx} }{2}G_g^2 L_f^2\mathbb{E} \| \bz^t - g(\bx^{t-1})\|^2   +  2\eta_{\bx} L^2 D_{\cY} \sqrt{\E\| \bz^{t}  -g(\bx^{t-1 })\|^2 } \\
    &+ 2\eta_{\bx} L \mathbb{E}\left( (2t-1-2s) \eta_{\bx}2(G_h+ G_f G_g)+ \frac{   \|\by^*(\bx^s) - \by^{t-1}  \|^2 -   \|\by^*(\bx^s) - \by^t\|^2 }{2\eta_{\by}}  +  \eta_{\by}  \frac{\sigma^2}{M}  + F(\bx^{t},\by^{t}) - F(\bx^{t-1},\by^{t-1}) \right)  \\
    &\quad- \frac{\eta_{\bx} }{8} \mathbb{E} \|\nabla \Phi_{1/2L}(\bx^{t-1})\|^2+  2L\eta_{\bx}^2(G_h^2+ G_f^2G_g^2) \\
     &\quad + C_1 \mathbb{E} \| \bz^{t+1} - g(\bx^{t})\|^2 + C_2 \sqrt{\mathbb{E} \| \bz^{t+1} - g(\bx^{t})\|^2 } - \pare{C_1 \mathbb{E} \| \bz^{t} - g(\bx^{t-1})\|^2 + C_2 \sqrt{\mathbb{E} \| \bz^{t} - g(\bx^{t-1})\|^2 }}\\
    &\leq \pare{ \frac{\eta_{\bx} }{2}G_g^2 L_f^2 -C_1}\mathbb{E} \| \bz^t - g(\bx^{t-1})\|^2   +  \pare{2\eta_{\bx} L^2 D_{\cY}-C_2 }\sqrt{\E\| \bz^{t}  -g(\bx^{t-1 })\|^2 } \\
    &+ 2\eta_{\bx} L \mathbb{E}\left( (2t-1-2s) \eta_{\bx}2(G_h+ G_f G_g)+ \frac{   \|\by^*(\bx^s) - \by^{t-1}  \|^2 -   \|\by^*(\bx^s) - \by^t\|^2 }{2\eta_{\by}}  +  \eta_{\by} \frac{\sigma^2}{M}  + F(\bx^{t},\by^{t}) - F(\bx^{t-1},\by^{t-1}) \right)  \\
    &\quad- \frac{\eta_{\bx} }{8} \mathbb{E} \|\nabla \Phi_{1/2L}(\bx^{t-1})\|^2+  2L\eta_{\bx}^2(G_h^2+ G_f^2G_g^2)   + C_1 \mathbb{E} \| \bz^{t+1} - g(\bx^{t})\|^2 + C_2 \sqrt{\mathbb{E} \| \bz^{t+1} - g(\bx^{t})\|^2 } .
\end{align*}}
Evoking Lemma~\ref{lem:tracking error} to replace $\mathbb{E} \| \bz^{t+1} - g(\bx^{t})\|^2$ yields:
{\small\begin{align*}
     &\E[\Psi^t - \Psi^{t-1}] \leq \pare{ \frac{\eta_{\bx} }{2}G_g^2 L_f^2 -C_1}\mathbb{E} \| \bz^t - g(\bx^{t-1})\|^2   +  \pare{2\eta_{\bx} L^2 D_{\cY}-C_2 }\sqrt{\E\| \bz^{t}  -g(\bx^{t-1 })\|^2 } \\
    &\ + 2\eta_{\bx} L \mathbb{E}\left( (2t-1-2s) \eta_{\bx}2(G_h+ G_f G_g)+ \frac{   \|\by^*(\bx^s) - \by^{t-1}  \|^2 -   \|\by^*(\bx^s) - \by^t\|^2 }{2\eta_{\by}}  +  \eta_{\by} \frac{\sigma^2}{M}  + F(\bx^{t},\by^{t}) - F(\bx^{t-1},\by^{t-1}) \right)  \\
    &\quad- \frac{\eta_{\bx} }{8} \mathbb{E} \|\nabla \Phi_{1/2L}(\bx^{t-1})\|^2+  2L\eta_{\bx}^2(G_h^2+ G_f^2G_g^2) \\
    &\quad + C_1\pare{ (1-\beta)^2 \E \norm{\bz^{t} - g(\by^{t-1}) }^2  + 4(1-\beta)^2 G_g^2 \norm{\bx^t - \bx^{t-1}}^2 + 2\beta^2 \frac{\sigma^2}{B}} \\
    &\quad + C_2 \sqrt{ (1-\beta)^2 \E \norm{\bz^{t} - g(\by^{t-1}) }^2  + 4(1-\beta)^2 G_g^2 \norm{\bx^t - \bx^{t-1}}^2 + 2\beta^2 \frac{\sigma^2}{B} } \\
    &\leq \underbrace{\pare{ \frac{\eta_{\bx} }{2}G_g^2 L_f^2 -C_1 +(1-\beta)^2C_1} }_{\clubsuit}\mathbb{E} \| \bz^t - g(\bx^{t-1})\|^2   + \underbrace{ \pare{2\eta_{\bx} L^2 D_{\cY}-C_2 + (1-\beta) C_2 }}_{\heartsuit}\sqrt{\E\| \bz^{t}  -g(\bx^{t-1 })\|^2 } \\
    &\ + 2\eta_{\bx} L \mathbb{E}\left( (2t-1-2s) \eta_{\bx}2(G_h+ G_f G_g)+ \frac{   \|\by^*(\bx^s) - \by^{t-1}  \|^2 -   \|\by^*(\bx^s) - \by^t\|^2 }{2\eta_{\by}}  +  \eta_{\by}  \frac{\sigma^2}{B}  + F(\bx^{t},\by^{t}) - F(\bx^{t-1},\by^{t-1}) \right)  \\
    &\quad- \frac{\eta_{\bx} }{8} \mathbb{E} \|\nabla \Phi_{1/2L}(\bx^{t-1})\|^2+  2L\eta_{\bx}^2(G_h^2+ G_f^2G_g^2) \\
    &\quad + C_1\pare{ 4(1-\beta)^2 G_g^2 \norm{\bx^t - \bx^{t-1}}^2 + 2\beta^2 \frac{\sigma^2}{B}}   + C_2 \sqrt{   4(1-\beta)^2 G_g^2 \norm{\bx^t - \bx^{t-1}}^2 + 2\beta^2 \frac{\sigma^2}{B} }
\end{align*}}
Since we choose $C_1 =\frac{\eta_{\bx} }{2(2\beta-\beta^2)}G_g^2 L_f^2 $ and $C_2 = \frac{2\eta_\bx L^2 D_{\by}}{\beta}$, we know $\clubsuit = 0$ and $\heartsuit=0$. Hence we have
\begin{align*}
     &\E[\Psi^t - \Psi^{t-1}]  \\
    &\leq   2\eta_{\bx} L \left( \frac{ \mathbb{E}\|\by^*(\bx^s) - \by^{t-1}  \|^2 -   \mathbb{E}\|\by^*(\bx^s) - \by^t\|^2 }{2\eta_{\by}}  +  \eta_{\by}  \frac{\sigma^2}{B}  + \mathbb{E}[F(\bx^{t},\by^{t}) - F(\bx^{t-1},\by^{t-1}) ]\right)  \\
    &\quad- \frac{\eta_{\bx} }{8} \mathbb{E} \|\nabla \Phi_{1/2L}(\bx^{t-1})\|^2+  2L\eta_{\bx}^2(G_h^2+ G_f^2G_g^2) + 2\eta_{\bx} L(2t-1-2s) \eta_{\bx}(G_h+ G_f G_g)\\
    &\quad + C_1\pare{ 4(1-\beta)^2 G_g^2 2\eta_{\bx}^2(G_h^2 + G_f^2G_g^2) + 2\beta^2 \frac{\sigma^2}{B}}   + C_2 \sqrt{   4(1-\beta)^2 G_g^2 2\eta_{\bx}^2(G_h^2 + G_f^2G_g^2) + 2\beta^2 \frac{\sigma^2}{B} }
\end{align*}
which concludes the proof.
    \end{proof}
\end{lemma}

\subsubsection{Proof of Theorem~\ref{thm:ncc}}
We first evoke Lemma~\ref{lemma: potential NC-C}:
 
\begin{align*}
   \mathbb{E} \|\nabla \Phi_{1/2L}&(\bx^{t-1})\|^2\leq   \frac{8\E[\Psi^{t-1}-\Psi^t  ]}{\eta_{\bx}} \\
   & + 
    16 L \left( \frac{ \mathbb{E}\|\by^*(\bx^s) - \by^{t-1}  \|^2 -   \mathbb{E}\|\by^*(\bx^s) - \by^t\|^2 }{2\eta_{\by}}    + \mathbb{E}[F(\bx^{t},\by^{t}) - F(\bx^{t-1},\by^{t-1}) ]\right)  \\
    &+  16L\eta_{\bx} (G_h^2+ G_f^2G_g^2) + 16 L(2t-1-2s) \eta_{\bx}(G_h+ G_f G_g) + 16 \eta_{\by}L \frac{\sigma^2}{B}\\
    &+ \frac{8C_1}{\eta_{\bx}}\pare{ 8(1-\beta)^2 G_g^2  \eta_{\bx}^2  (G_h^2 + G_f^2G_g^2) + 2\beta^2 \frac{\sigma^2}{B}}   +\frac{8C_2}{\eta_{\bx}} \sqrt{   8(1-\beta)^2 G_g^2 \eta_{\bx}^2   (G_h^2 + G_f^2G_g^2) + 2\beta^2 \frac{\sigma^2}{B} }.
\end{align*}
We sum above inequality for $t = jS$ to $(j+1)S-1$, with $\bx^s = \bx^{jS}$:

\begin{align*}
 &\sum_{t=jS}^{(j+1)S-1}  \mathbb{E} \|\nabla \Phi_{1/2L}(\bx^{t})\|^2\\
 \leq  &\sum_{t=jS}^{(j+1)S-1} \left[\frac{8\E[\Psi^{t}-\Psi^{t+1}  ]}{\eta_{\bx}} + 
    16   L \left( \frac{ \mathbb{E}\|\by^*(\bx^{jS}) - \by^{t}  \|^2 -   \mathbb{E}\|\by^*(\bx^{jS}) - \by^{t+1}\|^2 }{2\eta_{\by}}    + \mathbb{E}[F(\bx^{t+1},\by^{t+1}) - F(\bx^{t},\by^{t}) ]\right) \right] \\
    &\quad +  16SL\eta_{\bx} (G_h^2+ G_f^2G_g^2) + \sum_{t=jS}^{(j+1)S-1} 16 L(2t+1-2jS) \eta_{\bx}(G_h+ G_f G_g) + 16S \eta_{\by}L \frac{\sigma^2}{B}\\
    &\quad + 8SC_1\pare{ 8(1-\beta)^2 G_f^2  \eta_{\bx}^2 (G_h^2 + G_f^2G_g^2) + 2\beta^2 \frac{\sigma^2}{B}}   + 8SC_2 \sqrt{   8(1-\beta)^2 G_f^2 \eta_{\bx}^2   (G_h^2 + G_f^2G_g^2) + 2\beta^2 \frac{\sigma^2}{B} }\\
     \leq  &  \frac{8\E[\Psi^{jS}-\Psi^{(j+1)S)}  ]}{\eta_{\bx}} \\
   &\quad +  16   L \left( \frac{ \mathbb{E}\|\by^*(\bx^{jS}) - \by^{jS}  \|^2 -   \mathbb{E}\|\by^*(\bx^{jS}) - \by^{(j+1)S}\|^2 }{2\eta_{\by}}    + \mathbb{E}[F(\bx^{(j+1)S},\by^{(j+1)S}) - F(\bx^{jS},\by^{jS}) ]\right)  \\
    &\quad +  16SL\eta_{\bx} (G_h^2+ G_f^2G_g^2) +   32S^2L  \eta_{\bx}(G_h+ G_f G_g) + 16S \eta_{\by}L \frac{\sigma^2}{B}\\
    &\quad + 8S\frac{C_1}{\eta_{\bx}}\pare{ 8(1-\beta)^2 G_f^2  \eta_{\bx}^2 (G_h^2 + G_f^2G_g^2) + 2\beta^2 \frac{\sigma^2}{B}}   + 8S\frac{C_2}{\eta_{\bx}} \sqrt{   8(1-\beta)^2 G_f^2 \eta_{\bx}^2   (G_h^2 + G_f^2G_g^2) + 2\beta^2 \frac{\sigma^2}{B} }.
\end{align*}

Now we further sum over $j = 0$ to $(T+1)/S-1$

\begin{align*}
    \frac{1}{T+1}\sum_{t=0}^{T} &\mathbb{E} \|\nabla \Phi_{1/2L}(\bx^{t})\|^2 =  \frac{1}{T+1}\sum_{j=0}^{(T+1)/S-1}  \sum_{t=jS}^{(j+1)S-1}\mathbb{E} \|\nabla \Phi_{1/2L}(\bx^{t})\|^2 \\
   & \leq   \frac{1}{T+1}  \frac{8\E[\Psi^{0}-\Psi^{T}  ]}{\eta_{\bx}} + 
     {16   L}  \left( \frac{ D^2_{\cY}}{2S\eta_{\by}}    +  \frac{\E[F(\bx^T,\by^T)- F(\bx^0, \by^0)]}{T+1} \right)  \\
    & +  16L\eta_{\bx} (G_h^2+ G_f^2G_g^2) +   32SL  \eta_{\bx}(G_h+ G_f G_g) + 16 \eta_{\by}L \frac{\sigma^2}{B}\\
    &  + 8\frac{C_1}{\eta_{\bx}}\pare{ 8(1-\beta)^2 G_g^2  \eta^2_{\bx} (G_h^2 + G_f^2G_g^2) + 2\beta^2 \frac{\sigma^2}{M}}   + 8\frac{C_2}{\eta_{\bx}} \sqrt{   8(1-\beta)^2 G_g^2 \eta_{\bx}^2   (G_h^2 + G_f^2G_g^2) + 2\beta^2 \frac{\sigma^2}{M} },
\end{align*}
where $C_1 =\frac{\eta_{\bx} }{2(2\beta-\beta^2)}G_g^2 L_f^2 $ and $C_2 = \frac{2\eta_\bx L^2 D_{\by}}{\beta}$. Notice that $\E[F(\bx^T,\by^T)- F(\bx^0, \by^0)] \leq \E[F(\bx^T,\by^T)-F(\bx^0,\by^T)+F(\bx^0,\by^T)- F(\bx^0, \by^0)] \leq T \eta_{\bx} (G_h + G_f G_g) + \Delta_{\Phi}$.
 
Now, We choose $S = \sqrt{\frac{  D^2_{\cY}}{4 \eta_{\bx}\eta_{\by} (G_h+ G_f G_g)}}$, $\eta_{\by} = \Theta\pr{  \frac{\epsilon^2}{L\sigma^2}  }$ and $\eta_{\bx} = \Theta \pr{ \frac{\epsilon^6}{ LD_{\cY} \sigma^2 (G_h + G_f G_g )}}$, $B = \Theta\pr{1}$, $M = \Theta\pr{\frac{\sigma^2 L^4 D^2_{\cY}}{\epsilon^4}}$ and $T =\Theta(\max\cbr{\frac{\Delta_{\hat\Phi}}{\eta_{\bx}\epsilon^2}, \frac{L\Delta_{\Phi}}{\epsilon^2}})$, which yields total gradient complexity of
\begin{align*}
    O\pare{  \max\left\{ \frac{ L^3 \sigma^2D^2_{\cY}(G_h + G_fG_g) \Delta_{\hat{\Phi}}}{\epsilon^8},  \frac{L\Delta_{\Phi}}{\epsilon^2} \right\} }.
\end{align*}

\section{Proof of Dual Composition Setting}\label{app: proof dual}
In this section we provide the proof of results in primal compostion setting (Theorem~\ref{thm:SCNC} and Theorem~\ref{thm: CNC}).

\paragraph{Proof Sketch:} Before delving into the detailed proof, we first provide a sketch. We construct the following two-level potential function:
\begin{align*}
    \hat{F}^{t+1} := & F(\bx^{t+1},\by^{t+1}) + s^{t+1} - \pare{\frac{1}{4\eta_{\by}} + 4 L_g^2 G_f^2 G_g^2 \eta_{\by}  + \frac{\eta_{\bx} L^2}{2} + \frac{96  G_g^2L_f^2}{\mu^2 \eta_{\bx} \beta} + \frac{48  G_g^2 L_f^2\beta}{\mu^2 \eta_{\bx}}} \norm{\by^{t+1} - \by^t}^2\\
    &  -\pare{\frac{1}{8\eta_{\by}}+\frac{48 G_g^2L_f^2}{\mu^2 \eta_{\bx} \beta}}\norm{\by^t - \by^{t-1}}^2 + \pare{\frac{7}{2\eta_{\bx}} + \mu - \frac{\eta_{\bx} L^2}{2} - \frac{2L^2_f}{\mu}}\norm{\bx^{t+1} - \bx^t}^2
\end{align*}
where $s^{t+1} :=  -\frac{2}{\eta^2_{\bx} \mu}\|\bx^{t+1}-\bx^t\|^2$ and
\begin{align*}
     \tilde{F}^{t+1} :=  \hat{F}^{t+1} -   \E \norm{\bz^{t+1}  -  g(\by^t))  }^2 -\frac{(1-\frac{\beta}{2})^2}{1-(1-\frac{\beta}{2})^2}\frac{4L_f^2}{\mu^2\eta_{\bx}}   \E \norm{  \bz^{t+1 } - \bz^{t}  }^2.
\end{align*}

$\E \norm{  \bz^{t+1 } - \bz^{t}  }^2$ is the iterate difference of our auxilliary variables, and it can be controlled by
  \begin{align*}
        \E \norm{  \bz^{t+1} - \bz^{t}  }^2 &\leq \pare{1-\frac{\beta}{2}}^2 \E \norm{  \bz^{t } - \bz^{t-1}  }^2+ 4\pare{1+\frac{2}{\beta}} G_g^2\pare{\norm{    \by^t - \by^{t-1}   }^2 +  \norm{    \by^{t-1} - \by^{t-2}   }^2 } \\
        & + 2\pare{1+\frac{2}{\beta}}\beta^2 G_g^2 \E \norm{    \by^t - \by^{t-1}   }^2.
    \end{align*}
Putting pieces together we can derive the following descent property:
 \begin{align*}  
  \E [ \tilde{F}^{t+1} - \tilde{F}^{t}]   & \geq  \min \cbr{\frac{1}{\eta_{\by}} , \frac{1}{\eta_{\bx}}  }\cdot \E \norm{\nabla G(\bx^t,\by^t)}^2   \\ 
    & \quad  -  O\pare{  \eta_{\by} (G_g^2+G_f^2)  +   \eta_{\bx} } \frac{\sigma^2}{B} -   O\pare{ \beta^2 L_f^2 G_g^2 {\eta_{\by}}+  \beta^2  \eta_{\by} L_f^2 G_g^2 + \beta^2 \eta_\bx L_f^2  }\frac{\sigma^2}{M} .   
\end{align*}
Unrolling above recursion will give desired rate. The converegcne of merely-convex-nonconcave setting follows similar technique. 

\subsection{Proof of Strongly-convex-nonconcave setting}\label{app: proof SCC}
\vspace{-0.25cm}

 

In nonconcave setting, besides the Lemma~\ref{lem:tracking error}, we also need the following bound on tracking error between two iterates.

\begin{lemma}[Convergence of Iterates difference] \label{lem:second order tracking}
For Algorithm~\ref{algorithm: CODA-Dual}, under the assumptions of Theorem~\ref{thm:SCNC}, the following statement holds true:
    \begin{align*}
        \E \norm{  \bz^{t+1} - \bz^{t}  }^2 &\leq \pare{1-\frac{\beta}{2}}^2 \E \norm{  \bz^{t } - \bz^{t-1}  }^2+ 4\pare{1+\frac{2}{\beta}} G_g^2\pare{\norm{    \by^t - \by^{t-1}   }^2 +  \norm{    \by^{t-1} - \by^{t-2}   }^2 } \\
        &\quad + \pare{1+\frac{2}{\beta}}\beta^2 \pare{3G_g^2\E \norm{    \by^t - \by^{t-1}   }^2 + 6 \frac{\sigma^2}{M}}.
    \end{align*}
    \begin{proof}
For the ease of presentation, we define the following two auxiliary variables:
\begin{align*}
&g^t  = \frac{1}{B}\sum_{\xi \in \cM^t} g(\by^t;\xi) , \ g^{t\mapsto t-1} =\frac{1}{M}\sum_{\xi \in  \cM^t} \pr{ g(\by^t;\xi)   - g(\by^{t-1};\xi) }.
\end{align*}
According to updating rule of $\bz$, we have:
        \begin{align*}
            \bz^{t+1} - \bz^{t} = (1-\beta) ( \bz^{t} - \bz^{t-1}) + (1-\beta) (g^{t\mapsto t-1}      -   g^{t-1\mapsto t-2}   )  + \beta (g^t - g^{t-1}).
        \end{align*}
       Taking expectation w.r.t. $\cM^t$, and $\cM^{t-1}$ yields:
        \begin{align*}
           \E\norm{ \bz^{t+1} - \bz^{t}}^2 &= \E\norm{  (1-\beta) ( \bz^{t} - \bz^{t-1}) + (1-\beta) ( g^{t\mapsto t-1}  - g^{t-1\mapsto t-2}   )  + \beta (g^t - g^{t-1})}^2\\
           &{\leq} \pare{1+\frac{\beta}{2-2\beta} }(1-\beta)^2\E\norm{   \bz^{t} - \bz^{t-1} +( g^{t\mapsto t-1}      -    g^{t-1\mapsto t-2}   ) }^2   +\pare{1+\frac{2-2\beta}{\beta}}\norm{  \beta (g^t - g^{t-1})}^2\\
           &{\leq} \pare{1-\frac{\beta}{2}}(1-\beta)\E\norm{  \bz^{t} - \bz^{t-1} +( g^{t\mapsto t-1}      -    g^{t-1\mapsto t-2}   ) }^2    +  \pare{1+\frac{2 }{\beta}}\beta^2 \E\norm{g^t - g^{t-1}}^2 \\
            &{\leq} \pare{1-\frac{\beta}{2}}(1-\beta) \pare{1+\frac{\beta}{2-2\beta} }\E\norm{  \bz^{t} - \bz^{t-1}   }^2   \\
            & \quad + \pare{1-\frac{\beta}{2}}(1-\beta) \pare{1+\frac{2-2\beta}{\beta} }\E\norm{   g^{t\mapsto t-1}      -    g^{t-1\mapsto t-2}  }^2   +   \pare{1+\frac{2 }{\beta}}\beta^2 \E\norm{g^t - g^{t-1}}^2 \\
            &{\leq} \pare{1-\frac{\beta}{2}}^2 \E\norm{  \bz^{t} - \bz^{t-1}   }^2   +   \underbrace{\pare{1+\frac{2}{\beta} }\E\norm{   g^{t\mapsto t-1}      -    g^{t-1\mapsto t-2}  }^2}_{T_1}   + \underbrace{ \pare{1+\frac{2 }{\beta}}\beta^2 \E\norm{g^t - g^{t-1}}^2}_{T_2}
        \end{align*}
        where in the first and third steps we use Young's inequality that $\norm{\ba+\bb}^2 \leq (1+a)\norm{\ba}^2 + (1+\frac{1}{a})\norm{\bb}^2  $. 
        Now we bound $T_1$ as follows:
        \begin{align*}
            T_1 \leq & 2\pare{1+\frac{2}{\beta} }\pare{\frac{1}{M}\sum_{\xi^{t} \in   \cM^{t}}\E\norm{  g(\by^t;\xi^t)      -  g(\by^{t-1};\xi^{t}) }^2} \\
            & + 2\pare{1+\frac{2}{\beta} }\pare{\frac{1}{M}\sum_{\xi^{t-1} \in \cM^{t-1}}\E\norm{ g(\by^{t-1};\xi^{t-1})      -  g(\by^{t-2};\xi^{t-1})}^2}\\
            \leq &4\pare{1+\frac{2 }{\beta}}G_g^2\pare{\E\norm{\by^t - \by^{t-1}   }^2 + \E \norm{    \by^{t-1} - \by^{t-2}   }^2 },
        \end{align*}
        and for $T_2$:
        \begin{align*}
            T_2 \leq  \pare{1+\frac{2}{\beta}}\beta^2 \pare{3G_g^2\E \norm{    \by^t - \by^{t-1}   }^2 + 6 \frac{\sigma^2}{M}}.
        \end{align*}
        Putting pieces together will conclude the proof.
    \end{proof}
\end{lemma}

The following Lemma establishes the difference of the exact gradient computed on $F$ versus the gradients we actually used in Algorithm~\ref{algorithm: CODA-Dual}. 
\begin{lemma}[Gradient difference]\label{lem:grad bound}
For Algorithm~\ref{algorithm: CODA-Dual}, under the assumptions of Theorem~\ref{thm:SCNC}, the following statement holds true:
\begin{align*}
&\E\norm{\bg^t_{\by}- \nabla_{\by} F(\bx^t,\by^t)}^2 \leq \frac{ 4(G_f^2+G_g^2 )\sigma^2}{B} + 2G_g^2L_f^2 \E\norm{ \bz^{t+1}  -  g(\by^t)    }^2.\\
&\E\norm{\bg^t_{\bx} - \nabla_{\bx} F(\bx^t,\by^t)}^2 \leq 2L_f^2 \E\norm{\bz^{t+1} -   g(\by^t)  }^2 + 2\frac{\sigma^2}{B}.
\end{align*}

\begin{proof}
    By the definition of $\bg_{\by}^t$, we have
    \begin{align*}
       \E \norm{\bg^t_{\by}- \nabla_{\by} F(\bx^t,\by^t)}^2 &= \E\norm{ \frac{1}{B}\sum_{(\zeta^t,\xi^t) \in \cB^t}  \nabla_2 f(\bx^t, \bz^{t+1};\zeta^t)  \nabla g(\by^t;\xi^t) - \nabla_{\by} F(\bx^t,\by^t)}^2\\
       & \leq 2\E\norm{\frac{1}{B}\sum_{(\zeta^t,\xi^t) \in \cB^t}  \nabla_2 f(\bx^t, \bz^{t+1};\zeta^t)  \nabla g(\by^t;\xi^t) - \nabla_2 f(\bx^t, \bz^{t+1})  \nabla g(\by^t ) }^2\\
       & \quad + 2\E\norm{\nabla_2 f(\bx^t, \bz^{t+1})  \nabla g(\by^t ) - \nabla_2 f(\bx^t, g(\by^{t}))  \nabla g(\by^t )   }^2\\
       & \leq   \frac{ 4(G_f^2+G_g^2 )\sigma^2}{B} + 2G_g^2L_f^2 \E\norm{ \bz^{t+1}  -  g(\by^t)    }^2,
    \end{align*}
    where at the last step we use the bounded variance assumption. Similarly by the definition of $\bg_{\bx}^t$ we have:
    \begin{align*}
      \E  \norm{\bg^t_{\bx}- \nabla_{\bx} F(\bx^t,\by^t)}^2 &= \E\norm{ \frac{1}{B}\sum_{(\zeta^t,\xi^t) \in \cB^t} \nabla_1 f(\bx^t,\bz^{t+1};\zeta^t)- \nabla_1 f(\bx^t,g(\by^{t}))}^2  \leq  2L_f^2 \E\norm{\bz^{t+1} -   g(\by^t)  }^2 + 2\frac{\sigma^2}{B}.
    \end{align*}

\end{proof}

\end{lemma}
 
\begin{lemma}[Connection between stationary measure and iterates] \label{lem:stationary measure}
For Algorithm~\ref{algorithm: CODA-Dual}, under the assumptions of Theorem~\ref{thm:SCNC}, if we define \[\hat{\bg}_{\bx}(\bx,\by) = \frac{1}{\eta_{\bx}}  \pare{ \bx - \cP_{\cX}\pare{\bx-\eta_{\bx} \nabla_{\bx} F(\bx,\by)}},
 \hat{\bg}_{\by}(\bx,\by)=  \frac{1}{\eta_{\by}}  \pare{ \by - \cP_{\cY}\pare{\by+\eta_{\by} \nabla_{\by} F(\bx,\by)}  }\], 
  then the following statement holds:
\begin{align*}
   &  \E\norm{\hat{\bg}_{\bx}(\bx^t,\by^t) }^2 \leq \frac{2}{\eta_{\bx}^2} \E\norm{  \bx^t - \bx^{t+1}}^2 +  4L_f^2 \E\norm{\bz^{t+1} - g(\by^t)}^2 + 4\frac{\sigma^2}{B},\\
    & \E\norm{\hat{\bg}_{\by}(\bx^t,\by^t) }^2 \leq \frac{2}{\eta_{\by}^2}  \E\norm{  \by^t - \by^{t+1}}^2 +4 L_f^2 G_g^2\E\norm{  \bz^{t+1}  - g(\by^t   )  }^2 +   \frac{8 (G_f^2+G_g^2)\sigma^2}{B}.
\end{align*}
    \begin{proof}
        We begin with proving the first statement. According to updating rule we have
        \begin{align*}
           \E \norm{\hat{\bg}_{\bx}(\bx^t,\by^t)) }^2 &= \E\norm{  \frac{1}{\eta_{\bx}}  \pare{ \bx^t - \cP_{\cX}\pare{\bx^t-\eta_{\bx} \nabla_{\bx} F(\bx^t,\by^t)} }  }^2 \\
            &\leq \frac{2}{\eta_{\bx}^2} \E\norm{  \bx^t - \cP_{\cX}\pare{\bx^t-\eta_{\bx} \bg_{\bx}^t }}^2 + \frac{2}{\eta_{\bx}^2} \E\norm{  \cP_{\cX}\pare{\bx^t-\eta \bg_{\bx}^t}- \cP_{\cX}\pare{\bx^t-\eta \nabla_{\bx} F(\bx^t,\by^t)}}^2\\
            &\leq \frac{2}{\eta_{\bx}^2} \E\norm{  \bx^t - \bx^{t+1}}^2 +  2\E\norm{     \bg_{\bx}^t-  \nabla_{\bx} F(\bx^t,\by^t)}^2\\
             &\leq \frac{2}{\eta_{\bx}^2} \E\norm{  \bx^t - \bx^{t+1}}^2 +  4\E\norm{ \nabla_1 f(\bx^t, \bz^{t+1};\zeta^t)-  \nabla_1 f(\bx^t, \bz^{t+1})}^2 \\
             &\quad +  4\E\norm{  \nabla_1 f(\bx^t, \bz^{t+1})-  \nabla_1 f(\bx^t, g(\by^{t}))}^2\\
                &\leq \frac{2}{\eta_{\bx}^2} \E\norm{  \bx^t - \bx^{t+1} }^2 +  4L_f^2 \E\norm{\bz^{t+1} - g(\by^t)}^2 + 4\frac{\sigma^2}{B},
        \end{align*}
        where at last step we apply Lemma~\ref{lem:grad bound}.
        Similarly, for the second statement we have:
           \begin{align*}
            \E\norm{\hat{\bg}_{\by}(\bx^t,\by^t)}^2 &=  \E\norm{  \frac{1}{\eta_{\by}}  \pare{ \by^t - \cP_{\cY}\pare{\by^t+\eta \nabla_{\by} F(\bx^t,\by^t)}}  }^2 \\
            &\leq \frac{2}{\eta_{\by}^2}  \E\norm{  \by^t - \cP_{\cY}\pare{\by^t+\eta_{\by} \bg_{\by}^t }}^2 + \frac{2}{\eta_{\by}^2}  \E\norm{  \cP_{\cY}\pare{\by^t+\eta_{\by} \bg_{\by}^t}- \cP_{\cY}\pare{\by^t+\eta_{\by} \nabla_{\by} F(\bx^t,\by^t)}}^2\\
            &\leq \frac{2}{\eta_{\by}^2}  \E\norm{  \by^t - \by^{t+1}} +  2 \E\norm{     \bg_{\by}^t-  \nabla_{\by} F(\bx^t,\by^t)}^2\\
            &= \frac{2}{\eta_{\by}^2}  \E\norm{  \by^t - \by^{t+1}} +  2 \E\norm{ \frac{1}{B}\sum_{(\zeta^t,\xi^t) \in \cB^t} \nabla_2 f(\bx^t, \bz^{t+1};\zeta^t) \nabla g(\by^t;\xi^t)   -  \nabla_2 f(\bx^t, g(\by^t )) \nabla g(\by^t) }^2\\
            &\leq \frac{2}{\eta_{\by}^2}  \E\norm{  \by^t - \by^{t+1}}   +  4 \E\norm{ \frac{1}{B}\sum_{(\zeta^t,\xi^t) \in \cB^t} \nabla_2 f(\bx^t, \bz^{t+1};\zeta^t) \nabla g(\by^t;\xi^t)   -  \nabla_2 f(\bx^t, \bz^{t+1}) \nabla g(\by^t) }^2 \\
            &\quad +4\E\norm{ \nabla_2 f(\bx^t, \bz^{t+1}) \nabla g(\by^t)  -  \nabla_2 f(\bx^t, g(\by^t )) \nabla g(\by^t) }^2\\ 
            &\leq \frac{2}{\eta_{\by}^2}  \E\norm{  \by^t - \by^{t+1}}^2 +4 L_f^2 G_g^2\E\norm{  \bz^{t+1}  - g(\by^t   )  }^2 +   \frac{8 (G_f^2+G_g^2)\sigma^2}{B}.
        \end{align*}
    \end{proof}
\end{lemma}

The following lemma connects $F(\bx^{t+1}, \by^{t+1})$ and $ F(\bx^{t+1}, \by^t)  $.
\begin{lemma}[Dual Variable One Iteration Analysis]\label{lem:dual descent}
For Algorithm~\ref{algorithm: CODA-Dual}, under the assumptions of Theorem~\ref{thm:SCNC}, the following statement holds true:
 \begin{align*}
    \E[F(\bx^{t+1}, \by^{t+1})  - F(\bx^{t+1}, \by^t)   ] &\geq \pare{ \frac{1}{2\eta_{\by}} - \frac{L}{2} - \frac{L^2\eta_{\bx}}{2}}\E\|  \by^{t+1}- \by^t \|^2 -  \frac{1}{2\eta_{\bx}}\E \norm{\bx^{t+1} - \bx^{t}}^2  \\
       &-  \pare{ {\eta_{\by}} 2(G_g^2+ G_f^2)\frac{\sigma^2}{B} +2\beta^2 L_f^2 G_g^2 {\eta_{\by}}\frac{\sigma^2}{M}  } \\
    &-   L_f^2 G_g^2 {\eta_{\by}}    (1-\beta)^2 \E \norm{\bz^{t} - g(\by^{t-1}) }^2 - 4(1-\beta)^2   L_f^2   G_g^4 {\eta_{\by}} \E\norm{\by^t - \by^{t-1}}^2  .
 \end{align*}
 \begin{proof}
 
By the optimality of projection and updating rule of $\w$:
 \begin{align*}
     \left\langle  \eta_{\by} \bg_{\by}^t +(\by^t - \by^{t+1}) , \by^t - \by^{t+1}\right\rangle \leq 0
 \end{align*}
 Re-arranging terms we have:
 \begin{align*}
      &    \left\langle  \bg_{\by}^t , \by^t -\by^{t+1}\right\rangle \leq -\frac{1}{\eta_{\by}}\| \by^{t+1}-\by^t \|^2
 \end{align*}
Adding and subtracting $\nabla_{\by} F(\bx^t, \by^t)$ yields:
\begin{align*}
 \left\langle \nabla_{\by} F(\bx^t, \by^t)  , \by^t -\by^{t+1}\right\rangle +  \left\langle  \bg_{\by}^t -\nabla_{\by} F(\bx^t, \by^t), \by^t -\by^{t+1}\right\rangle  \leq -\frac{1}{\eta_{\by}}\| \by^{t+1}-\by^t \|^2
\end{align*}



Applying Cauchy-Schwartz inequality yields:
 \begin{align*}
   & \frac{1}{\eta_{\by}}\|  \by^{t+1}- \by^t \|^2 \leq  \left\langle    \nabla_{\by} F(\bx^t, \by^t)  ,  \by^{t+1} - \by^t\right\rangle+\frac{1}{2\eta_{\by}} \left\|\by^{t+1} -\by^t\right\|^2   +   \frac{1}{2}{\eta_{\by}}  \left\|  \bg_{\by}^t -\nabla_{\by} F(\bx^t, \by^t) \right\|^2 
 \end{align*}
Now, we take expectation over the randomness of $\cB^t$, and apply Lemma~\ref{lem:grad bound}:
  \begin{align*} 
  & \frac{1}{2\eta_{\by}}\E\|  \by^{t+1}- \by^t \|^2 \leq \E \left\langle    \nabla_{\by} F(\bx^t, \by^t)  ,  \by^{t+1} - \by^t\right\rangle    +   \eta_{\by} \frac{ 2(G_f^2+G_g^2 )\sigma^2}{B} + \eta_{\by} G_g^2L_f^2 \E\norm{ \bz^{t+1}  -  g(\by^t)    }^2.
 \end{align*}
 
 Plugging in Lemma~\ref{lem:tracking error} yields:
 \begin{align*}
    \frac{1}{2\eta_{\by}}\E\|  \by^{t+1}- \by^t \|^2 &\leq \E \left\langle    \nabla_{\by} F(\bx^t, \by^t)  ,  \by^{t+1} - \by^t\right\rangle    +    +   \eta_{\by} \frac{ 2(G_f^2+G_g^2 )\sigma^2}{B}  \\
   &\quad + \eta_{\by} G_g^2L_f^2 \pr{  (1-\beta)^2 \E \norm{\bz^{t} - g(\by^{t-1}) }^2  + 4(1-\beta)^2 G_g^2 \norm{\by^t - \by^{t-1}}^2 + 2\beta^2 \frac{\sigma^2}{M}  } \\
 \end{align*}

On the other hand, since $F$ is $L$ smooth, we have:
  \begin{align*}
    & F(\bx^{t+1}, \by^{t+1})  - F(\bx^{t+1}, \by^t)  \geq \inprod{\nabla_{\by} F(\bx^{t+1}, \by^t)}{ \by^{t+1} -  \by^t} - \frac{L }{2}   \norm{ \by^{t+1} -  \by^t}^2\\
     &\quad= \inprod{\nabla_{\by} F(\bx^{t}, \by^t)}{ \by^{t+1} -  \by^t} - \frac{L }{2}   \norm{ \by^{t+1} -  \by^t}^2 + \inprod{\nabla_{\by} F(\bx^{t+1} , \by^t)- \nabla_{\by} F(\bx^{t}, \by^t)}{ \by^{t+1} -  \by^t} \\ 
     &\quad \geq  \inprod{\nabla_{\by} F(\bx^{t}, \by^t)}{ \by^{t+1} -  \by^t} - \frac{L }{2}   \norm{ \by^{t+1} -  \by^t}^2 -    \frac{1}{L\sqrt{\eta_{\bx} }}\norm{\nabla_{\by} F(\bx^{t+1} , \by^t)- \nabla_{\by} F(\bx^{t}, \by^t)}  {L\sqrt{\eta_{\bx}}}\norm{ \by^{t+1} -  \by^t}  \\
     &\quad \geq  \inprod{\nabla_{\by} F(\bx^{t}, \by^t)}{ \by^{t+1} -  \by^t} - \frac{L }{2}   \norm{ \by^{t+1} -  \by^t}^2  -    \frac{1}{2 {\eta_{\bx} }}\norm{ \bx^{t+1}  -  \bx^{t} }^2 -  \frac{L^2 \eta_{\bx}}{2}\norm{ \by^{t+1} -  \by^t}^2  \\
     &\quad\geq  \inprod{\nabla_{\by} F(\bx^{t}, \by^t)}{ \by^{t+1} -  \by^t} -\pare{ \frac{L }{2} + \frac{L^2\eta_{\bx}}{2}   }\norm{ \by^{t+1} -  \by^t}^2  -    \frac{1}{2\eta_{\bx}}\norm{\bx^{t+1} - \bx^{t}}^2   \\
 \end{align*}

 Putting pieces together will conclude the proof:
 \begin{align*}
       \E[F(\bx^{t+1}, \by^{t+1})  - F(\bx^{t+1}, \by^t)   ] &\geq \pare{ \frac{1}{2\eta_{\by}} - \frac{L}{2} - \frac{L^2\eta_{\bx}}{2}}\E\|  \by^{t+1}- \by^t \|^2 -  \frac{1}{2\eta_{\bx}}\E \norm{\bx^{t+1} - \bx^{t}}^2  \\
       &-  \pare{ {\eta_{\by}} 2(G_g^2+ G_f^2)\frac{\sigma^2}{B} +2\beta^2 L_f^2 G_g^2 {\eta_{\by}}\frac{\sigma^2}{M}  } \\
    &-   L_f^2 G_g^2 {\eta_{\by}}    (1-\beta)^2 \E \norm{\bz^{t} - g(\by^{t-1}) }^2 - 4(1-\beta)^2   L_f^2   G_g^4 {\eta_{\by}} \E\norm{\by^t - \by^{t-1}}^2  .
 \end{align*}
 \end{proof}
\end{lemma}

 The following lemma connects $F( \bx^{t+1},\by^{t+1})  $ and $ F(\bx^{t},\by^t)$.
 
 \begin{lemma}[One Iteration Descent Lemma]\label{lm:descent lemma}
 For Algorithm~\ref{algorithm: CODA-Dual}, under the assumptions of Theorem~\ref{thm:SCNC}, the following statement holds true:
\begin{align*}
    \E[F(\bx^{t+1}, \by^{t+1}) & - F(\bx^{t}, \by^t)   ] \geq \pare{ \frac{1}{2\eta_{\by}} - \frac{L}{2} -  \frac{L^2\eta_{\bx}}{2} }\E\|  \by^{t+1}- \by^t \|^2 \\
    &- \pare{ 4(1-\beta)^2   L_f^2 G_f^2 G_g^2 {\eta_{\by}} + \frac{L^2\eta_{\bx}}{2} } \E\norm{\by^t - \by^{t-1}}^2 \\
    &+ \left (\frac{\mu}{2} - \frac{3}{ 2\eta_{\bx}}-\frac{1}{4\eta_{\bx} (1-\beta)^2}\right)\E \norm{\bx^{t+1} - \bx^{t}}^2 +\left(\mu - \frac{1}{2\eta_{\bx}}- \frac{\eta_{\bx} L^2}{2}\right)\E \norm{\bx^{t} - \bx^{t-1}}^2   \\
    &-  \pr{L_f^2 G_g^2 {\eta_{\by}} + \eta_{\bx} L_f^2  }(1-\beta)^2\E \norm{\bz^{t} - g(\by^{t-1}) }^2   -   \pare{ 2{\eta_{\by}}  (G_g^2+ G_f^2)\frac{\sigma^2}{B} +2\beta^2 L_f^2 G_g^2 {\eta_{\by}}\frac{\sigma^2}{M}  }.
 \end{align*}

 \begin{proof}
  According to property of projection we have:
 \begin{align}
     \left \langle  \bg_{\bx}^t +\frac{1}{\eta} (\bx^{t+1}-\bx^{t} ) , \bx -\bx^{t+1} \right \rangle \geq 0 \label{eq:smpe_update_1}\\
     \left \langle  \bg_{\bx}^t +\frac{1}{\eta} (\bx^{t+1}-\bx^{t} ) , \bx^t -\bx^{t+1} \right \rangle \geq 0\label{eq:smpe_update_2}\\
     \left \langle  \bg_{\bx}^{t-1} +\frac{1}{\eta} (\bx^{t}-\bx^{t-1} ) , \bx^{t+1} -\bx^{t} \right \rangle \geq 0 \label{eq:smpe_update_3}
 \end{align}
  
 Since $F(\cdot,\by)$ is strongly convex for any $\by\in\mathcal{Y}$, we have:
 \begin{align*}
  F( \bx^{t+1},\by^t) -F(\bx^{t},\by^t) &\geq  \left \langle  \nabla_{\bx} F(\bx^{t},\by^t), \bx^{t+1} -\bx^{t} \right \rangle + \frac{\mu}{2}  \|\bx^{t+1} -\bx^{t}\|^2\\
     & \geq   \left \langle  \nabla_{\bx} F(\bx^{t},\by^t) -\bg_{\bx}^{t-1}, \bx^{t+1} -\bx^{t} \right \rangle + \frac{\mu}{2}  \|\bx^{t+1} -\bx^{t}\|^2 - \frac{1}{\eta_{\bx}} \left \langle   \bx^{t}-\bx^{t-1}, \bx^{t+1} -\bx^{t} \right \rangle\\
     & \geq   \left \langle  \nabla_{\bx} F(\bx^{t},\by^t) - \nabla_{\bx} F(\bx^{t-1},\by^{t-1}), \bx^{t+1} -\bx^{t} \right \rangle + \frac{\mu}{2}  \|\bx^{t+1} -\bx^{t}\|^2\\
      & \quad + \left \langle \bg_{\bx}^{t-1}- \nabla_{\bx} F(\bx^{t-1},\by^{t-1}), \bx^{t+1} -\bx^{t} \right \rangle  - \frac{1}{\eta_{\bx}} \left \langle   \bx^{t}-\bx^{t-1}, \bx^{t+1} -\bx^{t} \right \rangle  
 \end{align*}
 where at second inequality we use the fact of (\ref{eq:smpe_update_3}) .
 For the first dot product term, we observe that:
 \begin{align}
      \left \langle  \nabla_{\bx} F(\bx^{t},\by^t) - \nabla_{\bx} F(\bx^{t-1}, \by^{t-1}), \bx^{t+1} -\bx^{t} \right \rangle &=  \left \langle  \nabla_{\bx} F(\bx^{t},\by^t) - \nabla_{\bx} F(\bx^{t},\by^{t-1}), \bx^{t+1} -\bx^{t} \right \rangle \nonumber\\
     &+ \left \langle  \nabla_{\bx} F(\bx^{t},\by^{t-1}) - \nabla_{\bx} F(\bx^{t-1}, \by^{t-1}), \Delta  \right \rangle \nonumber\\
     &+ \left \langle  \nabla_{\bx} F(\bx^{t},\by^{t-1}) - \nabla_{\bx} F(\bx^{t-1},\by^{t-1}), \bx^{t} -\bx^{t-1} \right \rangle\nonumber\\
     & \geq -\frac{L^2 \eta_{\bx}}{2}  \|\by^t - \by^{t-1} \|^2 - \frac{1}{2\eta_{\bx}} \|\bx^{t+1} -\bx^{t}\|^2 - \frac{\eta_{\bx} L^2}{2}  \|\bx^{t} -\bx^{t-1}\|^2\nonumber\\
     &\quad - \frac{1}{2\eta_{\bx}}  \|\Delta \|^2 + \mu  \|\bx^{t} -\bx^{t-1}\|^2\label{eq:one iteration 1},
 \end{align}
 where $ \Delta = (\bx^{t+1}- \bx^{t})-(\bx^{t}-\bx^{t-1}).  $
 
 For the second dot product, we have:
 \begin{align*}
   \E \left \langle \bg_{\bx}^{t-1}- \nabla_{\bx} F(\bx^{t-1},\by^{t-1}), \bx^{t+1} -\bx^{t} \right \rangle&\geq - \eta_{\bx} (1-\beta)^2 \E \norm{ \bg_{\bx}^{t-1}- \nabla_{\bx} F(\bx^{t-1},\by^{t-1})}^2 - \frac{1}{4\eta (1-\beta)^2}\E\norm{ \bx^{t+1} -\bx^{t}}^2\\ 
   & \geq -\eta_{\bx} (1-\beta)^2 L_f^2 \E\norm{\bz^{t} - g(\by^{t-1}) }^2  - \frac{1}{4\eta_{\bx} (1-\beta)^2} \E\norm{ \bx^{t+1} -\bx^{t}}^2,
 \end{align*}
 where we plug in Lemma~\ref{lem:grad bound} at last step.
 
For the third dot product, we use the following identity:
 \begin{align}
     \frac{1}{\eta_{\bx}}\left \langle   \bx^{t}-\bx^{t-1}, \bx^{t+1} -\bx^{t} \right \rangle = \frac{1}{\eta_{\bx}} \left( \frac{1}{2}\|\bx^{t}-\bx^{t-1}\|^2 + \frac{1}{2}\|\bx^{t+1}-\bx^{t}\|^2 - \frac{1}{2}\|\Delta \|^2  \right) \label{eq:dot_product_id}
 \end{align}
 By putting pieces together we have:
 \begin{align}
     \E[ F( \bx^{t+1},\by^t) -F(\bx^{t},\by^t)] &\geq -\frac{L^2 \eta_{\bx}}{2} \E\|\by^t -  \by^{t-1} \|^2 + \left (\frac{\mu}{2} - \frac{1}{ \eta_{\bx}}-\frac{1}{4\eta_{\bx} (1-\beta)^2}\right) \E\|\bx^{t+1} -\bx^{t}\|^2 \nonumber \\
     & \quad +\left(\mu - \frac{1}{2\eta_{\bx}}- \frac{\eta_{\bx} L^2}{2}\right)  \E\|\bx^{t} -\bx^{t-1}\|^2  -\eta_{\bx} (1-\beta)^2  L_f^2 \E\norm{\bz^{t} - g(\by^{t-1}) }^2. \label{eq:descent lem 1}
 \end{align}
 Recall Lemma~\ref{lem:dual descent} gives the following lower bound of $ \E[F(\bx^{t+1}, \by^{t+1})  - F(\bx^{t+1}, \by^t)   ] $:
 \begin{align}
        \E[F(\bx^{t+1}, \by^{t+1})  - F(\bx^{t+1}, \by^t)   ] &\geq \pare{ \frac{1}{2\eta_{\by}} - \frac{L}{2} - \frac{L^2\eta_{\bx}}{2}}\E\|  \by^{t+1}- \by^t \|^2 -  \frac{1}{2\eta_{\bx}}\E \norm{\bx^{t+1} - \bx^{t}}^2 \nonumber \\
       &-   \pare{ 2{\eta_{\by}}  (G_g^2+ G_f^2)\frac{\sigma^2}{B} +2\beta^2 L_f^2 G_g^2 {\eta_{\by}}\frac{\sigma^2}{M}  } \nonumber \\
    &-   L_f^2 G_g^2 {\eta_{\by}}    (1-\beta)^2 \E \norm{\bz^{t} - g(\by^{t-1}) }^2 - 4(1-\beta)^2   L_f^2 G_f^2 G_g^2 {\eta_{\by}} \E\norm{\by^t - \by^{t-1}}^2  . \label{eq:descent lem 2}
 \end{align}
 Hence, Adding (\ref{eq:descent lem 1}) and (\ref{eq:descent lem 2}) yields:
 \begin{align*}
    \E[F(\bx^{t+1}, \by^{t+1}) & - F(\bx^{t}, \by^t)   ] \geq \pare{ \frac{1}{2\eta_{\by}} - \frac{L}{2} -  \frac{L^2\eta_{\bx}}{2} }\E\|  \by^{t+1}- \by^t \|^2 \\
    & - \pare{ 4(1-\beta)^2   L_f^2 G_f^2 G_g^2 {\eta_{\by}} + \frac{L^2\eta_{\bx}}{2} } \E\norm{\by^t - \by^{t-1}}^2 \\
    &+ \left (\frac{\mu}{2} - \frac{3}{ 2\eta_{\bx}}-\frac{1}{4\eta_{\bx} (1-\beta)^2}\right)\E \norm{\bx^{t+1} - \bx^{t}}^2 +\left(\mu - \frac{1}{2\eta_{\bx}}- \frac{\eta_{\bx} L^2}{2}\right)\E \norm{\bx^{t} - \bx^{t-1}}^2   \\
    &-  \pr{ L_f^2 G_g^2 {\eta_{\by}} + \eta_{\bx} L_f^2  }(1-\beta)^2\E \norm{\bz^{t} - g(\by^{t-1}) }^2    -   \pare{ 2{\eta_{\by}}  (G_g^2+ G_f^2)\frac{\sigma^2}{B} +2\beta^2 L_f^2 G_g^2 {\eta_{\by}}\frac{\sigma^2}{M}  },
 \end{align*}
which concludes the proof.
 
 \end{proof}
\end{lemma}

 \begin{lemma}\label{lem:potential}
  For Algorithm~\ref{algorithm: CODA-Dual}, under the assumptions of Theorem~\ref{thm:SCNC},  define the following auxillary function $\hat F^{t+1}$
\begin{align*}
    \hat{F}^{t+1} := & F(\bx^{t+1},\by^{t+1}) + s^{t+1} - \pare{\frac{1}{4\eta_{\by}} + 4 L_g^2 G_f^2 G_g^2 \eta_{\by}  + \frac{\eta_{\bx} L^2}{2} + \frac{168  G_g^2L_f^2}{\mu^2 \eta_{\bx} \beta} + \frac{48  G_g^2 L_f^2\beta}{\mu^2 \eta_{\bx}}} \norm{\by^{t+1} - \by^t}^2\\
    &  -\pare{\frac{1}{8\eta_{\by}}+\frac{84 G_g^2L_f^2}{\mu^2 \eta_{\bx} \beta}}\norm{\by^t - \by^{t-1}}^2 + \pare{\frac{7}{2\eta_{\bx}} + \mu - \frac{\eta_{\bx} L^2}{2} - \frac{2L^2_f}{\mu}}\norm{\bx^{t+1} - \bx^t}^2
\end{align*}
where $s^{t+1} :=  -\frac{2}{\eta^2_{\bx} \mu}\|\bx^{t+1}-\bx^t\|^2$
then the following statement holds true:
\begin{align*}
  \E [ \hat{F}^{t+1} - \hat{F}^{t}] &\geq C_1\E \norm{\by^{t+1} - \by^t}^2 + \frac{1}{8\eta_{\by}}\E \norm{\by^t - \by^{t-1}}^2 + \frac{1}{8\eta_{\by}}\E \norm{\by^{t-1} - \by^{t-2}}^2 + C_2\E \norm{\bx^{t+1} - \bx^t}^2\\
      &-   (1-\beta)^2 \pare{  L_g^2 G_f^2 {\eta_{\by}} + \eta_{\bx} L_f^2}  \E \norm{\bz^{t} -g(\by^{t-1}) }^2 \\
    & -  \pare{ 2{\eta_{\by}}  (G_g^2+ G_f^2)\frac{\sigma^2}{B} +2\beta^2 L_f^2 G_g^2 {\eta_{\by}}\frac{\sigma^2}{M}  +\frac{72L_f^2}{\mu^2\eta_\bx}  \beta \frac{\sigma^2}{M}}  -\pare{1-\frac{\beta}{2}}^2  \frac{4L_f^2}{\mu^2 \eta_{\bx}}    \E \norm{  \bz_j^{t } - \bz^{t-1}  }^2 .
\end{align*}
where
\begin{align}
    &C_1 = \frac{1}{4\eta_{\by}} - \frac{L}{2} -\eta_{\bx} L^2 - 4L_g^2G_f^2G_g^2\eta_{\by} - \frac{168  G_g^2L_f^2}{\mu^2\eta_{\bx} \beta} - \frac{48   G_g^2 L_f^2 \beta}{\mu^2 \eta_{\bx}}, \label{eq: def C1} \\
    &C_2 = \frac{1}{\eta_{\bx}} + \frac{3}{2}\mu - \frac{\eta_{\bx} L^2}{2} - \frac{1}{4\eta_{\bx} (1-\beta)^2}-\frac{2L^2_f}{\mu}.\label{eq: def C2}
\end{align}

\begin{proof}
Define $\Delta^{t+1}:= \bx^{t+1} - \bx^{t} - (\bx^{t} - \bx^{t-1})$.    According to (\ref{eq:smpe_update_2}) and (\ref{eq:smpe_update_3}):

\begin{align*}
    \frac{1}{\eta}  \left \langle \Delta^{t+1}, \bx^{t} - \bx^{t+1}  \right \rangle &\geq \inprod{\bg_{\bx}^t - \bg_{\bx}^{t-1} }{\bx^{t+1} - \bx^{t} }. 
\end{align*}

If we define $\bg^t_{\bx}(\bx,\bz) = \nabla_1 f(\bx,\bz;\cB^t) + \nabla h(\bx)$.  
\begin{align*}
    \frac{1}{\eta}  \left \langle \Delta^{t+1}, \bx^{t} - \bx^{t+1}  \right \rangle &\geq \inprod{\bg_{\bx}^t - \bg^t_{\bx}(\bx^{t},\bz^t;\xi^t) }{\bx^{t+1} - \bx^{t} }\\ 
   &\quad  +  \inprod{ \bg^t_{\bx}(\bx^{t},\bz^t;\xi^t )- \bg_{\bx}^{t-1}}{\bx^{t+1} - \bx^{t} }\\
     &\geq \inprod{\bg_{\bx}^t - \bg_{\bx}( \bx^{t},\bz^{t}) }{\bx^{t+1} - \bx^{t} }\\ 
   &\quad +  \inprod{ \bg_{\bx}(\bx^{t},\bz^{t})- \bg_{\bx}^{t-1}}{\bx^{t+1} - \bx^{t} - (\bx^{t} - \bx^{t-1}) }\\
   &\quad +  \inprod{ \bg_{\bx}( \bx^{t},\bz^{t})- \bg_{\bx}^{t-1}}{  (\bx^{t} - \bx^{t-1}) }\\
   &\geq -\frac{L_f^2}{\mu}  \norm{  \bz^{t+1} -  \bz^{t}  }^2 -  \frac{\mu}{4} \norm{\bx^{t+1} - \bx^{t} }^2\\ 
   &\quad - \frac{\eta_{\bx} L_f^2 }{2}  \norm{   \bx^{t} - \bx^{t-1} }^2 - \frac{1}{2\eta_{\bx}} \norm{\Delta^{t+1} }^2  +  \mu \norm{  \bx^{t} - \bx^{t-1}  }^2\\
\end{align*}
where we use the Lipschitzness of $\nabla f(\bx,\bz;\xi)$ at the last step.
Since $\frac{1}{\eta_{\bx} }\left \langle \Delta^{t+1}, \bx^{t} - \bx^{t+1}  \right \rangle = \frac{1}{2\eta_{\bx}}\norm{\bx^{t} - \bx^{t-1}}^2 - \frac{1}{2\eta_{\bx}}\norm{\bx^{t+1} - \bx^{t}}^2 - \frac{1}{2\eta_{\bx}} \norm{\Delta^{t+1}}^2$
 we have:
\begin{align*}
   \frac{1}{2\eta_{\bx}}\norm{\bx^{t} - \bx^{t-1}}^2 - \frac{1}{2\eta_{\bx}}\norm{\bx^{t+1} - \bx^{t}}^2  
   &\geq -\frac{L_f^2}{\mu} \norm{  \bz^{t+1} -  \bz^{t}  }^2 -  \frac{\mu}{4} \norm{\bx^{t+1} - \bx^{t} }^2\\ 
   &\quad - \frac{\eta_{\bx} L^2_f }{2}  \norm{   \bx^{t} - \bx^{t-1} }^2    +  \mu \norm{  \bx^{t} - \bx^{t-1}  }^2\\
\end{align*}
Multiplying both sides with $\frac{4}{\mu \eta_{\bx}}$ yields:
\begin{align}
   \frac{2}{\mu \eta_{\bx}^2}\norm{\bx^{t} - \bx^{t-1}}^2 -  \frac{2}{\mu \eta_{\bx}^2}\norm{\bx^{t+1} - \bx^{t}}^2  
   &\geq -\frac{4L_f^2}{ \mu^2 \eta_{\bx}}  \norm{  \bz^{t+1} -  \bz^{t}  }^2 -  \frac{1}{\eta_{\bx}} \norm{\bx^{t+1} - \bx^{t} }^2 \nonumber \\ 
   &\quad - \frac{2L^2_f}{\mu}  \norm{   \bx^{t} - \bx^{t-1} }^2    +  \frac{4}{\eta_{\bx}} \norm{  \bx^{t} - \bx^{t-1}  }^2 \label{eq: potential1 1}
\end{align} 
Recall that in Lemma~\ref{lm:descent lemma} we have:
 \begin{align*}
    \E[F(\bx^{t+1}, \by^{t+1}) & - F(\bx^{t}, \by^t)   ] \geq \pare{ \frac{1}{2\eta_{\by}} - \frac{L}{2} -  \frac{L^2\eta_{\bx}}{2} }\E\|  \by^{t+1}- \by^t \|^2 \\
    & - \pare{ 4(1-\beta)^2   L_f^2 G_f^2 G_g^2 {\eta_{\by}} + \frac{L^2\eta_{\bx}}{2} } \E\norm{\by^t - \by^{t-1}}^2 \\
    &+ \left (\frac{\mu}{2} - \frac{3}{ 2\eta_{\bx}}-\frac{1}{4\eta_{\bx} (1-\beta)^2}\right)\E \norm{\bx^{t+1} - \bx^{t}}^2 +\left(\mu - \frac{1}{2\eta_{\bx}}- \frac{\eta_{\bx} L^2}{2}\right)\E \norm{\bx^{t} - \bx^{t-1}}^2   \\
    &-  \pr{ L_f^2 G_g^2 {\eta_{\by}} + \eta_{\bx} L_f^2  }(1-\beta)^2\E \norm{\bz^{t} - g(\by^{t-1}) }^2   -   \pare{ 2{\eta_{\by}}  (G_g^2+ G_f^2)\frac{\sigma^2}{B} +2\beta^2 L_f^2 G_g^2 {\eta_{\by}}\frac{\sigma^2}{M}  }.
 \end{align*}
Combining Lemma~\ref{lm:descent lemma} together with (\ref{eq: potential1 1}) yields:
 \begin{align*}
   \E&[F(\bx^{t+1}, \by^{t+1}) + s^{t+1}  -  (F(\bx^{t}, \by^t) +s^{t})  ] \\
    &\geq \pare{ \frac{1}{2\eta_{\by}} - \frac{L}{2} -  \frac{L^2\eta_{\bx}}{2} }\E\|  \by^{t+1}- \by^t \|^2- \pare{ 4(1-\beta )^2 L_f^2 G_f^2G_g^2 \eta_{\by} + \frac{L^2\eta_{\bx}}{2} } \E\norm{\by^t - \by^{t-1}}^2 \\
    &+ \left (\frac{\mu}{2} - \frac{3}{ 2\eta_{\bx}}-\frac{1}{4\eta_{\bx} (1-\beta )^2} - \frac{1}{\eta_{\bx}}\right)\E \norm{\bx^{t+1} - \bx^{t}}^2 +\left(\mu - \frac{1}{2\eta_{\bx}}- \frac{\eta_{\bx} L^2}{2} + \frac{4}{\eta_{\bx}} - \frac{2L^2_f}{\mu}\right)\E \norm{\bx^{t} - \bx^{t-1}}^2   \\
    &-    (1-\beta)^2 \pare{ L_f^2 G_g^2 {\eta_{\by}} + \eta_{\bx} L_f^2} \E \norm{\bz^{t} -g(\by^{t-1}) }^2    -   \pare{ 2{\eta_{\by}}  (G_g^2+ G_f^2)\frac{\sigma^2}{B} +2\beta^2 L_f^2 G_g^2 {\eta_{\by}}\frac{\sigma^2}{M}  }   - \frac{4L_f^2}{\mu^2 \eta_{\bx}}  \E\norm{\bz^{t+1} -  \bz^{t}}^2.
 \end{align*}

 Now we plug in Lemma~\ref{lem:second order tracking} to replace $\E\norm{\bz^{t+1} -  \bz^{t}}^2$ together with the fact that $\frac{1}{\beta} \geq 1$ and get:

 \begin{align*}
   \E[F(\bx^{t+1}, \by^{t+1}) &+ s^{t+1}  -  (F(\bx^{t}, \by^t) +s^{t})  ] \\
    &\geq \pare{ \frac{1}{2\eta_{\by}} - \frac{L}{2} -  \frac{L^2\eta_{\bx}}{2} }\E\|  \by^{t+1}- \by^t \|^2\\
    &  \quad- \pare{ 4  L_f^2 G_f^2G_g^2 \eta_{\by} + \frac{L^2\eta_{\bx}}{2} +\frac{84 G_g^2L_f^2}{\mu^2 \eta_{\bx} \beta}    +\frac{48G_g^2L_f^2\beta}{\mu^2 \eta_{\bx}}  }    \E\norm{\by^t - \by^{t-1}}^2  -\frac{48G_g^2L_f^2}{\mu^2 \eta_{\bx} \beta}   \E\norm{\by^{t-1} - \by^{t-2}}^2\\
    & \quad + \left (\frac{\mu}{2} - \frac{5}{ 2\eta_{\bx}}-\frac{1}{4\eta_{\bx} (1-\beta)^2}  \right)\E \norm{\bx^{t+1} - \bx^{t}}^2 +\left(\mu + \frac{7}{2\eta_{\bx}}- \frac{\eta_{\bx} L^2}{2} - \frac{2L^2_f}{\mu} \right)\E \norm{\bx^{t} - \bx^{t-1}}^2   \\
    & \quad  -   (1-\beta)^2 \pare{  L_f^2 G_g^2 {\eta_{\by}} + \eta_{\bx} L_f^2}   \E \norm{\bz^{t} -g(\by^{t-1}) }^2 \\
    &  \quad -  \pare{ 2{\eta_{\by}}  (G_g^2+ G_f^2)\frac{\sigma^2}{B} +2\beta^2 L_f^2 G_g^2 {\eta_{\by}}\frac{\sigma^2}{M}  }  -\pare{1-\frac{\beta}{2}}^2  \frac{4L_f^2}{\mu^2 \eta_{\bx}}    \E \norm{  \bz^{t } - \bz^{t-1}  }^2 + \frac{4L_f^2}{\mu^2\eta_\bx} \cdot 18\beta \frac{\sigma^2}{M} .
 \end{align*}
  
Recall our definition of potential function $\hat F^{t+1}$
\begin{align*}
    \hat{F}^{t+1} := & F(\bx^{t+1},\by^{t+1}) + s^{t+1} - \pare{\frac{1}{4\eta_{\by}} + 4 L_f^2 G_f^2 G_g^2 \eta_{\by}  + \frac{\eta_{\bx} L^2}{2} + \frac{168  G_g^2L_f^2}{\mu^2 \eta_{\bx} \beta} + \frac{48  G_g^2 L_f^2\beta}{\mu^2 \eta_{\bx}}} \norm{\by^{t+1} - \by^t}^2\\
    &  -\pare{\frac{1}{8\eta_{\by}}+\frac{84 G_g^2L_f^2}{\mu^2 \eta_{\bx} \beta}}\norm{\by^t - \by^{t-1}}^2 + \pare{\frac{7}{2\eta_{\bx}} + \mu - \frac{\eta_{\bx} L^2}{2} - \frac{2L^2_f}{\mu}}\norm{\bx^{t+1} - \bx^t}^2 .
\end{align*}

We conclude that:
\begin{align*}
  \E [ &\hat{F}^{t+1} - \hat{F}^{t}] \geq \pare{\frac{1}{4\eta_{\by}} - \frac{L}{2} -\eta_{\bx} L^2 - 4L_f^2G_f^2G_g^2\eta_{\by} - \frac{168  G_g^2L_f^2}{\mu^2\eta_{\bx} \beta} - \frac{48   G_g^2 L_f^2 \beta}{\mu^2 \eta_{\bx}}}\E \norm{\by^{t+1} - \by^t}^2\\
  &+ \frac{1}{8\eta_{\by}}\E \norm{\by^t - \by^{t-1}}^2 + \frac{1}{8\eta_{\by}}\E \norm{\by^{t-1} - \by^{t-2}}^2 + \pare{\frac{1}{\eta_{\bx}} + \frac{3}{2}\mu - \frac{\eta_{\bx} L^2}{2} - \frac{1}{4\eta_{\bx} (1-\beta)^2}-\frac{2L^2_f}{\mu}} \E \norm{\bx^{t+1} - \bx^t}^2\\
      &-   (1-\beta)^2 \pare{  L_f^2 G_g^2 {\eta_{\by}} + \eta_{\bx} L_f^2}  \E \norm{\bz^{t} -g(\by^{t-1}) }^2 \\
    &  -   \pare{ 2{\eta_{\by}}  (G_g^2+ G_f^2)\frac{\sigma^2}{B} +2\beta^2 L_f^2 G_g^2 {\eta_{\by}}\frac{\sigma^2}{M}   + \frac{72L_f^2}{\mu^2\eta_\bx}  \beta \frac{\sigma^2}{M} }   -\pare{1-\frac{\beta}{2}}^2  \frac{4L_f^2}{\mu^2 \eta_{\bx}}    \E \norm{  \bz^{t} - \bz^{t-1}  }^2 .
\end{align*}

\end{proof}
\end{lemma}

 \begin{lemma}\label{lemma: potential 2 scnc}

Let $C_1, C_2$ be defined in (\ref{eq: def C1}) and (\ref{eq: def C2}). If the following conditions hold:
\begin{align}
   \frac{1}{8\eta_{\by}}-4(2 C_1\eta_{\by}^2 L_f^2G_g^2 + 2C_2\eta_{\bx}^2 L_f^2 + 1)(1-\beta)^2 G_f^2   -\frac{(1-\frac{\beta}{2})^2}{1-(1-\frac{\beta}{2})^2}\frac{4L_f^2}{\mu^2\eta_{\bx}} \pr{ \frac{12}{\beta} + 6\beta  } G_g^2  \geq 0, \label{eq:cond 1}\\
   1-\pare{2C_1\eta_{\by}^2 L_f^2G_g^2 + 2C_2\eta^2_{\bx} L_f^2 +   L_f^2 G_g^2 {\eta_{\by}} + \eta_{\bx} L_f^2+1}(1-\beta )^2 \geq 0, \label{eq:cond2}\\
   \frac{1}{8\eta_{\by}} - \frac{(1-\frac{\beta}{2})^2}{1-(1-\frac{\beta}{2})^2}\frac{4L_f^2}{\mu^2\eta_{\bx}}\frac{6}{\beta}G_g^2 \geq 0 \label{eq:cond3},
\end{align}
 then we have:
 \begin{align*}  
  \E [ \tilde{F}^{t+1} - \tilde{F}^{t}]   & \geq  \frac{C_1}{2}\eta_{\by}^2\E \norm{\hat{\bg}_{\by}(\bx^t,\by^t)}^2 +  \frac{C_2}{2} \eta^2_{\bx}\E \norm{\hat{\bg}_{\bx}(\bx^t,\by^t)}^2 \\ 
    & \quad   -  \pare{  4C_1 \eta_{\by}^2 (G_g^2+G_f^2)  + 2C_2 \eta_{\bx}^2 } \frac{\sigma^2}{B} -  \pare{2\beta^2 L_f^2 G_g^2 {\eta_{\by}}+ 4\beta^2C_1  \eta_{\by}^2 L_f^2 G_g^2 +4\beta^2 C_2\eta^2_\bx L_f^2  }\frac{\sigma^2}{M}-\frac{72L_f^2}{\mu^2\eta_\bx}  \beta \frac{\sigma^2}{M},   
\end{align*}
where \begin{align*}
     \tilde{F}^{t+1} :=  \hat{F}^{t+1} -   \E \norm{\bz^{t+1}  -  g(\by^t))  }^2 -\frac{(1-\frac{\beta}{2})^2}{1-(1-\frac{\beta}{2})^2}\frac{4L_f^2}{\mu^2\eta_{\bx}}   \E \norm{  \bz^{t+1 } - \bz^{t}  }^2.
\end{align*}
 \begin{proof}
     
 According to Lemma~\ref{lem:potential}:
\begin{align*}
  \E [ \hat{F}^{t+1} - \hat{F}^{t}] &\geq C_1\E \norm{\by^{t+1} - \by^t}^2 + \frac{1}{8\eta_{\by}}\E \norm{\by^t - \by^{t-1}}^2 + \frac{1}{8\eta_{\by}}\E \norm{\by^{t-1} - \by^{t-2}}^2 + C_2\E \norm{\bx^{t+1} - \bx^t}^2\\
      &-   (1-\beta)^2 \pare{  L_f^2 G_g^2 {\eta_{\by}} + \eta_{\bx} L_f^2}  \E \norm{\bz^{t} -g(\by^{t-1}) }^2 \\
    & -   \pare{ 2{\eta_{\by}}  (G_g^2+ G_f^2)\frac{\sigma^2}{B} +2\beta^2 L_f^2 G_g^2 {\eta_{\by}}\frac{\sigma^2}{M} +\frac{72L_f^2}{\mu^2\eta_\bx}  \beta \frac{\sigma^2}{M} }   -\pare{1-\frac{\beta}{2}}^2  \frac{4L_f^2}{\mu^2 \eta_{\bx}}    \E \norm{  \bz_j^{t } - \bz^{t-1}  }^2 .
\end{align*}
where
\begin{align*}
    &C_1 = \frac{1}{4\eta_{\by}} - \frac{L}{2} -\eta_{\bx} L^2 - 4L_f^2G_f^2G_g^2\eta_{\by} - \frac{168  G_g^2L_f^2}{\mu^2\eta_{\bx} \beta} - \frac{48   G_g^2 L_f^2 \beta}{\mu^2 \eta_{\bx}},\\
    &C_2 = \frac{1}{\eta_{\bx}} + \frac{3}{2}\mu - \frac{\eta_{\bx} L^2}{2} - \frac{1}{4\eta_{\bx} (1-\beta)^2}-\frac{2L^2_f}{\mu}.
\end{align*}
 
 Now we plug in Lemma~\ref{lem:stationary measure}:
 \begin{align*}
  \E [ \hat{F}^{t+1} - \hat{F}^{t}] \geq &C_1\pare{\frac{1}{2}\eta_{\by}^2\E \norm{\hat{\bg}_{\by}(\bx^t,\by^t)}^2  - 4\eta_{\by}^2(G_f^2+G_g^2) \frac{\sigma^2}{B} - 2\eta_{\by}^2G_g^2 L_f^2  \E\norm{ \bz^{t+1}  -  g(\bx^t))    }^2}   \\
  & + C_2 \pare{\frac{1}{2} \eta^2_{\bx}\E \norm{\hat{\bg}_{\bx}(\bx^t,\by^t)}^2 - 2\frac{\eta_{\bx}^2\sigma^2}{B}  -2\eta_{\bx}^2 L_f^2 \E\norm{\bz^{t+1}  -  g(\by^t)) }^2} \\ 
& -   (1-\beta)^2 \pare{  L_f^2 G_g^2 {\eta_{\by}} + \eta_{\bx} L_f^2}  \E \norm{\bz^{t} -g(\by^{t-1}) }^2  -   \pare{ 2{\eta_{\by}}  (G_g^2+ G_f^2)\frac{\sigma^2}{B} +2\beta^2 L_f^2 G_g^2 {\eta_{\by}}\frac{\sigma^2}{M} + \frac{72L_f^2}{\mu^2\eta_\bx}  \beta \frac{\sigma^2}{M} }\\
&-\pr{1-\frac{\beta}{2}}^2  \frac{4L_f^2}{\mu^2 \eta_{\bx}}   \E \norm{  \bz^{t } - \bz^{t-1}  }^2+ \frac{1}{8\eta_{\by}}\E \norm{\by^t - \by^{t-1}}^2 + \frac{1}{8\eta_{\by}}\E \norm{\by^{t-1} - \by^{t-2}}^2. 
\end{align*}
 
 Plugging in Lemma~\ref{lem:tracking error} yields:

 {  \begin{align*}
 & \E [ \hat{F}^{t+1} - \hat{F}^{t}] \geq \frac{C_1}{2}\eta_{\by}^2\E \norm{\hat{\bg}_{\by}(\bx^t,\by^t)}^2 + C_2 \frac{1}{2} \eta^2_{\bx}\E \norm{\hat{\bg}_\bx(\bx^t,\by^t)}^2 \\
  &-2C_1  \eta_{\by}^2 L_f^2 G_g^2 \pare{  (1-\beta)^2 \E \norm{\bz^{t} -g(\by^{t-1}) }^2 + 4(1-\beta)^2 G_g^2 \norm{\by^t - \by^{t-1}}^2 + 2\beta^2 \frac{\sigma^2}{M} }   \\ 
  &-2C_2\eta^2_\bx L_f^2   \pare{(1-\beta)^2 \E \norm{\bz^{t} -g(\by^{t-1}) }^2 + 4(1-\beta)^2 G_g^2 \norm{\by^t - \by^{t-1}}^2 + 2\beta^2 \frac{\sigma^2}{M} } \\
& -   (1-\beta)^2 \pare{  L_f^2 G_g^2 {\eta_{\by}} + \eta_{\bx} L_f^2}  \E \norm{\bz^{t} -g(\by^{t-1}) }^2   -  \pare{  4C_1 \eta_{\by}^2 (G_g^2+G_f^2)  + 2C_2 \eta_{\bx}^2 } \frac{\sigma^2}{B} -  2\beta^2 L_f^2 G_g^2 {\eta_{\by}}\frac{\sigma^2}{M}- \frac{72L_f^2}{\mu^2\eta_\bx}  \beta \frac{\sigma^2}{M}\\
&-\pr{1-\frac{\beta}{2}}^2  \frac{4L_f^2}{\mu^2 \eta_{\bx}}     \E \norm{  \bz^{t } - \bz^{t-1}  }^2+ \frac{1}{8\eta_{\by}}\E \norm{\by^t - \by^{t-1}}^2 + \frac{1}{8\eta_{\by}}\E \norm{\by^{t-1} - \by^{t-2}}^2 .
\end{align*}}
Rearranging the terms yields:
  \begin{align*}  
   &\E [ \hat{F}^{t+1} - \hat{F}^{t}] \\
   & \geq \frac{C_1}{2}\eta_{\by}^2\E \norm{\hat{\bg}_\by(\bx^t,\by^t)}^2 +  \frac{C_2}{2} \eta^2_{\bx}\E \norm{\hat{\bg}_\bx(\bx^t,\by^t)}^2 + \pare{\frac{1}{8\eta_{\by}}-( 2 C_1\eta_{\by}^2 L_f^2G_g^2 + 2C_2\eta_{\bx}^2 L_f^2 )4(1-\beta)^2 G_g^2 }  \E \norm{\by^t - \by^{t-1}}^2 \\
    &  -\pare{2C_1\eta_{\by}^2 L_f^2G_g^2 + 2C_2\eta^2_{\bx} L_f^2 +   L_f^2 G_g^2 {\eta_{\by}} + \eta_{\bx} L_f^2}(1-\beta )^2   \E \norm{\bz^{t} -g(\by^{t-1}) }^2 +  \frac{1}{8\eta_{\by}}\E \norm{\by^{t-1} - \by^{t-2}}^2 \\
    &    -  \pare{  4C_1 \eta_{\by}^2 (G_g^2+G_f^2)  + 2C_2 \eta_{\bx}^2 } \frac{\sigma^2}{B} -  \pare{2\beta^2 L_f^2 G_g^2 {\eta_{\by}}+ 4\beta^2C_1  \eta_{\by}^2 L_f^2 G_g^2 +4\beta^2 C_2\eta^2_\bx L_f^2  }\frac{\sigma^2}{M}-\frac{72L_f^2}{\mu^2\eta_\bx}  \beta \frac{\sigma^2}{M}\\
    & -\pr{1-\frac{\beta}{2}}^2  \frac{4L_f^2}{\mu^2 \eta_{\bx}}   \E \norm{  \bz^{t } - \bz^{t-1}  }^2.
\end{align*}

Recall our definition of potential function $\tilde F^{t+1}$:
\begin{align*}
     \tilde{F}^{t+1} :=  \hat{F}^{t+1} -   \E \norm{\bz^{t+1}  -  g(\by^t))  }^2 -\frac{(1-\frac{\beta}{2})^2}{1-(1-\frac{\beta}{2})^2}\frac{4L_f^2}{\mu^2\eta_{\bx}}   \E \norm{  \bz^{t+1 } - \bz^{t}  }^2.
\end{align*}
Hence we have:
  \begin{align*}  
   &\E [ \tilde{F}^{t+1} - \tilde{F}^{t}] \\
   & \geq \frac{C_1}{2}\eta_{\by}^2\E \norm{\hat{\bg}_\by(\bx^t,\by^t)}^2 +  \frac{C_2}{2} \eta^2_{\bx}\E \norm{\hat{\bg}_\bx(\bx^t,\by^t)}^2 + \pare{\frac{1}{8\eta_{\by}}-( 2 C_1\eta_{\by}^2 L_f^2G_g^2 + 2C_2\eta_{\bx}^2 L_f^2 )4(1-\beta)^2 G_g^2 }  \E \norm{\by^t - \by^{t-1}}^2  \\
   &  -  \E \norm{\bz^{t+1}  -  g(\by^t))  }^2+  \E \norm{\bz^{t} -g(\by^{t-1}) }^2\\
    &   -\pare{2C_1\eta_{\by}^2 L_f^2G_g^2 + 2C_2\eta^2_{\bx} L_f^2 +   L_f^2 G_g^2 {\eta_{\by}} + \eta_{\bx} L_f^2}(1-\beta )^2   \E \norm{\bz^{t} -g(\by^{t-1}) }^2 +  \frac{1}{8\eta_{\by}}\E \norm{\by^{t-1} - \by^{t-2}}^2  \\
    &     -  \pare{  4C_1 \eta_{\by}^2 (G_g^2+G_f^2)  + 2C_2 \eta_{\bx}^2 } \frac{\sigma^2}{B} -  \pare{2\beta^2 L_f^2 G_g^2 {\eta_{\by}}+ 4\beta^2C_1  \eta_{\by}^2 L_f^2 G_g^2 +4\beta^2 C_2\eta^2_\bx L_f^2  }\frac{\sigma^2}{M}-\frac{72L_f^2}{\mu^2\eta_\bx}  \beta \frac{\sigma^2}{M} \\
    &  -\pr{1-\frac{\beta}{2}}^2  \frac{4L_f^2}{\mu^2 \eta_{\bx}}   \E \norm{  \bz^{t } - \bz^{t-1}  }^2 -\frac{(1-\frac{\beta}{2})^2}{1-(1-\frac{\beta}{2})^2}\frac{4L_f^2}{\mu^2\eta_{\bx}}    \E \norm{  \bz^{t+1 } - \bz^{t}  }^2 +\frac{(1-\frac{\beta}{2})^2}{1-(1-\frac{\beta}{2})^2}\frac{4L_f^2}{\mu^2\eta_{\bx}}    \E \norm{  \bz^{t } - \bz^{t-1}  }^2.
\end{align*}

Now we plug in Lemma~\ref{lem:tracking error} and~\ref{lem:second order tracking}:

 {  \begin{align*}  
   \E [ \tilde{F}^{t+1}  &- \tilde{F}^{t}]  \geq \frac{C_1}{2}\eta_{\by}^2\E \norm{\hat{\bg}_{\by}(\bx^t,\by^t)}^2 +  \frac{C_2}{2} \eta^2_{\bx}\E \norm{\hat{\bg}_{\bx}(\bx^t,\by^t)}^2 \\
   & + \pare{\frac{1}{8\eta_{\by}}-4( 2 C_1\eta_{\by}^2 L_f^2G_g^2 + 2C_2\eta_{\bx}^2 L_f^2 + 1)(1-\beta)^2 G_g^2   -\frac{(1-\frac{\beta}{2})^2}{1-(1-\frac{\beta}{2})^2}\frac{4L_f^2}{\mu^2\eta_{\bx}} \pr{ \frac{12}{\beta} + 6\beta  } G_g^2   }  \E \norm{\by^t - \by^{t-1}}^2 \\ 
    &  +\pare{1-\pare{2C_1\eta_{\by}^2 L_f^2G_g^2 + 2C_2\eta^2_{\bx} L_f^2 +   L_f^2 G_g^2 {\eta_{\by}} + \eta_{\bx} L_f^2+1}(1-\beta )^2 }   \E \norm{\bz^{t} -g(\by^{t-1}) }^2\\
    &  +  \pare{\frac{1}{8\eta_{\by}} - \frac{(1-\frac{\beta}{2})^2}{1-(1-\frac{\beta}{2})^2}\frac{4L_f^2}{\mu^2\eta_{\bx}}\frac{6}{\beta}G_g^2  }\E \norm{\by^{t-1} - \by^{t-2}}^2 \\
    &    -  \pare{  4C_1 \eta_{\by}^2 (G_g^2+G_f^2)  + 2C_2 \eta_{\bx}^2 } \frac{\sigma^2}{B} -  \pare{2\beta^2 L_f^2 G_g^2 {\eta_{\by}}+ 4\beta^2C_1  \eta_{\by}^2 L_f^2 G_g^2 +4\beta^2 C_2\eta^2_\bx L_f^2  }\frac{\sigma^2}{M}-\frac{72L_f^2}{\mu^2\eta_\bx}  \beta \frac{\sigma^2}{M}  \\
    & + (1-\frac{\beta}{2})^2\frac{4L_f^2}{\mu^2\eta_{\bx}} \underbrace{\pare{\frac{1}{1-(1-\frac{\beta}{2})^2}-\frac{(1-\frac{\beta}{2})^2}{1-(1-\frac{\beta}{2})^2}-1} }_{=0}     \E \norm{  \bz^{t } - \bz^{t-1}  }^2.
\end{align*}}

By our choice of $\eta_{\by}$, we know that
\begin{align*}
   \frac{1}{8\eta_{\by}}-4( 2 C_1\eta_{\by}^2 L_f^2G_g^2 + 2C_2\eta_{\bx}^2 L_f^2 + 1)(1-\beta)^2 G_f^2   -\frac{(1-\frac{\beta}{2})^2}{1-(1-\frac{\beta}{2})^2}\frac{4L_f^2}{\mu^2\eta_{\bx}} \pr{ \frac{12}{\beta} + 6\beta  } G_g^2  \geq 0,\\
   1-\pare{2C_1\eta_{\by}^2 L_f^2G_g^2 + 2C_2\eta^2_{\bx} L_f^2 +   L_f^2 G_g^2 {\eta_{\by}} + \eta_{\bx} L_f^2+1}(1-\beta )^2 \geq 0,\\
   \frac{1}{8\eta_{\by}} - \frac{(1-\frac{\beta}{2})^2}{1-(1-\frac{\beta}{2})^2}\frac{4L_f^2}{\mu^2\eta_{\bx}}\frac{6}{\beta}G_g^2 \geq 0.
\end{align*}
 Now we can have the clean bound:
 
 \begin{align*}  
  \E [ \tilde{F}^{t+1} - \tilde{F}^{t}]   & \geq  \frac{C_1}{2}\eta_{\by}^2\E \norm{\hat{\bg}_{\by}(\bx^t,\by^t)}^2 +  \frac{C_2}{2} \eta^2_{\bx}\E \norm{\hat{\bg}_{\bx}(\bx^t,\by^t)}^2 \\ 
    & \quad  -  \pare{  4C_1 \eta_{\by}^2 (G_g^2+G_f^2)  + 2C_2 \eta_{\bx}^2 } \frac{\sigma^2}{B} -  \pare{2\beta^2 L_f^2 G_g^2 {\eta_{\by}}+ 4\beta^2C_1  \eta_{\by}^2 L_f^2 G_g^2 +4\beta^2 C_2\eta^2_\bx L_f^2  }\frac{\sigma^2}{M} -\frac{72L_f^2}{\mu^2\eta_\bx}  \beta \frac{\sigma^2}{M} .   
\end{align*}
  \end{proof}

 \end{lemma}

 \subsection{Proof of Theorem~\ref{thm:SCNC}}
Evoking Lemma~\ref{lemma: potential 2 scnc} and summing the inequality from $t=1$ to $T$ yields:
 \begin{align*}  
   & \frac{\E [ \tilde{F}^{T} - \tilde{F}^{0}]}{T} +  \pare{  4C_1 \eta_{\by}^2 (G_g^2+G_f^2)  + 2C_2 \eta_{\bx}^2 } \frac{\sigma^2}{B} +  \pare{2\beta^2 L_f^2 G_g^2 {\eta_{\by}}+ 4\beta^2C_1  \eta_{\by}^2 L_f^2 G_g^2 +4\beta^2 C_2\eta^2_\bx L_f^2  }\frac{\sigma^2}{M}+\frac{72L_f^2}{\mu^2\eta_\bx}  \beta \frac{\sigma^2}{M}  \\
    &\geq \frac{1}{T} \sum_{t=0}^{T-1} \frac{C_1}{2}\eta_{\by}^2\E \norm{\hat{\bg}_{\by}(\bx^t,\by^t)}^2 +  \frac{C_2}{2} \eta^2_{\bx}\E \norm{\hat{\bg}_{\bx}(\bx^t,\by^t)}^2 \\
    &\geq  \min\left\{ \frac{C_1}{2}\eta_{\by}^2, \frac{C_2}{2} \eta_{\bx}^2 \right\} \frac{1}{T} \sum_{t=0}^{T-1} \E \norm{\hat{\bg}_{\by}(\bx^t,\by^t)}^2 +  \E \norm{\hat{\bg}_{\bx}(\bx^t,\by^t)}^2
\end{align*}
We compute the upper and lower bound of $C_1$ and $C_2$. For $C_1$
\begin{align*}
    C_1 = \frac{1}{4\eta_{\by}} - \frac{L}{2} -\eta_{\bx} L^2 - 4L_f^2G_f^2G_g^2\eta_{\by} - \frac{168  G_g^2L_f^2}{\mu^2\eta_{\bx} \beta} - \frac{48   G_g^2 L_f^2 \beta}{\mu^2 \eta_{\bx}}
\end{align*}
The upper bound $C_1 \leq \frac{1}{4\eta_{\by}}$ holds trivially. For lower bound,
since we choose
\begin{align*}
    \eta_{\by} = \min\cbr{\frac{1}{20 L}, \frac{1}{40 \eta_{\bx} L^2}, \frac{1}{10 L_f^2 G_g^2 G_f^2}, \frac{\mu^2\eta_{\bx} \beta}{3840 G_g^2 L_f^2}}
\end{align*}
we can conclude that $C_1 \geq \frac{1}{8\eta_{\by}}$. 

For $C_2$:
\begin{align*}
      &C_2 = \frac{1}{\eta_{\bx}} + \frac{3}{2}\mu - \frac{\eta_{\bx} L^2}{2} - \frac{1}{4\eta_{\bx} (1-\beta)^2}-\frac{2L^2_f}{\mu}.
\end{align*}
The upper bound $C_2 \leq \frac{1}{\eta_{\bx}} + \frac{3\mu}{2}$ holds trivially. For lower bound, since we choose:
\begin{align*}
    \eta_{\bx} = \min \cbr{ \frac{1}{\sqrt{6}L^2}, \frac{\mu}{12 L^2_f}}, \beta =0.1 \leq 1-\frac{\sqrt{3}}{2}
\end{align*}
it holds that $C_2 \geq \frac{1}{3\eta_{\bx}}$.

 Since $\frac{1}{8\eta_{\by}} \leq C_1 \leq  \frac{1}{4\eta_{\by}}$ and $ \frac{1}{3\eta_{\bx}}\leq C_2 \leq \frac{1}{\eta_{\bx}} + \frac{3\mu}{2} \leq \frac{2}{\eta_{\bx}}$, we have:
  \begin{align*}  
   & \frac{\E [ \tilde{F}^{T} - \tilde{F}^{0}]}{T} +  \pare{  \eta_{\by} (G_g^2+G_f^2)  + 4 \eta_{\bx} } \frac{\sigma^2}{B} +  \pare{2\beta^2 L_f^2 G_g^2 {\eta_{\by}}+ \beta^2 \eta_{\by} L_f^2 G_g^2 +8\beta^2\eta_\bx L_f^2  }\frac{\sigma^2}{M}+\frac{72L_f^2}{\mu^2\eta_\bx}  \beta \frac{\sigma^2}{M} \\ 
    &\geq  \min\left\{ \frac{1}{16}\eta_{\by} , \frac{1}{6} \eta_{\bx}  \right\} \frac{1}{T} \sum_{t=0}^{T-1} \E \norm{\hat{\bg}_{\by}(\bx^t,\by^t)}^2 +  \E \norm{\hat{\bg}_{\bx}(\bx^t,\by^t)}^2
\end{align*}
We then need to verify the conditions ~\ref{eq:cond 1}, \ref{eq:cond2} and \ref{eq:cond3} in Lemma~\ref{lemma: potential 2 scnc} can hold under our choice of $\eta_{\bx}$, $\eta_{\by}$ and $\beta$. To guarantee (\ref{eq:cond 1}) holding, we need:
\begin{align*}
    \eta_{\by} \leq \min\cbr{ \frac{1}{64 G_f^2 (\eta_x L_f^2 G_g^2 + 4\eta_x L_f^2 + 1)}, 2\times 10^{-4}\times\frac{\mu^2 \eta_{\bx}}{L_f^2 G_g^2} }.
\end{align*}
To guarantee condition (\ref{eq:cond2}) holding, we need 
\begin{align*}
    \eta_\bx \leq \frac{23}{1700 L_f^2}, \eta_{\by} \leq \frac{23}{500 L_f^2 G_g^2}.
\end{align*}

To guarantee condition (\ref{eq:cond3}) holding, we need:
\begin{align*}
    \eta_{\by} \leq \frac{\mu^2 \eta_{\bx}}{14400 L_f^2 G_g^2}.
\end{align*}
Next we examine how large the $\E [ \tilde{F}^{T} - \tilde{F}^{0}]$ is.
By definition of potential function, we have:
\begin{align*}
        \E [ \tilde{F}^{T} - \tilde{F}^{0}] &=  \hat{F}^{T} -   \E \norm{\bz^{T} -g(\by^{T-1}) }^2 -\frac{(1-\frac{\beta}{2})^2}{1-(1-\frac{\beta}{2})^2}\frac{4L_f^2}{\mu^2\eta_{\bx}}   \E \norm{  \bz^{t} - \bz^{t-1}  }^2\\
        &\quad - \pare{\hat{F}^{0} -   \E \norm{\bz^{0} -g(\by^{-1}) }^2 -\frac{(1-\frac{\beta}{2})^2}{1-(1-\frac{\beta}{2})^2}\frac{4L_f^2}{\mu^2\eta_{\bx}}   \E \norm{  \bz^{0} - \bz^{-1}  }^2} \\
        &\leq \hat{F}^{T}  - \hat{F}^{0} +     \E \norm{\bz^{0} -g(\by^{-1})}^2 + \frac{(1-\frac{\beta}{2})^2}{\beta - \frac{\beta^2}{4}}\frac{4L_f^2}{\mu^2\eta_{\bx}}   \E \norm{  \bz_j^{0} - \bz_j^{-1}  }^2.
\end{align*}
By convention $\bx^{-1} = \bx^0$, $\by^{0} = \by^{-1}$, and our choice $\E\norm{\bz_j^{0} -g(\by^{0})}^2 \leq O(1)$, , we have
 \begin{align*}
        \E [ \tilde{F}^{T} - \tilde{F}^{0}] 
        \leq \hat{F}^{T}  - \hat{F}^{0} + O(1).
\end{align*} 
 Next we examine how large the $\E [ \hat{F}^{T} - \hat{F}^{0}]$ is.

 \begin{align*}
     \E[\hat{F}^{T}  - \hat{F}^{0}]& = F(\bx^{T},\by^T) + s^{T} - \pare{\frac{1}{4\eta_{\by}} + 4 L_f^2 G_f^2 G_g^2 \eta_{\by}  + \frac{\eta_{\bx} L^2}{2} + \frac{168  G_g^2L_f^2}{\mu^2 \eta_{\bx} \beta} + \frac{48  G_g^2 L_f^2\beta}{\mu^2 \eta_{\bx}}} \norm{\by^T - \by^{T-1}}^2\\
    &  \quad-\pare{\frac{1}{8\eta_{\by}}+\frac{84 G_g^2L_f^2}{\mu^2 \eta_{\bx} \beta}}\norm{\by^{T-1} - \by^{T-2}}^2 + \pare{\frac{7}{2\eta_{\bx}} + \mu - \frac{\eta_{\bx} L^2}{2} - \frac{2L^2_f}{\mu}}\norm{\bx^{T} - \bx^{T-1}}^2  \\
    &\quad-F(\bx^{0},\by^{0}) - s^{0} + \pare{\frac{1}{4\eta_{\by}} + 4 L_f^2 G_f^2 G_g^2 \eta_{\by}  + \frac{\eta_{\bx} L^2}{2} + \frac{168  G_g^2L_f^2}{\mu^2 \eta_{\bx} \beta} + \frac{48  G_g^2 L_f^2\beta}{\mu^2 \eta_{\bx}}} \norm{\by^{0} - \by^{-1}}^2\\
    &  \quad+\pare{\frac{1}{8\eta_{\by}}+\frac{84 G_g^2}{\mu^2 \eta_{\bx} \beta}}\norm{\by^{-1} - \by^{-2}}^2 - \pare{\frac{7}{2\eta_{\bx}} + \mu - \frac{\eta_{\bx} L^2}{2} -\frac{2L_f^2}{\mu}}\norm{\bx^{0} - \bx^{-1}}^2 \\
    & \leq F_{\max}     + \pare{\frac{7}{2\eta_{\bx}} + \mu - \frac{\eta_{\bx} L^2}{2} -\frac{2L_f^2}{\mu}}\norm{\bx^{T} - \bx^{T-1}}^2 ,
 \end{align*}
where we used the convention $\bx^{-1} = \bx^0$ and $\by^{-1} = \by^{-2}$. Notice that $\E\norm{\bx^{T} - \bx^{T-1}}^2 \leq \eta_{\bx}^2\E \norm{\nabla_1 f(\bx^{T-1}, \bz^T;\cB^{T-1}) + \nabla h(\bx^{T-1}) }^2 \leq 2\eta_{\bx}^2(G_f^2+G_h^2) $,
we have $\E[\hat{F}^{T}  - \hat{F}^{0}] \leq F_{\max} + 7\eta_{\bx}  (G_f^2+G_h^2) + 2\mu \eta_{\bx}^2 (G_f^2+G_h^2)$.
 

 
Using the fact that $\beta - \frac{\beta^2}{2} \geq \frac{1}{2}\beta$ when $\beta \leq 2$ yield:
 
 \begin{align*}
     &\frac{1}{T} \sum_{t=1}^{T} \norm{\nabla G(\bx^t,\by^t)}^2 
     \leq \max\cbr{\frac{16}{\eta_{\by}}, \frac{6}{\eta_{\bx}}} \cdot \frac{ F_{\max} + 7\eta_{\bx}  (G_f^2+G_h^2) + 2\mu \eta_{\bx}^2 (G_f^2+G_h^2) +O(1)  }{T} \\
     &\quad +\max\cbr{\frac{16}{\eta_{\by}}, \frac{6}{\eta_{\bx}}}\cdot\pare{  \pare{  \eta_{\by} (G_g^2+G_f^2)  + 4 \eta_{\bx} } \frac{\sigma^2}{B} +  \pare{2\beta^2 L_f^2 G_g^2 {\eta_{\by}}+ \beta^2 \eta_{\by} L_f^2 G_g^2 +8\beta^2\eta_\bx L_f^2  }\frac{\sigma^2}{M} +\frac{72L_f^2}{\mu^2\eta_\bx}  \beta \frac{\sigma^2}{M}  } \\
     &\leq  O\pr{ \frac{ F_{\max} +  \eta_{\bx}  (G_f^2+G_h^2) +  \mu \eta_{\bx}^2 (G_f^2+G_h^2) }{\eta_{\by} T} }   +   O\pare{  G_g^2 + G_f^2 +  L_f^2 G_g^2 + \frac{\eta_{\bx} L_f^2}{\eta_{\by}} + \frac{\kappa^2}{\eta_{\bx}} } \max\cbr{\frac{\sigma^2}{B},\frac{\sigma^2}{M}  }.
 \end{align*} 
To guarantee RHS is less than $\epsilon^2$, we need:
 $T = O\pr{\frac{F_{\max} }{\eta_{\by} \epsilon^2}}$, and 
 \begin{align*}
    M&=\Theta \pare{ \max\cbr{ \frac{\kappa^3 L \sigma^2}{\epsilon^2}   1} }, B= \Theta \pare{ \max\cbr{   \frac{\kappa^2 L_f^2\sigma^2}{\epsilon^2}, 1} }, \beta = 0.1, \\
     \eta_\bx &= \Theta\pare{\min \cbr{ \frac{1}{ L^2 }, \frac{\mu}{ L^2_f}}},  \eta_{\by} = \Theta\pare{ \min\cbr{\frac{1}{  L}, \frac{1}{ \eta_{\bx} L^2}, \frac{1}{  L_f^2 G_g^2 G_f^2}, \frac{\mu^2\eta_{\bx}  }{  G_g^2 L_f^2}}}
 \end{align*}
 which yields the total gradient complexity:
 \begin{align*}
     O\pr{ \max\cbr{   \frac{\kappa^2 L^2_f \sigma^2}{\epsilon^2}, 1} \cdot \frac{\kappa^3 F_{\max}}{\epsilon^2} }. 
 \end{align*}
 
\subsection{Proof of Convex-nonconcave Setting} \label{app: proof CNC}

\begin{lemma}[Connection between stationary measure and iterates] \label{lem:stationary measure cnc}
For Algorithm~\ref{algorithm: CODA-Dual}, under assumptions of Theorem~\ref{thm: CNC}, if we define \[\hat{\bg}_{\bx}(\bx^t,\by) = \frac{1}{\eta_{\bx}}  \pare{ \bx^t - \cP_{\cX}\pare{\bx^t-\eta_{\bx} \nabla_{\bx} F(\bx^t,\by) - \eta_{\bx}\alpha_t \bx^t}},
 \hat{\bg}_{\by}(\bx,\by)=  \frac{1}{\eta_{\by}}  \pare{ \by - \cP_{\cY}\pare{\by+\eta_{\by} \nabla_{\by} F(\bx,\by)}  }\], 
  then the following statement holds:
\begin{align*}
   &  \E\norm{\hat{\bg}_{\bx}(\bx^t,\by^t) }^2 \leq  \frac{2}{\eta_{\bx}^2} \E\norm{  \bx^t - \bx^{t+1} }^2 +  8L_f^2 \E\norm{\bz^{t+1} - g(\by^t)}^2 + 8\frac{\sigma^2}{B}  , \\
    & \E\norm{\hat{\bg}_{\by}(\bx^t,\by^t) }^2 \leq \frac{2}{\eta_{\by}^2}  \E\norm{  \by^t - \by^{t+1}}^2+4 L_f^2 G_g^2\E\norm{  \bz^{t+1}  - g(\by^t   )  }^2 + \frac{8(G_f^2+G_g^2)\sigma^2}{B} .
\end{align*}
    \begin{proof}
        We begin with proving the first statement. According to updating rule we have
        \begin{align*}
           \E \norm{\hat{\bg}_{\bx}(\bx^t,\by^t)) }^2 &= \E\norm{  \frac{1}{\eta_{\bx}}  \pare{ \bx^t - \cP_{\cX}\pare{\bx^t-\eta_{\bx} \nabla_{\bx} F(\bx^t,\by^t)- \eta_{\bx}\alpha_t \bx^t} }  }^2 \\
            &\leq \frac{2}{\eta_{\bx}^2} \E\norm{  \bx^t - \cP_{\cX}\pare{\bx^t-\eta_{\bx} \bg_{\bx}^t }}^2 + \frac{2}{\eta_{\bx}^2} \E\norm{  \cP_{\cX}\pare{\bx^t-\eta \bg_{\bx}^t}- \cP_{\cX}\pare{\bx^t-\eta \nabla_{\bx} F(\bx^t,\by^t)- \eta_{\bx}\alpha_t \bx^t}}^2\\
            &\leq \frac{2}{\eta_{\bx}^2} \E\norm{  \bx^t - \bx^{t+1}}^2 +  2\E\norm{  \frac{1}{B}\sum_{(\zeta,\xi)\in\cB^t} \nabla_1 f(\bx^t,\bz^{t+1};\zeta)   -  \nabla_{\bx} f(\bx^t,\bz^{t+1})}^2\\
             &\leq \frac{2}{\eta_{\bx}^2} \E\norm{  \bx^t - \bx^{t+1}}^2 +  8\E\norm{ \frac{1}{B}\sum_{(\zeta,\xi)\in\cB^t}\nabla_1 f(\bx^t, \bz^{t+1};\zeta )-  \nabla_1 f(\bx^t, \bz^{t+1})}^2  \\
             & \quad +  8\E\norm{  \nabla_1 f(\bx^t, \bz^{t+1})-  \nabla_1 f(\bx^t, g(\by^{t}))}^2\\
                &\leq \frac{2}{\eta_{\bx}^2} \E\norm{  \bx^t - \bx^{t+1} }^2 +  8L_f^2 \E\norm{\bz^{t+1} - g(\by^t)}^2 + 8\frac{\sigma^2}{B} \\
        \end{align*}
        where at last step we apply Lemma~\ref{lem:grad bound}.
        Similarly, for the second statement we have:
           \begin{align*}
            \E\norm{\hat{\bg}_{\by}(\bx^t,\by^t)}^2 &=  \E\norm{  \frac{1}{\eta_{\by}}  \pare{ \by^t - \cP_{\cY}\pare{\by^t+\eta \nabla_{\by} F(\bx^t,\by^t)}}  }^2 \\
            &\leq \frac{2}{\eta_{\by}^2}  \E\norm{  \by^t - \cP_{\cY}\pare{\by^t+\eta_{\by} \bg_{\by}^t }}^2 + \frac{2}{\eta_{\by}^2}  \E\norm{  \cP_{\cY}\pare{\by^t+\eta_{\by} \bg_{\by}^t}- \cP_{\cY}\pare{\by^t+\eta_{\by} \nabla_{\by} F(\bx^t,\by^t)}}^2\\
            &\leq \frac{2}{\eta_{\by}^2}  \E\norm{  \by^t - \by^{t+1}} +  2 \E\norm{     \bg_{\by}^t-  \nabla_{\by} F(\bx^t,\by^t)}^2\\
            &= \frac{2}{\eta_{\by}^2}  \E\norm{  \by^t - \by^{t+1}} +  2 \E\norm{\frac{1}{B}\sum_{(\zeta,\xi)\in\cB^t} \nabla_2 f(\bx^t, \bz^{t+1};\zeta ) \nabla g(\by^t;\xi )   -  \nabla_2 f(\bx^t, g(\by^t )) \nabla g(\by^t) }^2\\
            &= \frac{2}{\eta_{\by}^2}  \E\norm{  \by^t - \by^{t+1}}   +  4 \E\norm{ \frac{1}{B}\sum_{(\zeta,\xi)\in\cB^t}\nabla_2 f(\bx^t, \bz^{t+1};\zeta ) \nabla g(\by^t;\xi ) -\nabla_2 f(\bx^t, \bz^{t+1}) \nabla g(\by^t )}^2 \\
            &\quad +4\E\norm{\nabla_2 f(\bx^t, \bz^{t+1}) \nabla g(\by^t )   -  \nabla_2 f(\bx^t, g(\by^t )) \nabla g(\by^t) }^2\\
            &\leq \frac{2}{\eta_{\by}^2}  \E\norm{  \by^t - \by^{t+1}}^2+4 L_f^2 G_g^2\E\norm{  \bz^{t+1}  - g(\by^t   )  }^2 + \frac{8(G_f^2+G_g^2)\sigma^2}{B} .
        \end{align*}
    \end{proof}
\end{lemma}

 \begin{proposition}\label{prop:lower bound}
 For Algorithm~\ref{algorithm: CODA-Dual}, under assumptions of Theorem~\ref{thm: CNC}, the following statement holds:
 \begin{align*}
     &\inprod{\nabla F_{t-1}(\bx^t,\by^{t-1}) -\nabla F_{t-1}(\bx^{t-1},\by^{t-1}) }{\bx^{t} - \bx^{t-1}} \\
     &\geq \frac{1}{2(L'  +\alpha_{t-1})}\norm{\nabla F_{t-1}(\bx^t,\by^{t-1}) -\nabla F_{t-1}(\bx^{t-1},\by^{t-1}) }^2   +\frac{\alpha_{t-1} L'}{2(L' + \alpha_{t-1})}\norm{\bx^{t} - \bx^{t-1}}^2 .
 \end{align*}
     \begin{proof}
         
    \begin{align*}
       &\inprod{\nabla F_{t-1}(\bx^t,\by^{t-1}) -\nabla F_{t-1}(\bx^{t-1},\by^{t-1}) }{\bx^t - \bx^{t-1}}\\
       =  &\frac{\frac{L'}{2}}{L' + \alpha_{t-1}}\underbrace{\inprod{\nabla F_{t-1}(\bx^t,\by^{t-1}) -\nabla F_{t-1}(\bx^{t-1},\by^{t-1}) }{\bx^t - \bx^{t-1}} }_{A}\\
       &+\frac{\frac{L'}{2}+\alpha_{t-1}}{L'+ \alpha_{t-1}}\underbrace{\inprod{\nabla F_{t-1}(\bx^t,\by^{t-1}) -\nabla F_{t-1}(\bx^{t-1},\by^{t-1}) }{\bx^t - \bx^{t-1}}}_{B}\\
        \geq  &\frac{1}{2(L' +\alpha_{t-1})}\norm{\nabla F_{t-1}(\bx^t,\by^{t-1}) -\nabla F_{t-1}(\bx^{t-1},\by^{t-1}) }^2   +\frac{\alpha_{t-1} L'}{2(L'+ \alpha_{t-1})}\norm{\bx^{t} - \bx^{t-1}}^2 
    \end{align*}
    where the last inequality is due to \\
    (1) $F_{t-1}(\bx,\by^{t-1})$ is  $L'$ smooth and hence $A \geq \frac{1}{L'} \norm{\nabla F_{t-1}(\bx^t,\by^{t-1}) -\nabla F_{t-1}(\bx^{t-1},\by^{t-1}) }^2 $; and \\
    (2) $F_{t-1}(\bx,\by^{t-1})$ is $\alpha_{t-1}$ strongly convex and $B \geq \alpha_{t-1}\norm{\bx^{t} - \bx^{t-1}}^2  $,

         \end{proof}

  \end{proposition}

\begin{lemma}\label{lm:descent lemma cnc}
For Algorithm~\ref{algorithm: CODA-Dual}, under assumptions of Theorem~\ref{thm: CNC}, then the following statement holds:
 {\small \begin{align*}
      \E[ F(\bx^{t+1},\by^{t+1}) - F(\bx^{t},\by^t)] & \geq -\pr{\frac{L^2 \eta_{\bx}}{2} + 4(1-\beta)^2   L_f^2  G_g^4 {\eta_{\by}}  }\E \|\by^t - \by^{t-1} \|^2 +\pare{ \frac{1}{2\eta_{\by}} - \frac{L}{2} - \frac{L^2\eta_{\bx}}{2}}\E\|  \by^{t+1}- \by^t \|^2\\
            & \quad +\pr{\frac{\alpha_{t-1} L'}{2(L'+\alpha_{t-1})} - \frac{1}{2\eta_\bx}}\E\norm{\bx^{t} -\bx^{t-1}}^2 + \pare{\frac{\alpha_{t-1}}{2}   - \frac{3}{2\eta_{\bx}} -  \frac{1}{4\eta_{\bx} (1-\beta)^2} } \E \|\bx^{t+1} -\bx^{t}\|^2 \\
            &\quad - \frac{1}{2} \alpha_{t-1}  (\E\norm{\bx^{t+1}}^2 - \E\norm{\bx^t}^2) -\pr{2\eta_{\bx}L_f^2 +  L_f^2 G_g^2 {\eta_{\by}}   }(1-\beta)^2  \E\norm{\bz^{t} - g(\by^{t-1}) }^2   \\
            &\quad-2 \eta_{\bx} (1-\beta)^2  \frac{\sigma^2}{B}-  \pare{  2{\eta_{\by}}(G_g^2+ G_f^2)\frac{\sigma^2}{B} +2\beta^2 L_f^2 G_g^2 {\eta_{\by}}\frac{\sigma^2}{M}  }. 
        \end{align*}}
    \begin{proof}
    According to updating rule:
     \begin{align}
     \left \langle  \bg_{\bx}^t +\frac{1}{\eta_{\bx}} (\bx^{t+1}-\bx^{t} ) , \bx -\bx^{t+1} \right \rangle \geq 0 \\
     \left \langle  \bg_{\bx}^t +\frac{1}{\eta_{\bx}} (\bx^{t+1}-\bx^{t} ) , \bx^t -\bx^{t+1} \right \rangle \geq 0 \\
     \left \langle  \bg_{\bx}^{t-1} +\frac{1}{\eta_{\bx}} (\bx^{t}-\bx^{t-1} ) , \bx^{t+1} -\bx^{t} \right \rangle \geq 0 \label{eq: cnc optimality condition}. 
 \end{align}
        Due to the $ \alpha_t$-strong convexity of $F_t(\cdot, \by)$, together with (\ref{eq: cnc optimality condition}) we have
        \begin{align*}
            F_t(\bx^{t+1},\by^t) - F_t(\bx^{t},\by^t) &\geq \inprod{\nabla F_t(\bx^{t},\by^t)}{\bx^{t+1} - \bx^{t}} + \frac{1}{2}\alpha_t \norm{ \bx^{t+1} - \bx^{t} }^2 \\ 
            & \geq \inprod{\nabla F_t(\bx^{t},\by^t) - \bg_{\bx}^{t-1}}{\bx^{t+1} - \bx^{t}} - \frac{1}{\eta_{\bx}}\inprod{\bx^t - \bx^{t-1}}{\bx^{t+1} - \bx^t}+\frac{1}{2} \alpha_t \norm{ \bx^{t+1} - \bx^{t} }^2 \\
            & = \underbrace{\inprod{\nabla F_t(\bx^{t},\by^t) - \nabla F_{t-1}(\bx^{t-1},\by^{t-1})}{\bx^{t+1} - \bx^{t}}}_{A} - \underbrace{\frac{1}{\eta_{\bx}}\inprod{\bx^t - \bx^{t-1}}{\bx^{t+1} - \bx^t} }_{B}\\
            & \quad + \underbrace{\inprod{\nabla F_{t-1}(\bx^{t-1},\by^{t-1})-\bg_{\bx}^{t-1}}{\bx^{t+1} - \bx^t}}_{C}+ \frac{1}{2}\alpha_t \norm{ \bx^{t+1} - \bx^{t} }^2
        \end{align*}
        Similar to the proof of Lemma~\ref{lm:descent lemma}, for A, we observe that:
 \begin{align}
      \left \langle  \nabla_{\bx} F_t(\bx^{t},\by^t) - \nabla_{\bx} F_{t-1}(\bx^{t-1}, \by^{t-1}), \bx^{t+1} -\bx^{t} \right \rangle &=  \underbrace{\left \langle  \nabla_{\bx} F_t(\bx^{t},\by^t) - \nabla_{\bx} F_{t-1}(\bx^{t},\by^{t-1}), \bx^{t+1} -\bx^{t} \right \rangle}_{\text{\ding{172}}} \nonumber\\
     &+ \underbrace{\left \langle  \nabla_{\bx} F_{t-1}(\bx^{t},\by^{t-1}) - \nabla_{\bx} F_{t-1}(\bx^{t-1}, \by^{t-1}), \Delta  \right \rangle}_{\text{\ding{173}}} \nonumber\\
     &+ \underbrace{\left \langle  \nabla_{\bx} F_{t-1}(\bx^{t},\by^{t-1}) - \nabla_{\bx} F_{t-1}(\bx^{t-1},\by^{t-1}), \bx^{t} -\bx^{t-1} \right \rangle}_{\text{\ding{174}}}\nonumber
 \end{align}
 where $ \Delta = (\bx^{t+1}- \bx^{t})-(\bx^{t}-\bx^{t-1}).  $
For \ding{172}, recall the definition of $F_t$:
\begin{align*}
    \text{\ding{172}} &= \inprod{\nabla_{\bx} F(\bx^{t},\by^t) - \nabla_{\bx} F(\bx^{t},\by^{t-1}) + \alpha_t \bx^t - \alpha_{t-1}\bx^{t}}{\bx^{t+1} -\bx^{t} }\\
    &= \inprod{\nabla_{\bx} F(\bx^{t},\by^t) - \nabla_{\bx} F(\bx^{t},\by^{t-1})  }{\bx^{t+1} -\bx^{t} }+\inprod{  \alpha_t \bx^t - \alpha_{t-1}\bx^{t}}{\bx^{t+1} -\bx^{t} }\\
    & \geq -\frac{L^2 \eta_{\bx}}{2}  \|\by^t - \by^{t-1} \|^2 - \frac{1}{2\eta_{\bx}} \|\bx^{t+1} -\bx^{t}\|^2 + \pr{\alpha_t - \alpha_{t-1}}\inprod{   \bx^t  }{\bx^{t+1}} - \pr{\alpha_t - \alpha_{t-1}}\norm{\bx^t}^2\\
    & = -\frac{L^2 \eta_{\bx}}{2}  \|\by^t - \by^{t-1} \|^2 - \frac{1}{2\eta_{\bx}} \|\bx^{t+1} -\bx^{t}\|^2  - \frac{1}{2}\pr{\alpha_{t-1} - \alpha_t}(\norm{\bx^{t+1}}^2 - \norm{\bx^t}^2) + \frac{1}{2}\pr{\alpha_{t-1} -\alpha_t } \|\bx^{t+1} -\bx^{t}\|^2 .
\end{align*}
For \ding{173}, we use Cauchy-Schwartz:
\begin{align*}
    \text{\ding{173}}  \geq - \frac{\eta_{\bx} }{2}  \|\nabla F_{t-1}(\bx^t,\by^{t-1}) -\nabla F_{t-1}(\bx^{t-1},\by^{t-1})\|^2  - \frac{1}{2\eta_{\bx}}  \|\Delta \|^2 .
\end{align*}

For \ding{174}, we use Proposition~\ref{prop:lower bound}:
  \begin{align*}
     &\inprod{\nabla F_{t-1}(\bx^t,\by^{t-1}) -\nabla F_{t-1}(\bx^{t-1},\by^{t-1}) }{\bx^{t} - \bx^{t-1}} \\
     &\geq \frac{1}{2(L' +\alpha_{t-1})}\norm{\nabla F_{t-1}(\bx^t,\by^{t-1}) -\nabla F_{t-1}(\bx^{t-1},\by^{t-1}) }^2  + \frac{\alpha_{t-1} L'}{2(L'+\alpha_{t-1})} \|\bx^{t} -\bx^{t-1}\|^2
 \end{align*}
 Since $\frac{1}{2(L' +\alpha_{t-1})} - \frac{\eta_\bx}{2} \geq 0 $, putting pieces together yields:
 \begin{align*}
     A &\geq -\frac{L^2 \eta_{\bx}}{2}  \|\by^t - \by^{t-1} \|^2 - \frac{1}{2\eta_{\bx}} \|\bx^{t+1} -\bx^{t}\|^2  - \frac{1}{2}\pr{\alpha_{t-1} - \alpha_t}(\norm{\bx^{t+1}}^2 - \norm{\bx^t}^2) + \frac{1}{2}\pr{\alpha_{t-1} -\alpha_t } \|\bx^{t+1} -\bx^{t}\|^2 \\
     &\quad -  \frac{1}{2\eta_{\bx}}  \|\Delta \|^2.
 \end{align*}

For B, we have:
\begin{align*}
    \frac{1}{\eta_{\bx}}\inprod{\bx^t - \bx^{t-1}}{\bx^{t+1} - \bx^t} =   \frac{1}{2\eta_{\bx}} \norm{\bx^t - \bx^{t-1}}^2 + \frac{1}{2\eta_{\bx}} \norm{\bx^{t+1} - \bx^t}^2 - \frac{1}{2\eta_{\bx}} \norm{\Delta}^2.
\end{align*}
 
 For C, we have:
 \begin{align*}
   \E \left \langle \bg_{\bx}^{t-1}- \nabla_{\bx} F_{t-1}(\bx^{t-1},\by^{t-1}), \bx^{t+1} -\bx^{t} \right \rangle&\geq - \eta_{\bx} (1-\beta)^2 \E \norm{ \bg_{\bx}^{t-1}- \nabla_{\bx} F_{t-1}(\bx^{t-1},\by^{t-1})}^2\\
   &\quad - \frac{1}{4\eta_{\bx} (1-\beta)^2}\E\norm{ \bx^{t+1} -\bx^{t}}^2\\ 
   & \geq -2\eta_{\bx} (1-\beta)^2 L_f^2 \E\norm{\bz^{t} - g(\by^{t-1}) }^2 -2 \eta_{\bx} (1-\beta)^2  \frac{\sigma^2}{B}\\
   &\quad - \frac{1}{4\eta_{\bx} (1-\beta)^2} \E\norm{ \bx^{t+1} -\bx^{t}}^2,
 \end{align*}
 where we plug in Lemma~\ref{lem:grad bound} at last step.

Putting lower bound of A, B and C together yields:
 \begin{align*}
            F_t(\bx^{t+1},\by^t) - F_t(\bx^{t},\by^t) & \geq -\frac{L^2 \eta_{\bx}}{2}  \|\by^t - \by^{t-1} \|^2+ \pare{\frac{1}{2} \alpha_{t-1}  - \frac{1}{\eta_{\bx}} -  \frac{1}{4\eta_{\bx} (1-\beta)^2} } \|\bx^{t+1} -\bx^{t}\|^2 \\
            & \quad \pr{\frac{\alpha_{t-1} L'}{2(L'+\alpha_{t-1})}- \frac{1}{2\eta_\bx}}\norm{\bx^{t} -\bx^{t-1}}^2 - \frac{1}{2}\pr{\alpha_{t-1} - \alpha_t}(\norm{\bx^{t+1}}^2 - \norm{\bx^t}^2)   \\
            &\quad-2\eta_{\bx} (1-\beta)^2 L_f^2 \E\norm{\bz^{t} - g(\by^{t-1}) }^2 -2 \eta_{\bx} (1-\beta)^2  \frac{\sigma^2}{B}.
        \end{align*}

Since $F_t(\bx^{t+1},\by^t) - F_t(\bx^{t},\by^t) = F(\bx^{t+1},\by^t) - F(\bx^{t},\by^t) + \frac{1}{2}\alpha_t (\norm{\bx^{t+1}}^2 - \norm{\bx^{t}}^2)$, and also recall in Lemma~\ref{lem:dual descent},
 \begin{align*}
         \E[F(\bx^{t+1}, \by^{t+1})  - F(\bx^{t+1}, \by^t)   ] &\geq \pare{ \frac{1}{2\eta_{\by}} - \frac{L}{2} - \frac{L^2\eta_{\bx}}{2}}\E\|  \by^{t+1}- \by^t \|^2 -  \frac{1}{2\eta_{\bx}}\E \norm{\bx^{t+1} - \bx^{t}}^2  \\
       &-  \pare{ {\eta_{\by}} 2(G_g^2+ G_f^2)\frac{\sigma^2}{B} +2\beta^2 L_f^2 G_g^2 {\eta_{\by}}\frac{\sigma^2}{M}  } \\
    &-   L_f^2 G_g^2 {\eta_{\by}}    (1-\beta)^2 \E \norm{\bz^{t} - g(\by^{t-1}) }^2 - 4(1-\beta)^2   L_f^2   G_g^4 {\eta_{\by}} \E\norm{\by^t - \by^{t-1}}^2,
 \end{align*}
we can conclude the proof:
 {\small \begin{align*}
    \E[ F(\bx^{t+1},\by^{t+1}) - F(\bx^{t},\by^t)] & \geq -\pr{\frac{L^2 \eta_{\bx}}{2} + 4(1-\beta)^2   L_f^2  G_g^4 {\eta_{\by}}  }\E \|\by^t - \by^{t-1} \|^2 +\pare{ \frac{1}{2\eta_{\by}} - \frac{L}{2} - \frac{L^2\eta_{\bx}}{2}}\E\|  \by^{t+1}- \by^t \|^2\\
            & \quad +\pr{\frac{\alpha_{t-1} L'}{2(L'+\alpha_{t-1})} - \frac{1}{2\eta_\bx}}\E\norm{\bx^{t} -\bx^{t-1}}^2 + \pare{\frac{\alpha_{t-1}}{2}   - \frac{3}{2\eta_{\bx}} -  \frac{1}{4\eta_{\bx} (1-\beta)^2} } \E \|\bx^{t+1} -\bx^{t}\|^2 \\
            &\quad - \frac{1}{2} \alpha_{t-1}  (\E\norm{\bx^{t+1}}^2 - \E\norm{\bx^t}^2) -\pr{2\eta_{\bx}L_f^2 +  L_f^2 G_g^2 {\eta_{\by}}   }(1-\beta)^2  \E\norm{\bz^{t} - g(\by^{t-1}) }^2   \\
            &\quad-2 \eta_{\bx} (1-\beta)^2  \frac{\sigma^2}{B}-  \pare{  2{\eta_{\by}}(G_g^2+ G_f^2)\frac{\sigma^2}{B} +2\beta^2 L_f^2 G_g^2 {\eta_{\by}}\frac{\sigma^2}{M}  }. 
 \end{align*}}
    \end{proof}
\end{lemma}
 \begin{lemma}\label{lem:potential cnc}
 For Algorithm~\ref{algorithm: CODA-Dual}, under assumptions of Theorem~\ref{thm: CNC}, define the following auxiliary function  
\begin{align*}
  \hat{F}^{t+1} := & F(\bx^{t+1},\by^{t+1}) + s^{t+1} - \pare{\frac{1}{4\eta_{\by}} + 4 L_f^2  G_g^4 \eta_{\by}  + \frac{\eta_{\bx} L^2}{2} + \frac{768  G_g^2L_f^2}{\alpha_{t+1}^2 \eta_{\bx} \beta} + \frac{576  G_g^2 L_f^2\beta}{\alpha_{t+1}^2 \eta_{\bx}}  } \norm{\by^{t+1} - \by^t}^2\\
    &  -\pare{\frac{1}{8\eta_{\by}}+\frac{768 G_g^2L_f^2}{\alpha_{t+1}^2 \eta_{\bx} \beta}}\norm{\by^t - \by^{t-1}}^2 + \frac{7}{2\eta_{\bx}} \norm{\bx^{t+1} - \bx^t}^2+ \frac{\alpha_t}{2}\norm{\bx^{t+1}}^2 ,
\end{align*}
where $s^{t+1} :=  -\frac{8}{\eta^2_{\bx} \alpha_{t+1}}\|\bx^{t+1}-\bx^t\|^2$
then the following statement holds true:
\begin{align*}
  \E [ \hat{F}^{t+1} - \hat{F}^{t}] &\geq C_1^t\E \norm{\by^{t+1} - \by^t}^2 + \pare{\frac{1}{8\eta_{\by}}-\frac{768 G_g^2L_f^2}{\alpha_{t+1}^2 \eta_{\bx} \beta}}\E \norm{\by^t - \by^{t-1}}^2 + \frac{1}{8\eta_{\by}} \E \norm{\by^{t-1} - \by^{t-2}}^2 \\
      &+ C_2^t \E \norm{\bx^{t+1} - \bx^t}^2 -   (1-\beta)^2 \pare{  L_f^2 G_g^2 {\eta_{\by}} + 2\eta_{\bx} L_f^2}  \E \norm{\bz^{t} -g(\by^{t-1}) }^2   + \frac{1}{2} \pr{ \alpha_{t} - \alpha_{t-1}}   \norm{\bx^{t+1}}^2 \\
    &     -\pare{1-\frac{\beta}{2}}^2  \frac{64L_f^2}{\alpha_t^2 \eta_{\bx}}    \E \norm{  \bz^{t} - \bz^{t-1}  }^2 + \frac{8}{\eta_\bx}\pr{\frac{\alpha_{t}}{\alpha_{t+1}} -\frac{\alpha_{t-1}}{\alpha_t}} \norm{\bx^{t+1} }^2 \\
    & -2 \eta_{\bx} (1-\beta)^2  \frac{\sigma^2}{B}-  \pare{  2{\eta_{\by}}(G_g^2+ G_f^2)\frac{\sigma^2}{B} +2\beta^2 L_f^2 G_g^2 {\eta_{\by}}\frac{\sigma^2}{M} +  \frac{1152 L_f^2 \beta}{\alpha_t^2 \eta_{\bx}}    \frac{\sigma^2}{M} },
\end{align*}

where
\begin{align}
    &C_1^t = \frac{1}{4\eta_{\by}} - \frac{L}{2} - \eta_{\bx} L^2  -  4 L_f^2  G_g^4 \eta_{\by}   - \frac{768  G_g^2L_f^2}{\alpha_{t+1}^2 \eta_{\bx} \beta} - \frac{576  G_g^2 L_f^2\beta}{\alpha_{t+1}^2 \eta_{\bx}}, \label{eq: def C1 cnc} \\
    &C_2^t = \pare{\frac{1}{ \eta_{\bx}}   +\frac{\alpha_{t-1}}{2} -\frac{1}{4\eta_{\bx} (1-\beta )^2}  }.\label{eq: def C2 cnc}
\end{align}

\begin{proof}
Define $\Delta^{t+1}:= \bx^{t+1} - \bx^{t} - (\bx^{t} - \bx^{t-1})$.    According to (\ref{eq:smpe_update_2}) and (\ref{eq:smpe_update_3}):

\begin{align*}
    \frac{1}{\eta}  \left \langle \Delta^{t+1}, \bx^{t} - \bx^{t+1}  \right \rangle &\geq \inprod{\bg_{\bx}^t - \bg_{\bx}^{t-1} }{\bx^{t+1} - \bx^{t} }. 
\end{align*}

If we define $\bg^t_{\bx}(\bx,\bz) = \nabla_1 f(\bx,\bz;\zeta^t) + \nabla h(\bx) + \alpha_t \bx$.  
\begin{align*}
    \frac{1}{\eta}  \left \langle \Delta^{t+1}, \bx^{t} - \bx^{t+1}  \right \rangle &\geq \inprod{\bg_{\bx}^t - \bg^t_{\bx}(\bx^{t},\bz^t ) }{\bx^{t+1} - \bx^{t} }\\ 
   &\quad  +  \inprod{ \bg^t_{\bx}(\bx^{t},\bz^t)- \bg_{\bx}^{t-1}}{\bx^{t+1} - \bx^{t} }\\
     &\geq \inprod{\bg_{\bx}^t - \bg_{\bx}^{t-1}( \bx^{t},\bz^{t}) }{\bx^{t+1} - \bx^{t} }\\ 
   &\quad +  \inprod{ \bg_{\bx}^{t-1}(\bx^{t},\bz^{t})- \bg_{\bx}^{t-1}}{\bx^{t+1} - \bx^{t} - (\bx^{t} - \bx^{t-1}) }\\
   &\quad +  \inprod{ \bg_{\bx}^{t-1}( \bx^{t},\bz^{t})- \bg_{\bx}^{t-1}}{  (\bx^{t} - \bx^{t-1}) }\\
   &\geq -\frac{L_f^2}{2c_t}  \norm{  \bz^{t+1} -  \bz^{t}  }^2 -  \frac{c_t}{2} \norm{\bx^{t+1} - \bx^{t} }^2 + \frac{\alpha_t - \alpha_{t-1}}{2} \pr{\norm{\bx^{t+1}}^2 - \norm{\bx^t}^2}\\ 
   &\quad   - \frac{\eta_{\bx}  }{2}  \norm{   \bg_{\bx}^{t-1}( \bx^{t},\bz^{t})- \bg_{\bx}^{t-1} }^2 - \frac{1}{2\eta_{\bx}} \norm{\Delta^{t+1} }^2\\
   & \quad + \frac{1}{2(L' +\alpha_{t-1})}\norm{\bg_{\bx}^{t-1}( \bx^{t},\bz^{t})- \bg_{\bx}^{t-1} }^2   +\frac{\alpha_{t-1} L'}{2(L'+ \alpha_{t-1})}\norm{\bx^{t} - \bx^{t-1}}^2    \\
\end{align*}
where we use Proposition~\ref{prop:lower bound} at the last step.
Since $\frac{1}{\eta_{\bx} }\left \langle \Delta^{t+1}, \bx^{t} - \bx^{t+1}  \right \rangle = \frac{1}{2\eta_{\bx}}\norm{\bx^{t} - \bx^{t-1}}^2 - \frac{1}{2\eta_{\bx}}\norm{\bx^{t+1} - \bx^{t}}^2 - \frac{1}{2\eta_{\bx}} \norm{\Delta^{t+1}}^2$, and $\frac{1}{2(L' +\alpha_{t-1})} - \frac{\eta_\bx}{2} \geq 0$, and $L' \geq \alpha_{t-1}$
 we have:
\begin{align*}
   \frac{1}{2\eta_{\bx}}\norm{\bx^{t} - \bx^{t-1}}^2 - \frac{1}{2\eta_{\bx}}\norm{\bx^{t+1} - \bx^{t}}^2  
   &\geq  -\frac{L_f^2}{2c_t}  \norm{  \bz^{t+1} -  \bz^{t}  }^2 -  \frac{c_t}{2} \norm{\bx^{t+1} - \bx^{t} }^2 + \frac{\alpha_t - \alpha_{t-1}}{2} \pr{\norm{\bx^{t+1}}^2 - \norm{\bx^t}^2}\\ 
   &\quad +\frac{\alpha_{t-1}  }{4}\norm{\bx^{t} - \bx^{t-1}}^2 .
\end{align*}

Rearranging terms yields:
\begin{align*}
      \frac{1}{2\eta_{\bx}}\norm{\bx^{t}  - \bx^{t-1}}^2 &+ \frac{\alpha_t - \alpha_{t-1}}{2}  \norm{\bx^{t}}^2    - \frac{1}{2\eta_{\bx}}\norm{\bx^{t+1} - \bx^{t}}^2  - \frac{\alpha_t - \alpha_{t-1}}{2}  \norm{\bx^{t+1}}^2 \\
   &\geq  -\frac{L_f^2}{2c_t}  \norm{  \bz^{t+1} -  \bz^{t}  }^2 -  \frac{c_t}{2} \norm{\bx^{t+1} - \bx^{t} }^2       +\frac{\alpha_{t}  }{4}\norm{\bx^{t} - \bx^{t-1}}^2,
\end{align*}
where we use the fact $\alpha_{t-1} \geq \alpha_t$.

Multiplying both sides with $\frac{16}{\alpha_{t} \eta_{\bx}}$ yields:
\begin{align}
     \frac{8}{\alpha_{t}\eta^2_{\bx}}\norm{\bx^{t}  - \bx^{t-1}}^2 &  + \frac{8}{\eta_\bx}\pr{1-\frac{  \alpha_{t-1}}{\alpha_t} } \norm{\bx^{t}}^2    - \frac{8}{\alpha_{t}\eta^2_{\bx}}\norm{\bx^{t+1} - \bx^{t}}^2  - \frac{8}{\eta_\bx}\pr{1-\frac{  \alpha_{t-1}}{\alpha_t} }  \norm{\bx^{t+1}}^2 \\
   &\geq  -\frac{8 L_f^2}{\alpha_{t} c_t \eta_\bx}  \norm{  \bz^{t+1} -  \bz^{t}  }^2 -  \frac{8 c_t}{\alpha_{t}  \eta_\bx } \norm{\bx^{t+1} - \bx^{t} }^2       +\frac{4  }{   \eta_\bx}\norm{\bx^{t} - \bx^{t-1}}^2 .\label{eq: potential1 1}
\end{align} 
Recall the definition of $s^t$, and choose $c_t = \frac{\alpha_t}{8}$:
\begin{align*}
    s^{t+1} - s^t \geq& -\frac{64 L_f^2}{\alpha^2_{t}   \eta_\bx}  \norm{  \bz^{t+1} -  \bz^{t}  }^2 -  \frac{1}{  \eta_\bx } \norm{\bx^{t+1} - \bx^{t} }^2     +\frac{4  }{   \eta_\bx}\norm{\bx^{t} - \bx^{t-1}}^2 + \frac{8}{\eta_\bx}\pr{\frac{\alpha_{t}}{\alpha_{t+1}} -\frac{\alpha_{t-1}}{\alpha_t}} \norm{\bx^{t+1} }^2.
\end{align*}
Recall that Lemma~\ref{lm:descent lemma cnc} gives:
{\small \begin{align*}
      \E[ F(\bx^{t+1},\by^{t+1}) - F(\bx^{t},\by^t)] & \geq -\pr{\frac{L^2 \eta_{\bx}}{2} + 4(1-\beta)^2   L_f^2  G_g^4 {\eta_{\by}}  }\E \|\by^t - \by^{t-1} \|^2 +\pare{ \frac{1}{2\eta_{\by}} - \frac{L}{2} - \frac{L^2\eta_{\bx}}{2}}\E\|  \by^{t+1}- \by^t \|^2\\
            & \quad +\pr{\frac{\alpha_{t-1} L'}{2(L'+\alpha_{t-1})} - \frac{1}{2\eta_\bx}}\E\norm{\bx^{t} -\bx^{t-1}}^2 + \pare{\frac{\alpha_{t-1}}{2}   - \frac{3}{2\eta_{\bx}} -  \frac{1}{4\eta_{\bx} (1-\beta)^2} } \E \|\bx^{t+1} -\bx^{t}\|^2 \\
            &\quad - \frac{1}{2} \alpha_{t-1}  (\E\norm{\bx^{t+1}}^2 - \E\norm{\bx^t}^2) -\pr{2\eta_{\bx}L_f^2 +  L_f^2 G_g^2 {\eta_{\by}}   }(1-\beta)^2  \E\norm{\bz^{t} - g(\by^{t-1}) }^2   \\
            &\quad-2 \eta_{\bx} (1-\beta)^2  \frac{\sigma^2}{B}-  \pare{  2{\eta_{\by}}(G_g^2+ G_f^2)\frac{\sigma^2}{B} +2\beta^2 L_f^2 G_g^2 {\eta_{\by}}\frac{\sigma^2}{M}  }. 
        \end{align*}}
Evoking Lemma~\ref{lm:descent lemma cnc} together with (\ref{eq: potential1 1}) yields:
 \begin{align*}
    &\E[F(\bx^{t+1}, \by^{t+1}) + s^{t+1}  - (F(\bx^{t}, \by^t) +s^{t})  ] \\
    &\geq \pare{ \frac{1}{2\eta_{\by}} - \frac{L}{2} -  \frac{L^2\eta_{\bx}}{2} }\E\|  \by^{t+1}- \by^t \|^2- \pare{ 4(1-\beta )^2 L_f^2 G_g^4 \eta_{\by} + \frac{L^2\eta_{\bx}}{2} } \E\norm{\by^t - \by^{t-1}}^2 \\
    &+ \left (\frac{\alpha_{t-1}}{2} - \frac{3}{ 2\eta_{\bx}}-\frac{1}{4\eta_{\bx} (1-\beta )^2} - \frac{1}{  \eta_\bx}    \right)\E \norm{\bx^{t+1} - \bx^{t}}^2 +\pr{\frac{4}{\eta_\bx}+\frac{\alpha_{t-1} L'}{2(L'+\alpha_{t-1})} - \frac{1}{2\eta_\bx}}\E \norm{\bx^{t} - \bx^{t-1}}^2   \\
    &-    (1-\beta)^2 \pare{ L_f^2 G_g^2 {\eta_{\by}} + 2\eta_{\bx} L_f^2} \E \norm{\bz^{t} -g(\by^{t-1}) }^2       - \frac{64L_f^2}{\alpha_t^2 \eta_{\bx}}  \E\norm{\bz^{t+1} -  \bz^{t}}^2\\
    & - \frac{1}{2} \alpha_{t-1}  (\norm{\bx^{t+1}}^2 - \norm{\bx^t}^2)+ \frac{8}{\eta_\bx}\pr{\frac{\alpha_{t}}{\alpha_{t+1}} -\frac{\alpha_{t-1}}{\alpha_t}} \norm{\bx^{t+1} }^2\\
    & -2 \eta_{\bx} (1-\beta)^2  \frac{\sigma^2}{B}-  \pare{  2{\eta_{\by}}(G_g^2+ G_f^2)\frac{\sigma^2}{B} +2\beta^2 L_f^2 G_g^2 {\eta_{\by}}\frac{\sigma^2}{M}  }.
 \end{align*}

 Now we plug in Lemma~\ref{lem:second order tracking} to replace $\E\norm{\bz^{t+1} -  \bz^{t}}^2$ together with the fact that $\frac{1}{\beta} \geq 1$ and get:
 
 \begin{align*}
    &\E[F(\bx^{t+1}, \by^{t+1}) + s^{t+1}  - (F(\bx^{t}, \by^t) +s^{t})  ] \\
    &\geq \pare{ \frac{1}{2\eta_{\by}} - \frac{L}{2} -  \frac{L^2\eta_{\bx}}{2} }\E\|  \by^{t+1}- \by^t \|^2- \pare{ 4(1-\beta )^2 L_f^2 G_g^4 \eta_{\by} + \frac{L^2\eta_{\bx}}{2}+ \frac{768 L_f^2 G_g^2}{\alpha_t^2 \eta_{\bx} \beta}+ \frac{576 L_f^2 G_g^2 \beta}{\alpha_t^2 \eta_{\bx}} } \E\norm{\by^t - \by^{t-1}}^2 \\
    &+ \left (\frac{\alpha_{t-1}}{2} - \frac{5}{ 2\eta_{\bx}}-\frac{1}{4\eta_{\bx} (1-\beta )^2}    \right)\E \norm{\bx^{t+1} - \bx^{t}}^2 + \frac{7}{2\eta_\bx} \E \norm{\bx^{t} - \bx^{t-1}}^2   \\
    &-    (1-\beta)^2 \pare{ L_f^2 G_g^2 {\eta_{\by}} + 2\eta_{\bx} L_f^2} \E \norm{\bz^{t} -g(\by^{t-1}) }^2   - \frac{64L_f^2}{\alpha_t^2 \eta_{\bx}}    \pare{1-\frac{\beta}{2}}^2 \E \norm{  \bz^{t } - \bz^{t-1}  }^2-  \frac{768L_f^2G_g^2}{\alpha_t^2 \eta_{\bx}\beta}     \norm{    \by^{t-1} - \by^{t-2}   }^2  \\
    & - \frac{1}{2} \alpha_{t-1}  (\norm{\bx^{t+1}}^2 - \norm{\bx^t}^2) + \frac{8}{\eta_\bx}\pr{\frac{\alpha_{t}}{\alpha_{t+1}} -\frac{\alpha_{t-1}}{\alpha_t}} \norm{\bx^{t+1} }^2 \\
    & -2 \eta_{\bx} (1-\beta)^2  \frac{\sigma^2}{B}-  \pare{  2{\eta_{\by}}(G_g^2+ G_f^2)\frac{\sigma^2}{B} +2\beta^2 L_f^2 G_g^2 {\eta_{\by}}\frac{\sigma^2}{M} +  \frac{64L_f^2}{\alpha_t^2 \eta_{\bx}}  \cdot 18\beta \frac{\sigma^2}{M} }.
 \end{align*}  
 
Recall our definition of potential function $\hat F^{t+1}$
\begin{align*}
    \hat{F}^{t+1} := & F(\bx^{t+1},\by^{t+1}) + s^{t+1} - \pare{\frac{1}{4\eta_{\by}} + 4 L_f^2  G_g^4 \eta_{\by}  + \frac{\eta_{\bx} L^2}{2} + \frac{768  G_g^2L_f^2}{\alpha_{t+1}^2 \eta_{\bx} \beta} + \frac{576  G_g^2 L_f^2\beta}{\alpha_{t+1}^2 \eta_{\bx}}  } \norm{\by^{t+1} - \by^t}^2\\
    &  -\pare{\frac{1}{8\eta_{\by}}+\frac{768 G_g^2L_f^2}{\alpha_{t+1}^2 \eta_{\bx} \beta}}\norm{\by^t - \by^{t-1}}^2 + \frac{7}{2\eta_{\bx}} \norm{\bx^{t+1} - \bx^t}^2+ \frac{\alpha_t}{2}\norm{\bx^{t+1}}^2 .
\end{align*} 
We conclude that:
{\small \begin{align*}
  \E [ \hat{F}^{t+1} - \hat{F}^{t}] &\geq \pare{\frac{1}{4\eta_{\by}} - \frac{L}{2} - \eta_{\bx} L^2  -  4 L_f^2  G_g^4 \eta_{\by}   - \frac{768  G_g^2L_f^2}{\alpha_{t+1}^2 \eta_{\bx} \beta} - \frac{576  G_g^2 L_f^2\beta}{\alpha_{t+1}^2 \eta_{\bx}}}\E \norm{\by^{t+1} - \by^t}^2\\
  &+ \pare{\frac{1}{8\eta_{\by}}-\frac{768 G_g^2L_f^2}{\alpha_{t+1}^2 \eta_{\bx} \beta}}\E \norm{\by^t - \by^{t-1}}^2 + \frac{1}{8\eta_{\by}}\E \norm{\by^{t-1} - \by^{t-2}}^2 + \pare{\frac{1}{ \eta_{\bx}}   +\frac{\alpha_{t-1}}{2} -\frac{1}{4\eta_{\bx} (1-\beta )^2}  } \E \norm{\bx^{t+1} - \bx^t}^2\\
      &-   (1-\beta)^2 \pare{  L_f^2 G_g^2 {\eta_{\by}} + 2\eta_{\bx} L_f^2}  \E \norm{\bz^{t} -g(\by^{t-1}) }^2   + \frac{1}{2} \pr{ \alpha_{t} - \alpha_{t-1}}   \norm{\bx^{t+1}}^2 \\
    &     -\pare{1-\frac{\beta}{2}}^2  \frac{64L_f^2}{\alpha_t^2 \eta_{\bx}}    \E \norm{  \bz^{t} - \bz^{t-1}  }^2 + \frac{8}{\eta_\bx}\pr{\frac{\alpha_{t}}{\alpha_{t+1}} -\frac{\alpha_{t-1}}{\alpha_t}} \norm{\bx^{t+1} }^2 \\
    & -2 \eta_{\bx} (1-\beta)^2  \frac{\sigma^2}{B}-  \pare{  2{\eta_{\by}}(G_g^2+ G_f^2)\frac{\sigma^2}{B} +2\beta^2 L_f^2 G_g^2 {\eta_{\by}}\frac{\sigma^2}{M}+  \frac{1152 L_f^2 \beta}{\alpha_t^2 \eta_{\bx}}    \frac{\sigma^2}{M}  }.
\end{align*}}

\end{proof}
\end{lemma}

 \begin{lemma}\label{lemma: potential 2 cnc}

Let $C_1, C_2$ be defined in (\ref{eq: def C1 cnc}) and (\ref{eq: def C2 cnc}). If the following conditions hold:
 
\begin{align}
    \frac{1}{8\eta_{\by}}-\frac{768 G_g^2L_f^2}{\alpha_{t+1}^2 \eta_{\bx} \beta}-16(1-\beta)^2G_g^2 (C^t_1\eta_\by L_f^2 G_g^2 + C^t_2 \eta_\bx^2 L_f^2)  -\frac{(1-\frac{\beta}{2})^2}{1-(1-\frac{\beta}{2})^2}\frac{64L_f^2}{\alpha_t^2\eta_\bx} \pr{ \frac{12}{\beta} + 6\beta  } G_g^2   \geq 0, \label{eq:cond 1 cnc}\\
   1- (1-\beta)^2 \pare{  L_f^2 G_g^2 {\eta_{\by}} + 2\eta_{\bx} L_f^2 +C^t_1 4 \eta_{\by}^2 L_f^2 G_g^2 + 4C^t_2\eta^2_\bx L_f^2 + 1}  \geq 0, \label{eq:cond 2 cnc}\\
  \frac{1}{8\eta_{\by}} - \frac{(1-\frac{\beta}{2})^2}{1-(1-\frac{\beta}{2})^2}\frac{64L_f^2}{\alpha_T^2\eta_{\bx}}\frac{6}{\beta}G_g^2   \geq 0 \label{eq:cond 3 cnc}.
\end{align}
 then for Algorithm~\ref{algorithm: CODA-Dual}, under assumptions of Theorem~\ref{thm: CNC}, then the following statement holds:
  \begin{align*}  
  \E [ \tilde{F}^{t+1} - \tilde{F}^{t}]   & \geq  \frac{C^t_1}{2}\eta_{\by}^2\E \norm{\hat{\bg}_{\by}(\bx^t,\by^t)}^2 +  \frac{C^t_2}{2} \eta^2_{\bx}\E \norm{\hat{\bg}_{\bx}(\bx^t,\by^t)}^2 - 2C^t_2 \eta_\bx^2 \alpha_t^2 \norm{\bx^t}^2 \\
  & \quad+ \frac{1}{2} \pr{ \alpha_{t} - \alpha_{t-1}}   \norm{\bx^{t+1}}^2 +\frac{8}{\eta_\bx}\pr{\frac{\alpha_{t}}{\alpha_{t+1}} -\frac{\alpha_{t-1}}{\alpha_t}} \norm{\bx^{t+1} }^2  \\ 
    & \quad -\pare{2 \eta_{\bx}  +4C^t_2 \eta_{\bx}^2 +4C^t_1 \eta_{\by}^2  (G_f^2 + G_g^2) +2{\eta_{\by}}(G_g^2+ G_f^2)  } \frac{\sigma^2}{B} \\
    &\quad - \pare{2\beta^2 L_f^2 G_g^2 {\eta_{\by}}  + 4C^t_1 \beta^2 \eta_{\by}^2 L_f^2 G_g^2 + 8C^t_2\beta^2 \eta_{\bx}^2 L_f^2  +  \frac{1152 L_f^2 \beta}{\alpha_t^2 \eta_{\bx}}     }\frac{\sigma^2}{M} ,
\end{align*}
where 
\begin{align*}
     \tilde{F}^{t+1} :=  \hat{F}^{t+1} -   \E \norm{\bz^{t+1}  -  g(\by^t))  }^2 -\frac{(1-\frac{\beta}{2})^2}{1-(1-\frac{\beta}{2})^2}\frac{64L_f^2}{\alpha_T^2\eta_{\bx}}   \E \norm{  \bz^{t+1 } - \bz^{t}  }^2.
\end{align*}
 \begin{proof}
     
 According to Lemma~\ref{lem:potential cnc}:
\begin{align*}
  \E [ \hat{F}^{t+1} - \hat{F}^{t}] &\geq C_1^t\E \norm{\by^{t+1} - \by^t}^2 + \pare{\frac{1}{8\eta_{\by}}-\frac{768 G_g^2L_f^2}{\alpha_{t+1}^2 \eta_{\bx} \beta}}\E \norm{\by^t - \by^{t-1}}^2+ \frac{1}{8\eta_{\by}} \E \norm{\by^{t-1} - \by^{t-2}}^2 \\
      &+ C_2^t \E \norm{\bx^{t+1} - \bx^t}^2 -   (1-\beta)^2 \pare{  L_f^2 G_g^2 {\eta_{\by}} + 2\eta_{\bx} L_f^2}  \E \norm{\bz^{t} -g(\by^{t-1}) }^2   + \frac{1}{2} \pr{ \alpha_{t} - \alpha_{t-1}}   \norm{\bx^{t+1}}^2 \\
    &     -\pare{1-\frac{\beta}{2}}^2  \frac{64L_f^2}{\alpha_t^2 \eta_{\bx}}    \E \norm{  \bz^{t} - \bz^{t-1}  }^2 + \frac{8}{\eta_\bx}\pr{\frac{\alpha_{t}}{\alpha_{t+1}} -\frac{\alpha_{t-1}}{\alpha_t}} \norm{\bx^{t+1} }^2 \\
    & -2 \eta_{\bx} (1-\beta)^2  \frac{\sigma^2}{B}-  \pare{  2{\eta_{\by}}(G_g^2+ G_f^2)\frac{\sigma^2}{B} +2\beta^2 L_f^2 G_g^2 {\eta_{\by}}\frac{\sigma^2}{M} +  \frac{1152 L_f^2 \beta}{\alpha_t^2 \eta_{\bx}}    \frac{\sigma^2}{M} },
\end{align*}

where
\begin{align}
    &C_1^t = \frac{1}{4\eta_{\by}} - \frac{L}{2} - \eta_{\bx} L^2  -  4 L_f^2  G_g^4 \eta_{\by}   - \frac{768  G_g^2L_f^2}{\alpha_{t+1}^2 \eta_{\bx} \beta} - \frac{576  G_g^2 L_f^2\beta}{\alpha_{t+1}^2 \eta_{\bx}},   \\
    &C_2^t = \pare{\frac{1}{ \eta_{\bx}}   +\frac{\alpha_{t-1}}{2} -\frac{1}{4\eta_{\bx} (1-\beta )^2}  }. 
\end{align} 
 Now we plug in Lemma~\ref{lem:stationary measure cnc}:
 \begin{align*}
  \E [ \hat{F}^{t+1} - \hat{F}^{t}] \geq &C^t_1\pare{\frac{1}{2}\eta_{\by}^2\E \norm{\hat{\bg}_{\by}(\bx^t,\by^t)}^2  - 4\eta_{\by}^2  \frac{(G_f^2 + G_g^2)\sigma^2}{B} - 2\eta_{\by}^2G_g^2 L_f^2  \E\norm{ \bz^{t+1}  -  g(\by^t))    }^2}   \\
  & + C^t_2 \pare{\frac{1}{2} \eta^2_{\bx}\E \norm{\hat{\bg}_{\bx}(\bx^t,\by^t)}^2 -  4\eta_\bx^2L_f^2 \E\norm{\bz^{t+1} - g(\by^t)}^2- 4\eta_\bx^2\frac{\sigma^2}{B}  } \\ 
& -   (1-\beta)^2 \pare{  L_f^2 G_g^2 {\eta_{\by}} +2 \eta_{\bx} L_f^2}  \E \norm{\bz^{t} -g(\by^{t-1}) }^2  \\
&+ \frac{1}{2} \pr{ \alpha_{t} - \alpha_{t-1}}  \E \norm{\bx^{t+1}}^2 +\frac{8}{\eta_\bx}\pr{\frac{\alpha_{t}}{\alpha_{t+1}} -\frac{\alpha_{t-1}}{\alpha_t}} \E\norm{\bx^{t+1} }^2\\
&-\pr{1-\frac{\beta}{2}}^2  \frac{64L_f^2}{\alpha_t^2 \eta_{\bx}}   \E \norm{  \bz^{t } - \bz^{t-1}  }^2+ \pare{\frac{1}{8\eta_{\by}}-\frac{768 G_g^2L_f^2}{\alpha_{t+1}^2 \eta_{\bx} \beta}}\E \norm{\by^t - \by^{t-1}}^2 + \frac{1}{8\eta_{\by}}\E \norm{\by^{t-1} - \by^{t-2}}^2 \\
&  -2 \eta_{\bx} (1-\beta)^2  \frac{\sigma^2}{B}-  \pare{  2{\eta_{\by}}(G_g^2+ G_f^2)\frac{\sigma^2}{B} +2\beta^2 L_f^2 G_g^2 {\eta_{\by}}\frac{\sigma^2}{M}  +  \frac{1152 L_f^2 \beta}{\alpha_t^2 \eta_{\bx}}    \frac{\sigma^2}{M} }.
\end{align*}

Plugging in Lemma~\ref{lem:tracking error} yields:

  \begin{align*}
  \E [ \hat{F}^{t+1} - \hat{F}^{t}] \geq& \frac{C^t_1}{2}\eta_{\by}^2\E \norm{\hat{\bg}_{\by}(\bx^t,\by^t)}^2 +  \frac{C_2^t}{2} \eta^2_{\bx}\E \norm{\hat{\bg}_\bx(\bx^t,\by^t)}^2  \\
  &-2C^t_1   \eta_{\by}^2 L_f^2 G_g^2 \pare{  (1-\beta)^2 \E \norm{\bz^{t} -g(\by^{t-1}) }^2 + 4(1-\beta)^2 G_g^2 \E\norm{\by^t - \by^{t-1}}^2 + 2\beta^2 \frac{\sigma^2}{M} }   \\ 
  &-4C^t_2\eta^2_\bx L_f^2   \pare{(1-\beta)^2 \E \norm{\bz^{t} -g(\by^{t-1}) }^2 + 4(1-\beta)^2 G_g^2 \E\norm{\by^t - \by^{t-1}}^2 + 2\beta^2 \frac{\sigma^2}{M} } \\
& -   (1-\beta)^2 \pare{  L_f^2 G_g^2 {\eta_{\by}} + 2\eta_{\bx} L_f^2}  \E \norm{\bz^{t} -g(\by^{t-1}) }^2    \\
&+ \frac{1}{2} \pr{ \alpha_{t} - \alpha_{t-1}}   \norm{\bx^{t+1}}^2 +\frac{8}{\eta_\bx}\pr{\frac{\alpha_{t}}{\alpha_{t+1}} -\frac{\alpha_{t-1}}{\alpha_t}} \norm{\bx^{t+1} }^2\\
&-\pr{1-\frac{\beta}{2}}^2  \frac{64L_f^2}{\alpha_t^2 \eta_{\bx}}     \E \norm{  \bz^{t } - \bz^{t-1}  }^2+ \pare{\frac{1}{8\eta_{\by}}-\frac{768 G_g^2L_f^2}{\alpha_{t+1}^2 \eta_{\bx} \beta}}\E \norm{\by^t - \by^{t-1}}^2 + \frac{1}{8\eta_{\by}}\E \norm{\by^{t-1} - \by^{t-2}}^2\\
&  -2 \eta_{\bx} (1-\beta)^2  \frac{\sigma^2}{B}-  \pare{  2{\eta_{\by}}(G_g^2+ G_f^2)\frac{\sigma^2}{B} +2\beta^2 L_f^2 G_g^2 {\eta_{\by}}\frac{\sigma^2}{M} +  \frac{1152 L_f^2 \beta}{\alpha_t^2 \eta_{\bx}}    \frac{\sigma^2}{M}  } - 4C_1 \eta_{\by}^2 \frac{(G_f^2 + G_g^2)\sigma^2}{B} - 4C_2 \eta_{\bx}^2 \frac{\sigma^2}{B}.
\end{align*} 
Rearranging the terms yields:

 { \small \begin{align*}
  \E [ \hat{F}^{t+1} - \hat{F}^{t}] \geq &\frac{C^t_1}{2}\eta_{\by}^2\E \norm{\hat{\bg}_{\by}(\bx^t,\by^t)}^2 + C^t_2 \frac{1}{2} \eta^2_{\bx}\E \norm{\hat{\bg}_\bx(\bx^t,\by^t)}^2 + \frac{1}{8\eta_{\by}}\E \norm{\by^{t-1} - \by^{t-2}}^2  \\ 
 &+ \pr{\frac{1}{8\eta_{\by}}-\frac{768 G_g^2L_f^2}{\alpha_{t+1}^2 \eta_{\bx} \beta} - 16(1-\beta)^2G_g^2 (C^t_1\eta_\by L_f^2 G_g^2 + C^t_2 \eta_\bx^2 L_f^2) }\E \norm{\by^t - \by^{t-1}}^2\\
& -   (1-\beta)^2 \pare{  L_f^2 G_g^2 {\eta_{\by}} + 2\eta_{\bx} L_f^2 +C^t_1 4 \eta_{\by}^2 L_f^2 G_g^2 + 4C^t_2\eta^2_\bx L_f^2 }  \E \norm{\bz^{t} -g(\by^{t-1}) }^2 \\ 
&+ \frac{1}{2} \pr{ \alpha_{t} - \alpha_{t-1}}   \norm{\bx^{t+1}}^2 +\frac{8}{\eta_\bx}\pr{\frac{\alpha_{t}}{\alpha_{t+1}} -\frac{\alpha_{t-1}}{\alpha_t}} \norm{\bx^{t+1} }^2 -\pr{1-\frac{\beta}{2}}^2  \frac{64L_f^2}{\alpha_t^2 \eta_{\bx}}     \E \norm{  \bz^{t } - \bz^{t-1}  }^2\\
&  -\pare{2 \eta_{\bx}  +4C^t_2 \eta_{\bx}^2 +4C^t_1 \eta_{\by}^2  (G_f^2 + G_g^2) +2{\eta_{\by}}(G_g^2+ G_f^2)  } \frac{\sigma^2}{B} \\
&- \pare{2\beta^2 L_f^2 G_g^2 {\eta_{\by}}  + 4C^t_1 \beta^2 \eta_{\by}^2 L_f^2 G_g^2 + 8C^t_2\beta^2 \eta_{\bx}^2 L_f^2 +  \frac{1152 L_f^2 \beta}{\alpha_t^2 \eta_{\bx}}       }\frac{\sigma^2}{M}   .
\end{align*}}

Recall our definition of potential function $\tilde F^{t+1}$:   $ \tilde{F}^{t+1} :=  \hat{F}^{t+1} -   \E \norm{\bz^{t+1}  -  g(\by^t))  }^2 -\frac{(1-\frac{\beta}{2})^2}{1-(1-\frac{\beta}{2})^2}\frac{64L_f^2}{\alpha_T^2\eta_{\bx}}   \E \norm{  \bz^{t+1 } - \bz^{t}  }^2$.
Hence we have:
 { \small \begin{align*}
  \E [ \hat{F}^{t+1} - &\hat{F}^{t}]  \geq\frac{C^t_1}{2}\eta_{\by}^2\E \norm{\hat{\bg}_{\by}(\bx^t,\by^t)}^2 + C^t_2 \frac{1}{2} \eta^2_{\bx}\E \norm{\hat{\bg}_\bx(\bx^t,\by^t)}^2 + \frac{1}{8\eta_{\by}}\E \norm{\by^{t-1} - \by^{t-2}}^2    \\ 
 &+ \pr{\frac{1}{8\eta_{\by}}-\frac{768 G_g^2L_f^2}{\alpha_{t+1}^2 \eta_{\bx} \beta} - 16(1-\beta)^2G_g^2 (C^t_1\eta_\by L_f^2 G_g^2 + C^t_2 \eta_\bx^2 L_f^2) }\E \norm{\by^t - \by^{t-1}}^2\\
& -   (1-\beta)^2 \pare{  L_f^2 G_g^2 {\eta_{\by}} + 2\eta_{\bx} L_f^2 +C^t_1 4 \eta_{\by}^2 L_f^2 G_g^2 + 4C^t_2\eta^2_\bx L_f^2 }  \E \norm{\bz^{t} -g(\by^{t-1}) }^2 \\ 
&+ \frac{1}{2} \pr{ \alpha_{t} - \alpha_{t-1}}   \norm{\bx^{t+1}}^2 +\frac{8}{\eta_\bx}\pr{\frac{\alpha_{t}}{\alpha_{t+1}} -\frac{\alpha_{t-1}}{\alpha_t}} \norm{\bx^{t+1} }^2 -\pr{1-\frac{\beta}{2}}^2  \frac{64L_f^2}{\alpha_T^2 \eta_{\bx}}     \E \norm{  \bz^{t } - \bz^{t-1}  }^2\\
  &  -  \E \norm{\bz^{t+1}  -  g(\by^t))  }^2+  \E \norm{\bz^{t} -g(\by^{t-1}) }^2\\
    & -\frac{(1-\frac{\beta}{2})^2}{1-(1-\frac{\beta}{2})^2}\frac{64L_f^2}{\alpha_T^2\eta_{\bx}}    \E \norm{  \bz^{t+1 } - \bz^{t}  }^2 +\frac{(1-\frac{\beta}{2})^2}{1-(1-\frac{\beta}{2})^2}\frac{64L_f^2}{\alpha_T^2\eta_{\bx}}    \E \norm{  \bz^{t } - \bz^{t-1}  }^2 \\
    &  -\pare{2 \eta_{\bx}  +4C^t_2 \eta_{\bx}^2 +4C^t_1 \eta_{\by}^2  (G_f^2 + G_g^2) +2{\eta_{\by}}(G_g^2+ G_f^2)  } \frac{\sigma^2}{B} - \pare{2\beta^2 L_f^2 G_g^2 {\eta_{\by}}  + 4C^t_1 \beta^2 \eta_{\by}^2 L_f^2 G_g^2 + 8C^t_2\beta^2 \eta_{\bx}^2 L_f^2   +  \frac{1152 L_f^2 \beta}{\alpha_t^2 \eta_{\bx}}  }\frac{\sigma^2}{M}   .
\end{align*}}


Now we plug in Lemma~\ref{lem:tracking error} and~\ref{lem:second order tracking}:

 { \small \begin{align*}  
   \E [ \tilde{F}^{t+1} - &\tilde{F}^{t}]   \geq \frac{C^t_1}{2}\eta_{\by}^2\E \norm{\hat{\bg}_{\by}(\bx^t,\by^t)}^2 +  \frac{C^t_2}{2} \eta^2_{\bx}\E \norm{\hat{\bg}_{\bx}(\bx^t,\by^t)}^2   \\
   & + \pare{\frac{1}{8\eta_{\by}}-\frac{768 G_g^2L_f^2}{\alpha_{t+1}^2 \eta_{\bx} \beta}-16(1-\beta)^2G_g^2 (C^t_1\eta_\by L_f^2 G_g^2 + C^t_2 \eta_\bx^2 L_f^2)  -\frac{(1-\frac{\beta}{2})^2}{1-(1-\frac{\beta}{2})^2}\frac{64L_f^2}{\alpha_t^2\eta_{\bx}} \pr{ \frac{12}{\beta} + 6\beta  } G_g^2   }  \E \norm{\by^t - \by^{t-1}}^2 \\ 
    &  +\pare{1- (1-\beta)^2 \pare{  L_f^2 G_g^2 {\eta_{\by}} + 2\eta_{\bx} L_f^2 +C^t_1 4 \eta_{\by}^2 L_f^2 G_g^2 + 4C^t_2\eta^2_\bx L_f^2 + 1} }   \E \norm{\bz^{t} -g(\by^{t-1}) }^2\\
    &  +  \pare{\frac{1}{8\eta_{\by}} - \frac{(1-\frac{\beta}{2})^2}{1-(1-\frac{\beta}{2})^2}\frac{64L_f^2}{\alpha_T^2\eta_{\bx}}\frac{6}{\beta}G_g^2  }\E \norm{\by^{t-1} - \by^{t-2}}^2 \\ 
    & + (1-\frac{\beta}{2})^2\frac{64L_f^2}{\alpha_T^2\eta_{\bx}} \underbrace{\pare{\frac{1}{1-(1-\frac{\beta}{2})^2}-\frac{(1-\frac{\beta}{2})^2}{1-(1-\frac{\beta}{2})^2}-1} }_{=0}     \E \norm{  \bz^{t } - \bz^{t-1}  }^2\\
     &+\frac{1}{2} \pr{ \alpha_{t} - \alpha_{t-1}}   \norm{\bx^{t+1}}^2 +\frac{8}{\eta_\bx}\pr{\frac{\alpha_{t}}{\alpha_{t+1}} -\frac{\alpha_{t-1}}{\alpha_t}} \norm{\bx^{t+1} }^2 \\
      &  -\pare{2 \eta_{\bx}  +4C^t_2 \eta_{\bx}^2 +4C^t_1 \eta_{\by}^2  (G_f^2 + G_g^2) +2{\eta_{\by}}(G_g^2+ G_f^2)  } \frac{\sigma^2}{B}- \pare{2\beta^2 L_f^2 G_g^2 {\eta_{\by}}  + 4C^t_1 \beta^2 \eta_{\by}^2 L_f^2 G_g^2 + 8C^t_2\beta^2 \eta_{\bx}^2 L_f^2   +  \frac{1152 L_f^2 \beta}{\alpha_t^2 \eta_{\bx}}    }\frac{\sigma^2}{M}   .
\end{align*}}

By our choice of $\eta_{\by}$, we know that
\begin{align*}
  \frac{1}{8\eta_{\by}}-\frac{768 G_g^2L_f^2}{\alpha_{t+1}^2 \eta_{\bx} \beta}-16(1-\beta)^2G_g^2 (C^t_1\eta_\by L_f^2 G_g^2 + C^t_2 \eta_\bx^2 L_f^2)  -\frac{(1-\frac{\beta}{2})^2}{1-(1-\frac{\beta}{2})^2}\frac{64L_f^2}{\alpha_t^2\eta_\bx} \pr{ \frac{12}{\beta} + 6\beta  } G_g^2   \geq 0,\\
  1- (1-\beta)^2 \pare{  L_f^2 G_g^2 {\eta_{\by}} + 2\eta_{\bx} L_f^2 +C^t_1 4 \eta_{\by}^2 L_f^2 G_g^2 + 4C^t_2\eta^2_\bx L_f^2 + 1}  \geq 0,\\
  \frac{1}{8\eta_{\by}} - \frac{(1-\frac{\beta}{2})^2}{1-(1-\frac{\beta}{2})^2}\frac{64L_f^2}{\alpha_T^2\eta_{\bx}}\frac{6}{\beta}G_g^2   \geq 0
\end{align*}
 Now we can have the clean bound:
 
{\small \begin{align*}  
  \E [ \tilde{F}^{t+1} - \tilde{F}^{t}]   & \geq  \frac{C^t_1}{2}\eta_{\by}^2\E \norm{\hat{\bg}_{\by}(\bx^t,\by^t)}^2 +  \frac{C^t_2}{2} \eta^2_{\bx}\E \norm{\hat{\bg}_{\bx}(\bx^t,\by^t)}^2  \\
  & \quad+ \frac{1}{2} \pr{ \alpha_{t} - \alpha_{t-1}}   \norm{\bx^{t+1}}^2 +\frac{8}{\eta_\bx}\pr{\frac{\alpha_{t}}{\alpha_{t+1}} -\frac{\alpha_{t-1}}{\alpha_t}} \norm{\bx^{t+1} }^2  \\ 
    & \quad -\pare{2 \eta_{\bx}  +4C^t_2 \eta_{\bx}^2 +4C^t_1 \eta_{\by}^2  (G_f^2 + G_g^2) +2{\eta_{\by}}(G_g^2+ G_f^2)  } \frac{\sigma^2}{B}\\
    &\quad- \pare{2\beta^2 L_f^2 G_g^2 {\eta_{\by}}  + 4C^t_1 \beta^2 \eta_{\by}^2 L_f^2 G_g^2 + 8C^t_2\beta^2 \eta_{\bx}^2 L_f^2 +  \frac{1152 L_f^2 \beta}{\alpha_t^2 \eta_{\bx}}   }\frac{\sigma^2}{M}   .
\end{align*}}
  \end{proof}

 \end{lemma}

\subsubsection{Proof Theorem~\ref{thm: CNC}} 
Evoking Lemma~\ref{lemma: potential 2 cnc} and re-arranging terms yields:
 \begin{align*}  
   \frac{C^t_1}{2}&\eta_{\by}^2\E \norm{\hat{\bg}_{\by}(\bx^t,\by^t)}^2 +  \frac{C^t_2}{2} \eta^2_{\bx}\E \norm{\hat{\bg}_{\bx}(\bx^t,\by^t)}^2  \leq  \E [ \tilde{F}^{t+1} - \tilde{F}^{t}]  - \frac{1}{2} \pr{ \alpha_{t} - \alpha_{t-1}}   \norm{\bx^{t+1}}^2   \\
   &\quad -\frac{8}{\eta_\bx}\pr{\frac{\alpha_{t}}{\alpha_{t+1}} -\frac{\alpha_{t-1}}{\alpha_t}} \norm{\bx^{t+1} }^2  \\ 
    & \quad   +\pare{2 \eta_{\bx}  +4C^t_2 \eta_{\bx}^2 +4C^t_1 \eta_{\by}^2  (G_f^2 + G_g^2) +2{\eta_{\by}}(G_g^2+ G_f^2)  } \frac{\sigma^2}{B} \\
    &\quad + \pare{2\beta^2 L_f^2 G_g^2 {\eta_{\by}}  + 4C^t_1 \beta^2 \eta_{\by}^2 L_f^2 G_g^2 + 8C^t_2\beta^2 \eta_{\bx}^2 L_f^2 +  \frac{1152 L_f^2 \beta}{\alpha_t^2 \eta_{\bx}}    }\frac{\sigma^2}{M}   .
\end{align*}

We compute the upper and lower bound of $C^t_1$ and $C^t_2$. For $C^t_1$
\begin{align*}
     &C_1^t = \frac{1}{4\eta_{\by}} - \frac{L}{2} - \eta_{\bx} L^2  -  4 L_f^2  G_g^4 \eta_{\by}   - \frac{768  G_g^2L_f^2}{\alpha_{t+1}^2 \eta_{\bx} \beta} - \frac{576  G_g^2 L_f^2\beta}{\alpha_{t+1}^2 \eta_{\bx}}.    
\end{align*}
The upper bound $C^t_1 \leq \frac{1}{4\eta_{\by}}$ holds trivially. For lower bound,
since we choose
\begin{align*}
    \eta_{\by} \leq  c \cdot \min\cbr{\frac{1}{  L}, \frac{1}{  \eta_{\bx} L^2}, \frac{1}{  L_f  G_g^2  }, \frac{\alpha_t^2 \eta_{\bx} \beta}{  G_g^2 L_f^2}} ,
\end{align*}
for some $c\geq 0$,
we can choose sufficiently small $c$ such that $C^t_1 \geq \frac{1}{8\eta_{\by}}$. 

For $C^t_2$:
\begin{align*}
      &C^t_2 =  \frac{1}{ \eta_{\bx}}   +\frac{\alpha_{t-1}}{2} -\frac{1}{4\eta_{\bx} (1-\beta )^2}.
\end{align*}
The upper bound $C^t_2 \leq \frac{1}{\eta_{\bx}} +\frac{\alpha_t}{2} $ holds trivially. For lower bound, since we choose:
\begin{align*}
   \beta =0.1 \leq 1-\frac{\sqrt{2}}{2}
\end{align*}
it holds that $C^t_2 \geq \frac{1}{2\eta_{\bx}}$.

 Since $\frac{1}{8\eta_{\by}} \leq C^t_1 \leq  \frac{1}{4\eta_{\by}}$ and $ \frac{1}{2\eta_{\bx}}\leq C^t_2 \leq \frac{1}{\eta_{\bx}}+\frac{\alpha_t}{2} \leq \frac{2}{\eta_{\bx}}$, we have:
 \begin{align*}  
   \frac{\eta_{\by}}{16}& \E \norm{\hat{\bg}_{\by}(\bx^t,\by^t)}^2 +  \frac{ \eta_{\bx}}{4} \E \norm{\hat{\bg}_{\bx}(\bx^t,\by^t)}^2  \leq  \E [ \tilde{F}^{t+1} - \tilde{F}^{t}]  - \frac{1}{2} \pr{ \alpha_{t} - \alpha_{t-1}}   \norm{\bx^{t+1}}^2   \\
   &\quad -\frac{8}{\eta_\bx}\pr{\frac{\alpha_{t}}{\alpha_{t+1}} -\frac{\alpha_{t-1}}{\alpha_t}} \norm{\bx^{t+1} }^2  \\ 
    & \quad   +\pare{2 \eta_{\bx}  +8 \eta_{\bx}  +  \eta_{\by}  (G_f^2 + G_g^2) +2{\eta_{\by}}(G_g^2+ G_f^2)  } \frac{\sigma^2}{B} + \pare{2\beta^2 L_f^2 G_g^2 {\eta_{\by}}  +  \eta_{\by}  L_f^2 G_g^2 + 16  \eta_{\bx}  L_f^2  +  \frac{1152 L_f^2 \beta}{\alpha_t^2 \eta_{\bx}}   }\frac{\sigma^2}{M}   .
\end{align*}
 Summing the above inequality from $t=0$ to $T-1$ yields:

  \begin{align*}  
     \frac{\E [ \tilde{F}^{T} - \tilde{F}^{0}]}{T} + \frac{\pr{\alpha_{T-1} - \alpha_0}D^2_{\cX}}{2T} &-  \frac{8\pr{\frac{\alpha_{T-1}}{\alpha_{T}} -\frac{\alpha_{0}}{\alpha_1}}D^2_{\cX}}{ T} + 4\eta_\bx\frac{1}{T} \sum_{t=0}^{T-1}  \alpha_t^2 D_\cX^2\\
   &+O\pare{  \pare{   \eta_{\bx}   + {\eta_{\by}}(G_g^2+ G_f^2)  } \frac{\sigma^2}{B} + \pare{  L_f^2 G_g^2 {\eta_{\by}}  + \eta_{\bx}  L_f^2  +  \frac{  L_f^2 \beta}{\alpha_t^2 \eta_{\bx}} }\frac{\sigma^2}{M} } \\  &\geq  \min\left\{ \frac{1}{16}\eta_{\by} , \frac{1}{4} \eta_{\bx}  \right\} \pr{\frac{1}{T} \sum_{t=0}^{T-1} \E \norm{\hat{\bg}_{\by}(\bx^t,\by^t)}^2 +  \E \norm{\hat{\bg}_{\bx}(\bx^t,\by^t)}^2}.
\end{align*}


We then need to verify the conditions ~\ref{eq:cond 1 cnc}, \ref{eq:cond 2 cnc} and \ref{eq:cond 3 cnc} in Lemma~\ref{lemma: potential 2 cnc} can hold under our choice of $\eta_{\bx}$, $\eta_{\by}$ and $\beta$. To guarantee (\ref{eq:cond 1 cnc}) holding, we need: 
\begin{align*}
    \eta_{\by} \leq \Theta\pare{\min\cbr{ \frac{1}{   \eta_x L_f^2 G_g^2  },\frac{1}{   L_f  G_g^2  },  \frac{\alpha_t^2 \eta_{\bx}}{L_f^2 G_g^2} } }, \forall \ 0 \leq  t \leq T-1.
\end{align*}
To guarantee condition (\ref{eq:cond 2 cnc}) holding, we need 
\begin{align*}
    \eta_\bx \leq \Theta\pare{  \frac{1}{  L_f^2}}, \eta_{\by} \leq \Theta\pare{  \frac{1}{L_f^2 G_g^2}} .
\end{align*}

To guarantee condition (\ref{eq:cond 3 cnc}) holding, we need:
\begin{align*}
    \eta_{\by} \leq \Theta\pare{\frac{\alpha_T^2 \eta_{\bx}}{  L_f^2 G_g^2}}.
\end{align*}
Next we examine how large the $\E [ \tilde{F}^{T} - \tilde{F}^{0}]$ is.
By definition of potential function, we have:
\begin{align*}
        \E [ \tilde{F}^{T} - \tilde{F}^{0}] &=  \hat{F}^{T} -   \E \norm{\bz^{T} -g(\by^{T-1}) }^2 -\frac{(1-\frac{\beta}{2})^2}{1-(1-\frac{\beta}{2})^2}\frac{64L_f^2}{\alpha_t^2\eta_{\bx}}   \E \norm{  \bz^{t} - \bz^{t-1}  }^2\\
        &\quad - \pare{\hat{F}^{0} -   \E \norm{\bz^{0} -g(\by^{-1}) }^2 -\frac{(1-\frac{\beta}{2})^2}{1-(1-\frac{\beta}{2})^2}\frac{64L_f^2}{\alpha_t^2\eta_{\bx}}   \E \norm{  \bz^{0} - \bz^{-1}  }^2} \\
        &\leq \hat{F}^{T}  - \hat{F}^{0} +     \E \norm{\bz^{0} -g(\by^{-1})}^2 + \frac{(1-\frac{\beta}{2})^2}{\beta - \frac{\beta^2}{4}}\frac{64L_f^2}{\alpha_t^2\eta_{\bx}}   \E \norm{  \bz_j^{0} - \bz_j^{-1}  }^2 \\ 
\end{align*}
By convention $\by^{0} = \by^{-1}$, $\bz^{0} = \bz^{-1}$, and our choice $\E\norm{\bz^{0} - g(\by^{0})}^2 \leq O(1)$,  we have
 \begin{align*}
        \E [ \tilde{F}^{T} - \tilde{F}^{0}] 
        \leq \hat{F}^{T}  - \hat{F}^{0} + O(1).
\end{align*} 
 Next we examine how large the $\E [ \hat{F}^{T} - \hat{F}^{0}]$ is. 
 \begin{align*} 
     \E[\hat{F}^{T}  - \hat{F}^{0}]& = F(\bx^{T},\by^T) + s^{T} - \pare{\frac{1}{4\eta_{\by}} + 4 L_f^2  G_g^4 \eta_{\by}  + \frac{\eta_{\bx} L^2}{2} + \frac{768  G_g^2L_f^2}{\alpha_{T}^2 \eta_{\bx} \beta} + \frac{576  G_g^2 L_f^2\beta}{\alpha_{T}^2 \eta_{\bx}}  } \norm{\by^T - \by^{T-1}}^2\\
    &  \quad-\pare{\frac{1}{8\eta_{\by}}+\frac{768 G_g^2L_f^2}{a_T^2 \eta_{\bx} \beta}}\norm{\by^{T-1} - \by^{T-2}}^2 +  \frac{7}{2\eta_{\bx}}  \norm{\bx^{T} - \bx^{T-1}}^2+ \frac{\alpha_{T-1}}{2}\norm{\bx^{T}}^2  \\
    &\quad-F(\bx^{0},\by^{0}) - s^{0} + \pare{\frac{1}{4\eta_{\by}} + 4 L_f^2   G_g^4 \eta_{\by}  + \frac{\eta_{\bx} L^2}{2} + \frac{768  G_g^2L_f^2}{a_0^2 \eta_{\bx} \beta} + \frac{576  G_g^2 L_f^2\beta}{a_0^2 \eta_{\bx}}} \norm{\by^{0} - \by^{-1}}^2\\
    &  \quad+\pare{\frac{1}{8\eta_{\by}}+\frac{768 G_g^2}{a_0^2 \eta_{\bx} \beta}}\norm{\by^{-1} - \by^{-2}}^2 -  \frac{7}{2\eta_{\bx}}  \norm{\bx^{0} - \bx^{-1}}^2- \frac{\alpha_{0}}{2}\norm{\bx^{0}}^2   \\
    & \leq F_{\max}     + \frac{7}{2\eta_{\bx}}  \norm{\bx^{T} - \bx^{T-1}}^2  ,
 \end{align*}
 
where we used the convention $\bx^{-1} = \bx^0$ and $\by^{-1} = \by^{-2}$. Notice that $\E\norm{\bx^{T} - \bx^{T-1}}^2 \leq \eta_{\bx}^2\E \norm{ \frac{1}{B}\sum_{(\zeta,\xi) \in \cB^{T-1}}\nabla_1 f(\bx^{T-1}, \bz^T;\zeta ) + \nabla h(\bx^{T-1}) }^2 \leq 2\eta_{\bx}^2(G_f^2+G_h^2) $,
we have $\E[\hat{F}^{T}  - \hat{F}^{0}] \leq F_{\max} + 7\eta_{\bx} (G_f^2+G_h^2)$.
 

 
Due to our choice that $\eta_{\by} \leq 8\eta_{\bx}$, $\alpha_0 = ... = \alpha_T = \alpha$ we have:
 
{\small \begin{align*}
     \frac{1}{T} \sum_{t=1}^{T} \pare{\E \norm{\hat{\bg}_{\by}(\bx^t,\by^t)}^2 +  \E \norm{\hat{\bg}_{\bx}(\bx^t,\by^t)}^2 } 
     &\leq  O\pr{   \frac{F_{\max} +  \eta_{\bx} (G_f^2+G_h^2)  }{\eta_{\by}T}   } +O\pare{     ( \frac{\eta_{\bx}}{\eta_{\by}} G_f^2 + G_g^2)    \frac{\sigma^2}{B} + \pare{  L_f^2 G_g^2    + \frac{\eta_{\bx}}{\eta_{\by}}   L_f^2 +  \frac{  L_f^2 \beta}{\alpha^2 \eta_{\bx}\eta_{\by}}  }\frac{\sigma^2}{M} } .   
 \end{align*} }
By the definition of stationary measure, we have
\begin{align*}
   \E  \norm{\nabla G(\bx^t,\by^t)}^2 &=  \frac{1}{\eta^2_{\bx}}\E \norm{ \bx^t - \cP_{\cX} \pare{\bx^t - \eta \nabla_{\bx} F(\bx^t,\by^t)} }^2 +  \E \norm{\hat{\bg}_{\by}(\bx^t,\by^t)}^2 \\
  & = \frac{2}{\eta^2_{\bx}}\E \norm{ \bx^t - \cP_{\cX} \pare{\bx^t - \eta_{\bx} \nabla_{\bx} F(\bx^t,\by^t) - \eta_{\bx} \alpha \bx^t} }^2\\
   & \quad + \frac{2}{\eta^2_{\bx}}\E \norm{   \cP_{\cX} \pare{\bx^t - \eta_{\bx} \nabla_{\bx} F(\bx^t,\by^t)} - \cP_{\cX} \pare{\bx^t - \eta_{\bx} \nabla_{\bx} F(\bx^t,\by^t) - \eta_{\bx} \alpha \bx^t}  }^2  +  \E \norm{\hat{\bg}_{\by}(\bx^t,\by^t)}^2 \\
    & = \frac{2}{\eta^2_{\bx}}\E \norm{ \bx^t - \cP_{\cX} \pare{\bx^t - \eta_{\bx} \nabla_{\bx} F(\bx^t,\by^t) - \eta_{\bx} \alpha \bx^t} }^2  + \frac{2}{\eta^2_{\bx}}\E \norm{    \alpha \bx^t   }^2  +  \E \norm{\hat{\bg}_{\by}(\bx^t,\by^t)}^2 \\
    & \leq  2\E \norm{ \hat{\bg}_{\bx}(\bx^t,\by^t)}^2+ 2 \E \norm{\hat{\bg}_{\by}(\bx^t,\by^t)}^2   + 2\E \norm{    \alpha \bx^t   }^2 .
\end{align*}

To guarantee $\frac{1}{T} \sum_{t=0}^T \E  \norm{\nabla G(\bx^t,\by^t)}^2$ is less than $\epsilon^2$, we need:
 $T = O\pr{\frac{F_{\max} }{\eta_{\by} \epsilon^2}}$, $\alpha = \Theta\pare{\frac{ \epsilon }{\cD_{\cX}}}$ and 
 \begin{align*}
    M&= \Theta \pare{ \max\cbr{   \frac{ L^6 L_f^2 D_{\cX}^4 \sigma^2}{\epsilon^4}, 1} },  B= \Theta \pare{ \max\cbr{   \frac{ L^3 L_f^2 D_{\cX}^2 \sigma^2}{\epsilon^4}, 1} },  \\ 
    \eta_{\by} &= \Theta\pare{\min\cbr{ \frac{1}{   \eta_x L_f^2 G_g^2  },\frac{1}{   L_f  G_g^2  },  \frac{\alpha^2 \eta_{\bx}}{L_f^2 G_g^2} } } , \eta_{\bx} = \Theta\pare{\frac{1}{L^2}} ,
 \end{align*}
 which yields the total gradient complexity:
 \begin{align*}
     O\pr{  \max\cbr{   \frac{ L^3 L_f^2 D_{\cX}^2 \sigma^2}{\epsilon^4}, 1} \cdot  \frac{F_{\max} D_{\cX}^2 L^4}{\epsilon^4} }. 
 \end{align*}
 
 \section{Proof of Both Sides Composition} \label{app: proof PD}
 In this section we provide the proof of results in primal and dual composition setting (Theorem~\ref{thm:wcwc}).

We first define some notations. Let $\bw = [\bx,\by] $ be stacked variable. $\bg_{\bx}(\bw) :=  \nabla f(g(\bx,\by)) \nabla_{\bx} g(\bx,\by) + \nabla h(\bx) + \frac{2}{\gamma}\pare{\bx-\bx^0} $, $\bg_{\by}(\bw) :=  \nabla f(g(\bx,\by)) \nabla_{\by} g(\bx,\by) - \nabla r(\by) - \frac{2}{\gamma} \pare{\by-\by^0} $. Hence we define stacked variable $\bg(\bw) := [\bg_{\bx}(\bw), \bg_{\by}(\bw)]$. 
Similarly, we define stacked variable $\bg^t := [\bg^t_{\bx}, \bg^t_{\by}]$ for the gradients we used in updating rule.
 \begin{proposition}\label{lem: tracking error ncnc}
 For Algorithm~\ref{algorithm: CODA-SCSC}, under the assumptions of Theorem~\ref{thm:wcwc}, the following statement holds true:
     \begin{align*}
    \E \norm{\bz^{t+1}  - g(\bx^t,\by^t) }^2 &\leq (1-\beta)^2 \E \norm{\bz^{t} - g(\bx^{t-1},\by^{t-1}) }^2  + 4(1-\beta)^2 G_g^2 \pare{\norm{\bx^t - \bx^{t-1}}^2+\norm{\by^t - \by^{t-1}}^2 }+ 2\beta^2 \frac{\sigma^2}{B}.
\end{align*}
 \end{proposition}
\begin{lemma}[Convergence of Iterates difference] \label{lem:second order tracking pd} 
For Algorithm~\ref{algorithm: CODA-SCSC}, under the assumptions of Theorem~\ref{thm:wcwc}, the following statement holds true:
    \begin{align*}
        \E \norm{  \bz^{t+1} - \bz^{t}  }^2 &\leq \pare{1-\frac{\beta}{2}}^2 \E \norm{  \bz^{t } - \bz^{t-1}  }^2+ 4\pare{1+\frac{2}{\beta}} G_g^2\pare{\norm{    \bw^t - \bw^{t-1}   }^2 +  \norm{    \bw^{t-1} - \bw^{t-2}   }^2 } \\
        & + 2\pare{1+\frac{2}{\beta}}\beta^2 G_g^2 \E \norm{    \bw^t - \bw^{t-1}   }^2.
    \end{align*}
    \begin{proof}
For the ease of presentation, we define the following two auxiliary variables:
\begin{align*}
&g^t  = g(\bx^t, \by^t;\cM^t)  , g^{t\mapsto t-1} = g(\bx^t, \by^t;\cM^t)   - g(\bx^{t-1},\by^{t-1};\cM^{t}) .
\end{align*}
According to updating rule of $\bz$, we have:
        \begin{align*}
            \bz^{t+1} - \bz^{t} = (1-\beta) ( \bz^{t} - \bz^{t-1}) + (1-\beta) (g^{t\mapsto t-1}      -   g^{t-1\mapsto t-2}   )  + \beta (g^t - g^{t-1}).
        \end{align*}
       Taking expectation w.r.t. $\cM^t$, $\cM^{t-1}$ yields:
        \begin{align*}
           \E\norm{ \bz^{t+1} - \bz^{t}}^2 &= \E\norm{  (1-\beta) ( \bz^{t} - \bz^{t-1}) + (1-\beta) ( g^{t\mapsto t-1}  - g^{t-1\mapsto t-2}   )  + \beta (g^t - g^{t-1})}^2\\
           &\stackrel{(1)}{\leq} \pare{1+\frac{\beta}{2-2\beta} }(1-\beta)^2\E\norm{   \bz^{t} - \bz^{t-1} +( g^{t\mapsto t-1}      -    g^{t-1\mapsto t-2}   ) }^2   +\pare{1+\frac{2-2\beta}{\beta}}\norm{  \beta (g^t - g^{t-1})}^2\\
           &\stackrel{(2)}{\leq} \pare{1-\frac{\beta}{2}}(1-\beta)\E\norm{  \bz^{t} - \bz^{t-1} +( g^{t\mapsto t-1}      -    g^{t-1\mapsto t-2}   ) }^2    +  \pare{1+\frac{2 }{\beta}}\beta^2 \E\norm{g^t - g^{t-1}}^2 \\
            &\stackrel{(3)}{\leq} \pare{1-\frac{\beta}{2}}(1-\beta) \pare{1+\frac{\beta}{2-2\beta} }\E\norm{  \bz^{t} - \bz^{t-1}   }^2   \\
            & \quad + \pare{1-\frac{\beta}{2}}(1-\beta) \pare{1+\frac{2-2\beta}{\beta} }\E\norm{   g^{t\mapsto t-1}      -    g^{t-1\mapsto t-2}  }^2   +   \pare{1+\frac{2 }{\beta}}\beta^2 \E\norm{g^t - g^{t-1}}^2 \\
            &\stackrel{(4)}{\leq} \pare{1-\frac{\beta}{2}}^2 \E\norm{  \bz^{t} - \bz^{t-1}   }^2   +   \underbrace{\pare{1+\frac{2}{\beta} }\E\norm{   g^{t\mapsto t-1}      -    g^{t-1\mapsto t-2}  }^2}_{T_1}   + \underbrace{ \pare{1+\frac{2 }{\beta}}\beta^2 \E\norm{g^t - g^{t-1}}^2}_{T_2},
        \end{align*}
        where in (1) and (3) we use Young's inequality that $\norm{\ba+\bb}^2 \leq (1+a)\norm{\ba}^2 + (1+\frac{1}{a})\norm{\bb}^2  $. 
        Now we bound $T_1$ as follows:
        \begin{align*}
            T_1 \leq &2\pare{1+\frac{2}{\beta} }\pare{\E\norm{  g(\bx^t, \by^t;\cM^t)      -  g(\bx^{t-1}, \by^{t-1};\cM^{t}) }^2}\\
            &+ 2\pare{1+\frac{2}{\beta} }\pare{\E\norm{ g(\bx^{t-1},\by^{t-1};\cM^{t-1})      -  g(\bx^{t-2},\by^{t-2};\cM^{t-1})}^2}\\
            \leq & 4\pare{1+\frac{2 }{\beta}}G_g^2\pare{\E\norm{\bx^t - \bx^{t-1}   }^2 +\E\norm{\by^t - \by^{t-1}   }^2+ \E \norm{    \bx^{t-1} - \bx^{t-2}   }^2  + \E \norm{    \by^{t-1} - \by^{t-2}   }^2 },
        \end{align*}
        and for $T_2$:
        \begin{align*}
            T_2 \leq  \pare{1+\frac{2}{\beta}}\beta^2 G_g^2\pare{\E \norm{\bx^t - \bx^{t-1}   }^2+\E \norm{    \by^t - \by^{t-1}   }^2 }.
        \end{align*}
        Putting pieces together will conclude the proof.
    \end{proof}
\end{lemma}
\begin{lemma}\label{lem:conv to optimal}
 
Let $\bw^*$ be the solution of $MVI(F^{\gamma},\cX\times \cY)$. For Algorithm~\ref{algorithm: CODA-SCSC}, under the assumptions of Theorem~\ref{thm:wcwc}, the following statement holds true:
    \begin{align*}
         A_{t+1}  \leq &\pare{1- \frac{\eta \mu}{4} }A_{t}  + \pare{8\eta^2 + 32C_2\frac{\eta}{\mu}} \frac{(G_f^2+G_g^2)\sigma^2}{B} + \pare{2\beta^2 C_1 +   C_2  \frac{8 \eta}{\mu }    G_g^2 L_f^2  \beta^2}\frac{\sigma^2}{M},
    \end{align*}
    where $A_t = \E\norm{\bw^{t}  -\bw^* }^2+C_1\E\norm{ \bz^{t+1} - g(\bx^{t},\by^{t}) }^2+ C_2 \E\norm{\bw^{t+1} - \bw^t}^2$, $C_1 =  220 \eta G_g^2 L_f^2 \max\cbr{1,\frac{1}{\mu}, G_g^2}$, and $C_2   = \frac{\mu }{16 \eta G_g^2 L_f^2}   C_1  = 14 \mu \max\cbr{1,\frac{1}{\mu}, G_g^2}$.
    \begin{proof}
    Recall that Algorithm~\ref{algorithm: CODA-SCSC} is optimizing a $\mu := \frac{1}{\gamma}-\rho$ strongly convex strongly concave function $F^{\gamma}(\bx,\by) = F(\bx,\by) + \frac{1}{2\gamma}\norm{\bx - \bx^0}^2 - \frac{1}{2\gamma}\norm{\by - \by^0}^2$, and it is $\tilde{L}:= L+ \frac{1}{\gamma}$ smooth.   According to updating rule we have:
        \begin{align*}
            \E\norm{\bw^{t+1} - \bw^t}^2  
            \leq & \E\norm{\bw^{t} - \bw^{t-1} - \eta (\bg^t - \bg^{t-1})     }^2\\\
            \leq & \pare{1+ \frac{1}{2(\frac{1}{2\mu\eta}-1)}} \E\norm{\bw^{t} - \bw^{t-1} - \eta (\nabla F^{\gamma}(\bx^t,\by^t) - \nabla F^{\gamma}(\bx^{t-1},\by^{t-1}))     }^2  \\
            & +  \pare{1+2(\frac{1}{2\mu\eta}-1)} \eta^2\E\norm{\bg^t- \nabla F^{\gamma}(\bx^t,\by^t) - (\nabla F^{\gamma}(\bx^{t-1},\by^{t-1}) -\bg^{t-1})}^2 \\
            \leq & \pare{1+ \frac{1}{2(\frac{1}{2\mu\eta}-1)}} \pare{1-2\eta \mu} \E\norm{\bw^{t} - \bw^{t-1}}^2 +\E\norm{ \eta (\nabla F^{\gamma}(\bx^t,\by^t) -\nabla F^{\gamma}(\bx^{t-1},\by^{t-1}))     }^2  \\
            & +  \frac{ 1}{\mu\eta} \eta^2\E\norm{\bg^t- \nabla F^{\gamma}(\bx^{t},\by^{t}) - (\nabla F^{\gamma}(\bx^{t-1},\by^{t-1}) -\bg^{t-1})}^2 \\
            \leq &   \pare{1-{\eta \mu}  + \eta^2 \tilde{L}^2} \E\norm{\bw^{t} - \bw^{t-1}}^2 \\
            &+  \frac{4 \eta}{\mu }    G_g^2 L_f^2\pare{ \E\norm{ \bz^{t+1} - g(\bx^t,\by^t) }^2 +  \E\norm{ \bz^{t} - g(\bx^{t-1},\by^{t-1}) }^2 } + 32 \frac{ \eta (G_g^2 + G_f^2 ) \sigma^2}{\mu B},
        \end{align*}
    where the last step is due to 
    \begin{align*}
        \E\norm{ \bg^t - \bar{\bg}^t  }^2=& \E\norm{ \bg_{\bx}^t - \bar{\bg}_{\bx}^t  }^2 + \E\norm{ \bg_{\by}^t - \bar{\bg}_{\by}^t  }^2\\ 
        = &\underbrace{\E\norm{ \frac{1}{B} \sum_{(\zeta,\xi)\in\cB^t} \nabla f(\bz^{t+1};\zeta) \nabla_{\bx} g(\bx^t,\by^t; \xi) -  \nabla f(\bz^{t+1} ) \nabla_{\bx} g(\bx^t,\by^t )  }^2}_{A} \\
        &+ \underbrace{\E\norm{  \frac{1}{B} \sum_{(\zeta,\xi)\in\cB^t} \nabla f(\bz^{t+1};\zeta) \nabla_{\by} g(\bx^t,\by^t; \xi) -  \nabla f(\bz^{t+1} ) \nabla_{\by} g(\bx^t,\by^t )}^2}_{B}. 
    \end{align*}
    For $A$,
    \begin{align*}
        A \leq& 2 \E\norm{ \frac{1}{B} \sum_{(\zeta,\xi)\in\cB^t} \nabla f(\bz^{t+1};\zeta) \nabla_{\bx} g(\bx^t,\by^t; \xi) -  \frac{1}{B} \sum_{(\zeta,\xi)\in\cB^t} \nabla f(\bz^{t+1};\zeta) \nabla_{\bx} g(\bx^t,\by^t )  }^2\\
        &+ 2 \E\norm{\frac{1}{B} \sum_{(\zeta,\xi)\in\cB^t} \nabla f(\bz^{t+1};\zeta) \nabla_{\bx} g(\bx^t,\by^t )  -   \nabla f(\bz^{t+1} ) \nabla_{\bx} g(\bx^t,\by^t )  }^2 \\
        \leq& 2\frac{1}{B^2} \sum_{(\zeta,\xi)\in\cB^t}\E\norm{   \nabla f(\bz^{t+1};\zeta) \nabla_{\bx} g(\bx^t,\by^t; \xi) -   \nabla f(\bz^{t+1};\zeta) \nabla_{\bx} g(\bx^t,\by^t )  }^2\\
        &+ 2\frac{1}{B^2} \sum_{(\zeta,\xi)\in\cB^t} \E\norm{ \nabla f(\bz^{t+1};\zeta) \nabla_{\bx} g(\bx^t,\by^t )  -   \nabla f(\bz^{t+1} ) \nabla_{\bx} g(\bx^t,\by^t )  }^2 \\
        = & 2 \frac{(G_f^2+G_g^2)\sigma^2}{B}.
    \end{align*}
    Similarly $B \leq 2 \frac{(G_f^2+G_g^2)\sigma^2}{B}$.
    Now we examine the convergence to optimal point:    
        \begin{align*}
            \E\norm{\bw^{t+1}-\bw^*}^2 &\leq \E\norm{\bw^{t}  -\bw^*- \eta (\bg^t- \bg(\bw^*) )}^2 \\
            &= \E\norm{\bw^{t}  -\bw^* }^2 + \eta^2\E\norm{ \bg^t- \bg(\bw^*)  }^2 - 2\eta \E\inprod{\bg^t- \bg(\bw^*)  }{\bw^{t}  -\bw^*} \\
            &= \E\norm{\bw^{t}  -\bw^* }^2 + \eta^2\E\norm{ \bg^t- \bg(\bw^*)  }^2 - 2\eta \E\inprod{\bg(\bw^t)- \bg(\bw^*)  }{\bw^{t}  -\bw^*} \\
            &\quad -  2\eta \E\inprod{\bg^t-\bg(\bw^t)  }{\bw^{t}  -\bw^*}  .
    \end{align*}
    Due to $\mu$ strong monotoncity we have:
 
    \begin{align*} 
          \E\norm{\bw^{t+1}-\bw^*}^2    &\leq (1-2\eta\mu) \E\norm{\bw^{t}  -\bw^* }^2 + \eta^2\E\norm{  \bar{\bg}^t- \bg(\bw^*)  }^2 +8\eta^2 \frac{(G_f^2+G_g^2)\sigma^2}{B}  \\
            &\quad -  2\eta \E\inprod{\bg^t-\bg(\bw^t)  }{\bw^{t}  -\bw^*} \\
             &\leq (1-2\eta\mu) \E\norm{\bw^{t}  -\bw^* }^2 + \eta^2\E\norm{   \bar{\bg}^t- \bg(\bw^*)  }^2 +8\eta^2 \frac{(G_f^2+G_g^2)\sigma^2}{B}   \\
            &\quad +   \pare{\frac{\eta}{\mu}\norm{ \bar{\bg}^t-\bg(\bw^t) }^2+ \eta \mu\norm{ \bw^{t}  -\bw^*}^2} \\
            &\leq (1-2\eta\mu + 2\eta^2 \tilde{L}^2) \E\norm{\bw^{t}  -\bw^* }^2 + 2\eta^2\E\norm{  \bar{\bg}^t- \bg(\bw^t)  }^2 +8\eta^2 \frac{(G_f^2+G_g^2)\sigma^2}{B} \\
            &\quad +   \pare{\frac{\eta}{\mu}\norm{ \bar{\bg}^t-\bg(\bw^t) }^2+ \eta \mu\norm{ \bw^{t}  -\bw^*}^2} \\
             &\leq (1- \frac{\eta\mu}{2} + 2\eta^2 \tilde{L}^2) \E\norm{\bw^{t}  -\bw^* }^2 + \pare{2\eta^2G_g^2 L_f^2+ \frac{\eta}{\mu}G_g^2 L_f^2}\E\norm{ \bz^{t+1} - g(\bx^t,\by^t) }^2 \\
             &\quad +8\eta^2 \frac{(G_f^2+G_g^2)\sigma^2}{B} .
        \end{align*}

        Now we add $C_1\E\norm{ \bz^{t+2} - g(\bx^{t+1},\by^{t+1}) }^2$ on both side, and $C_1(1-\frac{\eta\mu}{2} + 2\eta^2 \tilde{L}^2)\E\norm{ \bz^{t+1} - g(\bx^{t},\by^{t}) }^2$
        \begin{align*}
            &\E\norm{\bw^{t+1}-\bw^*}^2  +C_1\E\norm{ \bz^{t+2} - g(\bx^{t+1},\by^{t+1}) }^2\\
             &\leq (1-\frac{\eta\mu}{2} + 2\eta^2 \tilde{L}^2) (\E\norm{\bw^{t}  -\bw^* }^2+C_1\E\norm{ \bz^{t+1} - g(\bx^{t},\by^{t}) }^2 ) \\
             & \quad + \pare{2\eta^2G_g^2 L_f^2+ \frac{\eta}{\mu}G_g^2 L_f^2}\E\norm{ \bz^{t+1} - g(\bx^t,\by^t) }^2 +8\eta^2 \frac{(G_f^2+G_g^2)\sigma^2}{B}   \\
             & \quad +C_1\E\norm{ \bz^{t+2} - g(\bx^{t+1},\by^{t+1}) }^2 - C_1(1-\frac{\eta\mu}{2} + 2\eta^2 \tilde{L}^2)\E\norm{ \bz^{t+1} - g(\bx^t,\by^t) }^2.
        \end{align*} 
        We apply Proposition~\ref{lem: tracking error ncnc} to replace $\E\norm{ \bz^{t+2} - g(\bx^{t+1},\by^{t+1}) }^2$:
        \begin{align*}
            &\E\norm{\bw^{t+1}-\bw^*}^2  +C_1\E\norm{ \bz^{t+2} - g(\bx^{t+1},\by^{t+1}) }^2\\
             &\leq (1-\frac{\eta\mu}{2} + 2\eta^2 \tilde{L}^2) (\E\norm{\bw^{t}  -\bw^* }^2+C_1\E\norm{ \bz^{t+1} - g(\bx^{t},\by^{t}) }^2 ) +8\eta^2 \frac{(G_f^2+G_g^2)\sigma^2}{B} \\
             &\quad  + \pare{2\eta^2G_g^2 L_f^2+ \frac{\eta}{\mu}G_g^2 L_f^2}\E\norm{ \bz^{t+1} - g(\bx^t,\by^t) }^2  \\
             & \quad +C_1 \pare{ (1-\beta)^2\E\norm{ \bz^{t+1} - g(\bx^t,\by^t) }^2  + 4(1-\beta)^2 G_g^2 \E\norm{\bw^{t+1}   - \bw^{t}}^2  + 2\beta^2 \frac{\sigma^2}{M} }   \\
             & \quad  - C_1(1-\frac{\eta\mu}{2} + 2\eta^2 \tilde{L}^2)\E\norm{ \bz^{t+1} - g(\bx^t,\by^t) }^2\\
             &\leq (1-\frac{\eta\mu}{2} + 2\eta^2 \tilde{L}^2) (\E\norm{\bw^{t}  -\bw^* }^2+C_1\E\norm{ \bz^{t+1} - g(\bx^{t},\by^{t}) }^2 ) +8\eta^2 \frac{(G_f^2+G_g^2)\sigma^2}{B}  + 2\beta^2 C_1 \frac{\sigma^2}{M}\\
             & \quad + 4(1-\beta)^2 C_1 G_g^2 \E\norm{\bw^{t+1} - \bw^t}^2\\
             & \quad + \pare{2\eta^2G_g^2 L_f^2+ \frac{\eta}{\mu}G_g^2 L_f^2+ (1-\beta)^2C_1 -(1-\frac{\eta\mu}{2} + 2\eta^2 \tilde{L}^2)C_1} \E\norm{ \bz^{t+1} - g(\bx^t,\by^t) }^2.
        \end{align*}
   
    Adding $C_2 \E\norm{\bw^{t+2} - \bw^{t+1}}^2$ and $(1-\frac{\eta\mu}{2} + 2\eta^2 \tilde{L}^2) C_2 \E\norm{\bw^{t+1} - \bw^t}^2$ on both sides yields:
   {\small \begin{align*}
         &\E\norm{\bw^{t+1}-\bw^*}^2  +C_1\E\norm{ \bz^{t+2} - g(\bx^{t+1},\by^{t+1}) }^2 + C_2 \E\norm{\bw^{t+2} - \bw^{t+1}}^2\\ 
             &\leq (1-\frac{\eta\mu}{2} + 2\eta^2 \tilde{L}^2) \pr{\E\norm{\bw^{t}  -\bw^* }^2+C_1\E\norm{ \bz^{t+1} - g(\bx^{t},\by^{t}) }^2+ C_2 \E\norm{\bw^{t+1} - \bw^t}^2} +8\eta^2 \frac{(G_f^2+G_g^2)\sigma^2}{B} + 2\beta^2 C_1 \frac{\sigma^2}{M}\\
             & \quad + 4(1-\beta)^2 C_1 G_g^2\norm{\bw^{t+1} - \bw^t}^2  + \pare{2\eta^2G_g^2 L_f^2+ \frac{\eta}{\mu}G_g^2 L_f^2+ (1-\beta)^2C_1 -(1-\frac{\eta\mu}{2} + 2\eta^2 \tilde{L}^2)C_1} \E\norm{ \bz^{t+1} - g(\bx^t,\by^t) }^2\\
             & \quad + C_2 \E\norm{\bw^{t+2} - \bw^{t+1}}^2 - (1-\frac{\eta\mu}{2} + 2\eta^2 \tilde{L}^2) C_2 \E\norm{\bw^{t+1} - \bw^t}^2.
    \end{align*}}
    We now plug in the bound for $\E\norm{\bw^{t+2} - \bw^{t+1}}^2$ and define $A_t = \E\norm{\bw^{t}  -\bw^* }^2+C_1\E\norm{ \bz^{t+1} - g(\bx^{t},\by^{t}) }^2+ C_2 \E\norm{\bw^{t+1} - \bw^t}^2$:
   {\small \begin{align*}
       & A_{t+1} \leq (1-\frac{\eta\mu}{2} + 2\eta^2 \tilde{L}^2) A_t +8\eta^2 \frac{(G_f^2+G_g^2)\sigma^2}{B} + 2\beta^2 C_1 \frac{\sigma^2}{M} \\
       & +  \pare{2\eta^2G_g^2 L_f^2+ \frac{\eta}{\mu}G_g^2 L_f^2+ (1-\beta)^2C_1 -(1-\frac{\eta\mu}{2} + 2\eta^2 \tilde{L}^2)C_1} \E\norm{ \bz^{t+1} - g(\bx^t,\by^t) }^2\\
       &  + \pare{  4 (1-\beta)^2 C_1 G_g^2 - (1-\frac{\eta\mu}{2} + 2\eta^2 \tilde{L}^2)C_2 } \E\norm{\bw^{t+1} - \bw^t}^2 \\
       & + C_2\\
       &\times \pare{ \pare{1- {\eta \mu}  + \eta^2 \tilde{L}^2} \E\norm{\bw^{t+1} - \bw^{t}}^2  +  \frac{4 \eta}{\mu }    G_g^2 L_f^2\pare{ \E\norm{ \bz^{t+2} - g(\bx^{t+1},\by^{t+1}) }^2 +  \E\norm{ \bz^{t+1} - g(\bx^{t},\by^{t}) }^2 }+32 \frac{ \eta (G_g^2 + G_f^2 ) \sigma^2}{\mu B} } \\
       & = (1-\frac{\eta\mu}{2} + 2\eta^2 \tilde{L}^2) A_t +8\eta^2 \frac{(G_f^2+G_g^2)\sigma^2}{B} + 2\beta^2 C_1 \frac{\sigma^2}{M}+32 \frac{C_2 \eta (G_g^2 + G_f^2 ) \sigma^2}{\mu B} \\
       & +  \pare{2\eta^2G_g^2 L_f^2+ \frac{\eta}{\mu}G_g^2 L_f^2+ (1-\beta)^2C_1 -(1-\frac{\eta\mu}{2} + 2\eta^2 \tilde{L}^2)C_1} \E\norm{ \bz^{t+1} - g(\bx^t,\by^t) }^2\\
       &  + \pare{  4 (1-\beta)^2 C_1 G_g^2 + C_2  \pare{1- {\eta \mu}  + \eta^2 \tilde{L}^2} - (1-\frac{\eta\mu}{2} + 2\eta^2 \tilde{L}^2)C_2 } \E\norm{\bw^{t+1} - \bw^t}^2 \\
       & + C_2  \frac{4 \eta}{\mu }    G_g^2 L_f^2\pare{ \E\norm{ \bz^{t+2} - g(\bx^{t+1},\by^{t+1}) }^2 +  \E\norm{ \bz^{t+1} - g(\bx^{t},\by^{t}) }^2 }. 
    \end{align*}}
    Again applying the bound for $\E\norm{ \bz^{t+2} - g(\bx^{t+1},\by^{t+1}) }^2$ yields:
    {\small \begin{align*}
        A_{t+1} &\leq  (1-\frac{\eta\mu}{2} + 2\eta^2 \tilde{L}^2) A_t+8\eta^2 \frac{(G_f^2+G_g^2)\sigma^2}{B} + 2\beta^2 C_1 \frac{\sigma^2}{M}+32 \frac{C_2 \eta (G_g^2 + G_f^2 ) \sigma^2}{\mu B} \\
       & +  \pare{2\eta^2G_g^2 L_f^2+ \frac{\eta}{\mu}G_g^2 L_f^2+ (1-\beta)^2C_1 -(1-\frac{\eta\mu}{2} + 2\eta^2 \tilde{L}^2)C_1} \E\norm{ \bz^{t+1} - g(\bx^t,\by^t) }^2\\
       &  + \pare{  4 (1-\beta)^2 C_1 G_g^2 + C_2  \pare{1-{\eta \mu} + \eta^2 \tilde{L}^2} - (1-\frac{\eta\mu}{2} + 2\eta^2 \tilde{L}^2)C_2 } \E\norm{\bw^{t+1} - \bw^t}^2 \\
       & + C_2  \frac{4 \eta}{\mu }    G_g^2 L_f^2\pare{(1-\beta)^2 \E \norm{\bz^{t+1} - g(\bx^{t},\by^{t}) }^2  + 4(1-\beta)^2 G_g^2  \E\norm{\bw^{t+1} - \bw^{t}}^2 + 2\beta^2 \frac{\sigma^2}{M} +  \E\norm{ \bz^{t+1} - g(\bx^{t},\by^{t}) }^2 } \\ 
       & =  (1-\frac{\eta\mu}{2} + 2\eta^2 \tilde{L}^2) A_t+ 8\eta^2 \frac{(G_f^2+G_g^2)\sigma^2}{B} + 2\beta^2 C_1 \frac{\sigma^2}{M} +   C_2  \frac{8 \eta}{\mu }    G_g^2 L_f^2  \beta^2 \frac{\sigma^2}{M}+32 \frac{C_2 \eta (G_g^2 + G_f^2 ) \sigma^2}{\mu B}  \\
       & +  \pare{2\eta^2G_g^2 L_f^2+ \frac{\eta}{\mu}G_g^2 L_f^2+ (1-\beta)^2C_1 -(1-\frac{\eta\mu}{2} + 2\eta^2 \tilde{L}^2)C_1 + C_2  \frac{4 \eta}{\mu }G_g^2 L_f^2  \pare{1+(1-\beta)^2 }} \E\norm{ \bz^{t+1} - g(\bx^t,\by^t) }^2\\
       &  +  \pare{  4 (1-\beta)^2 C_1 G_g^2 -\pare{\frac{\eta \mu}{2}   -  \frac{16 \eta}{\mu }    G_g^4  L_f^2  (1-\beta)^2  }C_2} \E\norm{\bw^{t+1} - \bw^t}^2 .
    \end{align*}}
    We choose $\eta \leq  \min\cbr{\frac{\mu^3}{512}, \frac{\mu}{8\tilde{L}^2} }$ and $ (1-\beta)^2  \leq \frac{\mu^2}{256 G^4_g  L^2_f}$, we know 
     \begin{align*}
        A_{t+1} &\leq    (1-\frac{\eta\mu}{2} + 2\eta^2 \tilde{L}^2) A_t  + 8\eta^2 \frac{(G_f^2+G_g^2)\sigma^2}{B} + 2\beta^2 C_1 \frac{\sigma^2}{M} +   C_2  \frac{8 \eta}{\mu }    G_g^2 L_f^2  \beta^2 \frac{\sigma^2}{M}+32 \frac{C_2 \eta (G_g^2 + G_f^2 ) \sigma^2}{\mu B}\\
       & +  \underbrace{\pare{2\eta^2G_g^2 L_f^2+ \frac{\eta}{\mu}G_g^2 L_f^2+ (1-\beta)^2C_1 -(1-\frac{\eta\mu}{4}  )C_1 + C_2  \frac{8\eta}{\mu } G_g^2 L_f^2  }}_{T_1} \E\norm{ \bz^{t+1} - g(\bx^t,\by^t) }^2\\
       &  + \underbrace{\pare{  \frac{\mu^2}{64G_g^2 L_f^2} C_1 - \frac{\eta \mu}{4}    C_2}}_{T_2} \E\norm{\bw^{t+1} - \bw^t}^2 .
    \end{align*}
    Since we choose $ C_2 = \frac{\mu }{16 \eta G_g^2 L_f^2} C_1$, we know $T_2 \leq 0$, together with $\eta \leq \frac{\mu^3}{512}$ which yields:
     \begin{align*}
        A_{t+1} &\leq    (1-\frac{\eta\mu}{4}  ) A_t  + \pare{8\eta^2 + 32C_2\frac{\eta}{\mu}} \frac{(G_f^2+G_g^2)\sigma^2}{B} + \pare{2\beta^2 C_1 +   C_2  \frac{8 \eta}{\mu }    G_g^2 L_f^2  \beta^2}\frac{\sigma^2}{M}    \\
       & +  \underbrace{\pare{2\eta^2G_g^2 L_f^2+ \frac{\eta}{\mu}G_g^2 L_f^2+ (1-\beta)^2C_1 -(\frac{1}{2}-\frac{\eta\mu}{4}  )C_1   }}_{T_1} \E\norm{ \bz^{t+1} - g(\bx^t,\by^t) }^2 .
    \end{align*} 
To ensure $T_1 \leq 0$, we need  
    \begin{align*}  
        C_1 &\leq \frac{2\eta^2G_g^2L_f^2 + \frac{\eta}{\mu} G_g^2 L_f^2}{\frac{1}{2}-\frac{\eta\mu}{4}  - (1-\beta)^2}, 
    \end{align*}
   It can be satisfied If we choose 
     \begin{align*} 
        (1-\beta)^2 &\leq \frac{1}{2}-\frac{\eta\mu}{4} - \frac{\mu}{220 \max\cbr{1,\frac{1}{\mu},G_f^2}} \\ 
        C_1 &=      220 \eta G_g^2 L_f^2 \max\cbr{1,\frac{1}{\mu}, G_g^2}
    \end{align*}
     and we know $T_1 \leq 0$, which concludes the proof. 
     
    \end{proof} 
\end{lemma}

\begin{proposition}[Three points]\label{prop: 3 points}
If $ \tilde\bw = \cP_{\cW}(\bv - \eta \bg  )$, then the following statements holds:
    \begin{align*}
        \inprod{\eta \bg}{  \tilde\bw} + \frac{1}{2} \norm{\bv - \tilde\bw}^2 + \frac{1}{2} \norm{\bw- \tilde\bw }^2 \leq \inprod{\eta \bg}{ \bw}  + \frac{1}{2} \norm{\bw - \bv}^2, \forall \bw \in \cW.
    \end{align*}
\end{proposition}

 \begin{lemma} \label{lem: bound itr difference}
 Define $\bu^{t+1} = \cP_{\cW}(\bw^t - \eta \bg^{t+1})$. Then for Algorithm~\ref{algorithm: CODA-SCSC}, under the assumptions of Theorem~\ref{thm:wcwc}, the following statement holds:
   \begin{align*}
         \frac{1}{8}\E \norm{\bw^t - \bw^{t+1}}^2  
         \leq  &   \frac{1}{2} \E\norm{\bw^* - \bw^t}^2  + \pare{ \frac{2\eta}{\mu} + 72\eta }G_g^2 L_f^2 \beta^2 \frac{\sigma^2}{M}+ \frac{48\eta (G_g^2+G_f^2)\sigma^2}{B} \\
          & + \pare{72\eta G_g^2 L_f^2 + 144  \eta G_g^2 G_f^2 L_f^2 + \frac{\eta}{\mu} G_g^2 L_f^2 }\E\norm{\bz^{t+1} - g(\bx^t,\by^t) }^2 .
     \end{align*} 
 \begin{proof}
     In Proposition~\ref{prop: 3 points}, we set $ \tilde\bw = \bw^{t+1}$, $\bv = \bw^t$ and $\bw = \bu^{t+1}$
     \begin{align}
         &\inprod{\eta \bg^t}{ \bw^{t+1}} + \frac{1}{2} \norm{\bw^t - \bw^{t+1}}^2 + \frac{1}{2} \norm{\bu^{t+1}- \bw^{t+1}}^2 \leq \inprod{\eta \bg^t}{\bu^{t+1}}  + \frac{1}{2} \norm{\bu^{t+1} - \bw^t}^2 \nonumber\\
          &\Longleftrightarrow \nonumber\\
          & \inprod{\eta \bg^t}{ \bw^{t+1} - \bu^{t+1}} + \frac{1}{2} \norm{\bw^t - \bw^{t+1}}^2 + \frac{1}{2} \norm{\bu^{t+1}- \bw^{t+1}}^2 \leq   \frac{1}{2} \norm{\bu^{t+1} - \bw^t}^2 . 
         \label{eq: iterate diff 1}
     \end{align}
    Again in Proposition~\ref{prop: 3 points}, we set $ \tilde\bw = \bu^{t+1}$, $\bv = \bw^t$, $\bg = \bg^{t+1}$
     \begin{align}
         &\inprod{\eta \bg^{t+1}}{ \bu^{t+1}} + \frac{1}{2} \norm{\bw^t - \bu^{t+1}}^2 + \frac{1}{2} \norm{\bw- \bu^{t+1}}^2 \leq \inprod{\eta \bg^{t+1}}{ \bw}  + \frac{1}{2} \norm{\bw - \bw^t}^2 \nonumber\\
         &\Longleftrightarrow \nonumber\\
         & \inprod{\eta \bg^{t+1}}{ \bu^{t+1} - \bw} + \frac{1}{2} \norm{\bw^t - \bu^{t+1}}^2 + \frac{1}{2} \norm{\bw- \bu^{t+1}}^2 \leq  \frac{1}{2} \norm{\bw - \bw^t}^2 \nonumber\\
         &\Longleftrightarrow \nonumber\\
         & \inprod{\eta \bg^{t+1}}{ \bw^{t+1} - \bw}-\inprod{\eta \bg^{t+1}}{ \bw^{t+1} - \bu^{t+1} } + \frac{1}{2} \norm{\bw^t - \bu^{t+1}}^2 + \frac{1}{2} \norm{\bw- \bu^{t+1}}^2 \leq  \frac{1}{2} \norm{\bw - \bw^t}^2.
         \label{eq: iterate diff 2}
     \end{align}
     Adding (\ref{eq: iterate diff 1}) and (\ref{eq: iterate diff 2}) yields:
     \begin{align*}
           &\eta\inprod{ \bg^t -\bg^{t+1} }{ \bw^{t+1} - \bu^{t+1}} +  \inprod{\eta \bg^{t+1}}{ \bw^{t+1} - \bw } \\
           &\qquad \qquad  \qquad + \frac{1}{2} \norm{\bw^t - \bw^{t+1}}^2 + \frac{1}{2} \norm{\bu^{t+1}- \bw^{t+1}}^2 + \frac{1}{2} \norm{\bw- \bu^{t+1}}^2 \leq \frac{1}{2} \norm{\bw - \bw^t}^2\\
            &\Longleftrightarrow\\
           &\eta\inprod{ \bg^t -\bg^{t+1} }{ \bw^{t+1} - \bu^{t+1}} +  \inprod{\eta \bg(\bw^{t+1})}{ \bw^{t+1} - \bw }+  \eta\inprod{ \bg^{t+1} - \bg(\bw^{t+1})}{ \bw^{t+1} - \bw }  \\
           &\qquad \qquad  \qquad + \frac{1}{2} \norm{\bw^t - \bw^{t+1}}^2 + \frac{1}{2} \norm{\bu^{t+1}- \bw^{t+1}}^2 + \frac{1}{2} \norm{\bw- \bu^{t+1}}^2 \leq \frac{1}{2} \norm{\bw - \bw^t}^2.
     \end{align*} 
     Setting $\bw = \bw^*$, according to $\frac{1}{\gamma}$-strongly monotone of $\bg(\cdot)$, we have
     \begin{align*} 
       &\eta\inprod{ \bg^t -\bg^{t+1} }{ \bw^{t+1} - \bu^{t+1}} +   \eta\mu\norm{ \bw^{t+1} - \bw^* }^2 +  \eta\inprod{ \bg^{t+1} - \bg(\bw^{t+1})}{ \bw^{t+1} - \bw^* }  \\
           &\qquad \qquad  \qquad + \frac{1}{2} \norm{\bw^t - \bw^{t+1}}^2 + \frac{1}{2} \norm{\bu^{t+1}- \bw^{t+1}}^2 + \frac{1}{2} \norm{\bw^*- \bu^{t+1}}^2 \leq \frac{1}{2} \norm{\bw^* - \bw^t}^2.
     \end{align*}
     For the two inner products, we apply Cauchy-Schwartz:
     \begin{align*}
         \eta\E\inprod{ \bg^{t+1} - \bg(\bw^{t+1})}{ \bw^{t+1} - \bw^* } &=  \eta\E\inprod{ \bar{\bg}^{t+1} - \bg(\bw^{t+1})}{ \bw^{t+1} - \bw^* } \\
         &\geq -\frac{1}{2}\eta\pare{\frac{1}{\mu}\E\norm{ \bar{\bg}^{t+1} - \bg(\bw^{t+1})}^2 +\mu\E\norm{ \bw^{t+1} - \bw^* }^2} \\
          &\geq -\frac{1}{2}\eta\pare{\frac{1}{\mu}2G_g^2 L_f^2\E\norm{ \bz^{t+2} -  g(\bx^{t+1},\by^{t+1}) }^2 +\mu\E\norm{ \bw^{t+1} - \bw^* }^2} ,
     \end{align*}
     and
     \begin{align*}
         \eta\inprod{ \bg^t -\bg^{t+1} }{ \bw^{t+1} - \bu^{t+1}} &\geq - \frac{1}{2} \eta \pare{ \norm{\bg^t -\bg^{t+1}}^2 + \norm{\bw^{t+1} - \bu^{t+1}}^2 }\\
      &   = - \frac{1}{2} \eta \pare{ \norm{\bg_{\bx}^t -\bg_{\bx}^{t+1}}^2 + \norm{\bg_{\by}^t -\bg_{\by}^{t+1}}^2 + \norm{\bw^{t+1} - \bu^{t+1}}^2 } .
     \end{align*}
     For $\norm{\bg_{\bx}^t -\bg_{\bx}^{t+1}}^2$ we have:
     \begin{align*}
         \norm{\bg_{\bx}^t -\bg_{\bx}^{t+1}}^2  
         \leq &6\underbrace{\E\norm{\nabla f(\bz^{t+1} )\nabla_{\bx} g(\bx^t,\by^t ) -\nabla f(\bz^{t+2})\nabla_{\bx} g(\bx^{t+1},\by^{t+1}) }^2}_{A} + \frac{4}{\gamma^2} \E\norm{\bx^t - \bx^{t+1}}^2\\ 
       &+ 6\underbrace{\E\norm{\nabla f(\bz^{t+1} )\nabla_{\bx} g(\bx^t,\by^t ) -\frac{1}{B}\sum_{(\zeta,\xi)\in\cB^{t}}\nabla f(\bz^{t+1};\zeta)\nabla_{\bx} g(\bx^t,\by^t;\xi) }^2}_{B} \\
       &+ 6\underbrace{\E\norm{\nabla f(\bz^{t+2})\nabla_{\bx} g(\bx^{t+1},\by^{t+1})  - \frac{1}{B}\sum_{(\zeta,\xi)\in\cB^{t+1}} \nabla f(\bz^{t+2};\zeta)\nabla_{\bx} g(\bx^{t+1},\by^{t+1};\xi) }^2}_{C}. 
     \end{align*}
     We bound above terms as follows. For A:
     \begin{align*}
         A \leq 2G_g^2 L_f^2 \E\norm{\bz^{t+1} - \bz^{t+2}}^2 + 2 G_f^2 L_g^2 \norm{(\bx^t,\by^t) - (\bx^{t+1},\by^{t+1})}^2
     \end{align*}
    For B and C:
    \begin{align*}
        B, C \leq  \frac{2(G_g^2+G_f^2)\sigma^2}{B}  
    \end{align*}
    Putting pieces together we know 
    \begin{align*}
        \norm{\bg_\bx^t - \bg_\bx^{t+1}}^2 & \leq 12G_g^2 L_f^2 \E\norm{\bz^{t+1} - \bz^{t+2}}^2 + 12 G_f^2 L_g^2 \norm{(\bx^t,\by^t) - (\bx^{t+1},\by^{t+1})}^2 + \frac{4}{\gamma^2} \E\norm{\bx^t - \bx^{t+1}}^2\\
        &\quad +\frac{48(G_g^2+G_f^2)\sigma^2}{B} .
    \end{align*} 
    Similarly, we have
    \begin{align*}
           \norm{\bg_\by^t - \bg_\by^{t+1}}^2 &\leq 12G_g^2 L_f^2 \E\norm{\bz^{t+1} - \bz^{t+2}}^2 + 12 G_f^2 L_g^2 \norm{(\bx^t,\by^t) - (\bx^{t+1},\by^{t+1})}^2 + \frac{4}{\gamma^2} \E\norm{\by^t - \by^{t+1}}^2 \\
           &\quad + \frac{48(G_g^2+G_f^2)\sigma^2}{B} .
    \end{align*}
     Hence 
    {\small \begin{align*}
        & \eta\inprod{ \bg^t -\bg^{t+1} }{ \bw^{t+1} - \bu^{t+1}} \geq - \frac{1}{2} \eta \pare{ 24G_g^2 L_f^2 \E\norm{\bz^{t+1} - \bz^{t+2}}^2 + \pare{24 G_f^2 L_g^2 +\frac{4}{\gamma^2}   }\norm{ \bw^t -  \bw^{t+1} }^2  + \frac{96(G_g^2+G_f^2)\sigma^2}{B} }. 
     \end{align*}}
     Plugging back yields:
\begin{align*} 
       &- \frac{1}{2} \eta \pare{ 24G_g^2 L_f^2 \E\norm{\bz^{t+1} - \bz^{t+2}}^2 + \pare{24 G_f^2 L_g^2 +\frac{4}{\gamma^2}   }\norm{ \bw^t -  \bw^{t+1} }^2  + \frac{96(G_g^2+G_f^2)\sigma^2}{B}  }\\
       &-\frac{1}{2}\eta\pare{\frac{1}{\mu}2G_g^2 L_f^2\E\norm{ \bz^{t+2} -  g(\bx^{t+1},\by^{t+1}) }^2 +\frac{1}{\gamma}\E\norm{ \bw^{t+1} - \bw^* }^2}  \\
           &+   \eta\frac{1}{\gamma}\norm{ \bw^{t+1} - \bw^* }+ \frac{1}{2} \norm{\bw^t - \bw^{t+1}}^2 + \frac{1}{2} \norm{\bu^{t+1}- \bw^{t+1}}^2 + \frac{1}{2} \norm{\bw^*- \bu^{t+1}}^2 \leq \frac{1}{2} \norm{\bw^* - \bw^t}^2\\
             &\Longleftrightarrow\\
        &-    12\eta G_g^2 L_f^2 \E\norm{\bz^{t+1} - \bz^{t+2}}^2   -  \frac{\eta}{\mu} G_g^2 L_f^2\E\norm{ \bz^{t+2} -  g(\bx^{t+1},\by^{t+1}) }^2 + \frac{\mu \eta}{2 }\norm{ \bw^{t+1} - \bw^* }^2  - \frac{48\eta (G_g^2+G_f^2)\sigma^2}{B}  \\
           &+ \frac{1}{2}\pare{1 -\eta\pare{24 G_f^2 L_g^2 +\frac{4}{\gamma^2}   }}\norm{\bw^t - \bw^{t+1}}^2 + \pr{\frac{1}{2}-\frac{\eta}{2}} \norm{\bu^{t+1}- \bw^{t+1}}^2 + \frac{1}{2} \norm{\bw^*- \bu^{t+1}}^2 \leq \frac{1}{2} \norm{\bw^* - \bw^t}^2.
     \end{align*}
     By our choice $\eta  \leq \frac{1}{2}\pare{24 G_f^2 L_g^2 +\frac{4}{\gamma^2} }^{-1} $, we have $\pare{1 -\eta\pare{24 G_f^2 L_g^2 +\frac{4}{\gamma^2}   }} \geq \frac{1}{2}$, which yields:
   {\small  \begin{align*}
         \frac{1}{4} \norm{\bw^t - \bw^{t+1}}^2   &  \leq \frac{1}{2} \norm{\bw^* - \bw^t}^2 + 12\eta G_g^2 L_f^2 \E\norm{\bz^{t+1} - \bz^{t+2}}^2   + \frac{\eta}{\mu} G_g^2 L_f^2\E\norm{ \bz^{t+2} -  g(\bx^{t+1},\by^{t+1}) }^2 + \frac{48\eta (G_g^2+G_f^2)\sigma^2}{B}\\
          &\leq \frac{1}{2} \norm{\bw^* - \bw^t}^2 + \frac{48\eta (G_g^2+G_f^2)\sigma^2}{B}\\
           & +  12\eta G_g^2 L_f^2 \pare{ 3\E\norm{\bz^{t+1} - g(\bx^t,\by^t)}^2 + 3\E\norm{  g(\bx^t,\by^t) -g(\bx^{t+1},\by^{t+1}) }^2 + 3\E\norm{g(\bx^{t+1},\by^{t+1})-\bz^{t+2}}^2}  \\
          &+ \frac{\eta}{\mu} G_g^2 L_f^2
          \pare{ (1-\beta)^2 \E \norm{\bz^{t+1} - g(\bx^{t},\by^{t}) }^2  + 4(1-\beta)^2 G_g^2 \norm{\bw^{t+1} - \bw^{t}}^2 + 2\beta^2 \frac{\sigma^2}{M} }\\
          &\leq \frac{1}{2} \norm{\bw^* - \bw^t}^2 + \pare{36\eta G_g^4 L_f^2  + 4(1-\beta)^2 \frac{\eta}{\mu} G_g^4 L_f^2 }\norm{\bw^t - \bw^{t+1}}^2 + 2\frac{\eta}{\mu} G_g^2 L_f^2 \beta^2 \frac{\sigma^2}{M} + \frac{48\eta (G_g^2+G_f^2)\sigma^2}{B}\\
          & + \pare{36\eta G_g^2 L_f^2 +(1-\beta)^2\frac{\eta}{\mu} G_g^2 L_f^2 }\E\norm{\bz^{t+1} - g(\bx^t,\by^t) }^2   +  36\eta G_g^2 L_f^2 \E\norm{ \bz^{t+2} -  g(\bx^{t+1},\by^{t+1}) }^2.
     \end{align*}}
     Plugging bound for $\E\norm{ \bz^{t+2} -  g(\bx^{t+1},\by^{t+1}) }^2$ yields:
    {\small \begin{align*}
         \frac{1}{4} \norm{\bw^t - \bw^{t+1}}^2  
          &\leq \frac{1}{2} \norm{\bw^* - \bw^t}^2 + \pare{36\eta G_g^4 L_f^2  + 4(1-\beta)^2  \frac{\eta}{\mu} G_g^4 L_f^2 }\norm{\bw^t - \bw^{t+1}}^2 + 2\frac{\eta}{\mu} G_g^2 L_f^2 \beta^2 \frac{\sigma^2}{M} + \frac{48\eta (G_g^2+G_f^2)\sigma^2}{B}\\
          & + \pare{36\eta G_g^2 L_f^2 +(1-\beta)^2\frac{\eta}{\mu} G_g^2 L_f^2 }\E\norm{\bz^{t+1} - g(\bx^t,\by^t) }^2 \\
          &+  36\eta G_g^2 L_f^2 \pare{ (1-\beta)^2 \E \norm{\bz^{t+1} - g(\bx^{t},\by^{t}) }^2  + 4(1-\beta)^2 G_g^2 \norm{\bw^{t+1} - \bw^{t}}^2 + 2\beta^2 \frac{\sigma^2}{M} } \\
          & =  \frac{1}{2} \norm{\bw^* - \bw^t}^2 + \pare{36\eta G_g^4 L_f^2  + 4(1-\beta)^2  \frac{\eta}{\mu} G_g^4 L_f^2 }\norm{\bw^t - \bw^{t+1}}^2 + \pare{2\frac{\eta}{\mu} + 72\eta }G_g^2 L_f^2 \beta^2 \frac{\sigma^2}{M}\\
          & + \pare{72\eta G_g^2 L_f^2 + 144 (1-\beta)^2 \eta G_g^4  L_f^2 +(1-\beta)^2\frac{\eta}{\mu} G_g^2 L_f^2 }\E\norm{\bz^{t+1} - g(\bx^t,\by^t) }^2 + \frac{48\eta (G_g^2+G_f^2)\sigma^2}{B} .
     \end{align*}}
     Since we choose $\eta \leq \min\cbr{\frac{1}{576 G_g^4 L_f^2}, \frac{\mu}{64 G_g^4   L_f^2}}$, we know $36\eta G_g^4 L_f^2  + 4(1-\beta)^2 \frac{\eta}{\mu} G_g^4 L_f^2 \leq \frac{1}{8}$, which yields:
     \begin{align*}
         \frac{1}{8}\E \norm{\bw^t - \bw^{t+1}}^2  
          &\leq    \frac{1}{2} \E\norm{\bw^* - \bw^t}^2  + \pare{ \frac{2\eta}{\mu} + 72\eta }G_g^2 L_f^2 \beta^2 \frac{\sigma^2}{M}+ \frac{48\eta (G_g^2+G_f^2)\sigma^2}{B} \\
          & + \pare{72\eta G_g^2 L_f^2 + 144  \eta G_g^4  L_f^2 + \frac{\eta}{\mu} G_g^2 L_f^2 }\E\norm{\bz^{t+1} - g(\bx^t,\by^t) }^2 .
     \end{align*}
 \end{proof}
 \end{lemma}

\subsection{Proof of Lemma~\ref{lem:conv PD SMV}}
\begin{proof}
    First according to~[Eq.(23) of~\cite{liu2021first}],
    \begin{align}
        \max_{\bw \in \cW} \inprod{\nabla F_k(\bw^{t+1})}{\bw^{t+1} - \bw} \leq (D+\frac{1}{\eta}) \E \norm{\bw^{t+1} - \bw^t}^2. \label{eq: pf conv SMV 0}
    \end{align}
    Revoking Lemma~\ref{lem: bound itr difference}
    \begin{align}
         \frac{1}{8} \norm{\bw^t - \bw^{t+1}}^2  
          &\leq    \frac{1}{2} \norm{\bw^* - \bw^t}^2   + \pare{ \frac{2\eta}{\mu} + 72\eta }G_g^2 L_f^2 \beta^2 \frac{\sigma^2}{M}+ \frac{48\eta (G_g^2+G_f^2)\sigma^2}{B} \nonumber\\
          &\quad  + \pare{72\eta G_g^2 L_f^2 + 144  \eta G_g^4  L_f^2 + \frac{\eta}{\mu} G_g^2 L_f^2 }\E\norm{\bz^{t+1} - g(\bx^t,\by^t) }^2 \nonumber\\
            &\leq    \frac{1}{2} \E\norm{\bw^* - \bw^t}^2  + \pare{ \frac{2\eta}{\mu} + 72\eta }G_g^2 L_f^2 \beta^2 \frac{\sigma^2}{M}+ \frac{48\eta (G_g^2+G_f^2)\sigma^2}{B} \nonumber\\
          &\quad  +  220 \eta G_g^2 L_f^2 \max\cbr{1,\frac{1}{\mu}, G_g^2}\E\norm{\bz^{t+1} - g(\bx^t,\by^t) }^2 \label{eq: pf conv SMV 1}.
     \end{align}

     To bound $\E\norm{\bw^* - \bw^t}^2$, we revoke Lemma~\ref{lem:conv to optimal}:
       \begin{align*}
        &\norm{\bw^{t+1}-\bw^*}^2 + C_1\E\norm{\bz^{t+2} - g(\bx^{t+1},\by^{t+1})}^2 + C_2 \norm{ \bw^{t+1}-\bw^t  }^2\\
        \leq &\pare{1- \frac{\eta \mu}{2} } \pare{ \norm{\bw^{t}-\bw^*}^2 + C_1\E\norm{\bz^{t+1} - g(\bx^{t},\by^{t})}^2 + C_2 \norm{ \bw^{t}-\bw^{t-1}  }^2} \\
        &+ \pare{8\eta^2 + 32C_2\frac{\eta}{\mu}} \frac{(G_f^2+G_g^2)\sigma^2}{B} + \pare{2\beta^2 C_1 +   C_2  \frac{8 \eta}{\mu }    G_g^2 L_f^2  \beta^2}\frac{\sigma^2}{M} \\
        \leq & \pare{1- \frac{\eta \mu}{2} }^{t} \pare{ \norm{\bw^{0}-\bw^*}^2 + C_1\E\norm{\bz^{1} - g(\bx^{0},\by^{0})}^2 + C_2 \norm{ \bw^{0}-\bw^{-1}  }^2} \\
        & + \frac{2}{\eta \mu}\pare{ \pare{8\eta^2 + 32C_2\frac{\eta}{\mu}} \frac{(G_f^2+G_g^2)\sigma^2}{B} + \pare{2  C_1 +   C_2  \frac{8 \eta}{\mu }    G_g^2 L_f^2  }\frac{\sigma^2}{M} }.   
    \end{align*}
    where $C_1 =  220 \eta G_g^2 L_f^2 \max\cbr{1,\frac{1}{\mu}, G_g^2}$, and $C_2 =   14 \mu \max\cbr{1,\frac{1}{\mu}, G_g^2}$. Plugging above bound in (\ref{eq: pf conv SMV 1}) yields:
    \begin{align*}
         \frac{1}{8} \norm{\bw^t - \bw^{t+1}}^2   
            &\leq  \pare{1- \frac{\eta \mu}{2} }^{t} \pare{ \norm{\bw^{0}-\bw^*}^2 + C_1\E\norm{\bz^{1} - g(\bx^{0},\by^{0})}^2 } \\
            &  + \frac{2}{\eta \mu}\pare{ \pare{8\eta^2 + 32C_2\frac{\eta}{\mu}} \frac{(G_f^2+G_g^2)\sigma^2}{B} + \pare{2  C_1 +   C_2  \frac{8 \eta}{\mu }    G_g^2 L_f^2  }\frac{\sigma^2}{M} }  \\ 
            &+   \pare{ \frac{2\eta}{\mu} + 72\eta }G_g^2 L_f^2 \beta^2 \frac{\sigma^2}{M}+ \frac{48\eta (G_g^2+G_f^2)\sigma^2}{B}.
     \end{align*}
     Combining with (\ref{eq: pf conv SMV 0}) yields:
     \begin{align*}
          \max_{\bw \in \cW} \inprod{\nabla F_k(\bw^{t+1})}{\bw^{t+1} - \bw} & \leq \pare{D+\frac{1}{\eta}} \pare{1- \frac{\eta \mu}{2} }^{t} \pare{ \norm{\bw^{0}-\bw^*}^2 + C_1\E\norm{\bz^{1} - g(\bx^{0},\by^{0})}^2  }\\
        &  + \pare{D+\frac{1}{\eta}}\frac{2}{\eta \mu}\pare{ \pare{8\eta^2 + 32C_2\frac{\eta}{\mu}} \frac{(G_f^2+G_g^2)\sigma^2}{B} + \pare{2  C_1 +   C_2  \frac{8 \eta}{\mu }    G_g^2 L_f^2  }\frac{\sigma^2}{M} }  \\ 
            &+  \pare{D+\frac{1}{\eta}} \pare{ \frac{2\eta}{\mu} + 72\eta }G_g^2 L_f^2 \beta^2 \frac{\sigma^2}{M}+ \pare{D+\frac{1}{\eta}}\frac{48\eta (G_g^2+G_f^2)\sigma^2}{B}.
     \end{align*}
     To ensure RHS is less than $\epsilon$, we need $t \geq \Omega( 
\frac{1}{\eta\mu} \ln(\frac{D+1/\eta}{\epsilon}))$ and $M = B \geq \Omega\pare{(D+\frac{1}{\eta})\frac{G_g^2 L_f^2\sigma^2}{\mu \epsilon} \max\cbr{1,\frac{1}{\mu},G_g^2})}$, which yields the gradient complexity of
\begin{align*}
    O\pare{ \pare{D+\frac{1}{\eta}}\frac{L^2\sigma^2}{\mu^2 \eta\epsilon} \log \pare{\frac{D+1/\eta}{\epsilon}}  }.
\end{align*}

\end{proof}

\subsection{Proof of Theorem~\ref{thm:wcwc}}
\begin{proof}
Let $\bar\bw$ be the solution of SVI induced by $F^\gamma_{\bw}:= F(\bx',\by') + \frac{1}{2\gamma}(\bw' - \bw)$. Let $\bw_{k^*} = (\bx_{k^*}, \by_{k^*} ) $ be the solution returned by Algorithm~\ref{algorithm: CODA-PD}.
 According to [Eq.(21) ~\cite{liu2021first}], if Algorithm~\ref{algorithm: CODA-SCSC} returns a $\tilde \epsilon$-accurate solution, then we have: 
 \begin{align*}
     \E\norm{\bw_{k^*}  - \bar{\bw}_{k^*}}^2 \leq \frac{2D\theta_K}{\sum_{k=0}^K \theta_k} + \frac{4c\tilde{\epsilon}}{\rho}.
 \end{align*}
 Since we choose $\theta_k = (k+1)^\alpha$, we have
 \begin{align*}
     \E\norm{\bw_{k^*} - \bar{\bw}_{k^*}}^2 \leq \frac{2D(\alpha+1)}{ K} + \frac{4c\tilde{\epsilon}}{\rho}.
 \end{align*}
 According to [Lemma~6~\cite{liu2021first}], we know that if $\norm{\bw - \bar\bw} \leq \gamma \epsilon$, then
 \begin{align*}
      dist\pare{0, \partial (F(\bar\bx,\bar\by) + \mathbf{1}_{\cX\times \cY}(\bar\bx,\bar\by) )} \leq  \epsilon
 \end{align*}
 To ensure $ \E\norm{\bw_{k^*} - \bar{\bw}_{k^*}}^2$ is less than $\gamma^2 \epsilon^2$ where $\gamma = \frac{1}{2\rho}$, we need $K \geq \frac{16 \rho^2 D(\alpha+1)}{  \epsilon^2}$, and $\tilde \epsilon \leq \frac{  \epsilon^2}{32\rho c}$, together by our choice $\eta  = \Theta\pr{\frac{1}{L^2}}$  which yields the gradient complexity of
 \begin{align*}
      O\pare{ \pare{D+L^2}\frac{  \rho^2 D  L^4\sigma^2}{  \epsilon^4} \log \pare{\frac{1}{\epsilon}}   }. 
 \end{align*}

\end{proof}

 \section{Proof of Varaince Reduction Algorithm}\label{app: proof VR}
In this Section we will provide proofs for Theorem~\ref{thm:NCSC+} and~\ref{thm:ncc+}.

In the proof we will use the following notations:
\begin{align*}
   & F_k(\bx^t,\by^t) = F(\bx^t,\by^t) + \frac{\mu_x+L}{2}\norm{\bx^t - \bx_k}^2,\\
   & \bar{\bg}^t_{\bx} =   \nabla_1 f(\bz^{t+1},\by^{t+1} )   \nabla g(\bx^t)  + (\mu_x + L) (\bx^t - \bx^0),\\
  &  \bar{\bg}^t_{\by} =   \nabla_{\by} f(\bz^{t+1},\by^{t+1} ) ,\\
\end{align*}
For notational convenience, we define $\tilde \sigma^2 = \max\cbr{\sigma^2, 2(G_f^2+G_g^2)\sigma^2}$.

The following Lemma shows that the tracking error for inner function is reduced via our variance reduction technique.
\begin{lemma}[Tracking error with variance reduction]
For Algorithm~\ref{algorithm: CODA-Primal+}, under assumptions of Theorem~\ref{thm:NCSC+}, the following statement holds:
\begin{align*}
      \E \norm{  \bz^{t+1} - g(\bx^t)}^2  & \leq \pare{1-\beta} \E\norm{\bz^{t} - g(\bx^{t-1})}^2 \\
        &\quad + \pare{1+\frac{1}{\beta}} \pare{6(1-\beta)^2G_g^2\E\norm{  \bx^{t-1}  -\bx^t}^2 + 3\beta^2  \frac{\sigma^2}{B_{\tau}}  + \frac{G_g^2}{B}\sum_{t'=t_p}^{t-1} \E\norm{ \bx^{t'}  - \bx^{t'-1}   }^2 }.
\end{align*}
where $t_p$ is the last iteration $t$ such that $t \mod \tau = 0$.
\begin{proof}
    Recall the updating rule:
    \begin{align*}
     \bz^{t+1} = (1-\beta)(\bz^t + g(\bx^t; \cI^t_{\bz}) - g(\bx^{t-1}; \cI^t_{\bz})) + \beta g^{t-1}, 
    \end{align*}
    where $g^{t} = g^{t-1} + g(\bx^t;\cI^t_{\bz}) - g(\bx^{t-1};\cI^t_{\bz})$ and every $\tau$ iterations, we use large mini-batch to update $g^{t}$: $g^{t} =  g(\bx^t;\cB^t_{\bz})  $. We hence have:
    \begin{align*}
       \E \norm{g^t - g(\bx^t)}^2 &= \E\norm{g^{t-1} - g(\bx^{t-1}) + \frac{1}{B}\sum_{\xi \in \cI^t_{\bz}} g(\bx^{t};\xi ) -  g(\bx^{t-1};\xi ) + g(\bx^{t-1})   - g(\bx^t)}^2\\
       &=  \E\norm{g^{t-1} - g(\bx^{t-1})}^2 +\E\norm{\frac{1}{B}\sum_{\xi \in \cI^t_{\bz}} 
 g(\bx^{t};\xi ) -  g(\bx^{t-1};\xi ) + g(\bx^{t-1})   - g(\bx^t)}^2\\
       &\leq \E\norm{g^{t-1} - g(\bx^{t-1})}^2 + \frac{1}{B^2}\sum_{\xi \in \cI^t_{\bz}} \E\norm{g(\bx^{t};\xi^{t}) -  g(\bx^{t-1};\xi^{t})  }^2\\
        &\leq \E\norm{g^{t-1} - g(\bx^{t-1})}^2 + \frac{G_g^2}{B}\E\norm{ \bx^{t}  - \bx^{t-1}   }^2\\
        &\leq \E\norm{g^{t_p} - g(\bx^{t_p})}^2 + \frac{G_g^2}{B} \sum_{t'=t_p}^{t-1}  \E\norm{ \bx^{t'}  - \bx^{t'-1}   }^2\\
        &\leq \frac{\sigma^2}{B_{\tau}} + \frac{G_g^2}{B}\sum_{t'=t_p}^{t-1} \E\norm{ \bx^{t'}  - \bx^{t'-1}   }^2.
    \end{align*}
Again by the updating rule we have:
    \begin{align*}
        \bz^{t+1} - g(\bx^t) = (1-\beta) (\bz^{t} - g(\bx^{t-1})) + (1-\beta) (g(\bx^{t-1}) - g(\bx^t)) + \beta (g^t - g(\bx^t)) + (1-\beta)(g(\bx^t;\xi^t) - g(\bx^{t-1};\xi^t)).
    \end{align*}
    Taking expected norm on both sides and applying Young's inequality yields:
    \begin{align*}
        \E \norm{  \bz^{t+1} - g(\bx^t)}^2 &\leq \underbrace{(1+\beta)(1-\beta)^2}_{\leq 1-\beta} \E\norm{\bz^{t} - g(\bx^{t-1})}^2 \\
        &+ (1+\frac{1}{\beta}) \E\norm{  (1-\beta) (g(\bx^{t-1}) - g(\bx^t)) + \beta (g^t - g(\bx^t)) + (1-\beta)(g(\bx^t;\xi^t) - g(\bx^{t-1};\xi^t)) }^2\\
        &\leq \pare{1-\beta} \E\norm{\bz^{t} - g(\bx^{t-1})}^2 \\
        &+ (1+\frac{1}{\beta}) \pare{6(1-\beta)^2G_g^2\E\norm{  \bx^{t-1}  -\bx^t}^2 + 3\beta^2\E\norm{  g^t - g(\bx^t) }^2 }\\
        & \leq \pare{1-\beta} \E\norm{\bz^{t} - g(\bx^{t-1})}^2 \\
        &+ (1+\frac{1}{\beta}) \pare{6(1-\beta)^2G_g^2\E\norm{  \bx^{t-1}  -\bx^t}^2 + 3\beta^2\pare{ \frac{\sigma^2}{B_{\tau}}} + \frac{G_g^2}{B} \sum_{t'=t_p}^{t-1}  \E\norm{ \bx^{t'}  - \bx^{t'-1}   }^2 }.
    \end{align*}
    \end{proof}

\end{lemma}
 
\begin{lemma}[Gradient bounds]\label{lem:grad bounds VR}
For Algorithm~\ref{algorithm: CODA-Primal+}, under assumptions of Theorem~\ref{thm:NCSC+}, the following statement holds:
    \begin{align*}
        &\E\norm{\bg^t_{\bx} -  \bar{\bg}^t_{\bx} }^2 \leq  \frac{\tilde\sigma^2}{B_{\tau}}  + \frac{1}{B }\sum_{t'=t_p}^{t-1} \pare{2G_g^2 L^2\norm{ \bz^{t'+1}   -  \bz^{t'+2}    }^2 +2G_g^2 L^2\norm{  \by^{t'+1}  -  \by^{t'+2}  }^2 + 2G_f^2 L_g^2\norm{  \bx^{t'} -  \bx^{t'+1}    }^2}\\
      &  \E\norm{\bq^t - \bar{\bg}^t_{\by}}^2 \leq \frac{\tilde\sigma^2}{B_{\tau}}  +  \frac{1}{B} \sum_{t'=t_p}^{t}  L^2_f \pare{ \E\norm{ \bz^{t'}   -  \bz^{t'+1}  }^2 + \E\norm{  \by^{t'}  -   \by^{t'+1} }^2 }.
    \end{align*}
    \begin{proof}
    By updating rule in Algorithm~\ref{algorithm: CODA-Primal+} we have:
        \begin{align*}
           \E\norm{\bg^t_{\bx} -  \bar{\bg}^t_{\bx} }^2  & =  \E\norm{\bg^{t-1}_{\bx} +   G^k_{\bx}(\bx^t,\bz^{t+1}, \by^{t+1}; \cB^t_{\bx})- G^k_{\bx}(\bx^{t-1},\bz^{t}, \by^{t}; \cB^t_{\bx})  - \bar{\bg}^t_{\bx} }^2 \\
            &   = \E\norm{\bg^{t-1}_{\bx} - \bar{\bg}^{t-1}_{\bx} +  \bar{\bg}^{t-1}_{\bx} +  G^k_{\bx}(\bx^t,\bz^{t+1}, \by^{t+1}; \cB^t_{\bx})- G^k_{\bx}(\bx^{t-1},\bz^{t}, \by^{t}; \cB^t_{\bx})  -  \bar{\bg}^t_{\bx} }^2 \\
            &   = \E\norm{\bg^{t-1}_{\bx} - \bar{\bg}^{t-1}_{\bx} }^2 +  \E\norm{\bar{\bg}^{t-1}_{\bx}   - G^k_{\bx}(\bx^{t-1},\bz^{t}, \by^{t}; \cB^t_{\bx})  -  (\bar{\bg}^t_{\bx}- G^k_{\bx}(\bx^t,\bz^{t+1}, \by^{t+1}; \cB^t_{\bx}) )}^2 ,
        \end{align*}
        where the last step is due to $E [\bar{\bg}^{t-1}_{\bx}   - G^k_{\bx}(\bx^{t-1},\bz^{t}, \by^{t}; \cB^t_{\bx})  -  (\bar{\bg}^t_{\bx}- G^k_{\bx}(\bx^t,\bz^{t+1}, \by^{t+1}; \cB^t_{\bx}) )] = 0 $. We now turn to bounding the variance term.
        \begin{align*}
           & \E\norm{\bar{\bg}^{t-1}_{\bx}   - G^k_{\bx}(\bx^{t-1},\bz^{t}, \by^{t}; \cB^t_{\bx})  -  (\bar{\bg}^t_{\bx}- G^k_{\bx}(\bx^t,\bz^{t+1}, \by^{t+1}; \cB^t_{\bx}) )}^2\\
            = & \frac{1}{B^2} \sum_{i=1}^B \E\norm{\bar{\bg}^{t-1}_{\bx}   - G^k_{\bx}(\bx^{t-1},\bz^{t}, \by^{t}; \zeta_i,\xi_i)  -  (\bar{\bg}^t_{\bx}- G^k_{\bx}(\bx^t,\bz^{t+1}, \by^{t+1}; \zeta_i,\xi_i) )}^2 \\
             \leq & \frac{1}{B^2} \sum_{i=1}^B \E\norm{ G^k_{\bx}(\bx^{t-1},\bz^{t}, \by^{t}; \zeta_i,\xi_i)  -   G^k_{\bx}(\bx^t,\bz^{t+1}, \by^{t+1}; \zeta_i,\xi_i)  }^2 \\
            \leq &  \frac{1}{B^2} \sum_{i=1}^B  \norm{\nabla_1 f(\bz^{t},\by^{t};\zeta_i )   \nabla g(\bx^{t-1};\xi_i) - \nabla_1 f(\bz^{t+1},\by^{t+1};\zeta_i )   \nabla g(\bx^t;\xi_i)}^2.
        \end{align*}
        By simple variance decomposition we have:
     {\small \begin{align*}
            \norm{\nabla_1 f(\bz^{t},\by^{t};\zeta_i )   \nabla g(\bx^{t-1};\xi_i) - \nabla_1 f(\bz^{t+1},\by^{t+1};\zeta_i )   \nabla g(\bx^t;\xi_i)}^2 \leq 2G_g^2 L^2\norm{ (\bz^{t},\by^{t})  -  (\bz^{t+1},\by^{t+1} )   }^2 + 2G_f^2 L_g^2\norm{  \bx^{t-1} -  \bx^{t}    }^2.
        \end{align*}}
        Hence
       {\small \begin{align*}
            \E\norm{\bg^t_{\bx} -  \bar{\bg}^t_{\bx} }^2   &\leq \E\norm{\bg^{t-1}_{\bx} - \bar{\bg}^{t-1}_{\bx} }^2  + \frac{1}{B } \pare{2G_g^2 L^2\norm{ \bz^{t}   -  \bz^{t+1}    }^2 +2G_g^2 L^2\norm{  \by^{t}  -  \by^{t+1}  }^2 + 2G_f^2 L_g^2\norm{  \bx^{t-1} -  \bx^{t}    }^2}\\
            &\leq  \E\norm{\bg^{t_p}_{\bx} - \bar{\bg}^{t_p}_{\bx} }^2  + \frac{1}{B }\sum_{t'=t_p}^{t-1} \pare{2G_g^2 L^2\norm{ \bz^{t'+1}   -  \bz^{t'+2}    }^2 +2G_g^2 L^2\norm{  \by^{t'+1}  -  \by^{t'+2}  }^2 + 2G_f^2 L_g^2\norm{  \bx^{t'} -  \bx^{t'+1}    }^2},
        \end{align*}}
        where $t_p$ denotes the last iteration such that $t \mod \tau = 0$. From the updating rule we know $\E\norm{\bg^{t_p}_{\bx} - \bar{\bg}^{t_p}_{\bx} }^2 \leq \frac{\tilde\sigma^2}{B_{\tau}}$.
        Similarly,
        \begin{align*}
            \E\norm{\bq^{t+1} - \bar{\bg}^{t+1}_{\by}}^2  &=  \E\norm{\bq^{t} - \bar{\bg}^{t}_{\by}}^2 + \E\norm{ G_k(\bx^{t+1},\bz^{t+2}, \by^{t+1};\cB^t_{\by}) - \bar{\bg}^{t+1}_{\by}- ( G_k(\bx^{t},\bz^{t+1}, \by^{t };\cB^t_{\by}) - \bar{\bg}^{t}_{\by} )}^2\\
           &\leq \E\norm{\bq^{t} - \bar{\bg}^{t}_{\by}}^2 + \frac{1}{B^2} \sum_{i=1}^B  \E\norm{\nabla_{\by} f(\bz^{t},\by^{t};\zeta_i )  - \nabla_{\by} f(\bz^{t+1},\by^{t+1};\zeta_i ) }^2 \\
             &\leq \E\norm{\bq^{t} - \bar{\bg}^{t}_{\by}}^2  +  \frac{1}{B} L^2  \E\norm{ (\bz^{t},\by^{t}  )  -  (\bz^{t+1},\by^{t+1} ) }^2 \\
             &\leq \E\norm{\bq^{t_p} - \bar{\bg}^{t_p}_{\by}}^2  +  \frac{1}{B} \sum_{t'=t_p}^{t}  L^2_f \pare{ \E\norm{ \bz^{t'}   -  \bz^{t'+1}  }^2 + \E\norm{  \by^{t'}  -   \by^{t'+1} }^2 }\\
             &\leq \frac{\tilde\sigma^2}{B_{\tau}}  +  \frac{1}{B} \sum_{t'=t_p}^{t}  L^2_f \pare{ \E\norm{ \bz^{t'}   -  \bz^{t'+1}  }^2 + \E\norm{  \by^{t'}  -   \by^{t'+1} }^2 }.
        \end{align*}
    \end{proof}
\end{lemma}

\begin{lemma} \label{lem:descent lemma VR}
For Algorithm~\ref{algorithm: CODA-Primal+}, under assumptions of Theorem~\ref{thm:NCSC+}, the following statement holds:
         \begin{align*}
            F_k(\bx^{t+1},\by) - & F_k(\bx,\by^{t+1}) \\
            \leq &-\inprod{\nabla_{\by} f(g(\bx^{t+1}),\by^{t+1}) - \nabla_{\by} f(g(\bx^{t}),\by^{t})  }{\by^{t+1} - \by} \\
            &+  \inprod{\nabla_{\by} f(g(\bx^{t}),\by^{t}) - \nabla_{\by} f(g(\bx^{t-1}),\by^{t-1})  }{\by^{t} - \by}  \\
               &+  \inprod{   \nabla_{\by} f(g(\bx^t),\by^{t}) - \nabla_{\by} f(g(\bx^{t-1}),\by^{t-1})  }{ \by^{t+1} - \by^t } \\
               & + \frac{L_f'}{2}\norm{\bx^{t+1} - \bx^t}^2 + \frac{1}{2\eta_{\by}} \pare{ \norm{\by - \by^t}^2 - \norm{\by - \by^{t+1}}^2 - \norm{  \by^{t+1} - \by^t}^2 }  \\
               &+\frac{1}{2\eta_{\bx}} \pare{ \norm{\bx - \bx^t}^2 - \norm{\bx - \bx^{t+1}}^2 - \norm{  \bx^{t+1} - \bx^t}^2 }- \frac{\mu_x}{2} \norm{\bx- \bx^{t}}^2 +  \frac{\mu_x}{2}\norm{ \bx^{t+1} - \bx }^2\\
               & +\frac{1}{\mu_x} \norm{\bg_{\bx}^t - \bar\bg_{\bx}^t  }^2+ \pare{ \frac{L_f^2G^2_g}{\mu_x} +  \frac{8L^2_f }{\mu_y} }\norm{ \bz^{t+1}    -  g(\bx^t)   }^2  + \frac{2L^2_f  }{\mu_y} \norm{\bz^{t} - g(\bx^{t-1})}^2 \\
               & + \frac{8 }{\mu_y} \norm{\bq^{t} - \nabla_{\by} f(\bz^{t+1},\by^{t})}^2 + \frac{2  }{\mu_y} \norm{\bq^{t-1} - \nabla_{\by} f(\bz^{t},\by^{t-1})}^2, 
        \end{align*}
 where $L_f' := L_f+\mu_x+L$.
        
    \begin{proof}
    Define $\phi(\bx^t,\bz^{t+1},\by^{t+1}) = h(\bx^t) + \nabla_1 f(\bz^{t+1},\by^{t+1}) g(\bx^t) $, $f_k(\bx,\by) := f(g(\bx^t),\by^{t+1}) + \frac{\mu_x}{2}\norm{\bx - \bx_k}^2 $. According to updating rule we have:
        \begin{align*}
            \bx^{t+1} &= \bx^t - \eta_{\bx} \bg_{\bx}^t  = \arg\min_{\bx} h(\bx) + \phi(\bx) - \phi(\bx^t) - \inprod{\nabla_{\bx} \phi(\bx^t,\bz^{t+1},\by^{t+1})}{\bx- \bx^t}.
        \end{align*}
        Then according to [Lemma 7.1~\cite{hamedani2018primal}] we have:
        \begin{align*}
            h(\bx^{t+1}) + \inprod{\bg_{\bx}^t }{ \bx^{t+1} - \bx } \leq h(\bx) + \frac{1}{2\eta_{\bx}} \pare{ \norm{\bx - \bx^t}^2 - \norm{\bx - \bx^{t+1}}^2 - \norm{  \bx^{t+1} - \bx^t}^2 },\\
              r(\by^{t+1}) - \inprod{\bg_{\by}^t }{ \by^{t+1} - \by } \leq r(\by) + \frac{1}{2\eta_{\by}} \pare{ \norm{\by - \by^t}^2 - \norm{\by - \by^{t+1}}^2 - \norm{  \by^{t+1} - \by^t}^2 }.
        \end{align*}
        Re-arranging terms yields:
         \begin{align}
            h(\bx^{t+1}) + &\inprod{\nabla_{\bx} f_k(\bx^t, \by^{t+1}) }{ \bx^{t+1} - \bx } \leq h(\bx) + \frac{1}{2\eta_{\bx}} \pare{ \norm{\bx - \bx^t}^2 - \norm{\bx - \bx^{t+1}}^2 - \norm{  \bx^{t+1} - \bx^t}^2 } + \varepsilon^{t+1}_{\bx},\label{eq:descent VR 1} \\
           r(\by^{t+1}) -   & \inprod{\nabla_{\by} f(g(\bx^t),\by^{t}) + \theta (\nabla_{\by} f(g(\bx^t),\by^{t}) - \nabla_{\by} f(g(\bx^{t-1}),\by^{t-1}) ) }{ \by^{t+1} - \by }\nonumber \\
           &\quad \leq r(\by) + \frac{1}{2\eta_{\by}} \pare{ \norm{\by - \by^t}^2 - \norm{\by - \by^{t+1}}^2 - \norm{  \by^{t+1} - \by^t}^2 } + \varepsilon^{t+1}_{\by},\label{eq:descent VR 2}
        \end{align}
        where
        \begin{align*}
            &\varepsilon^{t+1}_{\bx} = \inprod{\bg_{\bx}^t - \nabla_{\bx} f(g(\bx^t),\by^{t+1})  }{ \bx^{t+1} - \bx } ,\\
            & \varepsilon^{t+1}_{\by} = \inprod{\nabla_{\by} f(g(\bx^t),\by^{t}) + \theta (\nabla_{\by} f(g(\bx^t),\by^{t}) - \nabla_{\by} f(g(\bx^{t-1}),\by^{t-1}) )-\bg_{\by}^t }{ \by^{t+1} - \by }.
        \end{align*}
        The rest of the proof follows [Lemma 13~\cite{zhang2022sapd+}]. Using $\mu_x$-convexity of $f_k$ w.r.t. $\bx$, we can lower bound the inner product in (\ref{eq:descent VR 1}) as
        \begin{align*}
            \inprod{\nabla_{\bx} f_k(\bx^t, \by^{t+1}) }{ \bx^{t+1} - \bx }  &= \inprod{\nabla_{\bx} f_k(\bx^t, \by^{t+1}) }{ \bx^{t } - \bx } + \inprod{\nabla_{\bx} f_k(\bx^t, \by^{t+1}) }{ \bx^{t+1} -\bx^{t } } \\
            &\geq f_k(\bx^t,\by^{t+1}) - f_k(\bx ,\by^{t+1}) + \frac{\mu_x}{2}\norm{\bx^t - \bx}^2+ \inprod{\nabla_{\bx} f_k(\bx^t, \by^{t+1}) }{ \bx^{t+1} -\bx^{t } }.
        \end{align*}
        Putting pieces together yields:
        \begin{align*}
           h(\bx^{t+1}) + f_k(\bx^t,\by^{t+1}) -& f_k(\bx ,\by^{t+1}) + r(\by^{t+1}) - r(\by) \\
           &\leq  - \frac{\mu_x}{2}\norm{\bx^t - \bx}^2-\inprod{\nabla_{\bx} f_k(\bx^t, \by^{t+1}) }{ \bx^{t+1} -\bx^{t } } \\
           & + \inprod{\nabla_{\by} f(g(\bx^t),\by^{t}) + \theta (\nabla_{\by} f(g(\bx^t),\by^{t}) - \nabla_{\by} f(g(\bx^{t-1}),\by^{t-1}) ) }{ \by^{t+1} - \by }\\
          & +h(\bx) + \frac{1}{2\eta_{\bx}} \pare{ \norm{\bx - \bx^t}^2 - \norm{\bx - \bx^{t+1}}^2 - \norm{  \bx^{t+1} - \bx^t}^2 } + \varepsilon^{t+1}_{\bx}\\
          & + \frac{1}{2\eta_{\by}} \pare{ \norm{\by - \by^t}^2 - \norm{\by - \by^{t+1}}^2 - \norm{  \by^{t+1} - \by^t}^2 } + \varepsilon^{t+1}_{\by}.
        \end{align*}
Adding $f_k(\bx^{t+1},\by^{t+1})$ on both sides yields:
\begin{align*}
  f_k(\bx^{t+1},\by^{t+1})+   h(\bx^{t+1})     -& f_k(\bx ,\by^{t+1}) + r(\by^{t+1}) - r(\by) -h(\bx) \\
 \leq  & f_k(\bx^{t+1},\by^{t+1}) -f_k(\bx^t,\by^{t+1}) -\inprod{\nabla_{\bx} f_k(\bx^t, \by^{t+1}) }{ \bx^{t+1} -\bx^{t } } - \frac{\mu_x}{2}\norm{\bx^t - \bx}^2 \\
  & + \inprod{\nabla_{\by} f(g(\bx^t),\by^{t}) +   (\nabla_{\by} f(g(\bx^t),\by^{t}) - \nabla_{\by} f(g(\bx^{t-1}),\by^{t-1}) ) }{ \by^{t+1} - \by }\\
          & +  \frac{1}{2\eta_{\bx}} \pare{ \norm{\bx - \bx^t}^2 - \norm{\bx - \bx^{t+1}}^2 - \norm{  \bx^{t+1} - \bx^t}^2 } + \varepsilon^{t+1}_{\bx}\\
          & + \frac{1}{2\eta_{\by}} \pare{ \norm{\by - \by^t}^2 - \norm{\by - \by^{t+1}}^2 - \norm{  \by^{t+1} - \by^t}^2 } + \varepsilon^{t+1}_{\by}.
\end{align*}
Due to $L_f' := L_f+\mu_x+L$ smoothness of $f_k$, we know $f_k(\bx^{t+1},\by^{t+1}) -f_k(\bx^t,\by^{t+1}) -\inprod{\nabla_{\bx} f_k(\bx^t, \by^{t+1}) }{ \bx^{t+1} -\bx^{t } } \leq \frac{L_f'}{2}\norm{\bx^t - \bx}^2  $. Hence we have

\begin{align*}
  f_k(\bx^{t+1},\by^{t+1})+   h(\bx^{t+1})     -& f_k(\bx ,\by^{t+1}) + r(\by^{t+1}) - r(\by) -h(\bx) \\
 \leq  & \frac{L_f'}{2}\norm{\bx^t - \bx}^2  - \frac{\mu_x}{2}\norm{\bx^t - \bx}^2 \\
  & + \inprod{\nabla_{\by} f(g(\bx^t),\by^{t}) +   (\nabla_{\by} f(g(\bx^t),\by^{t}) - \nabla_{\by} f(g(\bx^{t-1}),\by^{t-1}) ) }{ \by^{t+1} - \by }\\
          & +  \frac{1}{2\eta_{\bx}} \pare{ \norm{\bx - \bx^t}^2 - \norm{\bx - \bx^{t+1}}^2 - \norm{  \bx^{t+1} - \bx^t}^2 } + \varepsilon^{t+1}_{\bx}\\
          & + \frac{1}{2\eta_{\by}} \pare{ \norm{\by - \by^t}^2 - \norm{\by - \by^{t+1}}^2 - \norm{  \by^{t+1} - \by^t}^2 } + \varepsilon^{t+1}_{\by}.
\end{align*}
Adding $f_k(\bx^{t+1},\by)$ on both sides yields:
\begin{align*}
  &\underbrace{f_k(\bx^{t+1},\by)+   h(\bx^{t+1}) - r(\by)    -  (h(\bx)+f_k(\bx ,\by^{t+1}) - r(\by^{t+1})  )}_{=F_k(\bx^{t+1},\by) - F_k(\bx,\by^{t+1})} \\
 \leq  & f_k(\bx^{t+1},\by) -  f_k(\bx^{t+1},\by^{t+1})+ \frac{L_f'}{2}\norm{\bx^t - \bx}^2  - \frac{\mu_x}{2}\norm{\bx^t - \bx}^2 \\
  & + \inprod{\nabla_{\by} f(g(\bx^t),\by^{t}) +   (\nabla_{\by} f(g(\bx^t),\by^{t}) - \nabla_{\by} f(g(\bx^{t-1}),\by^{t-1}) ) }{ \by^{t+1} - \by }\\
          & +  \frac{1}{2\eta_{\bx}} \pare{ \norm{\bx - \bx^t}^2 - \norm{\bx - \bx^{t+1}}^2 - \norm{  \bx^{t+1} - \bx^t}^2 } + \varepsilon^{t+1}_{\bx}\\
          & + \frac{1}{2\eta_{\by}} \pare{ \norm{\by - \by^t}^2 - \norm{\by - \by^{t+1}}^2 - \norm{  \by^{t+1} - \by^t}^2 } + \varepsilon^{t+1}_{\by}.
\end{align*}
Due to $\mu_y$ strong concavity of $f_k$ w.r.t. $\by$, we have $f_k(\bx^{t+1},\by) -  f_k(\bx^{t+1},\by^{t+1}) \leq \inprod{ \nabla_{\by} f_k(\bx^{t+1},\by^{t+1})   }{\by - \by^{t+1}} - \frac{\mu_y}{2}\norm{ \by - \by^{t+1}}^2$. So we have
\begin{align*}
    F_k(\bx^{t+1},\by) - F_k(\bx,\by^{t+1})  
 \leq  & \inprod{ \nabla_{\by} f (g(\bx^{t+1}),\by^{t+1})   }{\by - \by^{t+1}} - \frac{\mu_y}{2}\norm{ \by - \by^{t+1}}^2+ \frac{L_f'}{2}\norm{\bx^t - \bx}^2  - \frac{\mu_x}{2}\norm{\bx^t - \bx}^2 \\
  & + \inprod{\nabla_{\by} f(g(\bx^t),\by^{t}) +   (\nabla_{\by} f(g(\bx^t),\by^{t}) - \nabla_{\by} f(g(\bx^{t-1}),\by^{t-1}) ) }{ \by^{t+1} - \by }\\
          & +  \frac{1}{2\eta_{\bx}} \pare{ \norm{\bx - \bx^t}^2 - \norm{\bx - \bx^{t+1}}^2 - \norm{  \bx^{t+1} - \bx^t}^2 } + \varepsilon^{t+1}_{\bx}\\
          & + \frac{1}{2\eta_{\by}} \pare{ \norm{\by - \by^t}^2 - \norm{\by - \by^{t+1}}^2 - \norm{  \by^{t+1} - \by^t}^2 } + \varepsilon^{t+1}_{\by}. 
\end{align*}
Notice the following identity:
\begin{align*}
    &\inprod{ \nabla_{\by} f (g(\bx^{t+1}),\by^{t+1})   }{\by - \by^{t+1}} +  \inprod{\nabla_{\by} f(g(\bx^t),\by^{t}) +   (\nabla_{\by} f(g(\bx^t),\by^{t}) - \nabla_{\by} f(g(\bx^{t-1}),\by^{t-1}) ) }{ \by^{t+1} - \by } \\
    =  &\inprod{\nabla_{\by} f(g(\bx^t),\by^{t}) -\nabla_{\by} f (g(\bx^{t+1}),\by^{t+1})    +   (\nabla_{\by} f(g(\bx^t),\by^{t}) - \nabla_{\by} f(g(\bx^{t-1}),\by^{t-1}) ) }{ \by^{t+1} - \by } \\
    = & - \inprod{\nabla_{\by} f (g(\bx^{t+1}),\by^{t+1})  -\nabla_{\by} f(g(\bx^t),\by^{t})   }{ \by^{t+1} - \by } + \inprod{   \nabla_{\by} f(g(\bx^t),\by^{t}) - \nabla_{\by} f(g(\bx^{t-1}),\by^{t-1})  }{ \by^{t } - \by } \\
    & +   \inprod{   \nabla_{\by} f(g(\bx^t),\by^{t}) - \nabla_{\by} f(g(\bx^{t-1}),\by^{t-1})  }{ \by^{t+1} - \by^t }.
\end{align*}
We thus arrive at
\begin{align*}
            F_k(\bx^{t+1},\by) & - F_k(\bx,\by^{t+1}) \\
            &\leq -\inprod{\nabla_{\by} f(g(\bx^{t+1}),\by^{t+1}) - \nabla_{\by} f(g(\bx^{t}),\by^{t})  }{\by^{t+1} - \by}  +  \inprod{\nabla_{\by} f(g(\bx^{t}),\by^{t}) - \nabla_{\by} f(g(\bx^{t-1}),\by^{t-1})  }{\by^{t} - \by}  \\
               &+  \inprod{   \nabla_{\by} f(g(\bx^t),\by^{t}) - \nabla_{\by} f(g(\bx^{t-1}),\by^{t-1})  }{ \by^{t+1} - \by^t } \\
               & + \frac{L_f'}{2}\norm{\bx^{t+1} - \bx^t}^2 + \frac{1}{2\eta_{\by}} \pare{ \norm{\by - \by^t}^2 - \norm{\by - \by^{t+1}}^2 - \norm{  \by^{t+1} - \by^t}^2 } - \frac{\mu_y}{2} \norm{\by- \by^{t+1}}^2\\
               &+\frac{1}{2\eta_{\bx}} \pare{ \norm{\bx - \bx^t}^2 - \norm{\bx - \bx^{t+1}}^2 - \norm{  \bx^{t+1} - \bx^t}^2 }- \frac{\mu_x}{2} \norm{\bx- \bx^{t}}^2 + \varepsilon^{t+1}_{\bx}+ \varepsilon^{t+1}_{\by}.
        \end{align*}
 
        Now it remains to bound $\varepsilon^{t+1}_{\bx}$ and $ \varepsilon^{t+1}_{\by}$. For $\varepsilon^{t+1}_{\bx}$:
        \begin{align*}
            \varepsilon^{t+1}_{\bx} \leq \frac{1}{2} \pare{ \frac{1}{\mu_x}\norm{\bg_{\bx}^t - \nabla_{\bx} f(g(\bx^t),\by^{t+1})  }^2 + \mu_x\norm{ \bx^{t+1} - \bx }^2  }.
        \end{align*}
        Notice that 
        \begin{align*}
            \norm{\bg_{\bx}^t - \nabla_{\bx} f(g(\bx^t),\by^{t+1})  }^2 & \leq 2\norm{\bg_{\bx}^t - \bar\bg_{\bx}^t  }^2+2\norm{\bar\bg_{\bx}^t - \nabla_{\bx} f(g(\bx^t),\by^{t+1})  }^2\\
            &\leq 2\norm{\bg_{\bx}^t - \bar\bg_{\bx}^t  }^2+2\norm{\nabla f(\bz^{t+1},\by^{t+1}) \nabla g(\bx^t) - \nabla  f(g(\bx^t),\by^{t+1}) \nabla g(\bx^t)  }^2\\
            &\leq 2\norm{\bg_{\bx}^t - \bar\bg_{\bx}^t  }^2+2L_f^2G^2_g\norm{ \bz^{t+1}    -  g(\bx^t)   }^2.
        \end{align*}
        Hence we have
        \begin{align*}
             \varepsilon^{t+1}_{\bx} \leq   \frac{1}{\mu_x} \norm{\bg_{\bx}^t - \bar\bg_{\bx}^t  }^2+ \frac{L_f^2G^2_g}{\mu_x}\norm{ \bz^{t+1}    -  g(\bx^t)   }^2  +  \frac{\mu_x}{2}\norm{ \bx^{t+1} - \bx }^2 .
        \end{align*}
        Similarly we have for $\varepsilon^{t+1}_{\by}$
        \begin{align*}
             \varepsilon^{t+1}_{\by} &\leq \frac{1}{2} \pare{ \frac{1}{\mu_y}\norm{2\bq^{t} -   \bq^{t-1} - \nabla_{\by} f(g(\bx^t),\by^{t}) -   (\nabla_{\by} f(g(\bx^t),\by^{t}) - \nabla_{\by} f(g(\bx^{t-1}),\by^{t-1}) )  }^2 + \mu_y\norm{ \by^{t+1} - \by }^2  } \\
             &\leq \frac{\mu_y}{2}\norm{ \by^{t+1} - \by }^2 + \frac{4 }{\mu_y} \norm{\bq^{t} - \nabla_{\by} f(g(\bx^t),\by^{t})}^2 + \frac{1 }{\mu_y} \norm{\bq^{t-1} - \nabla_{\by} f(g(\bx^{t-1}),\by^{t-1})}^2\\
             &\leq \frac{\mu_y}{2}\norm{ \by^{t+1} - \by }^2 + \frac{8 }{\mu_y} \norm{\bq^{t} - \nabla_{\by} f(\bz^{t+1},\by^{t})}^2 + \frac{2  }{\mu_y} \norm{\bq^{t-1} - \nabla_{\by} f(\bz^{t},\by^{t-1})}^2\\
             & \quad + \frac{8 }{\mu_y}L^2_f \norm{ \bz^{t+1} - g(\bx^t)}^2 + \frac{2  }{\mu_y}L^2_f \norm{\bz^{t} - g(\bx^{t-1})}^2.
        \end{align*}
        Putting pieces together will conclude the proof.
        \begin{align*}
            F_k(\bx^{t+1},\by) - & F_k(\bx,\by^{t+1}) \\
            \leq &-\inprod{\nabla_{\by} f(g(\bx^{t+1}),\by^{t+1}) - \nabla_{\by} f(g(\bx^{t}),\by^{t})  }{\by^{t+1} - \by} \\
            &+  \inprod{\nabla_{\by} f(g(\bx^{t}),\by^{t}) - \nabla_{\by} f(g(\bx^{t-1}),\by^{t-1})  }{\by^{t} - \by}  \\
               &+  \inprod{   \nabla_{\by} f(g(\bx^t),\by^{t}) - \nabla_{\by} f(g(\bx^{t-1}),\by^{t-1})  }{ \by^{t+1} - \by^t } \\
               & + \frac{L_f'}{2}\norm{\bx^{t+1} - \bx^t}^2 + \frac{1}{2\eta_{\by}} \pare{ \norm{\by - \by^t}^2 - \norm{\by - \by^{t+1}}^2 - \norm{  \by^{t+1} - \by^t}^2 }  \\
               &+\frac{1}{2\eta_{\bx}} \pare{ \norm{\bx - \bx^t}^2 - \norm{\bx - \bx^{t+1}}^2 - \norm{  \bx^{t+1} - \bx^t}^2 }- \frac{\mu_x}{2} \norm{\bx- \bx^{t}}^2 +  \frac{\mu_x}{2}\norm{ \bx^{t+1} - \bx }^2\\
               & +\frac{1}{\mu_x} \norm{\bg_{\bx}^t - \bar\bg_{\bx}^t  }^2+ \pare{ \frac{L_f^2G^2_g}{\mu_x} +  \frac{8L^2_f }{\mu_y} }\norm{ \bz^{t+1}    -  g(\bx^t)   }^2  + \frac{2L^2_f  }{\mu_y} \norm{\bz^{t} - g(\bx^{t-1})}^2 \\
               & + \frac{8 }{\mu_y} \norm{\bq^{t} - \nabla_{\by} f(\bz^{t+1},\by^{t})}^2 + \frac{2  }{\mu_y} \norm{\bq^{t-1} - \nabla_{\by} f(\bz^{t},\by^{t-1})}^2. 
        \end{align*}
    \end{proof}
\end{lemma}

\begin{lemma}[Primal-dual gap convergence]\label{lem:PD conevrgence VR}For Algorithm~\ref{algorithm: CODA-Primal+}, under assumptions of Theorem~\ref{thm:NCSC+}, the following statement holds:
    \begin{align*}
         \E[ F_k(\bx_{k+1},\by^*(\bx_{k+1})) - F_k(\bx^*(\by_{k+1}),\by_{k+1})]  
                \leq   & \frac{ \E\norm{\by^*_k(\bx_{k+1}) - \by_k}^2}{2\eta_{\by} T}  + \frac{ \E \norm{\bx^*_k(\by_{k+1}) - \bx_k}^2 }{4\eta_{\bx}  T}   +  \frac{C\E\norm{\bz^{0}- g(\bx^{0})}^2}{T}   \\  
                & +  \frac{\tilde\sigma^2}{\mu_x B_{\tau}}+  \frac{10\tilde\sigma^2}{\mu_y B_{\tau}} + \pare{\frac{24G_g^2 L^2_f \tau }{\mu_x B } +  \frac{120 \tau}{\mu_y B} + \frac{20L^2_f }{\mu_y} + \frac{2L_f^2G^2_g}{\mu_x}    }  9 \frac{\tilde\sigma^2}{B_\tau},
    \end{align*}
    where $C = \frac{12G_g^2 L^2 \tau }{\mu_x B } +  \frac{60 \tau}{\mu_y B} + \frac{12L^2_f }{\mu_y} + \frac{L_f^2G^2_g}{\mu_x} $, $\bx^*_k(\by ) = \arg\min_{\bx\in\cX} F_k(\bx,\by) $ and $\by^*_k(\bx ) = \arg\max_{\by\in\cY} F_k(\bx,\by) $.
    \begin{proof}

    Due to Lemma~\ref{lem:descent lemma VR} we have
  \begin{align*}
            F_k(\bx^{t+1},\by) - & F_k(\bx,\by^{t+1}) \\
            \leq &-\inprod{\nabla_{\by} f(g(\bx^{t+1}),\by^{t+1}) - \nabla_{\by} f(g(\bx^{t}),\by^{t})  }{\by^{t+1} - \by}\\
            &+  \inprod{\nabla_{\by} f(g(\bx^{t}),\by^{t}) - \nabla_{\by} f(g(\bx^{t-1}),\by^{t-1})  }{\by^{t} - \by}  \\
               &+  \inprod{   \nabla_{\by} f(g(\bx^t),\by^{t}) - \nabla_{\by} f(g(\bx^{t-1}),\by^{t-1})  }{ \by^{t+1} - \by^t } \\
               & + \frac{L_f'}{2}\norm{\bx^{t+1} - \bx^t}^2 + \frac{1}{2\eta_{\by}} \pare{ \norm{\by - \by^t}^2 - \norm{\by - \by^{t+1}}^2 - \norm{  \by^{t+1} - \by^t}^2 }  \\
               &+\frac{1}{2\eta_{\bx}} \pare{ \norm{\bx - \bx^t}^2 - \norm{\bx - \bx^{t+1}}^2 - \norm{  \bx^{t+1} - \bx^t}^2 }- \frac{\mu_x}{2} \norm{\bx- \bx^{t}}^2 +  \frac{\mu_x}{2}\norm{ \bx^{t+1} - \bx }^2\\
               & +\frac{1}{\mu_x} \norm{\bg_{\bx}^t - \bar\bg_{\bx}^t  }^2+ \pare{ \frac{L_f^2G^2_g}{\mu_x} +  \frac{8L^2_f }{\mu_y} }\norm{ \bz^{t+1}    -  g(\bx^t)   }^2  + \frac{2L^2_f  }{\mu_y} \norm{\bz^{t} - g(\bx^{t-1})}^2 \\
               & + \frac{8 }{\mu_y} \norm{\bq^{t} - \nabla_{\by} f(\bz^{t+1},\by^{t})}^2 + \frac{2  }{\mu_y} \norm{\bq^{t-1} - \nabla_{\by} f(\bz^{t},\by^{t-1})}^2. 
        \end{align*}
    Define $\cF_k^{t+1} = F_k(\bx^{t+1},\by) -   F_k(\bx,\by^{t+1})  + C\pare{\E\norm{ \bz^{t+1} - g(\bx^t)}^2 - \E\norm{\bz^{t} - g(\bx^{t-1})}^2}$, and $\Gamma^t =  \inprod{\nabla_{\by} f(g(\bx^{t}),\by^{t}) - \nabla_{\by} f(g(\bx^{t-1}),\by^{t-1})  }{\by^{t} - \by} $.  Adding $C\pare{\E\norm{ \bz^{t+1} - g(\bx^t)}^2 - \E\norm{\bz^{t} - g(\bx^{t-1})}^2}$  on both side yields:
    \begin{align*}
         \cF_k^{t+1} \leq -&\Gamma^{t+1} + \Gamma^t  +  \inprod{   \nabla_{\by} f(g(\bx^t),\by^{t}) - \nabla_{\by} f(g(\bx^{t-1}),\by^{t-1})  }{ \by^{t+1} - \by^t }\\
                & + \frac{L_f'}{2}\norm{\bx^{t+1} - \bx^t}^2 + \frac{1}{2\eta_{\by}} \pare{ \norm{\by - \by^t}^2 - \norm{\by - \by^{t+1}}^2 - \norm{  \by^{t+1} - \by^t}^2 }  \\
               &+\frac{1}{2\eta_{\bx}} \pare{ \norm{\bx - \bx^t}^2 - \norm{\bx - \bx^{t+1}}^2 - \norm{  \bx^{t+1} - \bx^t}^2 }- \frac{\mu_x}{2} \norm{\bx- \bx^{t}}^2 +  \frac{\mu_x}{2}\norm{ \bx^{t+1} - \bx }^2\\
               & +\frac{1}{\mu_x} \norm{\bg_{\bx}^t - \bar\bg_{\bx}^t  }^2+ \pare{ \frac{L_f^2G^2_g}{\mu_x} +  \frac{8L^2_f }{\mu_y} }\norm{ \bz^{t+1}    -  g(\bx^t)   }^2  + \frac{2L^2_f  }{\mu_y} \norm{\bz^{t} - g(\bx^{t-1})}^2 \\
               & + \frac{8 }{\mu_y} \norm{\bq^{t} - \nabla_{\by} f(\bz^{t+1},\by^{t})}^2 + \frac{2  }{\mu_y} \norm{\bq^{t-1} - \nabla_{\by} f(\bz^{t},\by^{t-1})}^2 \\ 
             &+C\pare{\E\norm{ \bz^{t+1} - g(\bx^t)}^2 - \E\norm{\bz^{t} - g(\bx^{t-1})}^2}.
    \end{align*}
    
    Notice that by Cauchy-schwartz and AM-GM inequality we have:
   {\small \begin{align*}
        \inprod{   \nabla_{\by} f(g(\bx^t),\by^{t}) - \nabla_{\by} f(g(\bx^{t-1}),\by^{t-1})  }{ \by^{t+1} - \by^t } &\leq \eta_{\by}\norm{ \nabla_{\by} f(g(\bx^t),\by^{t}) - \nabla_{\by} f(g(\bx^{t-1}),\by^{t-1})  }^2 + \frac{1}{4\eta_{\by}} \norm{\by^{t+1} - \by^t}^2  \\
        & \leq \eta_{\by}L_f^2 \norm{ (g(\bx^t),\by^{t}) -  (g(\bx^{t-1}),\by^{t-1})  }^2 + \frac{1}{4\eta_{\by}} \norm{\by^{t+1} - \by^t}^2  \\
         &\leq \eta_{\by}L_f^2 \pare{G_g^2\norm{  \bx^{t-1}  -  \bx^{t}   }^2 + \norm{\by^{t-1} - \by^t}^2} + \frac{1}{4\eta_{\by}} \norm{\by^{t+1} - \by^t}^2.
    \end{align*}}
    Plugging back yields:
         \begin{align*}
         \cF_k^{t+1}  
            \leq -&\Gamma^{t+1} + \Gamma^t   + \frac{1}{2\eta_{\by}} \pare{ \norm{\by - \by^t}^2 - \norm{\by - \by^{t+1}}^2 }  - \frac{1}{4\eta_{\by}} \norm{  \by^{t+1} - \by^t}^2 + \eta_{\by} L_f^2 \norm{\by^{t-1} - \by^t}^2 \\
               &+\pare{\frac{1}{2\eta_{\bx}} - \frac{\mu_x}{2}} \pare{ \norm{\bx - \bx^t}^2 - \norm{\bx - \bx^{t+1}}^2 }  -\pare{\frac{1}{2\eta_{\bx}} - \frac{L'_f}{2} } \norm{  \bx^{t+1} - \bx^t}^2 + \eta_{\by} L_f^2 G_g^2 \norm{\bx^{t-1} - \bx^t}^2 \\
               & + \frac{1}{\mu_x} \norm{\bg_{\bx}^t - \bar\bg_{\bx}^t }^2   + \frac{8 }{\mu_y} \norm{\bq^{t} -\bar\bg_{\by}^t }^2 + \frac{2  }{\mu_y} \norm{\bq^{t-1} - \bar\bg_{\by}^{t-1} }^2\\
             &  +\pare{ \frac{8L^2_f }{\mu_y} + \frac{L_f^2G^2_g}{\mu_x} }\norm{ \bz^{t+1} - g(\bx^t)}^2 + \frac{2L^2_f  }{\mu_y} \norm{\bz^{t} - g(\bx^{t-1})}^2\\
             &+C\pare{\E\norm{ \bz^{t+1} - g(\bx^t)}^2 - \E\norm{\bz^{t} - g(\bx^{t-1})}^2}.
    \end{align*}
 Plugging in the bound from Lemma~\ref{lem:grad bounds VR} for $\norm{\bg_{\bx}^t - \bar\bg_{\bx}^t }^2 $ and $\norm{\bq^{t} -\bar\bg_{\by}^t }^2$ yields:
    \begin{align*}
            \cF_k^{t+1}  
             \leq &-\Gamma^{t+1} +\Gamma^{t}   + \frac{1}{2\eta_{\by}} \pare{ \norm{\by - \by^t}^2 - \norm{\by - \by^{t+1}}^2 }   - \frac{1}{4\eta_{\by}} \norm{  \by^{t+1} - \by^t}^2 + \eta_{\by} L_f^2 \norm{\by^{t-1} - \by^t}^2 \\
               &+\pare{\frac{1}{2\eta_{\bx}} - \frac{\mu_x}{2}} \pare{ \norm{\bx - \bx^t}^2 - \norm{\bx - \bx^{t+1}}^2 } -\pare{\frac{1}{2\eta_{\bx}} - \frac{L'_f}{2} } \norm{  \bx^{t+1} - \bx^t}^2 + \eta_{\by} L_f^2 G_g^2 \norm{\bx^{t-1} - \bx^t}^2 \\
               & +  \frac{\tilde\sigma^2}{\mu_x B_{\tau}}  + \frac{1}{\mu_x B }\sum_{t'=t_p}^{t-1} \pare{2G_g^2 L^2_f\norm{ \bz^{t'+1}   -  \bz^{t'+2}    }^2 +2G_g^2 L^2_f\norm{  \by^{t'+1}  -  \by^{t'+2}  }^2 + 2G_f^2 L_g^2\norm{  \bx^{t'} -  \bx^{t'+1}    }^2} \\
               & +  \frac{10\tilde\sigma^2}{\mu_y B_{\tau}}  +  \frac{8}{\mu_y B} \sum_{t'=t_p}^{t-1}  L^2_f \pare{ \E\norm{ \bz^{t'}   -  \bz^{t'+1}  }^2 + \E\norm{  \by^{t'}  -   \by^{t'+1} }^2 } \\
               &+  \frac{2}{\mu_y B} \sum_{t'=t_p}^{t-2}  L^2_f \pare{ \E\norm{ \bz^{t'}   -  \bz^{t'+1}  }^2 + \E\norm{  \by^{t'}  -   \by^{t'+1} }^2 }\\
               &+ \pare{ \frac{8L^2_f }{\mu_y} + \frac{L_f^2G^2_g}{\mu_x} }\norm{ \bz^{t+1} - g(\bx^t)}^2 + \frac{2 L^2_f  }{\mu_y}  \norm{\bz^{t} - g(\bx^{t-1})}^2 +C\pare{\E\norm{ \bz^{t+1} - g(\bx^t)}^2 - \E\norm{\bz^{t} - g(\bx^{t-1})}^2}.
    \end{align*}
    Since $\sum_{t'=t_p}^{t-2}  L^2 \pare{ \E\norm{ \bz^{t'}   -  \bz^{t'+1}  }^2 + \E\norm{  \by^{t'}  -   \by^{t'+1} }^2 } \leq \sum_{t'=t_p}^{t-1}  L^2 \pare{ \E\norm{ \bz^{t'}   -  \bz^{t'+1}  }^2 + \E\norm{  \by^{t'}  -   \by^{t'+1} }^2 }$, and $\norm{\bz^t - \bz^{t+1}}^2 \leq 3\norm{\bz^t - g(\bx^{t-1})}^2 + 3 G_g^2\norm{\bx^{t-1} - \bx^t  }^2 + 3\norm{g(\bx^t) - \bz^{t+1}}^2$, we have
    {\small  \begin{align*}
            \cF_k^{t+1}  
             \leq -&\Gamma^{t+1} +\Gamma^{t}   + \frac{1}{2\eta_{\by}} \pare{ \norm{\by - \by^t}^2 - \norm{\by - \by^{t+1}}^2 }   - \frac{1}{4\eta_{\by}} \norm{  \by^{t+1} - \by^t}^2 + \eta_{\by} L_f^2 \norm{\by^{t-1} - \by^t}^2 \\
               &+\pare{\frac{1}{2\eta_{\bx}} - \frac{\mu_x}{2}} \pare{ \norm{\bx - \bx^t}^2 - \norm{\bx - \bx^{t+1}}^2 } -\pare{\frac{1}{2\eta_{\bx}} - \frac{L'_f}{2} } \norm{  \bx^{t+1} - \bx^t}^2 + \eta_{\by} L_f^2 G_g^2 \norm{\bx^{t-1} - \bx^t}^2 \\
               & +  \frac{\tilde\sigma^2}{\mu_x B_{\tau}}  + \frac{2G_g^2 L^2_f}{\mu_x B }\sum_{t'=t_p}^{t-1}   \pare{3\norm{\bz^{t'+1} - g(\bx^{t'})}^2 + 3 G_g^2\norm{\bx^{t'} - \bx^{t'+1}  }^2 + 3\norm{g(\bx^{t'+1}) - \bz^{t'+2}}^2} \\
               & +\frac{1}{\mu_x B }\sum_{t'=t_p}^{t-1} \pare{ 2G_g^2 L^2_f\norm{  \by^{t'+1}  -  \by^{t'+2}  }^2 + 2G_f^2 L_g^2\norm{  \bx^{t'} -  \bx^{t'+1}    }^2 }\\
               & +  \frac{10\tilde\sigma^2}{\mu_y B_{\tau}}  +  \frac{10}{\mu_y B} \sum_{t'=t_p}^{t-1}  L^2_f \pare{ 3\norm{\bz^{t'} - g(\bx^{t'-1})}^2 + 3 G_g^2\norm{\bx^{t'-1} - \bx^{t'}  }^2 + 3\norm{g(\bx^{t'}) - \bz^{t'+1}}^2 + \E\norm{  \by^{t'}  -   \by^{t'+1} }^2 } \\ 
               &+ \pare{ \frac{8L^2_f }{\mu_y} + \frac{L_f^2G^2}{\mu_x} }\norm{ \bz^{t+1} - g(\bx^t)}^2 + \frac{2L^2_f  }{\mu_y} \norm{\bz^{t} - g(\bx^{t-1})}^2 +C\pare{\E\norm{ \bz^{t+1} - g(\bx^t)}^2 - \E\norm{\bz^{t} - g(\bx^{t-1})}^2}.
    \end{align*}   }

Summing from $t=0$ to $T-1$ yields
{\small\begin{align*}
        \frac{1}{T}\sum_{t=0}^{T-1}\cF_k^{t+1}  
             &\leq \Gamma^{0} + \frac{1}{2\eta_{\by} T} \pare{ \norm{\by - \by^0}^2 - \norm{\by - \by^{T}}^2 } +\pare{\frac{1}{2\eta_{\bx}  T} - \frac{\mu_x}{2T}} \pare{ \norm{\bx - \bx^0}^2 - \norm{\bx - \bx^{T}}^2 } \\ 
                 &   -  \frac{1}{4\eta_{\by}}  \frac{1}{T}\sum_{t=0}^{T-1}\norm{  \by^{t+1} - \by^t}^2 +  \frac{1}{T}\sum_{t=0}^{T-1} \eta_{\by} L_f^2 \norm{\by^{t-1} - \by^t}^2\\
                 &-\pare{\frac{1}{2\eta_{\bx}} - \frac{L'_f}{2}  } \frac{1}{T}\sum_{t=0}^{T-1} \norm{  \bx^{t+1} - \bx^t}^2  +\eta_{\by} L_f^2 G_g^2\frac{1}{T}\sum_{t=0}^{T-1}  \norm{\bx^{t-1} - \bx^t}^2   +  \frac{\tilde\sigma^2}{\mu_x B_{\tau}}+  \frac{10\tilde\sigma^2}{\mu_y B_{\tau}} \\
              &  + \frac{2G_g^2 L^2}{\mu_x B }\frac{1}{T}\sum_{t=0}^{T-1}\sum_{t'=c(t)}^{t-1}   \pare{3\norm{\bz^{t'+1} - g(\bx^{t'})}^2 + 3 G_g^2\norm{\bx^{t'} - \bx^{t'+1}  }^2 + 3\norm{g(\bx^{t'+1}) - \bz^{t'+2}}^2} \\
               & +\frac{1}{\mu_x B }\frac{1}{T}\sum_{t=0}^{T-1}\sum_{t'=c(t)}^{t-1} \pare{ 2G_g^2 L^2\norm{  \by^{t'+1}  -  \by^{t'+2}  }^2 + 2G_f^2 L_g^2\norm{  \bx^{t'} -  \bx^{t'+1}    }^2 }\\
               &  +  \frac{10}{\mu_y B}  \frac{1}{T}\sum_{t=0}^{T-1}\sum_{t'=c(t)}^{t-1}  L^2_f \pare{ 3\norm{\bz^{t'} - g(\bx^{t'-1})}^2 + 3 G_g^2\norm{\bx^{t'-1} - \bx^{t'}  }^2 + 3\norm{g(\bx^{t'}) - \bz^{t'+1}}^2 + \E\norm{  \by^{t'}  -   \by^{t'+1} }^2 } \\ 
               &+ \pare{ \frac{8L^2_f }{\mu_y} + \frac{L_f^2G^2_g}{\mu_x} } \frac{1}{T}\sum_{t=0}^{T-1}\norm{ \bz^{t+1} - g(\bx^t)}^2 + \frac{2L^2_f  }{\mu_y}  \frac{1}{T}\sum_{t=0}^{T-1}\norm{\bz^{t} - g(\bx^{t-1})}^2\\
               &+C \frac{1}{T}\sum_{t=0}^{T-1} \pare{\E\norm{ \bz^{t+1} - g(\bx^t)}^2 - \E\norm{\bz^{t} - g(\bx^{t-1})}^2} .
    \end{align*}}
    Re-arranging terms yields:
    \begin{align*} 
          \frac{1}{T}\sum_{t=0}^{T-1}\cF_k^{t+1}          \leq &\Gamma^{0} + \frac{1}{2\eta_{\by} T} \pare{ \norm{\by - \by^0}^2 - \norm{\by - \by^{T}}^2 } +\pare{\frac{1}{2\eta_{\bx}  T} - \frac{\mu_x}{2T}} \pare{ \norm{\bx - \bx^0}^2 - \norm{\bx - \bx^{T}}^2 }+  \frac{\tilde\sigma^2}{\mu_x B_{\tau}}+  \frac{10\tilde\sigma^2}{\mu_y B_{\tau}} \\ 
                 &   - \pare{\frac{1}{4\eta_{\by}} -\eta_{\by}L_f^2 - \frac{2G_g^2 L^2_f\tau}{\mu_x B} -\frac{10\tau L^2}{\mu_y B}  } \frac{1}{T}\sum_{t=0}^{T-1}\norm{  \by^{t+1} - \by^t}^2 \\
                 &-\pare{\frac{1}{2\eta_{\bx}} - \frac{L'_f}{2} - \eta_{\by}L_f^2 G_g^2 -  \frac{6G_g^4 L^2_f \tau}{\mu_x B }   - \frac{2G_f^2 L_g^2
 \tau}{\mu_x B} - \frac{30G_g^2 L^2_f \tau}{\mu_y B} } \frac{1}{T}\sum_{t=0}^{T-1} \norm{  \bx^{t+1} - \bx^t}^2 \\
              &  + \pare{\frac{12G_g^2 L^2_f \tau }{\mu_x B } +  \frac{60 \tau}{\mu_y B} + \frac{8L^2_f }{\mu_y} + \frac{L_f^2G^2_g}{\mu_x} + C  }\frac{1}{T}\sum_{t=0}^{T-1}    \norm{\bz^{t+1} - g(\bx^{t})}^2  \\
              &- \pare{C-\frac{2 L^2_f }{\mu_y} }    \frac{1}{T}\sum_{t=0}^{T-1}   \E\norm{\bz^{t} - g(\bx^{t-1})}^2 .  
    \end{align*}
    Plugging the bound for $\E\norm{ \bz^{t+1} - g(\bx^t)}^2$ yields:
    \begin{align*}
        &\pare{\frac{12G_g^2 L^2 \tau }{\mu_x B } +  \frac{60 \tau}{\mu_y B} + \frac{8L^2_f }{\mu_y} + \frac{L_f^2G^2_g}{\mu_x} + C    }\frac{1}{T}\sum_{t=0}^{T-1}    \norm{\bz^{t+1} - g(\bx^{t})}^2    - \pare{C-\frac{2 L^2 }{\mu_y} }     \frac{1}{T}\sum_{t=0}^{T-1}   \E\norm{\bz^{t} - g(\bx^{t-1})}^2 \\
        \leq &\pare{\frac{12G_g^2 L^2 \tau }{\mu_x B } +  \frac{60 \tau}{\mu_y B} + \frac{8L^2_f }{\mu_y} + \frac{L_f^2G^2_g}{\mu_x} + C    }\frac{1}{T}\sum_{t=0}^{T-1}  \\
        & \times \pare{ \pare{1-\beta} \E\norm{\bz^{t} - g(\bx^{t-1})}^2 + \pare{1+\frac{1}{\beta}} \pare{6(1-\beta)^2G_g^2\E\norm{  \bx^{t-1}  -\bx^t}^2 + 3\beta^2  \frac{\sigma^2}{B_{\tau}}  + \frac{G_g^2}{B}\sum_{t'=t_p}^{t-1}  \E\norm{ \bx^{t'}  - \bx^{t'-1}   }^2 }}  \\
        &- \pare{C-\frac{2 L^2 }{\mu_y} }     \frac{1}{T}\sum_{t=0}^{T-1}   \E\norm{\bz^{t} - g(\bx^{t-1})}^2  \\
        \leq &\pare{\frac{12G_g^2 L^2 \tau }{\mu_x B } +  \frac{60 \tau}{\mu_y B} + \frac{8L^2_f }{\mu_y} + \frac{L_f^2G^2_g}{\mu_x} + C    }\frac{1}{T}\sum_{t=0}^{T-1} \pare{1-\beta} \E\norm{\bz^{t} - g(\bx^{t-1})}^2  \\
        & \pare{\frac{12G_g^2 L^2 \tau }{\mu_x B } +  \frac{60 \tau}{\mu_y B} + \frac{8L^2_f }{\mu_y} + \frac{L_f^2G^2_g}{\mu_x} + C    }  \pare{1+\frac{1}{\beta}} \pare{6 G_g^2 +   \frac{\tau G_g^2}{B}  }  \frac{1}{T}\sum_{t=0}^{T-1}  \E\norm{ \bx^{t}  - \bx^{t-1}   }^2  \\
        & + \pare{\frac{12G_g^2 L^2 \tau }{\mu_x B } +  \frac{60 \tau}{\mu_y B} + \frac{8L^2_f }{\mu_y} + \frac{L_f^2G^2_g}{\mu_x} + C    }   \pare{1+\frac{1}{\beta}} 3 \frac{\sigma^2}{B_\tau} - \pare{C-\frac{2 L^2_f }{\mu_y} }     \frac{1}{T}\sum_{t=0}^{T-1}   \E\norm{\bz^{t} - g(\bx^{t-1})}^2 .
    \end{align*}
    We choose $\beta = \frac{1}{2}$, and $C = \frac{12G_g^2 L^2 \tau }{\mu_x B } +  \frac{60 \tau}{\mu_y B} + \frac{12L^2_f }{\mu_y} + \frac{L_f^2G^2_g}{\mu_x} $, so we have:
    \begin{align*}
       \pare{\frac{12G_g^2 L^2 \tau }{\mu_x B } +  \frac{60 \tau}{\mu_y B} + \frac{8L^2_f }{\mu_y} + \frac{L_f^2G^2_g}{\mu_x} + C    }  \pare{1-\beta}  -\pare{ C-\frac{2 L^2_f }{\mu_y}} = 0. 
    \end{align*}

    Define $\Delta_{\bx} = \frac{1}{T}\sum_{t=0}^{T-1}\E\norm{\bx^{t+1} - \bx^t}^2$ and $\Delta_{\by} = \frac{1}{T}\sum_{t=0}^{T-1}\E\norm{\by^{t+1} - \by^t}^2$, and we have
   {\small \begin{align*}
       &  \frac{1}{T}\sum_{t=0}^{T-1}\cF_k^{t+1}   
               \leq \Gamma^{0} + \frac{1}{2\eta_{\by} T} \pare{ \norm{\by - \by^0}^2 - \norm{\by - \by^{T}}^2 } +\pare{\frac{1}{2\eta_{\bx}  T} - \frac{\mu_x}{2T}} \pare{ \norm{\bx - \bx^0}^2 - \norm{\bx - \bx^{T}}^2 }+  \frac{\tilde\sigma^2}{\mu_x B_{\tau}}+  \frac{10\tilde\sigma^2}{\mu_y B_{\tau}} \\ 
                 &   - \pare{\frac{1}{4\eta_{\by}} -\eta_{\by}L_f^2 - \frac{2G_g^2 L^2_f\tau}{\mu_x B} -\frac{10\tau L^2_f}{\mu_y B}  } \Delta_{\by} \\
                 &- \pare{\frac{1}{2\eta_{\bx}} - \frac{L'_f}{2} - \eta_{\by}L_f^2 G_g^2 -  \frac{6G_g^4 L^2 \tau}{\mu_x B }   - \frac{2G_f^2 L_g^2
 \tau}{\mu_x B} - \frac{30G_g^2 L^2_f \tau}{\mu_y B} - \pare{\frac{72G_g^2 L^2_f \tau }{\mu_x B } +  \frac{360 \tau}{\mu_y B} + \frac{60L^2_f }{\mu_y} + \frac{6L_f^2G^2_g}{\mu_x}    }  (6 +   \frac{\tau}{B}) G_g^2  } \Delta_{\bx} \\
 & + \pare{\frac{24G_g^2 L^2 \tau }{\mu_x B } +  \frac{120 \tau}{\mu_y B} + \frac{20L^2_f }{\mu_y} + \frac{2L_f^2G^2_g}{\mu_x}    }  9 \frac{\tilde\sigma^2}{B_\tau}.
    \end{align*} }
    We choose $\eta_{\by} \leq \Theta\pare{\min\cbr{ \frac{1}{L_f}}, \frac{\mu_y B}{\tau G_g^2 L_f^2}}$ and $\eta_{\bx} \leq \Theta\pare{\min\cbr{\frac{\mu_y B}{G_g^2 L^2_f \tau}, \frac{1}{L'_f}, \frac{\rho B^2 }{G_g^4 L^2_f \tau^2} , \frac{\mu_y B^2}{\tau^2} }}$ which yields:
    \begin{align*}
        \frac{1}{T}\sum_{t=0}^{T-1}\cF_k^{t+1}   
               \leq  &\Gamma^{0} + \frac{1}{2\eta_{\by} T} \pare{ \norm{\by - \by^0}^2 - \norm{\by - \by^{T}}^2 } +\pare{\frac{1}{2\eta_{\bx}  T} - \frac{\mu_x}{2T}} \pare{ \norm{\bx - \bx^0}^2 - \norm{\bx - \bx^{T}}^2 } \\
               &+  \frac{\tilde\sigma^2}{\mu_x B_{\tau}}+  \frac{10\tilde\sigma^2}{\mu_y B_{\tau}} + \pare{\frac{24G_g^2 L^2 \tau }{\mu_x B } +  \frac{120 \tau}{\mu_y B} + \frac{20L^2_f }{\mu_y} + \frac{2L_f^2G^2_g}{\mu_x}    }  9 \frac{\tilde\sigma^2}{B_\tau}.
               \end{align*}
        Notice that
        \begin{align*}
            \frac{1}{T}\sum_{t=0}^{T-1}\cF_k^{t+1} &=  \frac{1}{T}\sum_{t=0}^{T-1} F_k(\bx^{t+1},\bx) - F_k(\bx,\by^{t+1}) + C (\E\norm{\bz^{t+1}- g(\bx^t)}^2 - \E\norm{\bz^{t}- g(\bx^{t-1})}^2) \\
           & =\frac{1}{T}\sum_{t=0}^{T-1} F_k(\bx^{t+1},\bx) - F_k(\bx,\by^{t+1}) + \frac{1}{T} C (\E\norm{\bz^{T}- g(\bx^{T-1})}^2 - \E\norm{\bz^{0}- g(\bx^{-1})}^2)\\
           & \geq \frac{1}{T}\sum_{t=0}^{T-1} F_k(\bx^{t+1},\bx) - F_k(\bx,\by^{t+1})       -  \frac{1}{T} C\E\norm{\bz^{0}- g(\bx^{-1})}^2.
        \end{align*}
        Since we choose   $\eta_{\bx} \leq \frac{1}{2L}$ and by convention $\bx^{-1} = \bx^0$, we have
               \begin{align*}  
            F_k(\bx_{k+1},\by) - F_k(\bx,\by_{k+1}) \leq &  \frac{1}{T}\sum_{t=0}^{T-1} F_k(\bx^{t+1},\by) - F_k(\bx,\by^{t+1}) + O\pare{\frac{C\E\norm{\bz^{0}- g(\bx^{0})}^2}{T}}  \\
                \leq &  \frac{ \E\norm{\by - \by^0}^2}{2\eta_{\by} T}  + \frac{ \E \norm{\bx - \bx^0}^2 }{4\eta_{\bx}  T} + O\pare{\frac{C\E\norm{\bz^{0}- g(\bx^{0})}^2}{T}}     +  \frac{\tilde\sigma^2}{\mu_x B_{\tau}}+  \frac{10\tilde\sigma^2}{\mu_y B_{\tau}} \\ 
                & + \pare{\frac{24G_g^2 L^2_f \tau }{\mu_x B } +  \frac{120 \tau}{\mu_y B} + \frac{20L^2_f }{\mu_y} + \frac{2L_f^2G^2_g}{\mu_x}    }  9 \frac{\tilde\sigma^2}{B_\tau}.
    \end{align*}

    \end{proof}
\end{lemma}

\begin{lemma}\label{lem:gap relation}
For Algorithm~\ref{algorithm: CODA-Primal+}, under assumptions of Theorem~\ref{thm:NCSC+}, the following statement holds:
        \begin{align*}
        \pare{ Q -1 } & \mathrm{Gap}^k( \bx_{k+1},\by_{k+1})  \leq   \mathrm{Gap}^k(\bx_{k},\by_{k}) \\
       & + Q \pare{   \frac{\tilde\sigma^2}{\mu_x B_{\tau}}+ \frac{10\tilde\sigma^2}{\mu_y B_{\tau}}  + O\pare{\frac{C\E\norm{\bz^{0}- g(\bx^{0})}^2}{T}} + \pare{\frac{24G_g^2 L^2_f \tau }{\mu_x B } +  \frac{120 \tau}{\mu_y B} + \frac{20L^2_f }{\mu_y} + \frac{2L_f^2G^2_g}{\mu_x}    }  9 \frac{\tilde\sigma^2}{B_\tau} },
        \end{align*}
        where $Q = \frac{\min\cbr { \frac{\mu_x}{4},\frac{\mu_y}{4} } }{\max\cbr{\frac{1}{2\eta_{\by}T},\frac{1}{4\eta_{\bx}T}}}$.
    \begin{proof}
        Due to~[Lemma. 1~\cite{yan2020optimal}]
        \begin{align*}
            \frac{\mu_x}{4} \norm{\bx^*_k(\by) - \bx'}^2 +\frac{\mu_y}{4} \norm{\by^*_k(\bx) - \by'}^2  \leq \mathrm{Gap}^k(\bx,\by) + \mathrm{Gap}^k(\bx',\by') .
        \end{align*}
        We set $\bx = \bx_{k+1}, \by = \by_{k+1}$, $\bx' = \bx_{k}, \by' = \by_{k}$, which yields 
        \begin{align}
            \frac{\mu_x}{4} \norm{\bx^*_k(\by_{k+1}) - \bx_{k}}^2 +\frac{\mu_y}{4} \norm{\by^*_k( \bx_{k+1}) - \by_{k}}^2  \leq \mathrm{Gap}^k( \bx_{k+1},\by_{k+1}) + \mathrm{Gap}^k(\bx_{k},\by_{k}) . \label{eq:dist to gap}
        \end{align}

        Due to Lemma~\ref{lem:PD conevrgence VR} we have
          \begin{align}  
            \mathrm{Gap}^k( \bx_{k+1},\by_{k+1})  
                \leq &  \max\cbr{\frac{1}{2\eta_{\by}T},\frac{1}{4\eta_{\bx}T}} \pare{ \E\norm{\by^*_k( \bx_{k+1}) - \by_k}^2  +   \E \norm{\bx^*_k(\by_{k+1}) - \bx_k}^2}   +  \frac{C\E\norm{\bz^{0}- g(\bx^{0})}^2}{T}   \nonumber\\
               &   +  \frac{\tilde\sigma^2}{\mu_x B_{\tau}}+  \frac{10\tilde\sigma^2}{\mu_y B_{\tau}} + \pare{\frac{24G_g^2 L^2_f \tau }{\mu_x B } +  \frac{120 \tau}{\mu_y B} + \frac{20L^2_f }{\mu_y} + \frac{2L_f^2G^2_g}{\mu_x}    }  9 \frac{\tilde\sigma^2}{B_\tau}.\label{eq:gap to dist} 
    \end{align}
    Combining (\ref{eq:dist to gap}) and (\ref{eq:gap to dist}) yields:
    \begin{align*}
         &Q \pare{ \mathrm{Gap}^k( \bx_{k+1},\by_{k+1}) -    \frac{\sigma^2}{\mu_x B_{\tau}}- \frac{10\sigma^2}{\mu_y B_{\tau}}  - \frac{C\E\norm{\bz^{0}- g(\bx^{0})}^2}{T} - \pare{\frac{24G_g^2 L^2_f \tau }{\mu_x B } +  \frac{120 \tau}{\mu_y B} + \frac{20L^2_f }{\mu_y} + \frac{2L_f^2G^2_g}{\mu_x}    }  9 \frac{\tilde\sigma^2}{B_\tau}  } \\
         & \leq  \mathrm{Gap}^k( \bx_{k+1},\by_{k+1}) + \mathrm{Gap}^k(\bx_{k},\by_{k}) \\
         \Longleftrightarrow & \pare{Q -1 }  \mathrm{Gap}^k( \bx_{k+1},\by_{k+1})  \leq   \mathrm{Gap}^k(\bx_{k},\by_{k}) \\
         & \qquad + Q \pare{   \frac{\sigma^2}{\mu_x B_{\tau}}+ \frac{10\sigma^2}{\mu_y B_{\tau}}  +\frac{C\E\norm{\bz^{0}- g(\bx^{0})}^2}{T}+ \pare{\frac{24G_g^2 L^2_f \tau }{\mu_x B } +  \frac{120 \tau}{\mu_y B} + \frac{20L^2_f }{\mu_y} + \frac{2L_f^2G^2_g}{\mu_x}    }  9 \frac{\tilde\sigma^2}{B_\tau}  }.
        \end{align*}
    \end{proof}
\end{lemma}
\begin{lemma}[Lemma 8~\cite{yan2020optimal}] \label{lem:three inequalities}
If $F(\bx,\by)$ is $\rho$ weakly convex in $\bx$, and let $F_k(\bx,\by) = F(\bx,\by)+\frac{\mu_x + \rho}{2}\norm{\bx - \bx_k}^2$. Then the following statements hold true for any $\alpha_1,\alpha_2 \in (0,1)$:
    \begin{align}
      F_k( \bx_{k+1}, \by ) - F_k(\bx,  \by_{k+1} )  & \geq \pare{1- \frac{ \mu_x + \rho}{\rho} \pare{\frac{1}{\alpha_1} -1 }} \pare{ F_{k+1}( \bx_{k+1}, \by ) - F_{k+1}(\bx,  \by_{k+1} )} \nonumber \\
       &- \frac{\mu_x + \rho}{2} \frac{\alpha_1}{1-\alpha_1} \norm{ \bx_{k+1} - \bx_{k}}^2  \label{eq:three inequalities 1}   \\
        F_k( \bx_{k+1}, \by ) - F_k(\bx,  \by_{k+1} ) &\geq   \Phi( \bx_{k+1}) - \Phi(\bx_{k})  +    \frac{\mu_x + \rho}{2}\norm{  \bx_{k+1} - \bx_{k}}^2 \label{eq:three inequalities 2}\\
         F_k( \bx_{k+1}, \by ) - F_k(\bx,  \by_{k+1} ) & \geq   \frac{\rho \alpha_2}{2}\norm{   \bx_{k} - \bx_{k}^*}^2 - \frac{ \rho \alpha_2}{2(1-\alpha_2)} \norm{ \bx_k - \bx_{k+1}}^2. \label{eq:three inequalities 3}
    \end{align}
\end{lemma}

\begin{lemma}\label{lem:recursion}
For Algorithm~\ref{algorithm: CODA-Primal+}, under assumptions of Theorem~\ref{thm:NCSC+}, the following statement holds:
    \begin{align*}
     (Q-1)  &\frac{a_3 \rho \alpha_2}{2}   \norm{\bx_{k} - \bx_{k}^*}^2 \leq \mathrm{Gap}^k(\bx_{k},\by_{k})   - (Q-1) a_1 \pare{1- \frac{ \mu_x + \rho}{\rho} \pare{\frac{1}{\alpha_1} -1 }  }\mathrm{Gap}^{k+1}(\bx_{k+1},\by_{k+1}) \\
       & \quad - (Q-1) a_2 \pare{ \Phi(\bx_{k+1}) - \Phi(\bx_{k}) }\\
       & \quad + Q \pare{   \frac{\sigma^2}{\mu_x B_{\tau}}+ \frac{10\sigma^2}{\mu_y B_{\tau}}  +\frac{C\E\norm{\bz^{0}- g(\bx^{0})}^2}{T}+ \pare{\frac{24G_g^2 L^2_f \tau }{\mu_x B } +  \frac{120 \tau}{\mu_y B} + \frac{20L^2_f }{\mu_y} + \frac{2L_f^2G^2_g}{\mu_x}    }  9 \frac{\tilde\sigma^2}{B_\tau}  },  
    \end{align*}
    where $a_1, a_2, a_3$ are such that
    \begin{align*}
        a_1 + a_2 + a_3 &= 1,\\
        a_1 \frac{(\mu_x+\rho) \alpha_1}{2(1-\alpha_1)} - a_2\frac{\mu_x + \rho}{2} + a_3\frac{ \rho \alpha_2}{2(1-\alpha_2)} &\leq 0,\\
        (Q-1) a_1 \pare{1- \frac{ \mu_x + \rho }{\rho} \pare{\frac{1}{\alpha_1} -1 }  }  &\geq 1.
    \end{align*}
    \begin{proof}
    Due to Lemma~\ref{lem:three inequalities}, $a_1 \cdot (\ref{eq:three inequalities 1}) +a_2 \cdot (\ref{eq:three inequalities 2}) + a_3 \cdot (\ref{eq:three inequalities 3})$ yields:
        \begin{align*}
       \frac{a_3 \rho \alpha_2}{2}   \norm{  \bx_{k} - \bx_{k}^*}^2 &\leq  \mathrm{Gap}^k( \bx_{k+1},\by_{k+1}) - a_1 \pare{1- \frac{ \mu_x+\rho}{\rho} \pare{\frac{1}{\alpha_1} -1 }  } \mathrm{Gap}^{k+1}( \bx_{k+1},\by_{k+1}) \\
       & \quad - a_2 \pare{ \Phi( \bx_{k+1}) - \Phi( \bx_{k}) } + \pare{ a_1 \frac{\mu_x + \rho}{2} \frac{\alpha_1}{1-\alpha_1} +   a_3\frac{ \rho \alpha_2}{2(1-\alpha_2)} - a_2 \frac{\mu_x + \rho}{2} }\norm{  \bx_{k+1} - \bx_{k}}^2 , 
    \end{align*}
    Multiply both sides with $Q-1$, plugging result from Lemma~\ref{lem:gap relation},  then the proof immediately follows the condition $a_1 \frac{(\mu_x + \rho) \alpha_1}{2(1-\alpha_1)} - a_2\frac{\mu_x + \rho}{2} + a_3\frac{ \rho \alpha_2}{2(1-\alpha_2)}  \leq 0$ and $(Q-1) a_1 \pare{1- \frac{ \mu_x + \rho }{\rho} \pare{\frac{1}{\alpha_1} -1 }  }  \geq 1$.
    \end{proof}
\end{lemma}

\subsection{Proof of Theorem~\ref{thm:NCSC+}}\label{app: proof ncsc+}
\begin{proof}
    First in Lemma~\ref{lem:recursion}, we set $ \mu_x = \rho = L$, and we choose $\eta_{\bx} \leq \frac{1}{2}\eta_{\by}$, together with assumption that $L \geq \mu_y$, we know $Q = \eta_{\bx} T  \mu_y$. We set:
    \begin{align}
        a_1 = \frac{1}{12}, a_2 = \frac{10}{12}, a_3= \frac{1}{12}, \alpha_1= \frac{3}{4}, \alpha_2 = \frac{3}{4}, \eta_{\bx}T \mu_y \geq 37 \label{eq: cond NCSC+}
    \end{align} , the condition in Lemma~\ref{lem:recursion} can hold, and we have: 
    \begin{align*}
     (\eta_{\bx}T \mu_y -1)  \frac{  3L}{96}   \norm{\bx_{k} - \bx_{k}^*}^2 &\leq \mathrm{Gap}^k(\bx_{k},\by_{k})   - \mathrm{Gap}^{k+1}(\bx_{k+1},\by_{k+1}) \\
       & \quad - (\eta_{\bx}T \mu_y-1) \frac{5}{6}\pare{ \Phi(\bx_{k+1}) - \Phi(\bx_{k}) } \\
       & \quad + \eta_{\bx}T \mu_y \pare{   \frac{\tilde\sigma^2}{\mu_x B_{\tau}}+ \frac{10\tilde\sigma^2}{\mu_y B_{\tau}}+ \pare{\frac{24G_g^2 L^2_f \tau }{\mu_x B } +  \frac{120 \tau}{\mu_y B} + \frac{20L^2_f }{\mu_y} + \frac{2L_f^2G^2_g}{\mu_x}    }  9 \frac{\tilde\sigma^2}{B_\tau} },  
    \end{align*}
    Due to the fact $4L^2  \norm{\bx_{k} - \bx_{k}^*}^2 =\norm{\Phi_{1/2L}(\bx_k)}^2$, and summing over $k=0$ to $K-1$ we have
     \begin{align*}
    &  (\eta_{\bx}T \mu_y -1) \frac{1}{K} \sum_{k=0}^{K-1}\frac{ 3 L}{384L^2}   \norm{\nabla \Phi_{1/2L}(\bx_k)}^2  \leq \frac{\mathrm{Gap}^k(\bx_{0},\by_{0})}{K}  + (\eta_{\bx}T \mu_y-1) \frac{5(\Phi(\bx_{0}) - \min_{\bx\in\cX} \Phi(\bx)) }{6K} \\
      & \quad + \eta_{\bx}T \mu_y O \pare{ \frac{C\E\norm{\bz^{0}- g(\bx^{0})}^2}{T}+  \frac{  \tilde\sigma^2}{\mu_x B_{\tau}}+ \frac{  \tilde\sigma^2}{\mu_y B_{\tau}} + \pare{\frac{ G_g^2 L^2_f \tau }{\mu_x B } +  \frac{ \tau}{\mu_y B} + \frac{ L^2_f }{\mu_y} + \frac{ L_f^2G^2_g}{\mu_x}    }    \frac{\tilde\sigma^2}{B_\tau}  },  \\
       \Longleftrightarrow   &  \frac{1}{K} \sum_{k=0}^{K-1}    \norm{\nabla \Phi_{1/2L}(\bx_k)}^2 \leq \frac{384 L\mathrm{Gap}^k(\bx_{0},\by_{0})}{3(\eta_{\bx}T \mu_y -1) K}  +   \frac{1920L \Delta_{\Phi}}{18 K} \\
       & \quad +  O \pare{  \frac{C\E\norm{\bz^{0}- g(\bx^{0})}^2}{T}+  \frac{L\tilde \sigma^2}{\mu_x B_{\tau}}+ \frac{ L \tilde\sigma^2}{\mu_y B_{\tau}} + \pare{\frac{ G_g^2 L^2_f \tau }{\mu_x B } +  \frac{ \tau}{\mu_y B} + \frac{ L^2_f }{\mu_y} + \frac{ L_f^2G^2_g}{\mu_x}    }    \frac{L \tilde\sigma^2}{B_\tau}  }.
    \end{align*} 
    Recall that by our initialization method, $\E\norm{\bz^{0}- g(\bx^{0})}^2 \leq \delta$.
    To ensure RHS is less than $\epsilon^2$, we choose $K = \Theta\pare{\frac{L\Delta_{\Phi}}{\epsilon^2}}$, $B_{\tau} = \Theta\pare{\frac{L^2 \sigma^2}{\mu_y \epsilon^2}}$, $\tau = \sqrt{\frac{B_{\tau}}{\kappa}}$ and $B = \sqrt{ \kappa B_\tau}$.
    
    It remains to determine value of $T$. According to condition~\ref{eq: cond NCSC+}, we need $T \geq \frac{6}{\mu_y \eta_{\bx}}$. Since 
    \begin{align*}
        \eta_{\bx} &\leq \Theta\pare{\min\cbr{\frac{\mu_y B}{G_g^2 L^2_f \tau}, \frac{1}{L'_f}, \frac{\rho B^2 }{G_g^4 L^2_f \tau^2} , \frac{\mu_y B^2}{\tau^2} }} \\
        &\leq \Theta\pare{\min\cbr{\frac{\mu_y \kappa }{G_g^2 L^2_f  }, \frac{1}{L }, \frac{\rho \kappa^2 }{G_g^4 L^2_f  } ,  \mu_y \kappa^2     }} = \Theta \pare{\frac{1}{L}},
    \end{align*}
    we know that $T \geq \Omega (\kappa)$. And hence we should choose $\delta \leq O(\kappa\epsilon^2/C) = O( \epsilon^2/L)$, since   $C = \frac{12G_g^2 L^2 \tau }{\mu_x B } +  \frac{60 \tau}{\mu_y B} + \frac{12L^2_f }{\mu_y} + \frac{L_f^2G^2_g}{\mu_x}  = O(\kappa L)$.  The final gradient complexity is 
    \begin{align*}
        \frac{TKB_\tau}{\tau} + TKB = O\pare{ \frac{\kappa^2 \sqrt{L} \Delta_\Phi \sigma}{\epsilon^3}       }.
    \end{align*}
\end{proof}

\subsection{Proof of Theorem~\ref{thm:ncc+}}\label{app: proof ncc+}
Let $\tilde \Phi(\bx) = \max_{\by \in \cY} \tilde F(\bx,\by)$, and $\tilde F(\bx,\by)$ is defined in (\ref{eq:augmented ncc}) .
According to~\cite{zhang2022sapd+} Lemma~27, a $\epsilon/2\sqrt{6}$ stationary point of Moreau Envelope of $\tilde \Phi(\bx)$ is $\epsilon$ stationary point of that of $\Phi(\bx)$. Then the complexity result follows immediately by setting $\mu=\frac{\epsilon^2}{L D^2_{\cY}}$ in Theorem~\ref{thm:NCSC+}.

\bibliography{references}
\bibliographystyle{plain}

\end{document}